\documentclass[11pt]{article}

\usepackage[utf8]{inputenc}
\usepackage[T1]{fontenc}
\usepackage{lmodern}
\usepackage[table]{xcolor}
\usepackage{amssymb}
\usepackage{amsbsy}
\usepackage{amsmath}
\usepackage{amsthm}
\usepackage{amsfonts}
\usepackage{graphicx}
\usepackage{mathtools}
\usepackage{geometry}
\usepackage{MnSymbol}
\usepackage{bm}
\definecolor{darkblue}{rgb}{0,0,.75}
\usepackage[numbers]{natbib}
\usepackage[colorlinks=true,allcolors=darkblue]{hyperref}
\usepackage{booktabs}
\usepackage{microtype}
\usepackage{float}
\usepackage{dsfont}
\usepackage{tikz}
\usepackage{enumitem}
\usepackage{mathrsfs}
\usepackage{upgreek}
\usepackage{placeins}

\mathtoolsset{showonlyrefs,showmanualtags}
\AddToHook{env/subequations/begin}{%
    \let\savedsubequationseqref\eqref
    \let\savedsubequationsrefeq\refeq
    \mathtoolsset{showonlyrefs=false}%
    \let\eqref\savedsubequationseqref
    \let\refeq\savedsubequationsrefeq
}
\allowdisplaybreaks

\makeatletter
\long\def\@makecaption#1#2{
  \vskip 0.8ex
  \setbox\@tempboxa\hbox{\small {\bf #1:} #2}
  \parindent 1.5em%
  \dimen0=\hsize
  \advance\dimen0 by -3em
  \ifdim \wd\@tempboxa >\dimen0
  \hbox to \hsize{
    \parindent 0em
    \hfil
    \parbox{\dimen0}{\def\baselinestretch{0.96}\small
      {\bf #1.} #2
    }
    \hfil}
  \else \hbox to \hsize{\hfil \box\@tempboxa \hfil}
  \fi
}
\makeatother

\newtheorem{theorem}{Theorem}
\newtheorem{lemma}{Lemma}[section]
\newtheorem{remark}{Remark}
\newtheorem{proposition}{Proposition}
\newtheorem{assumption}{Assumption}
\renewenvironment{proof}{\noindent{\bf Proof}\hspace*{1em}}{\qed\\}

\newcommand{\ind}{\mathord{\text{\usefont{U}{dsrom}{m}{n}1}}}
\newcommand{\E}{\mathbb{E}}
\newcommand{\R}{\mathbb{R}}
\newcommand{\N}{\mathbb{N}}
\renewcommand{\P}{\mathbb{P}}
\newcommand{\brc}[1]{\left\{{#1}\right\}}
\newcommand{\prn}[1]{\left({#1}\right)}
\newcommand{\brk}[1]{\left[{#1}\right]}
\newcommand{\norm}[1]{\left\|{#1}\right\|}
\newcommand{\what}[1]{\widehat{#1}}
\newcommand{\wt}[1]{\widetilde{#1}}
\newcommand{\wb}[1]{\overline{#1}}
\newcommand{\dist}{\mathsf{dist}}
\newcommand{\norms}[1]{\|{#1}\|}
\newcommand{\<}{\langle}
\renewcommand{\>}{\rangle}
\DeclareMathOperator*{\argmin}{argmin}
\newcommand{\half}{\frac{1}{2}}
\newcommand{\mc}[1]{\mathcal{#1}}
\newcommand{\C}{\mathbb{C}}
\newcommand{\tvnorms}[1]{\norms{#1}_{\rm TV}}
\DeclareMathOperator{\Tr}{\mathrm{Tr}}
\newcommand{\supp}{\mathop{\rm supp}}
\newcommand{\rank}{\mathop{\rm rank}}
\newcommand{\bb}{\boldsymbol{b}}
\newcommand{\bh}{\boldsymbol{h}}
\newcommand{\bu}{\boldsymbol{u}}
\newcommand{\bv}{\boldsymbol{v}}
\newcommand{\bx}{\boldsymbol{x}}
\newcommand{\by}{\boldsymbol{y}}
\newcommand{\bz}{\boldsymbol{z}}
\newcommand{\bA}{\boldsymbol{A}}
\newcommand{\bB}{\boldsymbol{B}}
\newcommand{\bC}{\boldsymbol{C}}
\newcommand{\bD}{\boldsymbol{D}}
\newcommand{\bI}{\boldsymbol{I}}
\newcommand{\bM}{\boldsymbol{M}}
\newcommand{\bP}{\boldsymbol{P}}
\newcommand{\bQ}{\boldsymbol{Q}}
\newcommand{\bR}{\boldsymbol{R}}
\newcommand{\bS}{\boldsymbol{S}}
\newcommand{\bT}{\boldsymbol{T}}
\newcommand{\bU}{\boldsymbol{U}}
\newcommand{\bV}{\boldsymbol{V}}
\newcommand{\bW}{\boldsymbol{W}}
\newcommand{\bX}{\boldsymbol{X}}
\newcommand{\bY}{\boldsymbol{Y}}
\newcommand{\bZ}{\boldsymbol{Z}}
\newcommand{\bbeta}{\boldsymbol{\beta}}
\newcommand{\btheta}{\boldsymbol{\theta}}
\newcommand{\bSigma}{\boldsymbol{\Sigma}}
\newcommand{\const}{\mathsf{C}}
\newcommand{\sconst}{\mathsf{c}}
\newcommand{\Cov}{\mathrm{Cov}}
\newcommand{\de}{\mathrm{d}}
\renewcommand{\Im}{\operatorname{Im}}
\renewcommand{\Re}{\operatorname{Re}}
\newcommand{\Spec}{\operatorname{Spec}}
\newcommand{\bzeros}{\boldsymbol{0}}
\newcommand{\bones}{\boldsymbol{1}}
\newcommand{\lin}{^\textup{lin}}
\newcommand{\bvarepsilon}{\boldsymbol{\varepsilon}}
\newcommand{\iu}{\mathrm{i}}
\newcommand{\sampleRisk}{\mathscr{R}_{\bX}}
\newcommand{\sampleVariance}{\mathscr{V}_{\bX}}
\newcommand{\sampleBias}{\mathscr{B}_{\bX}}
\newcommand{\effectiveRisk}{\mathsf{R}_n}
\newcommand{\effectiveVariance}{\mathsf{V}_n}
\newcommand{\effectiveBias}{\mathsf{B}_n}
\title{Dimension-free spectral laws and the risk of spectral filters}

\author{Chen Cheng\thanks{Department of Statistics, University of Illinois Urbana-Champaign}}

\begin{document}

    \maketitle
    \begin{abstract}
    Spectral properties of the sample covariance are key to modern machine learning with high-dimensional data $\bx_1, \dots, \bx_n$. Random matrix theory (RMT) supplies powerful tools; however, most classical results assume proportional asymptotics, where $d/n$ tends to a finite limit. In functional data analysis and kernel methods, a fixed infinite-dimensional Hilbert space $\mc{H}$ is often the more natural setting for the observations, beyond the scope of proportional asymptotics.

    This calls for a dimension-free spectral theory for i.i.d.\ observations $\bx_1,\ldots,\bx_n$ in a Hilbert space $\mc H$ with trace-class population covariance $\bSigma$. Building on the dimension-free ridge regression analysis of~\citet{cheng2024dimensionfree}, we establish two deterministic equivalents with explicit nonasymptotic error bounds. The first approximates the empirical spectral law $\what\mu_n$ of the sample Gram matrix by a deterministic law $\mu_n$. The second characterizes the conditional prediction risk of analytic spectral filters in linear regression.

    Previous sharp approximations to prediction risk in this setting have focused on ridge regression and closely-related problems. Our analysis covers a broader class of filters at the bulk scale, with an explicit user-friendly spectral stability condition to check. We apply our theory to power and critical spectral laws, and we recover the proportional regime. For prediction, we study (negative) ridge estimators, iterated ridge, and gradient flow.

    Mathematically, building on the idea of stochastic interpolation~\citep{cheng2024dimensionfree}, we construct a deterministic backward path to interpolate the spectral parameter. The resulting bounds on resolvent and mixed-resolvent traces give mesoscopic approximations to the spectral law within the bulk and, by Cauchy’s formula, deterministic equivalents for the prediction variance and squared bias of analytic spectral filters.
\end{abstract}

    \section{Introduction}
The central problem in statistical learning is to infer structure from an i.i.d.\ sample $\bx_1, \dots, \bx_n$, with or without paired observations $y_1, \dots, y_n$. The accompanying central challenge of high-dimensionality in modern machine learning arises from increasingly complex datasets (e.g., ImageNet~\citep{deng2009imagenet}, LAION-5B~\citep{schuhmann2022laion}, and FineWeb~\citep{penedo2024fineweb}) and model architectures (e.g., diffusion~\citep{ho2020denoising,rombach2022latent} and LLMs~\citep{brown2020fewshot,grattafiori2024llama}), going beyond the classical understanding of large-sample asymptotics~\citep{vandervaartwellner2023weak} and sparsity~\citep{wainwright2019highdimensional}. Random matrix theory (RMT) has become popular in theoretical machine learning in those regimes where learning is possible but challenging. Not only useful for modeling frontier learning models such as kernel methods, random features, and neural networks~\citep{liang2020just,misiakiewicz2024nonasymptotic,defilippis2024dimensionfree,montanari2022interpolation}, it is also successful in understanding many empirical phenomena in machine learning, such as double-descent, benign-overfitting, memorization, and scaling laws~\citep{bartlett2020benign,mei2022generalization,cheng2022memorize,cheng2026memorization,paquette20244+,atanasov2026scaling}.

\begin{figure}[htbp!]
    \centering
    \includegraphics[width=.92\textwidth]{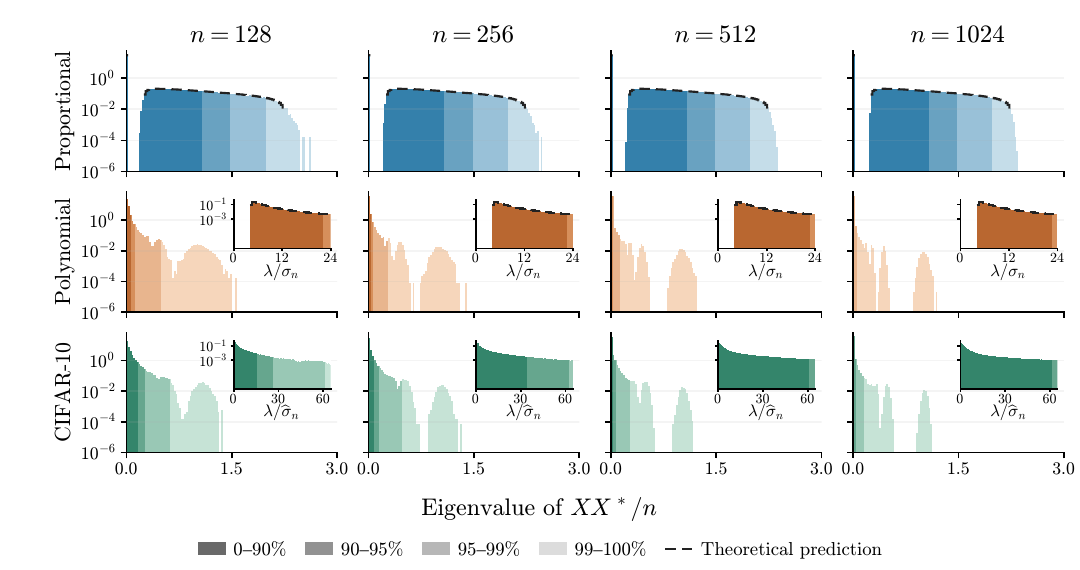}
    \caption{Spectra of $\bX\bX^*/n$ pooled over 2,000 trials, of (i) proportional isotropic Gaussian features with $n=4d$, (ii) Zipf's law with i.i.d.\ Gaussian coordinates, $\sigma_i=i^{-1.2}$ and (iii) CIFAR-10 ResNet-50 second-to-last layer features ($d=2,048$). The zoomed in view show rescaled histograms. For CIFAR-10, we estimate $\what{\sigma}_n$ from the separate training set.}
    \label{fig:intro-gram-spectra}
\end{figure}

Existing tools in RMT primarily focus on the proportional setup when $n/d$ converges to a constant, in which $\bx_i \in \R^d$ has population covariance $\bSigma$ and the sample covariance $\what{\bSigma} = \bX^* \bX/n$ does not concentrate~\citep{marchenko1967distribution}. This setup, including its extension to polynomial scalings when $n/d^\ell$ converges~\citep{bloemendal2014isotropic,lu2025equivalence}, still allows the input dimension to grow. In functional data analysis and in regression with a fixed kernel, it is natural to increase $n$ while keeping the underlying Hilbert space $\mc H$ and its covariance $\bSigma$ fixed, even when $\mc H$ is infinite-dimensional. Empirically, an infinite-dimensional $\bSigma$ whose population spectrum $\sigma_1 \geq \sigma_2 \geq \dots$ is trace-class ($\Tr \bSigma < \infty$), e.g., Zipf's law with polynomial decay, $\sigma_i = i^{-1-\alpha}$ for some $\alpha > 0$, is often a better model for real-world data, as shown in Fig.~\ref{fig:intro-gram-spectra}. We observe that the sample spectrum from a fixed Hilbert space with growing size $n$ is a better proxy for the sample spectrum of CIFAR-10 features.

The pioneering work of~\citet{cheng2024dimensionfree} addressed this setting for the prediction risk of ridge regression problems, with later works studying closely-related problems~\citep{misiakiewicz2024nonasymptotic,defilippis2024dimensionfree}. However, a more applicable and complete dimension-free spectral theory at the bulk scale, beyond ridge regression, is still missing. In this work, we develop nonasymptotic and deterministic approximations for (i) spectral parameters and laws and (ii) the prediction risk of a general class of spectral filters. Given i.i.d.\ $\bx_1, \cdots, \bx_n \in \mc{H}$, we study the two following problems:
\begin{enumerate}
    \item \emph{Unsupervised learning: the empirical spectrum.} How can we characterize the (top) empirical spectrum of the sample covariance $\what{\bSigma}$?

    Equivalently, we consider the Gram matrix $\bX \bX^*/n$. Let $\what{\sigma}_1 \geq \dots \geq \what{\sigma}_n \geq 0$ be the sample eigenvalues of $\bX \bX^*/n$. We seek a sharp characterization of the unnormalized empirical spectral distribution
    \begin{align}
        \what{\mu}_n := \frac{1}{n} \sum_{j=1}^n \delta_{n \what{\sigma}_j}\label{eq:intro-gram-esd}
    \end{align}
    at the bulk scale, i.e., on any interval $E \subset \R$ such $\what{\mu}_n(E) \in (c, 1-c)$ for some $c > 0$.
    \item \emph{Supervised learning: risk of spectral filters.} How can we characterize the prediction error of a spectral filter $\upphi$?

    We consider a linear model $\by=\bX\bbeta+\bvarepsilon$ with noise $\E[\bvarepsilon\mid\bX]=0$. The spectral filter $\upphi$ multiplies each sample eigendirection at eigenvalue $x$ by $\upphi(x)$. Let the residual function be $\uppsi(x)=1-x\upphi(x)$; this yields the estimator $\what{\bbeta}_{\upphi}=\upphi(\what{\bSigma})\bX^*\by/n$, which covers a broad class of spectral algorithms, e.g., $\uppsi(x)=\lambda/(x+\lambda)$ for ridge, $\uppsi(x)=(\lambda/(x+\lambda))^k$ for iterated ridge, $\uppsi(x)=e^{-tx}$ for gradient flow, and $\uppsi(x)=\ind\{x<\lambda\}$ for principal component regression (PCR). We seek sharp characterizations for the conditional generalization error
    \begin{align}
        \sampleRisk[\upphi]:=\E\big[|\<\bx_{\mathrm{new}},\what{\bbeta}_{\upphi}-\bbeta\>|^2 \mid \bX\big] =\sampleVariance[\upphi]+\sampleBias[\upphi].\label{eq:intro-filter-risk}
    \end{align}
    Here $\bx_{\mathrm{new}}$ denotes a test feature with the same law as $\bx_1$, drawn independently of the training data. The expectation averages over the test feature and training noise conditional on $\bX$. The term $\sampleVariance[\upphi]=\tau^2\Tr\prn{\bSigma\what{\bSigma}\upphi(\what{\bSigma})^2}/n$ is the conditional variance when $\E[\bvarepsilon\bvarepsilon^*\mid\bX]=\tau^2\bI_n$, and $\sampleBias[\upphi]=\|\bSigma^{1/2}\uppsi(\what{\bSigma})\bbeta\|^2$ is the squared bias.
\end{enumerate}

\subsection{Our contributions}
We give our solutions to the two problems, and summarize our contributions in the paper.
\begin{enumerate}
    \item \emph{Deterministic equivalents to the empirical spectrum.} In Section~\ref{sec:gram-bulk-law}, we track the empirical spectral law $\what{\mu}_n$ and its Stieltjes transform $\what{s}_n$ by their deterministic counterparts:
    \begin{align*}
        (s_n(\bSigma),\mu_n(\bSigma)) \;\mathrel{\substack{
        \textstyle\xleftarrow{\text{\footnotesize deterministic equivalent}}\\[-0.6ex]
        \textstyle\xrightarrow[\text{\footnotesize closely approximates}]{}
        }}\; (\what{s}_n(\bX),\what{\mu}_n(\bX)).
    \end{align*}
    At fixed height on the bulk scale, the Stieltjes transform bound captures the optimal $n^{-1}$ power for averaged fluctuations (cf.~\citep{knowles2017anisotropic}).

    \item \emph{Deterministic equivalents for prediction risk.}
    In Section~\ref{sec:spectral-filter-risk}, for a given analytic spectral filter, at fixed signal and noise level, we obtain
    \begin{align*}
        (\effectiveVariance[\upphi;\bSigma],\effectiveBias[\upphi;\bSigma]) \;\mathrel{\substack{
        \textstyle\xleftarrow{\text{\footnotesize deterministic equivalent}}\\[-0.6ex]
        \textstyle\xrightarrow[\text{\footnotesize closely approximates}]{}
        }}\; (\sampleVariance[\upphi],\sampleBias[\upphi]).
    \end{align*}
    Our bounds capture the optimal corresponding averaged and local fluctuation rate, $n^{-1}$ and $n^{-1/2}$, respectively.

    \item \emph{Technical innovation: the backward interpolation path.} In Section~\ref{sec:interpolation-path}, we present our technical core , where we achieve the preceding goals by approximating weighted traces of resolvents and mixed resolvents. For a deterministic weight $\bQ$, we have
    \begin{align*}
        \Tr\prn{\bQ\bSigma(\bSigma-\omega(z/n)\bI)^{-1}} \;\mathrel{\substack{
        \textstyle\xleftarrow{\text{\footnotesize deterministic equivalent}}\\[-0.6ex]
        \textstyle\xrightarrow[\text{\footnotesize closely approximates}]{}
        }}\; \Tr\prn{\bQ\bX^*(\bX\bX^*-z\bI_n)^{-1}\bX}.
    \end{align*}
    Similar deterministic equivalents exist for the mixed-resolvent traces. We obtain these comparisons along a deterministic backward path, inspired by and innovated upon the analysis of~\citet{cheng2024dimensionfree}. This adjustment cancels the leading drift at each step and connects the sample trace to its population counterpart while keeping each step deterministic.
\end{enumerate}
We emphasize that all our approximations are nonasymptotic and that all deterministic equivalents can vary in $n$, even if they could have a limit in $n$. For power law and critical decay, the spectral approximation remains valid at mesoscopic heights $h=n^{-\delta}$ in the interior of the bulk, with $h$ measured relative to the bulk width and $0<\delta<1/17$.

\subsection{Related work}

\begin{description}
    \item[Random matrix theory.] The study of random sample covariance matrices is classical and dates back to Wishart's characterization of their distribution under Gaussian sampling~\citep{wishart1928generalised}. The Marchenko-Pastur law subsequently described their limiting empirical spectrum when dimension and sample size grow proportionally~\citep{marchenko1967distribution}. More recently, these global descriptions have been sharpened to local and anisotropic laws, controlling resolvents and spectral statistics on increasingly fine scales~\citep{knowles2017anisotropic}. Very recent work also establishes quantitative approximations for mixed-resolvent traces in the proportional regime~\citep{stark2026quantitative}.

    Relevant works include extensions to polynomial dimensional scalings~\citep{ding2023global}. However, these settings do not directly address an infinite decaying population spectrum with a shrinking empirical bulk. Removing restrictions on the population condition number is essential here, since even finite-dimensional truncations become increasingly ill-conditioned. \citet{cheng2024dimensionfree} provided a methodological foundation through their leave-one-out interpolation strategy. Although they do not directly derive spectral laws and local laws over the bulk, their ridge formulas suggest the deterministic equivalent through analytic continuation.

    \item[Kernel random matrices.]
    Our work closely relates to kernel methods, which provide a natural Hilbert-space setting through a possibly infinite-dimensional feature map, even when the input space is fixed~\citep[Section~4.5]{steinwart2008support}. Despite having the Hilbert space embedding, much of the kernel random matrix literature instead considers growing input dimension.  The foundational work of \citet{elkaroui2010spectrum} studies smooth inner-product and distance kernels when input dimension grows proportionally to sample size, showing that their spectra admit a linear approximation involving the sample Gram matrix. Subsequent work treats nonsmooth kernels, universality and spectral edges~\citep{cheng2013spectrum,do2013spectrum,fan2019spectral}. Further extensions include polynomial dimensional scalings that keep higher-degree feature components~\citep{misiakiewicz2022spectrum,dubova2026universality}. Our framework instead permits a fixed kernel and input distribution, and studies the shrinking bulk of a decaying feature spectrum under quadratic concentration.

    \item[Spectral regularization and risk characterization.] Studying the risk of spectral regularization and spectral filters permits a common framework for ridge regression, iterative regularization and early stopping. Classical analyses in this direction give risk bounds and establish learning rates and statistical optimality~\citep{logerfo2008spectral,caponnetto2007optimal}.

    Precise risk characterizations are available under proportional asymptotics, including for ridge regression and gradient flow~\citep{dobriban2018highdimensional,ali2019continuous}. Beyond this regime, \citet{cheng2024dimensionfree} establish nonasymptotic, dimension-free deterministic equivalents for ridge prediction risk at the bulk scale, with nonnegative regularization. \citet{liganshilinrisk} establish asymptotic risk formulas for analytic spectral algorithms under power-law decay, with regularization polynomially above the bulk scale. Their formulas apply the filter directly to the population covariance, whereas ours preserve corrections from a self-consistent spectral parameter at the bulk scale. Notably, \citet{velikanov2024generalization} derive closely related risk formulas for an infinite Gaussian feature model, relying on concentration and stability approximations whose justification in that setting is left open. Our results give deterministic equivalents for analytic filters at the bulk scale under dimension-free concentration assumptions, including negative ridge and gradient flow, with explicit nonasymptotic errors for both variance and squared bias.

    \item[Dimension-free learning theory.] Beyond the ridge regression work~\citep{cheng2024dimensionfree}, dimension-free deterministic equivalents have been used to analyze kernel ridge regression and generalized cross-validation~\citep{misiakiewicz2024nonasymptotic} and random feature regression~\citep{defilippis2024dimensionfree}. Related random matrix approaches have also been used to study scaling laws~\citep{atanasov2026scaling}. These works relate prediction risk to spectral decay and signal alignment, making it possible to investigate how to choose regularization, how many random features are required to achieve statistical accuracy, and how learning rates vary with sample size.

    Our extension to analytic filters provides quantitative inputs for studying iterative regularization and early stopping through the same spectral quantities. In particular, the deterministic risk expressions maintain the source-dependent bias-variance trade-off over the course of training. As a result, they offer a foundation for comparing regularization methods and studying stopping rules, with direct control over the approximation error. More generally, our spectral laws connect these learning problems to the geometry of the empirical spectrum.
\end{description}

\paragraph{Notation}
For a positive integer $n$, let $[n]:=\{1,\ldots,n\}$ and let $\C_+:=\{z\in\C:\Im z>0\}$. For sets $A,B\subseteq\C$, define $\dist(A,B):=\inf\{|a-b|:a\in A,\ b\in B\}$, with points identified with singleton sets. We write $a\wedge b:=\min\{a,b\}$ and $a\vee b:=\max\{a,b\}$. Let $\mc H$ be a real or complex separable Hilbert space. We denote by $\mc L(\mc H)$ the Banach space of bounded linear operators acting on $\mc H$. The inner product on $\mc H$, denoted $\langle \bx,\by\rangle_{\mc H}$, is assumed to be conjugate-linear in its first argument; with respect to an orthonormal basis, it can be expressed in coordinates as $\bx^*\by$. For $\bA\in\mc{L}(\mc{H})$, we denote its adjoint by $\bA^*$, its operator norm by $\norm{\bA}$, its Hilbert-Schmidt (Frobenius) norm by $\norm{\bA}_{F}$, its spectrum by $\Spec(\bA)$, and its trace by $\Tr\prn{\bA}$ whenever the latter is well defined. We write $|\bA|=(\bA^*\bA)^{1/2}$ for the operator modulus and $\Tr|\bA|$ for the trace norm when $\bA$ is trace-class. The induced Hilbert-space norm is denoted by $\norm{\bx}_{\mc H}^2:=\<\bx,\bx\>_{\mc H}$. For any bounded positive semidefinite operator $\bA$ on $\mc{H}$, define
\begin{align*}
    \langle\bx,\by\rangle_{\bA}
    &:=\langle\bx,\bA\by\rangle_{\mc{H}},
    &
    \norm{\bx}_{\bA}^2
    &:=\langle\bx,\bx\rangle_{\bA}
    =\|\bA^{1/2}\bx\|_{\mc{H}}^2.
\end{align*}
When $\bA$ is not injective, $\norm{\cdot}_{\bA}$ is understood as a
seminorm. We use $\bI$ for the identity operator and $\bzeros$ and $\bones$ for the all-zero and all-one vectors, respectively, with the underlying dimensions clear from context. For self-adjoint operators $\bA,\bB$, the relation $\bA\succeq\bB$ means that $\bA-\bB$ is positive semidefinite. We write $\lambda_{\min}^+(\bA)$ for the smallest positive eigenvalue of a finite-rank positive semidefinite operator $\bA$. For $\bx,\by\in\mc H$, the rank-one operator $\bx\by^*$ acts by $(\bx\by^*)\bv=\langle\by,\bv\rangle_{\mc H}\bx$.

For a real or complex random variable $Z$ and $1\leq p<\infty$, write $\|Z\|_{L^p}:=(\E|Z|^p)^{1/p}$, with $\|Z\|_{L^\infty}$ its essential supremum in modulus. For a function of a random element $Y$, we use $L^p(Y)$ when specifying that the norm is taken under the law of $Y$. For square-integrable complex random variables, we use the bilinear covariance $\Cov(Z,W):=\E[(Z-\E Z)(W-\E W)]$.

For nonnegative quantities $f,g$, we write $f=O(g)$ if $f\leq Cg$ and $f=\Omega(g)$ if $f\geq cg$ on the stated range, for positive constants $C,c$. We write $f\asymp g$ or $f=\Theta(g)$, when both bounds hold. For signed or complex quantities the comparisons apply to their absolute values. We write $f=o(g)$ if $|f/g|\to0$ along the specified limit, with $n\to\infty$ understood for sequences unless stated otherwise. We use subscripts on $O$, $\Omega$, $\Theta$ and $o$ to indicate the parameters on which the implicit constants may depend.

    \section{Model and main results}\label{sec:main-results}

\subsection{Model setup}\label{sec:model-setup}

\paragraph{Dimension-free sample covariances.}
Let $\mc{H}$ be a real separable Hilbert space of dimension $d\in\N\cup\{\infty\}$. We observe independent and identically distributed feature vectors $\bx_1,\ldots,\bx_n\in\mc{H}$ with mean zero and population covariance operator $\bSigma:=\E[\bx_i\bx_i^*]$. For $k\in[n]$, define the partial data operator $\bX_k:\mc{H}\to\R^k$ and its unnormalized partial sample covariance operator by
\begin{align}
    (\bX_k\bv)_j:=\langle\bx_j,\bv\rangle_{\mc H},\quad j\in[k], \qquad \text{and} \qquad \bT_k:=\bX_k^*\bX_k=\sum_{j=1}^k\bx_j\bx_j^*. \label{eq:partial-sample-covariance}
\end{align}
Setting $\bT_0:=0$ and writing $\bX:=\bX_n$, we denote by $r_n:=\rank(\bX)$. Let the positive eigenvalues of the normalized sample covariance operator $\bX^*\bX/n$ be $\what{\sigma}_{1,n}\geq\cdots\geq\what{\sigma}_{r_n,n}>0$. In the following discussion, we let $\mc{F}_k:=\sigma(\bx_1,\ldots,\bx_k)$, with $\mc{F}_0$ trivial, and we abbreviate conditional expectation by $\E_k[\cdot]:=\E[\cdot\mid\mc{F}_k]$.

\paragraph{$\mc{H}$-valued feature distribution.} We impose the two distributional conditions used in earlier dimension-free ridge-regression analysis~\citep{cheng2024dimensionfree}, and~\citep{misiakiewicz2024nonasymptotic,defilippis2024dimensionfree} impose closely related concentration conditions on kernel eigenfunctions or random feature maps. Either condition below yields the same dimension-free Hanson-Wright inequality, while condition~(ii) permits dependence among the coordinates.

\begin{assumption}[Feature distribution]\label{assmp:features}
    The covariance operator $\bSigma$ is positive semi-definite and trace class. Without loss of generality, we assume $\norm{\bSigma}=1$. Let
    \begin{align*}
        \bSigma=\sum_{j\geq1}\sigma_j\bu_j\bu_j^*,
        \qquad
        1=\sigma_1\geq\sigma_2\geq\dots>0,
        \qquad
        \sum_{j\geq1}\sigma_j<\infty,
    \end{align*}
    where the sum is finite when $d<\infty$ and zero-covariance directions are omitted. The features admit
    \begin{align}
        \bx_i=\bSigma^{1/2}\bz_i
        :=\sum_{j\geq1}\sqrt{\sigma_j}\,z_{ij}\bu_j,                 \label{eq:feature-expansion}
    \end{align}
    where the coordinate sequences are centered and isotropic\footnote{When $d=\infty$, this notation is slightly informal and is used for presentational cleanliness: $\bz_i$ need not be an $\mc{H}$-valued random element. The statements $\E[\bz_i]=\bzeros$ and $\E[\bz_i\bz_i^*]=\bI$ mean that $\E[z_{ij}]=0$ and $\E[z_{ij}z_{i\ell}]=\ind\{j=\ell\}$ for every $j,\ell\geq1$. The $\mc{H}$-valued feature $\bx_i$ is defined by the almost surely convergent series in Eq.~\eqref{eq:feature-expansion}.}, i.e., $\E[\bz_i] = 0$ and $\E[\bz_i \bz_i^*] = \bI$. There exists a constant $C_{\bx}>0$ such that one of the following conditions holds:
    \begin{enumerate}
        \item[(i)] (Independent sub-Gaussian coordinates).
        For every $i$, the coordinates $\{z_{ij}\}_{j\geq1}$ are independent, not necessarily identically distributed, and
        \begin{align}
            \norm{z_{ij}}_{\psi_2}
            :=\sup_{p\geq1}p^{-1/2}\prn{\E|z_{ij}|^p}^{1/p}
            \leq C_{\bx}.                                                 \label{eq:subgaussian-feature-assumption}
        \end{align}

        \item[(ii)] (Convex concentration).
        For every $m\leq d$ when $d<\infty$, and every $m\geq1$ when $d=\infty$, every convex $1$-Lipschitz function $\varphi:\R^m\to\R$, and every $t>0$,
        \begin{align}
            \P\brc{
            \left|
            \varphi(z_{i1},\ldots,z_{im})
            -\E\varphi(z_{i1},\ldots,z_{im})
            \right|\geq t
            }
            \leq 2\exp\prn{-\frac{t^2}{C_{\bx}^2}}.                    \label{eq:convex-concentration-feature-assumption}
        \end{align}
    \end{enumerate}
    When $d=\infty$, the series in Eq.~\eqref{eq:feature-expansion} converges in $\mc{H}$ almost surely and in $L^2$.
\end{assumption}

The convex-concentration condition in (ii) includes several important cases, such as independent uniformly bounded coordinates and distributions satisfying a dimension-free log-Sobolev inequality; in particular, it includes uniformly strongly log-concave distributions~\citep{boucheron2013concentration,bakry2014analysis}. Under either condition, the following Hanson-Wright inequality holds in $\mc{H}$. For finite-dimensional real vectors, the bound under condition~(i) follows from the original Hanson-Wright inequality~\citep[Theorem~1.1]{rudelson2013hanson}, whereas the bound under condition~(ii) follows from its convex-concentration extension~\citep[Theorem~2.3]{adamczak2015note}. The Hilbert-space extension by finite-rank approximation is also used in dimension-free ridge-regression analysis~\citep{cheng2024dimensionfree}.

\begin{proposition}[Hanson-Wright inequality]\label{prop:hanson-wright}
    Let $\bx$ be an independent copy of $\bx_i$ satisfying Assumption~\ref{assmp:features}. For any bounded operator $\bA$ on the complexification of $\mc{H}$ and any $t\geq0$, there exists a universal constant $c>0$ such that
    \begin{align}
        \P\prn{\left|\bx^*\bA\bx-\Tr\prn{\bSigma\bA}\right|\geq t}
        \leq 4\exp\brc{
        -c\prn{
        \frac{t^2}{C_{\bx}^4
        \norm{\bSigma^{1/2}\bA\bSigma^{1/2}}_{F}^2}
        \wedge
        \frac{t}{C_{\bx}^2
        \norm{\bSigma^{1/2}\bA\bSigma^{1/2}}}
        }
        }.                                                               \label{eq:hanson-wright}
    \end{align}
    For real self-adjoint $\bA$, the prefactor $4$ can be replaced by $2$. If $\bA$ is measurable with respect to a $\sigma$-field $\mc{G}$ and $\bx$ is independent of $\mc{G}$, then the same bound holds conditionally on $\mc{G}$.
\end{proposition}

\paragraph{Main examples: power law and critical law decay.} To showcase how our results apply, we examine two infinite-dimensional populations as primary illustrative examples. For $\alpha,\beta>0$, their eigenvalues are
\begin{align}
    \text{(Power law)}\quad \sigma_i=\frac{1}{i^{\alpha+1}} \qquad \text{and} \qquad \text{(Critical law)}\quad \sigma_i=\frac{1}{i(1+\log i)^{\beta+1}}. \label{eq:spectrum-decay-examples}
\end{align}
In the critical case, the logarithmic factor makes the spectrum summable at the boundary exponent $i^{-1}$. Of the two spectra in Eq.~\eqref{eq:spectrum-decay-examples}, only the critical one satisfies the tail effective-rank criterion for benign overfitting in minimum-norm linear regression~\cite{bartlett2020benign}. \citet{tsigler2023benign} extended this analysis to ridge regression. For these same spectra, \citet{cheng2024dimensionfree} obtained dimension-free deterministic equivalents for ridge bias and variance. Subsequent work extended this approach to kernel ridge regression~\cite{misiakiewicz2024nonasymptotic} and to random feature ridge regression, where it also yielded excess risk rates under power-law decay~\cite{defilippis2024dimensionfree}. These analyses are building blocks for understanding scaling laws for prediction error in kernel and random feature regression~\cite{bordelon2020spectrum,bahri2024explaining,atanasov2026scaling}. Here we study the empirical Gram spectrum itself and the risk of general spectral filters in these two regimes.

\subsection{Main result I: Spectral parameter and bulk laws}\label{sec:gram-bulk-law}
Our first results characterize the empirical spectral distributions (e.s.d.) of the sample covariance operators $\bT_n=\bX^*\bX$ at scales adapted to the population covariance $\bSigma$ in a dimension-free manner, capturing the bulk behavior. The $r_n$ positive eigenvalues of $\bT_n$ are also the positive eigenvalues of $\bX\bX^*$. We therefore study the e.s.d.\ $\what{\mu}_n$ of the unscaled Gram matrix $\bX\bX^*$ in~Eq.~\eqref{eq:intro-gram-esd}, and its associated empirical Stieltjes transform $\what{s}_n$. As previewed in the introduction, the main result in this section is essentially
\begin{align*}
    (s_n(\bSigma), \mu_n(\bSigma)) \quad \mathrel{\scalebox{2}[1]{$\approx$}} \quad (\what{s}_n(\bX), \what{\mu}_n(\bX)),
\end{align*}
where the deterministic pair $(s_n,\mu_n)$ depends only on $n$ and $\bSigma$. Theorem~\ref{thm:gram-stieltjes-approx} quantifies the approximation for Stieltjes transforms and interval masses.

Key to our analysis, we consider the \emph{effective spectral parameter} $\omega(z)$, defined as the unique solution in $\C_+$ to
\begin{align}
    z \omega(z)^{-1} = 1 + \frac{1}{n} \Tr \prn{\bSigma(\omega(z) \bI - \bSigma)^{-1}}.       \label{def:effective-spectral-parameter-self-consistent-equation}
\end{align}
The map $z\mapsto\omega(z)$ is analytic on $\C_+$.
Written equivalently, it is clear that the effective spectral parameter $\omega(z)$ is closely connected to the notion of \textit{effective regularization} $\lambda_\star$ defined in~\citep{cheng2024dimensionfree} for dimension free ridge regression. This relationship can be interpreted in terms of the deterministic variance and bias quantities that arise in proportional-asymptotic linear regression~\cite{bartlett2021deep, hastie2022surprises}. In particular, for $\lambda \in (0, +\infty)$, one has
\begin{equation}\label{eq:ridge-effective-regularization}
    n \prn{\frac{z}{\omega(z)}-1} = \Tr \prn{\bSigma(\omega(z) \bI - \bSigma)^{-1}} \qquad \text{and} \qquad n \prn{1 - \frac{\lambda}{\lambda_\star(\lambda)}} = \Tr \prn{\bSigma(\bSigma + \lambda_\star(\lambda) \bI)^{-1}},
\end{equation}
and we have the one-to-one correspondence $\omega(-\lambda) = -\lambda_\star(\lambda)$. The deeper connection originates from the spectral analysis of random matrix theory under the finite-dimensional proportional regime. On the normalized covariance scale, $-1/\omega(z/n)$ is the deterministic companion Stieltjes transform evaluated at $z/n$, as used in classical support analysis~\citep{silverstein1995analysis,bai1999exact}. We next use $\omega(z)$ to study the finite-$n$ spectral law. For real $t \in \R$, we additionally let $\omega(t)$ be the boundary value whenever $ \lim_{\eta \downarrow 0} \omega(t+ \iu \eta)$ is finite.

For the unscaled Gram matrix $\bX\bX^*$, its $n$ eigenvalues, counted with multiplicity, consist of $n\what{\sigma}_{j,n}$ for $j\in[r_n]$ and $n-r_n$ zeros. Thus its empirical law and Stieltjes transform are
\begin{align*}
    \what{\mu}_n&=\frac1n\sum_{j=1}^{r_n}\delta_{n\what{\sigma}_{j,n}}+\frac{n-r_n}{n}\delta_0,
    &\what{s}_n(z)&:=\int\frac{\what{\mu}_n(\de x)}{x-z}=\frac1n\Tr\prn{(\bX\bX^*-z\bI_n)^{-1}}.
\end{align*}
On the same scale, set $s_n(z):=-1/(n\omega(z/n))$. Substituting $\omega(z/n)=-1/(ns_n(z))$ into Eq.~\eqref{def:effective-spectral-parameter-self-consistent-equation} at $z/n$ gives the self-consistent equation
\begin{align}
    s_n(z)^{-1}=-z+\sum_{j\geq1}\frac{\sigma_j}{1+n\sigma_j s_n(z)}. \label{eq:population-scalar-equation}
\end{align}
By Lemma~\ref{lem:population-scalar-stability}, this equation has a unique solution $s_n(z)\in\C_+$ for $z\in\C_+$ and a unique positive solution for real $z<0$. Denote by $\mu_n$ the unique probability measure on $[0,\infty)$ with Stieltjes transform $s_n$, whose existence follows from the same lemma. We will establish that the empirical pair $(\what{s}_n, \what{\mu}_n)$ ``concentrates'' around this deterministic pair $(s_n, \mu_n)$. The general strategy of recovering the Gram matrix's spectral law from the effective regularization parameter $\lambda_\star$ was sketched in the Ph.D. dissertation of \citet[Chapter~4]{cheng2025thesis}, and in this paper, we establish rigorously in Theorem~\ref{thm:gram-stieltjes-approx} with a finite-sample $L^p$ bound for $\what{s}_n(z)-s_n(z)$ and transfers it to a bound on interval masses.

To approximate $\what{s}_n(z)$, we introduce necessary notations to characterize stability of the population spectrum. When $d<\infty$, extend the population eigenvalues by zero beyond $d$. For $p\ge1$, define the tail trace and tail effective rank parameters by
\begin{align}
    \tau_p(k):=\sum_{j>k}\sigma_j^p \qquad \text{and} \qquad r_p(k):=\sum_{j>k}\prn{\frac{\sigma_j}{\sigma_{k+1}}}^p.
    \label{eq:tail-parameters-trace}
\end{align}
We set $r_p(k)=0$ when $\sigma_{k+1}=0$. Fix $z=t+\iu\eta\in\C_+\cup\R_{<0}$ and write $z_-:=(t\wedge0)+\iu\eta$, so that $|z_-|=\dist(z,\R_{\geq0})$. There is a constant $\rho=\rho(C_{\bx})\in(0,1/2)$, depending only on the feature constant in Assumption~\ref{assmp:features}, for which the results below hold. For each integer $1\leq m\leq d$, define
\begin{align*}
    \Lambda_n(z,m)&:=\frac{|z|+n\sigma_{m+1}+\tau_1(m)\log(n)\log\prn{n r_1(m)}}{|z_-|}.
\end{align*}
When $\sigma_{m+1}=0$, we set $\Lambda_n(z,m)=|z|/|z_-|$. The following assumption quantifies the required sample size $n$ and the Stieltjes parameter scale $z$ in $\C_+ \cup \R_{<0}$.
\begin{assumption}[Spectral stability]\label{assmp:spectral-stability}
    For $\nu,\kappa>0$, define the spectral stability parameter
    \begin{align}
        \Lambda_n(z)&:=\min_{1\leq m\leq d \wedge \rho(\lfloor n/2\rfloor-2)}\Lambda_n(z,m). \label{eq:path-spectral-envelope}
    \end{align}
    Assume $n=\Omega_{C_{\bx},\nu,\kappa}(\Lambda_n(z)^{16})$ and $|z_-|\ge n^{-\nu}$.
\end{assumption}
We point out that the same $m$ attains the minimum for every $z$, since the $m$-dependent part of the numerator of $\Lambda_n(z,m)$ is independent of $z$. The selection of $\kappa$ is primarily a technical parameter, as it specifies the admissible range $2 \le p \le \kappa \log n$ on which the $L^p$ bounds (presented below) for the approximation error are valid uniformly. For an interval $E=[a,b]$ and $0<h<(b-a)/2$, write $E_h^+=[a-h,b+h]$ and $E_h^-=[a+h,b-h]$, we have the following deterministic equivalents for the spectrum. See its proof in Appendix~\ref{proof:gram-bulk-law}.

\begin{theorem}[Deterministic equivalents of the spectra]\label{thm:gram-stieltjes-approx}
    Under Assumption~\ref{assmp:features}, fix $\nu,\kappa>0$.
    \begin{subequations}
        \begin{enumerate}
            \item[(i)] Let $z\in\C_+\cup\R_{<0}$ satisfy Assumption~\ref{assmp:spectral-stability}. For $2\leq p\leq\kappa\log n$, the Stieltjes transform $\what{s}_n$ of the e.s.d.\ of $\bX \bX^*$ and its deterministic equivalent $s_n$ satisfy
            \begin{align}
                \|\what s_n(z)-s_n(z)\|_{L^p}=O_{C_{\bx},\nu,\kappa}\prn{\frac{p^2\Lambda_n(z)^{17}}{n|z|}}.\label{eq:gram-stieltjes-approx}
            \end{align}
            \item[(ii)] Fix $E=[a,b]$, $0<h<(b-a)/2$, and $D,\eta>0$. Suppose $\kappa\log n\geq2$ and Assumption~\ref{assmp:spectral-stability} holds at $t+\iu\eta$ for every $t\in[a-h,b+h]$. Then with probability $1-O(n^{-D})$,
            \begin{align}
                |\what\mu_n(E)-\mu_n(E)|\leq O_{C_{\bx},\nu,\kappa,D}\prn{\frac{(\log n)^2}{n}\int_{a-h}^{b+h}\frac{\Lambda_n(t+\iu\eta)^{17}}{|t+\iu\eta|}\,\de t+\frac{\eta}{h}}+\mu_n(E_h^+\setminus E_h^-).\label{eq:gram-interval-approx}
            \end{align}
        \end{enumerate}
    \end{subequations}
\end{theorem}
This deterministic equivalent theorem provides a dimension-free, nonasymptotic, and quantitative description for the e.s.d.\ of the unnormalized Gram matrix $\bX \bX^*$. The user can adapt to specific scenarios, depending on the population $\bSigma$ and the sample size $n$, choosing an appropriate scale $z=b_n+a_n \zeta$ and verifying Assumption~\ref{assmp:spectral-stability} under that scale. Here $\zeta$ ranges over a fixed compact subset of $\C_+$, $a_n$ is the appropriate ``bulk width'' and $b_n$ is the appropriate ``bulk shift''.
\begin{table}[htb!]
    \centering
    \begin{tabular}{@{}lccc@{}}
        \hline
        & Power law & Critical law & Proportional regime \\
        \hline
        $z$ & $n\sigma_n\zeta$ & $\tau_1(n)+n\sigma_n\zeta$ & $n\sigma_1\zeta$ \\
        $\Lambda_n(z)$ & $O((\log n)^2)$ & $O((\log n)^3)$ & $O((\log n)^2)$ \\
        \hline
    \end{tabular}
    \vspace{.5em}
    \caption{Spectral scale and stability bounds for applying Theorem~\ref{thm:gram-stieltjes-approx}.}\label{tab:spectral-stability-scales}
    \begin{figure}[H]
        \centering
        \includegraphics[width=.95\linewidth]{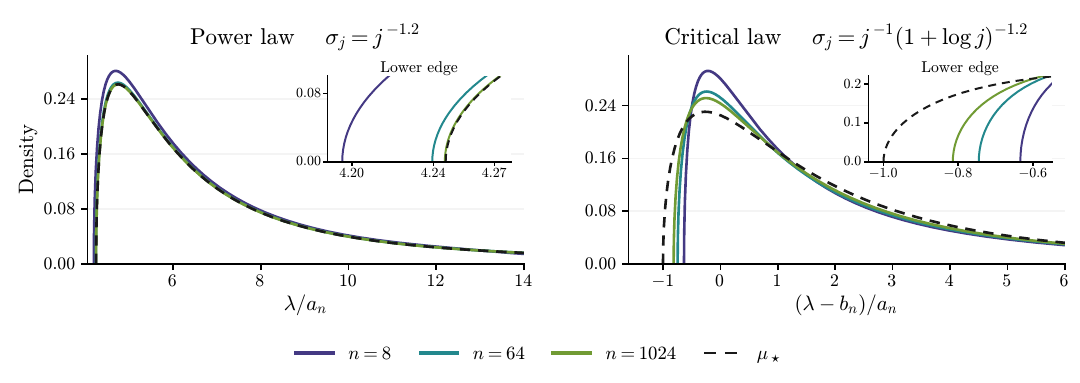}
        \caption{Densities of the deterministic laws $\mu_n$ (solid) and their limits $\mu_\star$ (dashed), in the bulk coordinates of Table~\ref{tab:spectral-stability-scales}. We show magnified views of each lower edge. $\mu_n$ and $\mu_\star$ are determined by Eqs.~\eqref{eq:population-scalar-equation} and~\eqref{eq:canonical-bulk-inverses}, respectively.}\label{fig:deterministic-bulk-laws}
    \end{figure}
\end{table}

We summarize the choices to apply our results to the main examples in the paper in Table~\ref{tab:spectral-stability-scales}. Under these rescalings, we plot the deterministic equivalent measure $\mu_n$ in $\zeta$ in Fig.~\ref{fig:deterministic-bulk-laws} against its asymptotic limit. We remark that the probability measures depend on $n$, and we further observe that, under the critical law, convergence to its limit is substantially slower in the critical case.

In the next Proposition~\ref{prop:canonical-bulk-laws}, we give quantitative nonasymptotic approximation bounds for the power and critical laws in Eq.~\eqref{eq:spectrum-decay-examples}. In both cases the bulk width is $a_n=n\sigma_n$. In the critical case the shift $b_n=\tau_1(n)$ is larger by a factor of order $\log n$.

\begin{proposition}[Power and critical bulk laws]\label{prop:canonical-bulk-laws}
    \begin{subequations}
        Under Assumption~\ref{assmp:features}, take either spectrum in Eq.~\eqref{eq:spectrum-decay-examples} and let $z=z_n(\zeta)$ be the corresponding entry of Table~\ref{tab:spectral-stability-scales}. Fix $\kappa>0$ and a compact $K\subset\C_+$. For all sufficiently large $n$ and $2\le p\le\kappa\log n$,
        \begin{align}
            \sup_{\zeta\in K}\|a_n(\what s_n(z)-s_n(z))\|_{L^p}
            =\frac{p^2}{n}\begin{cases}O_{C_{\bx}, \kappa, K, \alpha}((\log n)^{34}),&\text{power law;} \\ O_{C_{\bx}, \kappa, K, \beta}((\log n)^{50}),&\text{critical law.}\end{cases}\label{eq:canonical-bulk-rate}
        \end{align}
        After the corresponding centering and scaling, $\mu_n$ and $\what\mu_n$ converge weakly to a fixed probability measure $\mu_{\star}$, almost surely for $\what\mu_n$. Its Stieltjes transform $s(\zeta)$ is the locally uniform limit of $a_ns_n(z)$ on $\C_+$ and the unique solution in $\C_+$ of
        \begin{align}
            s(\zeta)^{-1}=-\zeta+\begin{cases}
            \dfrac{\pi\,s(\zeta)^{-\alpha/(1+\alpha)}}{(1+\alpha)\sin(\pi/(1+\alpha))},&\text{power law;} \vspace{.6em} \\
            -\log s(\zeta),&\text{critical law.}
        \end{cases}\label{eq:canonical-bulk-inverses}
    \end{align}
    Let $\kappa_{\star}:=\inf\supp(\mu_{\star})$. In each case, $\mu_{\star}$ has support $[\kappa_{\star},\infty)$ and a bounded density $\rho_{\star}$ satisfying
    \begin{align}
        \kappa_{\star}=
        \begin{cases}
            \dfrac{1}{\alpha} \cdot \left(\dfrac{\pi\alpha}{(1+\alpha)^2\sin(\pi/(1+\alpha))}\right)^{1+\alpha},&\text{power law;} \vspace{.6em} \\
            -1,&\text{critical law.}
        \end{cases}
        \qquad
        \|\rho_{\star}\|_\infty\le
        \begin{cases} \dfrac{1}{\pi\alpha\kappa_{\star}} ,&\text{power law;} \vspace{.6em} \\ \dfrac{1}{\pi} ,&\text{critical law.}\end{cases}\label{eq:canonical-power-edge-density}
    \end{align}
\end{subequations}
\end{proposition}
The proof is in Appendix~\ref{proof:canonical-bulk-laws}, which also contains the exact parameterizations of both limiting densities. We also point out that for the critical law, its limiting rescaled law $\mu_\star$ is independent of the exponent $\beta$, unlike the power law decay, whose $\mu_\star$ and $\kappa_\star$ rely on $\alpha$.

Intriguingly, the limiting laws identified by Eq.~\eqref{eq:canonical-bulk-inverses} are free stable laws characterized in~\citep{bercovici1999stable} in free probability. Power law decay gives a scaled positive free stable law of index $1/(1+\alpha)$, while critical decay gives the reflection of the exceptional free $1$-stable law studied in~\citep{hasebe2020properties}.

\begin{remark}[Gram matrix geometry and benign overfitting]
    In the critical case, $a_n/b_n\sim\beta/\log n\to0$, so Proposition~\ref{prop:canonical-bulk-laws} implies that the empirical spectral distribution of $b_n^{-1}\bX\bX^*$ converges almost surely to a point mass $\delta_1$, as in the classical fixed-dimensional asymptotics. For power-law decay, the bulk of $a_n^{-1}\bX\bX^*$ retains a nondegenerate spread, as in proportional asymptotics. This gives a spectral interpretation of the effective-rank criterion for benign overfitting in~\citet{bartlett2020benign,tsigler2023benign}, whose bounds show that, under their hypotheses, the prediction variance of ridgeless regression vanishes in the critical case and persists under power-law decay. We give a quantitative interpretation here by Gram matrix spectral geometry.
\end{remark}

To wrap up the discussion on our first set of main results, we point the readers to two additional extensions of our bulk law results in Appendix~\ref{app:additional-bulk-laws}. The first is Proposition~\ref{prop:mesoscopic-bulk}, which extends the spectral height from $O(1)$ in bulk scale to $\Im\zeta = n^{-\delta}$ for fixed $0<\delta<1/17$, providing its mesoscopic version (albeit, suboptimal as $\delta < 1$). The second is Proposition~\ref{prop:proportional-gram-law}, where we recover the classical generalized Marchenko-Pastur law, showing the applicability and adaptability of our master Theorem~\ref{thm:gram-stieltjes-approx}.

We next turn to the second set of our results, which provides the deterministic equivalents of prediction risk of analytic spectral filters, with explicit formulas of effective variance and bias.

\subsection{Main result II: Effective risk of spectral filters}\label{sec:spectral-filter-risk}

Our second main result characterizes the conditional excess prediction error of spectral filters in linear regression. For the design $\bX$ of Section~\ref{sec:model-setup} and a fixed signal $\bbeta\in\mc H$, we observe
\begin{align*}
    \by=\bX\bbeta+\bvarepsilon \qquad \text{where} \qquad \E[\bvarepsilon\mid\bX]=0, \quad \E[\bvarepsilon\bvarepsilon^*\mid\bX]=\tau^2\bI_n.
\end{align*}
This includes the regime of centered independent noise $\varepsilon_i$ for $1 \leq i \leq n$, with common variance $\tau^2$ and noise vector independent of $\bX$. We make the signal assumption $\bbeta = \bSigma^{1/2} \btheta$ with $\norm{\btheta} < \infty$, as is in~\citep{cheng2024dimensionfree}. For a given filter $\upphi$ and noise level $\tau$, we seek the deterministic approximation
\begin{align*}
    (\effectiveVariance[\upphi;\bSigma],\effectiveBias[\upphi;\bSigma]) \quad \mathrel{\scalebox{2}[1]{$\approx$}} \quad (\sampleVariance[\upphi],\sampleBias[\upphi])
\end{align*}
for the two terms in Eq.~\eqref{eq:intro-filter-risk}. Recall that with $\what{\bSigma}=\bX^*\bX/n$, the variance and squared bias are
\begin{align}
    \sampleVariance[\upphi]=\frac{\tau^2}{n}\Tr\prn{\bSigma\what{\bSigma}\upphi(\what{\bSigma})^2}
    \qquad\text{and}\qquad
    \sampleBias[\upphi]=\norm{\bSigma^{1/2}\uppsi(\what{\bSigma})\bbeta}^2.\label{eq:conditional-filter-risks}
\end{align}
In literature, instance-wise risk bounds have been used to compare gradient descent with ridge regression~\cite{wu2026riskcomparisons} and to establish the optimality, up to constant factors, of optimally tuned PCR among monotone spectral filters~\cite{kim2026monotone}. Our deterministic equivalents could provide sharper estimates of risk differences and be useful in extending the comparisons. In this paper, we focus on analytic filters, for which a characterization via Cauchy’s formula is available, directly connecting prediction risk to their resolvents on $\C$.

To incorporate spectral filters $\upphi$ that may be undefined or unbounded at or near zero, such as ridgeless and negative ridge estimators, we choose $\eta_n \geq 0$ such that $[n\eta_n,\infty)=:I_n$ contains $\supp\mu_n$. On a hard edge event, that is $\what{\sigma}_n\ge\eta_n$, the spectrum of $\bX^*\bX$ lies in $\{0\}\cup I_n$. Let $\uppsi(x)=1-x\upphi(x)$ for the filter residual.
\begin{assumption}[Analytic spectral filters]\label{assmp:analytic-spectral-filter}
    For $\supp\mu_n\subseteq I_n=[n\eta_n,\infty)$, let $\Omega$ be a conjugation-symmetric open neighborhood of $\{0\}\cup I_n$, with piecewise smooth boundary $\Gamma=\partial\Omega$. The set $\Omega$ may be disconnected. The restriction of the filter $z\mapsto\upphi(z/n)$ to $I_n$ extends holomorphically to a neighborhood of $\wb{\Omega}$, with $\upphi(\wb{z}/n)=\wb{\upphi(z/n)}$. For this extension, let $\uppsi(z/n)=1-(z/n)\upphi(z/n)$ and assume that
    \begin{align*}
        \lim_{|z|\to\infty}\uppsi(z/n)=0 \quad \text{for $z \in \wb{\Omega}$},
        \qquad\text{and}\qquad
        \int_\Gamma |\uppsi(z/n)| \cdot |\de z|<\infty \quad \text{when $\Gamma$ is unbounded}.
    \end{align*}
\end{assumption}

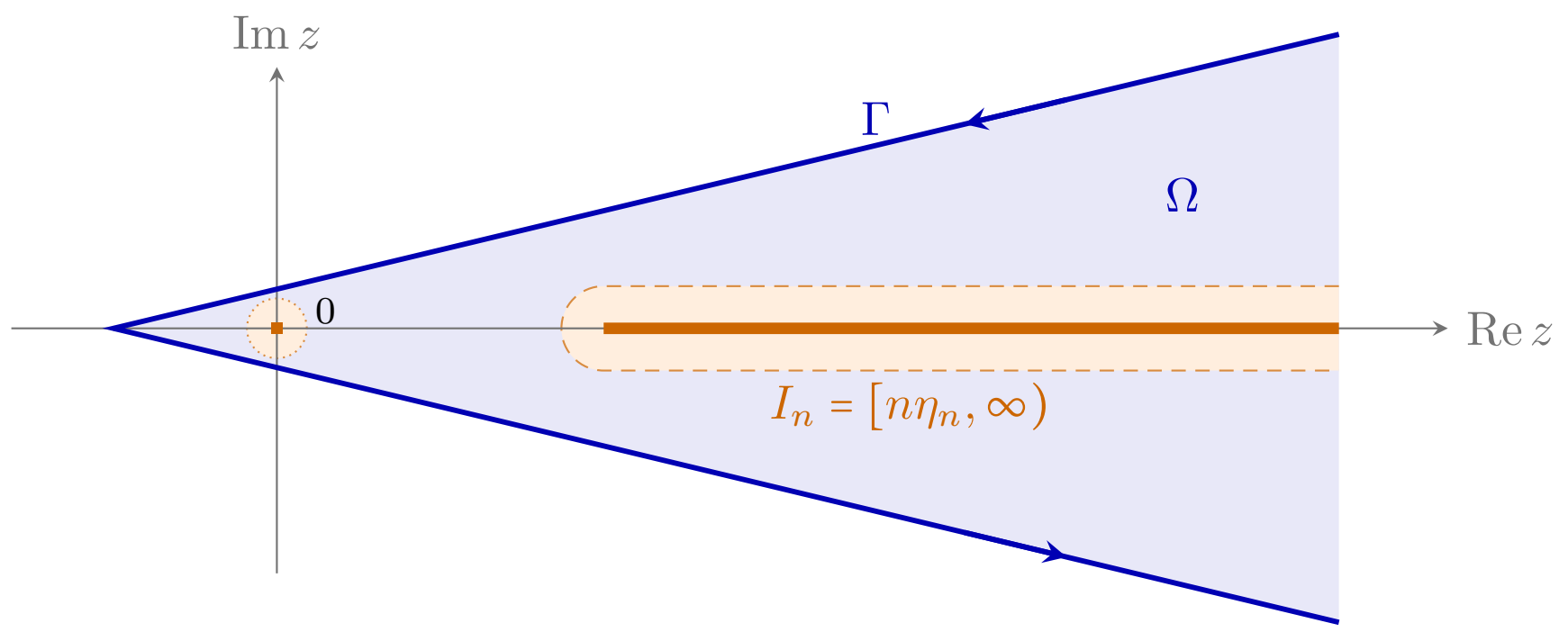
\begin{figure}[htbp]
    \centering
    \begin{tikzpicture}[x=1cm,y=1cm,font=\small,>=stealth]
        \colorlet{filterblue}{blue!70!black}
        \colorlet{filterorange}{orange!80!black}
        \fill[filterblue!9] (-1.20,0)--(7.8,2.16)--(7.8,-2.16)--cycle;
        \fill[orange!13] (7.8,.31)--(2.4,.31) arc[start angle=90,end angle=270,radius=.31]--(7.8,-.31)--cycle;
        \draw[filterorange!75,densely dashed] (7.8,.31)--(2.4,.31) arc[start angle=90,end angle=270,radius=.31]--(7.8,-.31);
        \filldraw[fill=orange!13,draw=filterorange!75,densely dotted] (0,0) circle[radius=.22];
        \draw[->,black!55] (-1.95,0)--(8.6,0) node[right] {$\Re z$};
        \draw[->,black!55] (0,-1.8)--(0,1.92) node[above] {$\Im z$};
        \node[right,inner sep=1pt,font=\scriptsize] at (.24,.13) {$0$};
        \draw[filterorange,line width=2.4pt] (2.4,0)--(7.8,0);
        \fill[filterorange] (-1.2pt,-1.2pt) rectangle (1.2pt,1.2pt);
        \node[below,text=filterorange,inner sep=2pt] at (4.65,-.32) {$I_n=[n\eta_n,\infty)$};
        \draw[filterblue,line width=1.1pt] (7.8,2.16)--(-1.20,0)--(7.8,-2.16);
        \draw[filterblue,line width=1.1pt,->] (5.80,1.68)--(5.05,1.50);
        \draw[filterblue,line width=1.1pt,->] (5.05,-1.50)--(5.80,-1.68);
        \node[above,text=filterblue,inner sep=2pt] at (4.40,1.344) {$\Gamma$};
        \node[text=filterblue] at (6.65,.98) {$\Omega$};
    \end{tikzpicture}
    \caption{Example of $\Omega$ and the contour $\Gamma$ in Assumption~\ref{assmp:analytic-spectral-filter}. The illustrated two-ray contour $\Gamma$ is chosen for evaluating the gradient flow estimator's effective terms $\effectiveVariance$ and $\effectiveBias$ in Section~\ref{sec:applications}.}
    \label{fig:analytic-filter-support}
\end{figure}

We remark that under Assumption~\ref{assmp:analytic-spectral-filter}, we need only specify $x\mapsto\upphi(x/n)$ on $I_n$, which means $\upphi$ need only be defined (analytically) on $[\eta_n, \infty)$. For any extension as in the assumption, the residual $\uppsi(z/n)=1-(z/n)\upphi(z/n)$ satisfies $\uppsi(0)=1$, since $0 \in \Omega$. Setting $\alpha_{\upphi}(\de z)=\uppsi(z/n) \de z/(2\pi \iu z)$ on the positively oriented boundary, we can apply Cauchy's formula at $z=nx\in I_n$ to $z\mapsto\uppsi(z/n)$ and $z\mapsto\upphi(z/n)$, which yields the following representations, respectively:
\begin{align}
    \uppsi(x)=\int_\Gamma \frac{z}{z-nx}\alpha_{\upphi}(\de z)
    \qquad\text{and}\qquad
    \upphi(x)=\int_\Gamma \frac{1}{x-z/n}\alpha_{\upphi}(\de z).\label{eq:analytic-filter-representation}
\end{align}
Since $\uppsi(0)=1$, it also holds by Cauchy's formula at $z=0$ that $\int_\Gamma \alpha_{\upphi}(\de z) = 1$. Writing $|\alpha_{\upphi}|(\de z)= |\uppsi(z/n)| /(2\pi|z|) \cdot |\de z|$, we let $\tvnorms{\alpha_{\upphi}}:=\int_\Gamma |\alpha_{\upphi}|(\de z)$. Applying the Cauchy representations in Eq.~\eqref{eq:analytic-filter-representation} to the two conditional risk formulas in Eq.~\eqref{eq:conditional-filter-risks}, we have the following contour-integral representations. Its proof is in Appendix~\ref{proof:conditional-risk-resolvents}.
\begin{lemma}\label{lem:conditional-risk-resolvents}
    Under Assumption~\ref{assmp:analytic-spectral-filter}, on $\{\what{\sigma}_n\ge\eta_n\}$, the conditional variance and squared bias satisfy
    \begin{align*}
        \sampleVariance[\upphi]&=\tau^2\iint_{\Gamma^2}\frac{z\Tr\prn{\bSigma(\bT_n-z\bI)^{-1}}-z'\Tr\prn{\bSigma(\bT_n-z'\bI)^{-1}}}{z-z'}\alpha_{\upphi}(\de z)\alpha_{\upphi}(\de z'),\\
        \sampleBias[\upphi]&=\iint_{\Gamma^2} zz'\bbeta^*(\bT_n-z\bI)^{-1}\bSigma(\bT_n-z'\bI)^{-1}\bbeta \alpha_{\upphi}(\de z)\alpha_{\upphi}(\de z').
    \end{align*}
    The variance integrand at $z'=z$ is understood to be the derivative of $z\Tr\prn{\bSigma(\bT_n-z\bI)^{-1}}$ in $z$.
\end{lemma}
To see how this contour-integral can be useful, we first rescale the effective spectral parameter in Eq.~\eqref{def:effective-spectral-parameter-self-consistent-equation} on the scale of the contour variables, i.e., $\wt{\omega}(z):=n\omega(z/n)=-1/s_n(z)$ for $z\notin I_n$. Following the notation in~\cite[Section~2.1]{tibshirani2022divided}, we write its first divided difference at $z,z'\notin I_n$ as
\begin{align*}
    \wt{\omega}[z,z']:=\frac{\wt{\omega}(z)-\wt{\omega}(z')}{z-z'}\quad(z\ne z'),\qquad\text{and}\qquad \wt{\omega}[z,z]:=\wt{\omega}'(z)=\omega'(z/n).
\end{align*}
We now invoke the resolvent and mixed-resolvent deterministic equivalents in Theorem~\ref{thm:terminal-weighted-traces} and Theorem~\ref{thm:terminal-mixed-resolvents}, whose details are in the later Section~\ref{sec:interpolation-path}. Informally, we have
\begin{align*}
    \Tr(\bSigma(\bT_n - z\bI)^{-1}) \quad &\stackrel{\text{Theorem~\ref{thm:terminal-weighted-traces}}}{\mathrel{\scalebox{2}[1]{$\approx$}}} \quad \frac{\wt{\omega}(z)}{z} - 1, \\
    \bbeta^*(\bT_n-z\bI)^{-1}\bSigma(\bT_n-z'\bI)^{-1}\bbeta \quad &\stackrel{\text{Theorem~\ref{thm:terminal-mixed-resolvents}}}{\mathrel{\scalebox{2}[1]{$\approx$}}} \quad \frac{\wt{\omega}[z,z']}{zz'}\bbeta^*\prn{\bI-\frac{n\bSigma}{\wt{\omega}(z)}}^{-1}\bSigma\prn{\bI-\frac{n\bSigma}{\wt{\omega}(z')}}^{-1}\bbeta.
\end{align*}
We can now immediately see that the effective risk formula $\effectiveRisk[\upphi; \bSigma]=\effectiveVariance[\upphi; \bSigma]+\effectiveBias[\upphi; \bSigma] = \effectiveVariance[\upphi]+\effectiveBias[\upphi]$ follows by substitution into Lemma~\ref{lem:conditional-risk-resolvents}.
\begin{subequations}\label{eq:analytic-filter-risk}
    \begin{align}
        \effectiveVariance[\upphi]&=\tau^2\iint_{\Gamma^2} (\wt{\omega}[z,z']-1)\,\alpha_{\upphi}(\de z)\alpha_{\upphi}(\de z'),\label{eq:analytic-filter-risk-variance}\\
        \effectiveBias[\upphi]&=\iint_{\Gamma^2}\wt{\omega}[z,z']\bbeta^*\prn{\bI-\frac{n\bSigma}{\wt{\omega}(z)}}^{-1}\bSigma\prn{\bI-\frac{n\bSigma}{\wt{\omega}(z')}}^{-1}\bbeta\,\alpha_{\upphi}(\de z)\alpha_{\upphi}(\de z').\label{eq:analytic-filter-risk-bias}
    \end{align}
\end{subequations}
The rigorous replacement requires more caution. As $\Gamma$ passes through the gap $(0,n\eta_n)$ on the real axis, Theorems~\ref{thm:terminal-weighted-traces} and~\ref{thm:terminal-mixed-resolvents} do not apply directly there. Using a contour argument allows us to bound the log-norm of the error at $z$ near the gap by its average over a small disk, by the submean inequality. The disk is large enough to apply Theorems~\ref{thm:terminal-weighted-traces} and~\ref{thm:terminal-mixed-resolvents} to control the function away from a narrow horizontal strip about the real axis that occupies only a small fraction of the disk. Using this strategy, the following theorem quantifies how closely these deterministic quantities approximate the conditional risks. See its proof in Appendix~\ref{proof:analytic-filter-risk}.
\begin{theorem}[Effective risk of analytic spectral filters]\label{thm:analytic-filter-risk}
    Under Assumptions~\ref{assmp:features} and~\ref{assmp:analytic-spectral-filter}, let $\kappa>0$ and $\alpha_{\upphi}$ be the measure in Eq.~\eqref{eq:analytic-filter-representation}. Let $\Delta_\Gamma:=\Lambda_n(-\dist(\Gamma,I_n))\vee\sup_{z\in\Gamma}(|z|/\dist(z,I_n))$ denote the spectral stability of the contour $\Gamma$. Assume $\bbeta=\bSigma^{1/2}\btheta$ with $\norm{\btheta} < \infty$, and for some constant $C_\Gamma > 0$,
    \begin{align*}
        -\log\!\big(1\wedge\dist(\Gamma,I_n)\big) \leq C_\Gamma \log n
        \qquad\text{and}\qquad
        \Delta_\Gamma \leq C_\Gamma n^{1/16} / \log n.
    \end{align*}
    Then for all large $n$ and $2\le p\le\kappa\log n$, the effective variance and squared bias in Eq.~\eqref{eq:analytic-filter-risk} satisfy
    \begin{align*}
        \|\ind\{\what{\sigma}_n\ge\eta_n\}(\sampleVariance[\upphi]-\effectiveVariance[\upphi])\|_{L^p}&=O_{C_{\bx},\kappa, C_\Gamma}\left(\frac{\tau^2p^2\Delta_\Gamma^{23}\tvnorms{\alpha_{\upphi}} }{n}\int_{\Gamma}\frac{|z|\vee n\eta_n}{\dist(z,I_n)}\cdot |\alpha_{\upphi}|(\de z)\right),\\
        \|\ind\{\what{\sigma}_n\ge\eta_n\}(\sampleBias[\upphi]-\effectiveBias[\upphi])\|_{L^p}&=O_{C_{\bx},\kappa, C_\Gamma}\left(\frac{p^{3/2}\Delta_\Gamma^{26}\|\btheta\|^2}{n^{5/2}}\prn{\int_{\Gamma} (|z| \vee n\eta_n) \cdot |\alpha_{\upphi}|(\de z)}^2\right).
    \end{align*}
\end{theorem}

We now establish that the variance and squared bias of spectral filters concentrate around their deterministic equivalents in Eq.~\eqref{eq:analytic-filter-risk}. However, it may occur to the reader that these quantities appear to depend on the contour $\Gamma$. On $\{\what{\sigma}_n\ge\eta_n\}$, all possible singularities of the integrands for $\sampleVariance[\upphi]$, $\sampleBias[\upphi]$ and their equivalents $\effectiveVariance[\upphi]$, $\effectiveBias[\upphi]$ lie in $\{0\}\cup I_n$, which is enclosed by $\Gamma$. By Cauchy's formula, we can deform the contour within the domain of analyticity of the filter, and the integrals depend only on the behavior of the integrands near the enclosed singularities. On the population side, for each $\sigma_j>0$ (pole at spectral parameter $\omega$), define
\begin{equation}\label{eq:filter-directional-measures}
    s_{\omega,j}(z):=\frac{\omega(z)}{z(\sigma_j-\omega(z))}=-\frac{1}{z(1+n\sigma_j s_n(nz))},\qquad z\in\C_+.
\end{equation}
As verified in Appendix~\ref{proof:filter-real-risk} (along with the proof of the following Proposition) using Lemma~\ref{lem:herglotz-nevanlinna}, $s_{\omega,j}$ is the Stieltjes transform of a probability measure $\mu_{\omega,j}$ on $\{0\}\cup[\eta_n,\infty)$, so $s_{\omega,j}(z)=\int (x-z)^{-1}\mu_{\omega,j}(\de x)$. Define the finite positive measure $\mu_\omega$ by the weighted sum $\mu_\omega(\de x):=(x/n)\sum_{j}\sigma_j\mu_{\omega,j}(\de x)$.

\begin{figure}[htbp!]
    \centering
    \includegraphics[width=\linewidth]{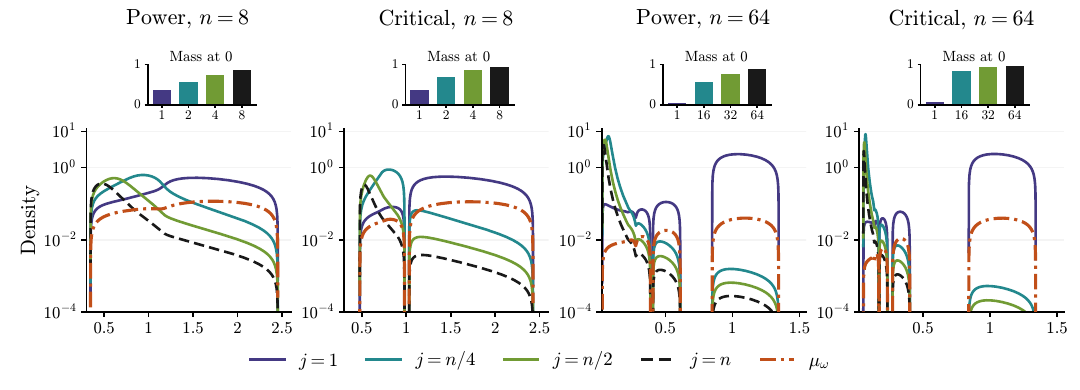}
    \caption{Densities of $\mu_{\omega,j}$ and their weighted measure $\mu_\omega$ for $\sigma_i=i^{-1.2}$ (power) and $\sigma_i=i^{-1}(1+\log i)^{-1.2}$ (critical) for $n \in \{8, 64\}$ and $j \in \{1, n/4, n/2, n\}$. The bar plots show $\mu_{\omega,j}(\{0\})$.}
    \label{fig:directional-measures}
\end{figure}

In Figure~\ref{fig:directional-measures}, we make plots of $\mu_{\omega,j}$ and $\mu_\omega$  for the power law and critical law, showing how the spectral mass shifts with the population eigenvalue index~$j$. Finally, we establish the following representation of $\effectiveRisk$.

\begin{proposition}\label{prop:filter-real-risk}
    Under the setup of Theorem~\ref{thm:analytic-filter-risk}, the effective risks in Eq.~\eqref{eq:analytic-filter-risk} satisfy
    \begin{subequations}\label{eq:filter-real-risk}
        \begin{align}
            \effectiveVariance[\upphi]&=\tau^2\int \upphi(x)^2\,\mu_\omega(\de x)=\frac{\tau^2}{n}\sum_j\sigma_j\int x\upphi(x)^2\,\mu_{\omega,j}(\de x),\label{eq:filter-real-variance}\\
            \effectiveBias[\upphi]&=\sum_j\sigma_j|\langle\bu_j,\bbeta\rangle|^2\brc{\left(\int \uppsi(x)\,\mu_{\omega,j}(\de x)\right)^2+
            \int\left(\int\frac{\uppsi(x)-\uppsi(x')}{x'-x}\,\mu_{\omega,j}(\de x)\right)^2\mu_\omega(\de x')}.\label{eq:filter-real-bias}
        \end{align}
    \end{subequations}
    At $x=x'>0$, the divided difference in Eq.~\eqref{eq:filter-real-bias} equals $-\uppsi'(x')$.
\end{proposition}

We note that Eqs.~\eqref{eq:analytic-filter-risk} and~\eqref{eq:filter-real-risk} provide equivalent deterministic risk expressions, each useful from a different viewpoint. The first contour representation captures the filter's analytic regularity, and its holomorphic extension determines the explicit approximation error through its behavior on $\Gamma$ in Theorem~\ref{thm:analytic-filter-risk}. On the other hand, Proposition~\ref{prop:filter-real-risk} provides real-line evaluation formulas, without separate approximation guarantees while avoiding working on the complex plane.

In practice, both formulas can be useful for evaluating the deterministic risk; while they should theoretically yield the same result, one may be simpler to compute than the other. We will see this in the next section.

\FloatBarrier

    \section{Applications} \label{sec:applications}
We will apply our main result II in Section~\ref{sec:spectral-filter-risk} to ridge-type filters and gradient flow under the power and critical laws in Eq.~\eqref{eq:spectrum-decay-examples}. Recall that in those results, we work under the hard edge event $\{\what{\sigma}_n \geq \eta_n\}$. This is true unconditionally setting $\eta_n=0$, by p.s.d.\ of $\bX \bX^*$. Under the working Assumption~\ref{assmp:features} of quadratic concentration, we can indeed take without loss of generality that
\begin{align*}
    n\eta_n = \Theta(\tau_1(n)) \qquad \text{such that} \qquad \P\prn{\what{\sigma}_n\geq\eta_n}\geq1-2e^{-n},
\end{align*}
for all $n=\Omega(1)$ (with dependence on the decay exponent and $C_{\bx}$ suppressed). We defer discussion on the hard edge event to Lemma~\ref{lem:smallest-eigenvalue}. At the same time, for our power law and critical law example, this is at the correct scale, since $\inf\supp\mu_n=\Theta(\tau_1(n))$ (cf.~Appendix~\ref{proof:application-spectral-bounds}).

It is possible to strengthen and take $n \eta_n = (1-c) \inf\supp\mu_n$ for arbitrarily small $c>0$, or even vanishing $c = c_n$ and still have $\what{\sigma}_n \geq \eta_n$ with overwhelmingly high probability, provided with stronger concentration of the features, e.g., finer Gaussian edge bounds (cf.~\citep[Section~3.3]{brailovskaya2024universality}). For both of our applications, having $n\eta_n = \Theta(\inf\supp\mu_n)$ is enough, so we leave this for future work.

\subsection{Ridge-type filters}\label{sec:negative-ridge}
We begin with a unifying framework for ridge-type filters outlined in the introduction. Indeed, with
\begin{align*}
    \uppsi_{k,\lambda}(x)=\left(\frac{\lambda}{x+\lambda}\right)^k\qquad\text{and}\qquad
    \upphi_{k,\lambda}(x)=\frac{1-\uppsi_{k,\lambda}(x)}{x},
\end{align*}
we can recover regular ridge ($k=1, \lambda > 0$), ridgeless ($\lambda = 0$), negative ridge ($k=1, \lambda < 0$) and their iterated version ($k>1$). We can also view these filters as estimators from an iterative algorithm, fits a ridge correction to the current residual at each step, initialized at $\what{\bbeta}_{0,\lambda}=0$ and defined recursively for $k=1,2,\dots$ by
\begin{align*}
    \what{\bbeta}_{k, \lambda} =\what{\bbeta}_{k-1, \lambda }+
    \argmin_{\bb\in \mathrm{ran}(\bX^*)}
    \left\{\frac{1}{2n}\big\|\by-\bX\what{\bbeta}_{k-1, \lambda}-\bX\bb\big\|^2+\frac{\lambda}{2}\|\bb\|^2\right\}.
\end{align*}
Here $\mathrm{ran}(\bX^*)$ denotes the range of $\bX^*$---it is important to have this constraint to allow for negative curvature in the regularization term, as the ambient space is in $\mc{H}$, which can have arbitrarily many descent directions when $\lambda < 0$. For $\lambda>0$ this is iterated Tikhonov regularization~\cite{logerfo2008spectral}.

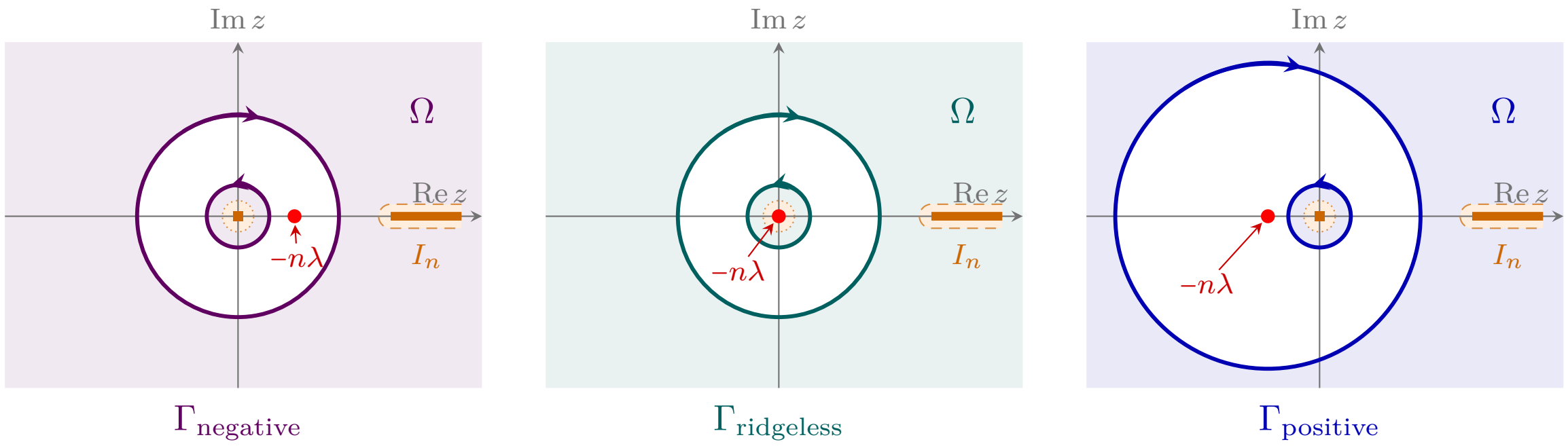
\begin{figure}[htbp]
    \centering
    \begin{tikzpicture}[scale=1.17,x=1cm,y=1cm,font=\small,>=stealth]
        \colorlet{ridgeorange}{orange!80!black}
        \foreach \offset/\shiftcenter/\pole/\name/\ink/\labelx/\labely in {0/0/.47/negative/violet!75!black/.49/-.34,4.50/0/0/ridgeless/teal!75!black/-.34/-.45,9/-.43/-.43/positive/blue!70!black/-.94/-.56}{
        \begin{scope}[xshift=\offset cm]
            \colorlet{ridgecontour}{\ink}
            \pgfmathsetmacro{\outerradius}{.84-\shiftcenter}
            \path[fill=ridgecontour,fill opacity=.09,even odd rule]
            (-1.94,-1.43) rectangle (2.03,1.45) (\shiftcenter,0) circle[radius=\outerradius];
            \fill[ridgecontour,fill opacity=.09] (0,0) circle[radius=.26];
            \fill[orange!13] (1.86,.10)--(1.27,.10) arc[start angle=90,end angle=270,radius=.10]--(1.86,-.10)--cycle;
            \filldraw[fill=orange!13,draw=ridgeorange!75,densely dotted] (0,0) circle[radius=.13];
            \draw[->,black!55] (-1.94,0)--(2.03,0) node[above left,font=\scriptsize] {$\Re z$};
            \draw[->,black!55] (0,-1.43)--(0,1.45) node[above,font=\scriptsize] {$\Im z$};
            \draw[ridgeorange!75,densely dashed] (1.86,.10)--(1.27,.10) arc[start angle=90,end angle=270,radius=.10]--(1.86,-.10);
            \draw[ridgeorange,line width=2.4pt] (1.27,0)--(1.86,0);
            \fill[ridgeorange] (-1.2pt,-1.2pt) rectangle (1.2pt,1.2pt);
            \node[below,text=ridgeorange,font=\scriptsize] at (1.57,-.12) {$I_n$};
            \draw[ridgecontour,line width=1.1pt] (0,0) circle[radius=.26];
            \draw[ridgecontour,line width=1.1pt,->] (32:.26) arc[start angle=32,end angle=105,radius=.26];
            \draw[ridgecontour,line width=1.1pt] (\shiftcenter,0) circle[radius=\outerradius];
            \draw[ridgecontour,line width=1.1pt,->] (\shiftcenter,0) ++(123:\outerradius) arc[start angle=123,end angle=77,radius=\outerradius];
            \node[text=ridgecontour] at (1.53,.88) {$\Omega$};
            \node[font=\scriptsize,text=red!80!black,inner sep=.8pt] (polelabel) at (\labelx,\labely) {$-n\lambda$};
            \draw[->,red!80!black,thin,shorten >=2.4pt] (polelabel)--(\pole,0);
            \fill[red] (\pole,0) circle[radius=1.65pt];
            \node[text=ridgecontour] at (0,-1.72) {$\Gamma_{\mathrm{\name}}$};
        \end{scope}
        }
    \end{tikzpicture}
    \caption{Two-circle contours $\Gamma$ for \textcolor{violet!75!black}{negative ridge}, \textcolor{teal!75!black}{ridgeless}, and \textcolor{blue!70!black}{positive ridge} and its enclosed $\Omega$, which is the union of the inner disk around $0$ and the exterior of the outer circle containing $I_n$. The \textcolor{red}{dot} marks $-n\lambda$.}
    \label{fig:ridge-filter-contours}
\end{figure}

To apply Theorem~\ref{thm:analytic-filter-risk}, it is instructive to construct the contour $\Gamma$ in Assumption~\ref{assmp:analytic-spectral-filter}. The two-circle contour in Figure~\ref{fig:ridge-filter-contours} suffices for all bounded $\lambda \geq -\eta_n/2$ and fixed $k$. We highlight that on the component around zero, we separately define $\upphi_{k,\lambda}(z/n)=0$ and hence $\uppsi_{k,\lambda}(z/n)=1$, while the rational formula is used on the component containing $I_n$. Since these components are disjoint, the two definitions give a holomorphic function on their union (cf.~\citep[Chapter~III]{conway1978functions}). In the ridge-type filters case, it is advantageous to directly use the contour formula in Eq.~\eqref{eq:analytic-filter-representation} and use the Cauchy's formula to obtain the exact forms of $\effectiveVariance$ and $\effectiveBias$. We write $\omega[x,x']:=\wt{\omega}[nx,nx']$ with the same divided-difference convention, summarize our results next, and defer its proof to Appendix~\ref{proof:ridge-type-risk}.

\begin{proposition}[Risk of ridge-type filters]\label{prop:ridge-type-risk}
    Under Assumption~\ref{assmp:features}, take either spectrum in Eq.~\eqref{eq:spectrum-decay-examples} and let $\bbeta=\bSigma^{1/2}\btheta$, where $\|\btheta\|<\infty$. Fix $\kappa>0$ and an integer $k\geq1$. There exists a deterministic $\eta_n=\Theta(\tau_1(n)/n)$ satisfying $\P\prn{\what\sigma_n\ge\eta_n}\ge1-2e^{-n}$, such that the following holds for all sufficiently large $n$ and $\lambda\geq-\eta_n/2$ with $\lambda=O(1)$.
    \begin{enumerate}
        \item[(i)] For the iterated ridge filter $\upphi_{k,\lambda}$, the effective variance and squared bias are
        \begin{align*}
            \effectiveVariance[\upphi_{k,\lambda}]&=\tau^2\brc{\sum_{m=1}^{2k-1}\frac{(k-|k-m|) \cdot \lambda^{m-1}\omega^{(m)}(-\lambda)}{m!}-1},\\
            \effectiveBias[\upphi_{k,\lambda}]&=\left.\sum_{a,b=0}^{k-1}\frac{\lambda^{a+b}}{a!b!}\frac{\de^{a+b}}{\de x^a\de {x'}^b}\brc{\omega[x,x']\,\bbeta^*\bSigma\prn{\bI-\frac{\bSigma}{\omega(x)}}^{-1}\prn{\bI-\frac{\bSigma}{\omega(x')}}^{-1}\bbeta}\right|_{x=x'=-\lambda}.
        \end{align*}
        We use $\wt{O}$ to extend the big-Oh notation by suppressing factors polynomial in $\log n$, and leave dependence on the decay exponent $\alpha$ or $\beta$ implicit. For $2\le p\le\kappa\log n$, we have
        \begin{align*}
            \big\|\ind\{\what\sigma_n\ge\eta_n\}(\sampleVariance[\upphi_{k,\lambda}]-\effectiveVariance[\upphi_{k,\lambda}])\big\|_{L^p}
            &=\wt{O}_{C_{\bx},\kappa,k}\left(\frac{\tau^2}{n}\right),\\
            \big\|\ind\{\what\sigma_n\ge\eta_n\}(\sampleBias[\upphi_{k,\lambda}]-\effectiveBias[\upphi_{k,\lambda}])\big\|_{L^p}
            &=\wt{O}_{C_{\bx},\kappa,k}\left(\frac{(\lambda+\eta_n)^2\|\btheta\|^2}{\sqrt n}\right).
        \end{align*}
        \item[(ii)] For $k=1$, let $\lambda_\star=-\omega(-\lambda)$ be the effective regularization in Eq.~\eqref{eq:ridge-effective-regularization}, and set $\mathsf{J}_n(\lambda)=n^{-1}\Tr(\bSigma^2(\bSigma+\lambda_\star\bI)^{-2})<1$. Then the formulas in (i) reduce to
        \begin{align*}
            \effectiveVariance[\upphi_{1,\lambda}]=\tau^2\frac{\mathsf{J}_n(\lambda)}{1-\mathsf{J}_n(\lambda)}
            \qquad\text{and}\qquad
            \effectiveBias[\upphi_{1,\lambda}]=\frac{\lambda_\star^2}{1-\mathsf{J}_n(\lambda)}\big\|\bSigma^{1/2}(\bSigma+\lambda_\star\bI)^{-1}\bbeta\big\|^2.
        \end{align*}
        For $\tau>0$ and $\btheta\neq0$, dividing the bounds in (i) by these quantities yields
        \begin{align*}
            \left\|\ind\{\what\sigma_n\geq\eta_n\}\left(\frac{\sampleVariance[\upphi_{1,\lambda}]}{\effectiveVariance[\upphi_{1,\lambda}]}-1\right)\right\|_{L^p}
            &=\wt{O}_{C_{\bx},\kappa}\left(\frac1n\cdot\frac{1-\mathsf{J}_n(\lambda)}{\mathsf{J}_n(\lambda)}\right),\\
            \left\|\ind\{\what\sigma_n\geq\eta_n\}\left(\frac{\sampleBias[\upphi_{1,\lambda}]}{\effectiveBias[\upphi_{1,\lambda}]}-1\right)\right\|_{L^p}
            &=\wt{O}_{C_{\bx},\kappa}\left(\frac1{\sqrt n}\cdot\frac{(\lambda+\eta_n)^2(1-\mathsf{J}_n(\lambda))\|\btheta\|^2}{\lambda_\star^2\|\bSigma(\bSigma+\lambda_\star\bI)^{-1}\btheta\|^2}\right).
        \end{align*}
    \end{enumerate}
\end{proposition}
\begin{remark}
    The relative (multiplicative) bounds in (ii) preserve the population spectrum and the source dependence explicitly, and can still vanish when these factors grow sufficiently slowly. For polynomially varying spectra, the ridgeless bias bounds of \citep{cheng2024dimensionfree} do not, in general, yield a vanishing relative error~\citep[Remark 4.2]{cheng2024dimensionfree}. We remove this barrier and recover optimal fluctuation rates $\wt{O}(n^{-1})$ and $\wt{O}(n^{-1/2})$ for the variance and squared bias. The proposition also covers negative ridge with $-\eta_n/2\leq\lambda<0$ and its iterated versions.
\end{remark}

\subsection{Gradient flow}\label{sec:gradient-flow}
We next consider gradient flow for the empirical squared loss in the Hilbert space $\mc{H}$ initialized at zero:
\begin{align*}
    \frac{\de}{\de t}\what{\bbeta}_t=\frac1n\bX^*(\by-\bX\what{\bbeta}_t)\qquad\text{and}\qquad \what{\bbeta}_0=0.
\end{align*}
Solving this equation yields $\what{\bbeta}_t=\upphi_t(\what{\bSigma})\bX^*\by/n$, with filter and residual
\begin{align*}
    \upphi_t(x)=\frac{1-e^{-tx}}{x}\qquad\text{and}\qquad\uppsi_t(x)=e^{-tx},
\end{align*}
where $\upphi_t(0)=t$ by continuity. The stopping time $t$ controls which directions are fitted, since the residual is small when $tx$ is large. This filter also arises from iterated ridge by taking $\lambda=k/t$ and $k\to\infty$, since $(1+tx/k)^{-k}\to e^{-tx}$.

To invoke Theorem~\ref{thm:analytic-filter-risk}, we consider the two-ray contour indicated in Figure~\ref{fig:analytic-filter-support}, with the origin of the rays at $-n/t$. To evaluate the effective variance and squared bias, it is in fact convenient to use the real-line formulas in Proposition~\ref{prop:filter-real-risk}, which reduce to the Laplace transforms of the probability measures in Eq.~\eqref{eq:filter-directional-measures} through
\begin{align*}
    \mathcal L_{\omega,j}(t):=\int_{[0,\infty)}e^{-tx}\mu_{\omega,j}(\de x),\qquad t\geq0.
\end{align*}
We summarize the resulting formulas and approximation bounds next, and defer the proof to Appendix~\ref{proof:gradient-flow-risk}.
\begin{proposition}[Risk of gradient flow]\label{prop:gradient-flow-risk}
    Under the same setup of Proposition~\ref{prop:ridge-type-risk}, for all sufficiently large $n$ and $t=t_n>0$, the effective variance and squared bias of $\upphi_t$ are
    \begin{align*}
        \effectiveVariance[\upphi_t]&=\frac{\tau^2}{n}\sum_j\sigma_j\int_0^t\bigl(\mathcal L_{\omega,j}(r)-\mathcal L_{\omega,j}(r+t)\bigr)\,\de r,\\
        \effectiveBias[\upphi_t]&=\sum_j\sigma_j|\langle\bu_j,\bbeta\rangle|^2\brc{\mathcal L_{\omega,j}(t)^2+\int_{[0,\infty)}\left(\int_0^t e^{-(t-r)x}\mathcal L_{\omega,j}(r)\,\de r\right)^2\mu_\omega(\de x)}.
    \end{align*}
    For $2\leq p\leq\kappa\log n$, using the $\wt{O}$ convention of Proposition~\ref{prop:ridge-type-risk}, we have the following bounds with implied constants independent of $t$:
    \begin{align*}
        \big\|\ind\{\what\sigma_n\geq\eta_n\}(\sampleVariance[\upphi_t]-\effectiveVariance[\upphi_t])\big\|_{L^p}
        &=\wt{O}_{C_{\bx},\kappa}\left(\frac{\tau^2}{n}\right),\\
        \big\|\ind\{\what\sigma_n\geq\eta_n\}(\sampleBias[\upphi_t]-\effectiveBias[\upphi_t])\big\|_{L^p}
        &=\wt{O}_{C_{\bx},\kappa}\left(\frac{(\eta_n+t^{-1})^2\|\btheta\|^2}{\sqrt n}\right).
    \end{align*}
\end{proposition}
The reader may notice that our bias approximation bounds can deteriorate near initialization $t \to 0+$ or under strong ridge regularization $\lambda\to\infty$. This is because our bounds are adapted to bulk scales. Indeed, strong regularization in those cases are easier regimes, and can be controlled separately by standard concentration arguments. We focus on growing stopping times $t\to\infty$ up to the bulk scale $t=O(\eta_n^{-1})$.

\begin{remark}
    Recent work on early stopping and gradient descent with Hilbert-space features uses concentration and nonasymptotic risk bounds to compare with ridge regression~\citep{wu2026riskcomparisons}, while existing asymptotic risk formulas exclude the power-law bulk scale~\citep{liganshilinrisk,li2026riskequivalence}. For both regimes considered here, the bulk stopping scale is $t=\Theta(\eta_n^{-1})=\Theta(n/\tau_1(n))$. Proposition~\ref{prop:gradient-flow-risk} supplies explicit deterministic equivalents at and beyond this scale, with quantitative errors at the $n^{-1}$ and $n^{-1/2}$ fluctuation rates after the displayed scale factors.
\end{remark}

An interesting immediate observation is that these formulas make the source-dependent variance-bias trade-off explicit through $|\langle\bu_j,\bbeta\rangle|^2$. Indeed, letting $(e^{-x\cdot}*\mathcal L_{\omega,j})(t)=\int_0^t e^{-(t-r)x}\mathcal L_{\omega,j}(r)\de r$ and differentiating the formulas in Proposition~\ref{prop:gradient-flow-risk} yields
\begin{align*}
    \frac{\de}{\de t}\effectiveVariance[\upphi_t]&=\frac{2\tau^2}{n}\sum_j\sigma_j\bigl(\mathcal L_{\omega,j}(t)-\mathcal L_{\omega,j}(2t)\bigr),\\
    \frac{\de}{\de t}\effectiveBias[\upphi_t]&=2\sum_j\sigma_j|\langle\bu_j,\bbeta\rangle|^2\brc{\mathcal L_{\omega,j}\mathcal L'_{\omega,j}+\int_{[0,\infty)}(e^{-x\cdot}*\mathcal L_{\omega,j})\bigl(\mathcal L_{\omega,j}-x(e^{-x\cdot}*\mathcal L_{\omega,j})\bigr)\mu_\omega(\de x)}.
\end{align*}
In the second line, all transforms and convolutions are evaluated at $t$, with
\begin{align*}
    \mathcal L'_{\omega,j}(t)=-\int_{[0,\infty)}x e^{-tx}\mu_{\omega,j}(\de x)
\end{align*}
This suggests an explicit stationary condition for the optimal (early) stopping time $t_{\mathrm{opt}}>0$, upon its existence, satisfying
\begin{align*}
    \left.\frac{\de}{\de t}\effectiveVariance[\upphi_t]\right|_{t=t_{\mathrm{opt}}}=-\left.\frac{\de}{\de t}\effectiveBias[\upphi_t]\right|_{t=t_{\mathrm{opt}}}.
\end{align*}
This could yield further insights into optimal stopping and how it relates to the source distribution over high and low frequency features. We defer the analysis of this trade-off and the resulting derivation of stopping rules to future work.

    \section{Numerical experiments} \label{sec:numerical}

We complement our theoretical results with three numerical experiments: (i) empirical spectral laws, (ii) ridge-type filters, and (iii) gradient flow. We compare the deterministic predictions with simulations under the power and critical laws of Section~\ref{sec:model-setup}. Codes are in the \href{https://github.com/Moriartycc/dimension-free-rmt}{public repo}.

Throughout these experiments, we use Gaussian covariates with diagonal covariance, taking decay exponents $0.2$ and $0.05$. We truncate the spectra at $d=2^{24}$ for computational efficiency and compute the predictions using the same truncated covariance. For the regression experiments, we take a Poisson source---to vary the source weights across population eigendirections---$\bbeta=\bSigma^{1/2}\btheta$, where
\begin{align*}
    \btheta\in\{\btheta_q:q>0\}, \qquad
    (\btheta_q)_j
    \propto \frac{q^{j-1}}{(j-1)!},
    \quad j=1,\ldots,d.
\end{align*}
$\btheta_q$ has coordinates proportional to the Poisson probabilities and we normalize so that $\norm{\btheta_q}_2=1$. The details of the sample sizes, filter parameters, and source parameter $q$ are in each experiment below.

\begin{figure}[htbp!]
    \centering
    \includegraphics[width=\linewidth]{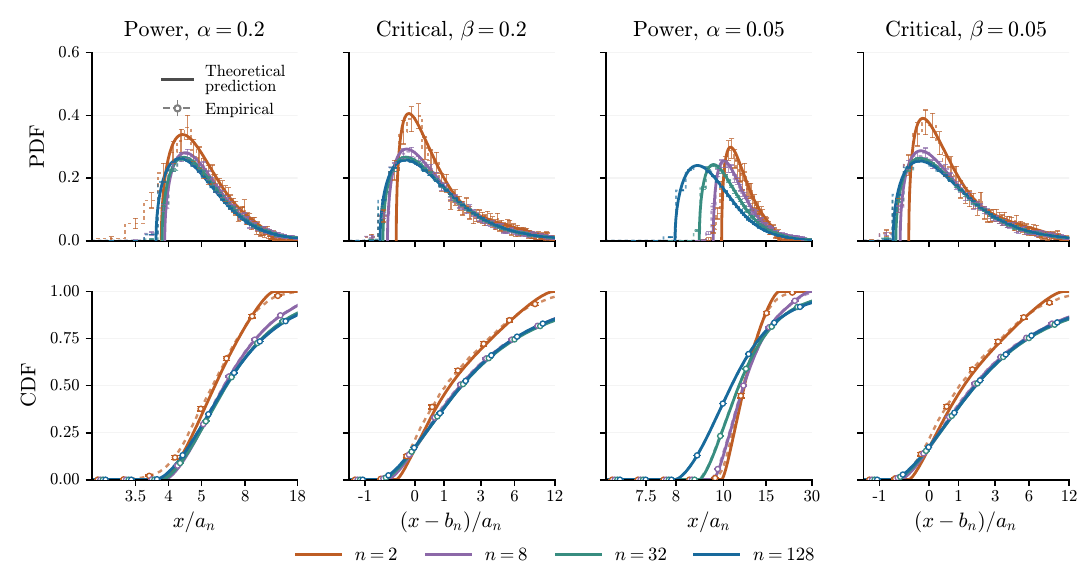}
    \caption{Empirical vs. predicted spectral laws under power and critical decay regimes with $n \in \{\textcolor[HTML]{BE5D24}{2},\textcolor[HTML]{8C68A8}{8},\textcolor[HTML]{388D80}{32},\textcolor[HTML]{186A9C}{128}\}$. We average over $2{,}000$ independent Gaussian designs. Solid curves show the theoretical predictions, and error bars give pointwise $95\%$ Monte Carlo intervals computed across matrices. For eigenvalues $x$ of $\bX\bX^*$, the horizontal coordinate is $x/a_n$ for the power law and $(x-b_n)/a_n$ for the critical law. We use shifted logarithmic x-axes for better visualization.}
    \label{fig:numerical-spectral-laws}
\end{figure}

\paragraph{Empirical and predicted spectral laws.}
We first test our main result I in Theorem~\ref{thm:gram-stieltjes-approx}, comparing the empirical spectral laws $\what{\mu}_n$ with their deterministic predictions $\mu_n$ under the four covariance settings above. For the purpose of nonasymptotic comparisons, we choose relatively small $n \in \{2,8,32,128\}$. Following Proposition~\ref{prop:canonical-bulk-laws}, we rescale the eigenvalues by $a_n=n\sigma_n$ in the power law case and center by $b_n=\tau_1(n)$ from Eq.~\eqref{eq:tail-parameters-trace} before rescaling by $a_n$ in the critical law case. Figure~\ref{fig:numerical-spectral-laws} presents the empirical/predicted PDFs and CDFs.

We observe that in four settings, the empirical histograms and CDFs agree increasingly well with the finite-sample predictions as $n$ grows. The predictions capture both the bulk shape and its dependence on $n$. The visible discrepancies occur only near the lower edge when $n$ is small.

\paragraph{Ridge-type filters.}
We next examine our main result II, comparing empirical results with the theoretical predictions in Eq.~\eqref{eq:analytic-filter-risk} for ridge-type filters first. We use the Poisson source $q=16$, noise standard deviation $\tau=0.01$, and sample sizes $n\in\{16,32,64,128,192,256,384,512\}$. We take $\lambda\in\{-0.6,0,0.5,1,1.5,2\}\kappa_n$, where $\kappa_n:=n^{-1}\inf\supp\mu_n$ is the predicted lower spectral edge of $\bX\bX^*/n$ for the truncated covariance.

Figure~\ref{fig:ridge-q16-k1-absolute} shows the unnormalized results for $k=1$. Figures~\ref{fig:ridge-q16-k1-relative} and~\ref{fig:ridge-q16-k2-relative} compare these quantities relative to ridgeless regression for $k=1$ and $k=2$, respectively. The normalized predictions are $\effectiveVariance{}/\mathsf{V}_0$, $\effectiveBias{}/\mathsf{B}_0$, and $\effectiveRisk{}/\mathsf{R}_0$, using (theoretical or empirical) ridgeless baselines at the same $n$ and $\bSigma$. The exact effective risk formulas for $k=2$, and additional corresponding results for a different shape source $q=24$ are in Appendix~\ref{app:additional-ridge-experiments}.

\begin{figure}[!htbp]
    \renewcommand{\thefigure}{\arabic{figure}a}
    \renewcommand{\theHfigure}{\arabic{figure}a}
    \centering
    \includegraphics[width=\linewidth]{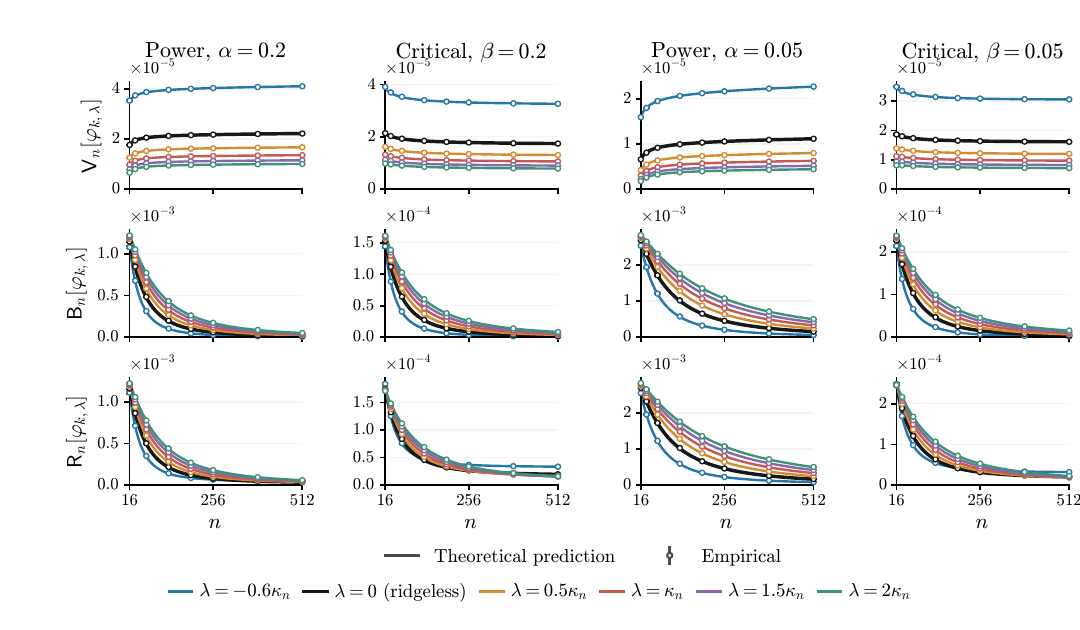}
    \caption{Ridge-type filters with $q=16$ and $k=1$. Rows show $\effectiveVariance{}$, $\effectiveBias{}$, and $\effectiveRisk{}$. Solid curves show the theoretical predictions. We show markers for the simulated results, with error bars by pointwise $95\%$ Monte Carlo intervals over $400$ designs.}
    \label{fig:ridge-q16-k1-absolute}
\end{figure}

From the plots, we see the empirical means closely match the theoretical predictions, with error bars that are nearly imperceptible. It is important to note that, unlike standard ridge regression---where $\effectiveVariance{}$ decreases monotonically with $\lambda$---iterated ridge does not generally exhibit monotonic behavior, as shown in Figure~\ref{fig:ridge-q16-k2-relative}. We also observe there is no universal dominance of negative, positive or ridgeless choices for ridge-type filters. The optimal choice depends on the spectrum, sample size and also the source shape, parameterized by $\bSigma, n ,q$ respectively in our experiments.

\begin{figure}[!htbp]
    \addtocounter{figure}{-1}
    \renewcommand{\thefigure}{\arabic{figure}b}
    \renewcommand{\theHfigure}{\arabic{figure}b}
    \centering
    \includegraphics[width=\linewidth]{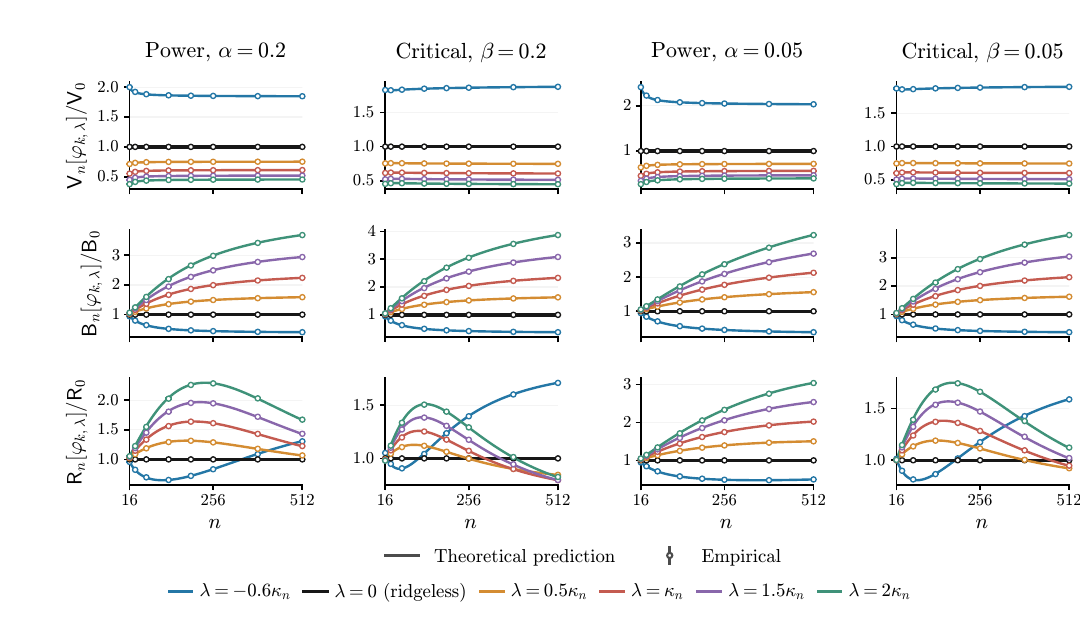}
    \caption{Ridge-type filters with $q=16$ and $k=1$, normalized by the respective ridgeless baselines.}
    \label{fig:ridge-q16-k1-relative}
\end{figure}

\begin{figure}[!htbp]
    \addtocounter{figure}{-1}
    \renewcommand{\thefigure}{\arabic{figure}c}
    \renewcommand{\theHfigure}{\arabic{figure}c}
    \centering
    \includegraphics[width=\linewidth]{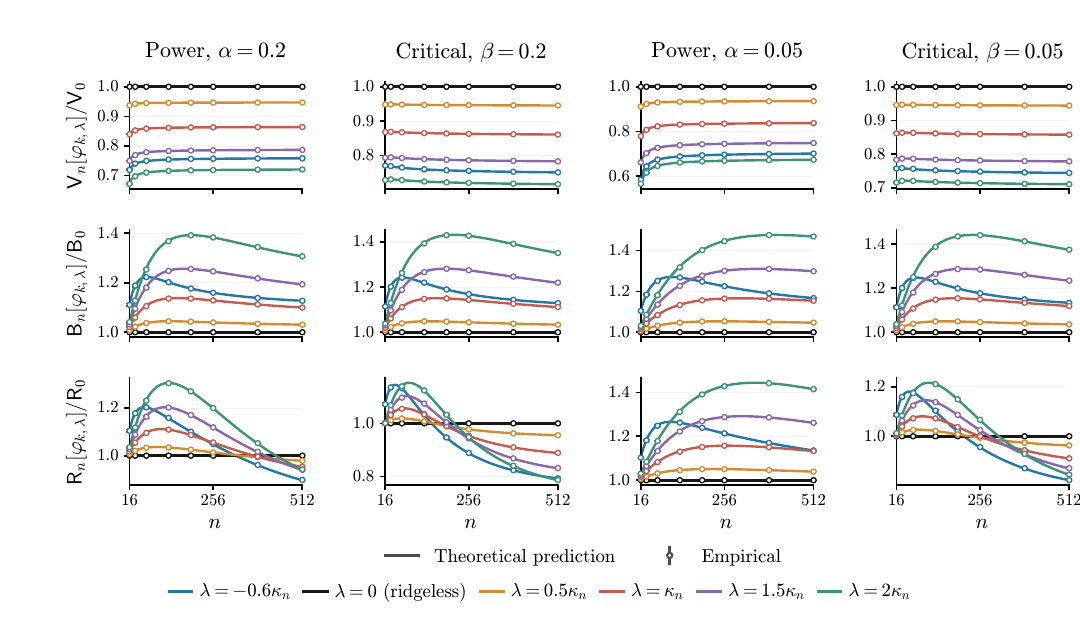}
    \caption{Iterated ridge with $q=16$ and $k=2$, normalized by the respective ridgeless baselines.}
    \label{fig:ridge-q16-k2-relative}
\end{figure}

\begin{figure}[!htbp]
    \centering
    \includegraphics[width=\linewidth]{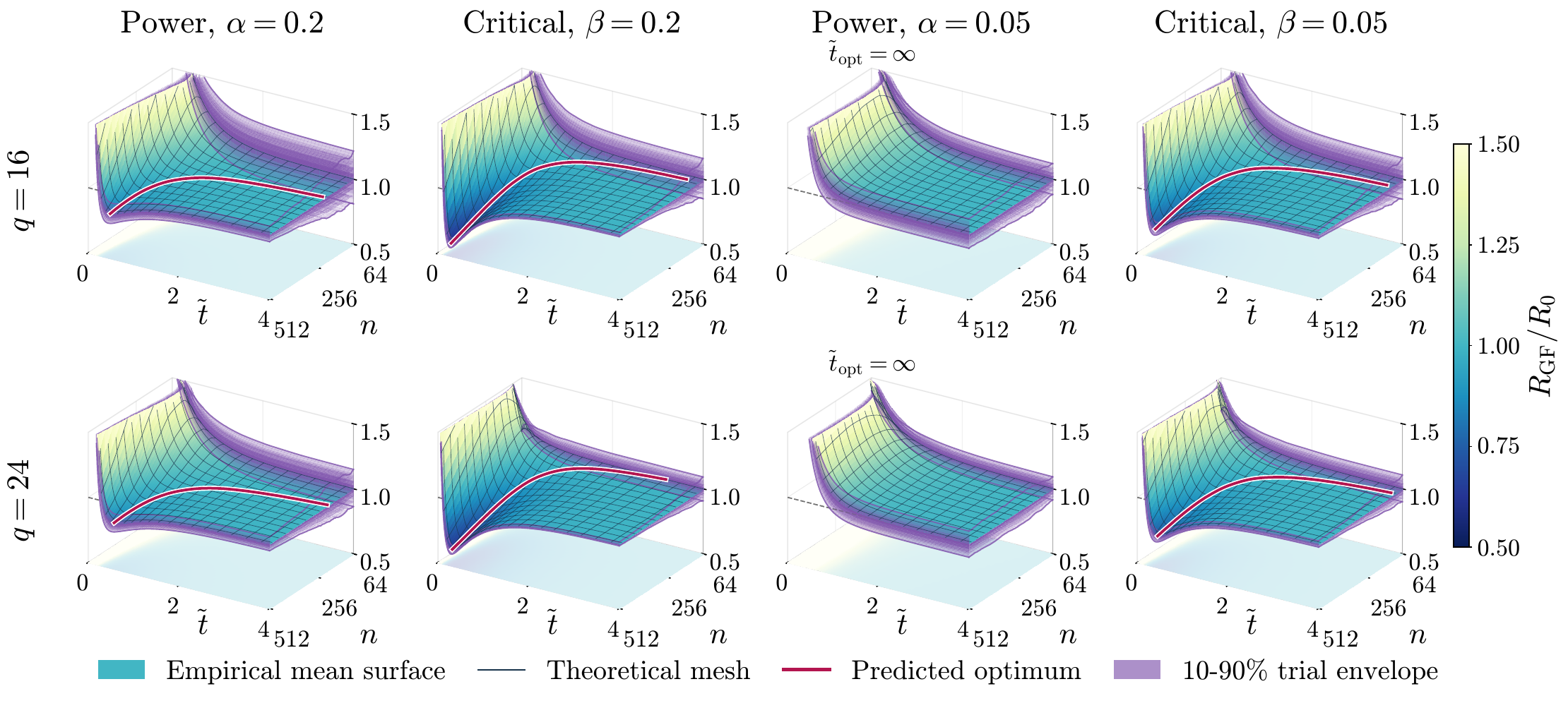}
    \caption{Gradient-flow risk relative to ridgeless regression. Surfaces show empirical means over $400$ designs, \textcolor[HTML]{18334C}{meshes} show theoretical predictions, and \textcolor[HTML]{B3154F}{curves} show predicted finite stopping optima. \textcolor[HTML]{8055AD}{Envelopes} show pointwise $10\%$ to $90\%$ quantiles across designs. See Appendix~\ref{app:additional-gradient-flow-experiments} for details.}
    \label{fig:numerical-gradient-flow}
\end{figure}

\paragraph{Gradient flow.}
In the last experiment, we investigate gradient flow risks with the theoretical predictions in Proposition~\ref{prop:gradient-flow-risk}. We use $q\in\{16,24\}$, $\tau=0.01$, and $n\in\{64,128,192,256,384,512\}$. Figure~\ref{fig:numerical-gradient-flow} compares empirical and predicted risks, normalized by their respective ridgeless baselines $(t=\infty)$, against $\wt{t}=t\tau_1(n)/n\in[0,4]$. We also make the plot with the predicted optimal stopping times from the stationary condition derived in Section~\ref{sec:gradient-flow}.

We once again see our theoretical predictions closely match the empirical mean risks across the four spectral settings and both sources. However, the benefit of early stopping is more nuanced, depending on the spectrum, sample size, and source distribution. For the power law with $\alpha=0.05$, our numerical optimization favors the ridgeless limit throughout the displayed $n$ range and therefore no optimal curve is shown---however, the numerical evaluation of $\effectiveRisk{}$ suggests an apparent phase transition at larger $n$. This transition lies beyond the empirical range $n\leq512$ in Figure~\ref{fig:numerical-gradient-flow}.

It is possible that there exists a finite $n_0$ such that $t_{\mathrm{opt}} = \infty$ when $n<n_0$, and $t_{\mathrm{opt}} < \infty$ when $n>n_0$. However, the discontinuities in the $t_{\mathrm{opt}}$ curve itself do not always imply such transition point. We give further numerical illustrations of this complicated behavior in Appendix~\ref{app:additional-gradient-flow-experiments}.

    \section{The deterministic self-consistent interpolation}
\label{sec:interpolation-path}
We now introduce the interpolation framework used throughout the paper. To this end, we first specify the appropriate domains on which we can analyze the matrix trace functionals of interest. Recall $z=t+\iu\eta$ and $z_-=(t\wedge0)+\iu\eta$ from Section~\ref{sec:gram-bulk-law}, with $|z_-|=\dist(z,\R_{\geq0})$. We write additionally $w=u+\iu v$ with $u,v\in\R$.
Denote by $t_+:=t \vee 0$ and $t_-:=(-t) \vee 0$, so $z_-=z-t_+=-t_-+\iu\eta$.

Throughout this section, $z$ denotes the spectral parameter on the
unnormalized scale of $\bT_k$, and thus $z/n$ is the corresponding spectral
parameter for the normalized sample covariance operator. Let their natural domain be $\mc{D} := \{(z,w) : \eta >0, v \geq 0\}$. We extend their domain to include subsets of the real line by
\begin{align*}
    \mc{D}_k & := \brc{(z, w): \eta, v \geq 0, 0 \notin \Spec(\bT_k - w \bSigma - z \bI) }.
\end{align*}
The domain $\mc{D}_k$ depends on the sample. Since $\Im(\bT_k - w\bSigma - z\bI) = -v \bSigma -\eta \bI \preceq -\eta \bI$, we have $\mc{D} \subset \mc{D}_k$ for all realizations $\bT_k$. For $(z,w)\in\mc{D}_k$ and $\bQ \in \mc{L}(\mc{H})$, define the relative-resolvent trace $H_k^{\bQ}(z,w)$; and for $(z,w) \in \mc{D}_{k-1}$, the leave-one-out trace $G_k^{\bQ}(z,w)$ by
\begin{subequations}
    \begin{align}
        \bR_k(z,w)
        &:=\prn{\bT_k-w\bSigma-z\bI}^{-1},   & &                              \label{eq:def-relative-resolvent}\\
        H_k^{\bQ}(z,w)
        &:=\Tr\prn{\bQ \prn{\bT_k-w\bSigma}\bR_k(z,w)}, & G_k^{\bQ}(z, w) &:=\Tr\prn{
        \bSigma^{\half} \bQ \bSigma^{\half}\prn{\bT_{k-1}-w\bSigma-z\bI}^{-1}
        }                       \label{eq:def-traces-interpolation}
    \end{align}
    The traces in
    Eq.~\eqref{eq:def-traces-interpolation} are finite because
    $\bSigma$  and  $\bT_k$ are trace-class, and we suppress their dependence in $\bQ$ when $\bQ=\bI$ by writing $H_k := H_k^{\bI}$ and $G_k := G_k^{\bI}$. We explain the terms relative-resolvent and leave-one-out traces by the following identities when $z \neq 0$:
    \begin{align}
        & \textrm{(Relative-resolvent trace)} &  H_k(z,w)
        & =z\Tr\prn{\bR_k(z,w)-(-z\bI)^{-1}};                        \label{eq:relative-resolvent-trace-identity} \\
        & \textrm{(Leave-one-out trace)} &  G_k(z,w)
        & =\E_{\bx_k}\brk{\left\langle
        \bx_k,
        \prn{\bT_{k-1}-w\bSigma-z\bI}^{-1}\bx_k
        \right\rangle_{\mc H}
        }.             \label{eq:leave-one-out-trace-identity}
    \end{align}
\end{subequations}
We point out that: (i) in Eq.~\eqref{eq:relative-resolvent-trace-identity} for $H_k$, neither resolvent within the right hand side need be trace
class, but their difference is, and $H_k$ measures (the rescaled trace of) $\bR_k$ relative to the
resolvent of the zero operator in $\mc{L}(\mc{H})$; (ii) in Eq.~\eqref{eq:leave-one-out-trace-identity},  $\E_{\bx_k}$ denotes expectation with respect to $\bx_k$ alone with
$\bT_{k-1}$ held fixed. We also identify the quantity $G_k^{\bQ}$ corresponds to $\mathscr{R}_{k-1}(\bQ)$ used by~\citeauthor{cheng2024dimensionfree} with real arguments. We highlight that these weights $\bQ$ in the resolvents allow us to probe selected directions, and their derivatives help the analysis of general spectral filters, especially for the bias term.

\subsection{Trace interpolation for spectral and regression analysis}
The preceding resolvent traces in Eqs.~\eqref{eq:relative-resolvent-trace-identity} and \eqref{eq:leave-one-out-trace-identity} are central to study both the spectrum of the sample covariance operator and the bias and variance of general regression estimators defined by spectral filters. We will construct a deterministic parameter sequence
$\{w_k\}_{k=0}^n$ with $w_n=0$ such that $H_k(z,w_k)$ changes only by a
second-order error in expectation. Its endpoints are
\begin{align}
    H_0(z,w_0)
    &=\Tr\prn{(-w_0\bSigma)(-w_0\bSigma-z\bI)^{-1}},
    &
    H_n(z,w_n)=H_n(z,0)
    &=\Tr\prn{\bT_n(\bT_n-z\bI)^{-1}}.                            \label{eq:relative-resolvent-endpoints}
\end{align}
We call $\{H_k(z,w_k)\}_{k=0}^n$ the \textit{relative-resolvent
interpolation process}, or simply the \textit{interpolation process}. For retrieval of the spectrum,
the terminal quantity $H_n(z,0)$ contains exactly
the spectral information carried by the empirical Stieltjes transform. Indeed,
since $\bT_n=\bX^*\bX$ and $\bX\bX^*$ has the same nonzero eigenvalues,
\begin{align}
    \what{s}_n(z)=\frac{1}{n}\Tr\prn{(\bX\bX^*-z\bI_n)^{-1}}
    =\frac{1}{z}\prn{\frac{1}{n}H_n(z,0)-1}.                     \label{eq:H-and-empirical-Stieltjes-transform}
\end{align}
Thus $H_n(z,0)$ determines the Stieltjes transform of the finite empirical
Gram matrix while remaining well-defined on an infinite-dimensional feature
space. By contrast, $\Tr\prn{\bR_n(z,0)}$ is not defined when
$\dim(\mc{H})=\infty$ as $\bR_n(z,0)$ equals $-z^{-1}\bI$ on the
infinite-dimensional kernel of $\bT_n$.

For $z=-\lambda<0$,
$H_n(-\lambda,0)=\Tr\prn{\bT_n(\bT_n+\lambda\bI)^{-1}}$ is the empirical
effective dimension familiar from kernel learning~\citep{zhang2002effective,caponnetto2007optimal,blanchard2020kernel}.
The subtraction of a reference resolvent, arranged so that the resulting difference belongs to the trace class, is a standard technique in the context of trace ideals and spectral shift theory~\citep{schatten1960norm,gohberg1969introduction,simon2005trace}.
\citet{misiakiewicz2024nonasymptotic} consider exactly
$n^{-1}H_n(-\lambda,0)$ as one of their matrix functionals and derive its
deterministic equivalent by a leave-one-out argument as an ingredient in their
analysis of kernel ridge regression training error and generalized
cross-validation. Our aim, more universal, is to develop a common interpolation framework for spectral analysis and regression. Along the same deterministic path, we control weighted resolvent traces and their products, which provide the quantities needed to analyze the variance and bias of general spectral filters.

The ingredients underlying the parameter update are classical.
Expected leave-one-out resolvents and their self-consistent equations are
standard tools for deterministic equivalents of sample covariance
matrices~\citep{hachem2007deterministic,chouard2022quantitative}, but they typically reduce to a simple linear interpolation.
In~\citep{cheng2024dimensionfree}, the authors construct a closely related stochastic interpolation process to cancel the first order drift terms. Our novel construction to be introduced next combines the stochastic interpolation with an expectation-level simplification--its parameter
is deterministic and is chosen to cancel the leading unconditional mean drift.
We call the resulting parameter sequence the \textit{deterministic
self-consistent path}, distinguishing it from both the adapted random path and
the standard linear path.

\subsection{The backward path and its stability}
\label{sec:expected-partial-interpolation}
We formally construct the alluded deterministic path. Fix $z\neq0$, and for $(z,w)\in\mc D\cup\{(z,w)\in\R^2:z<0,\ w\leq0\}$, define the expected leave-one-out trace
\begin{align}
    \wb{G}_k(z,w)
    &:=\E\left[G_k(z,w)\right] =\E\Tr\prn{
    \bSigma\prn{\bT_{k-1}-w\bSigma-z\bI}^{-1}
    }.                                                           \label{eq:def-expected-leave-one-out-trace}
\end{align}
Starting from $w_n=0$, we use this deterministic trace to define the backward updates for $k=n,\ldots,1$ as below (cf.~Theorem~\ref{thm:self-consistent-path-stability}.(i) guarantees that each update is well-defined)
\begin{subequations}
    \begin{align}
        \textrm{(Self-consistent backward path)} & & w_{k-1} :=w_k-\frac{1}{1+\wb{G}_k(z,w_k)}.
        \label{eq:self-consistent-expected-trace-update}
    \end{align}
    We highlight that the path itself depend only on the unweighted leave-one-out traces, and the same path will also be used for weighted traces. It is instructive to compare this path with the standard linear interpolation and the adapted forward construction of~\citet{cheng2024dimensionfree}, written after a change of variables, start from $w_0\lin=\mu_0=-z/\omega(z/n)$ and satisfy
    \begin{align}
        \textrm{(Linear interpolation path)} & & w_{k-1}\lin :=w_k\lin-\frac{1}{1+\wb{G}_1(z,w_0\lin)}; \label{eq:self-consistent-linear-interpolation} \\
        \textrm{(Self-consistent forward path)} & & \mu_{k-1} :=\mu_k-\frac{1}{1+G_k(z, \mu_k)}. \label{eq:self-consistent-forward-update}
    \end{align}
\end{subequations}
The defining equation for $\omega(z/n)$ implies $1+\wb G_1(z,w_0\lin)=-n/w_0\lin$, so the linear path has constant increments and $w_k\lin=(1-k/n)w_0\lin$, while the backward path instead updates its increment using the current expected trace.
\begin{figure}[!ht]
    \centering
    \includegraphics[width=\textwidth]{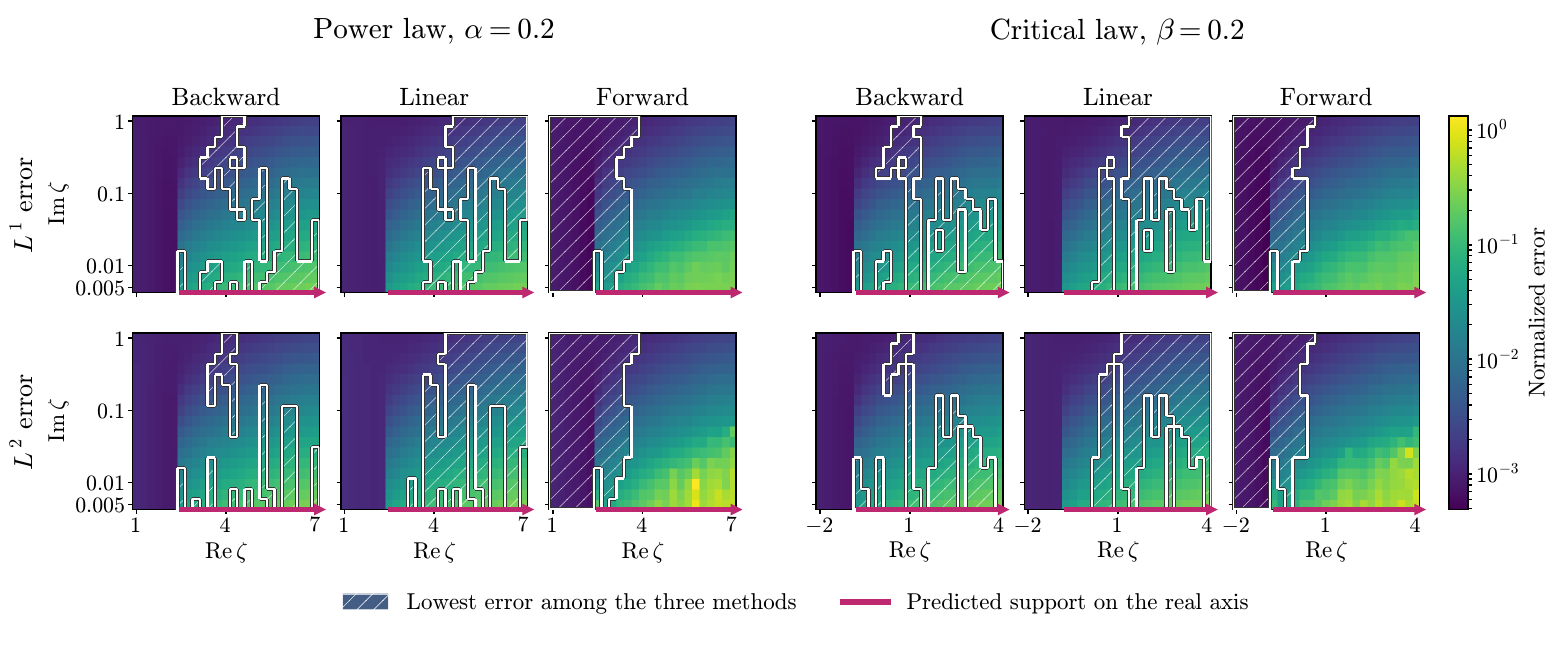}
    \caption{Endpoint trace errors for backward, linear, and forward interpolation under power law (left) and critical law (right) spectra. Rows show empirical $L^1$ and $L^2$ errors. We use brighter colors to indicate larger errors. See Appendix~\ref{app:numerical-details} for error definitions and numerical details.}
    \label{fig:interpolation-endpoint-errors}
\end{figure}

\begin{remark}
    We highlight the beneficial deterministic definition of the backward path, as we can impose the target endpoint $w_n=0$ and evaluate
    Eq.~\eqref{eq:self-consistent-expected-trace-update} explicitly. While the martingale updates in~Eq.~\eqref{eq:self-consistent-forward-update} used by~\citet{cheng2024dimensionfree} are generated forward from implicit self-consistent equations, which
    require inverting the map $w\mapsto w-(1+G_k(z,w))^{-1}$ at every step. This introduces questions of
    existence and possible instability. The backward definition circumvents solving the implicit iteration and fixes
    the terminal relative-resolvent trace $H_n(z,0)$ exactly. Figure~\ref{fig:interpolation-endpoint-errors} illustrates the resulting accuracy: backward interpolation has smaller endpoint trace errors than forward interpolation over much of the bulk at small spectral heights, particularly in $L^2$. Related linear-interpolation estimates appear in~\citet[Theorem~5]{hastie2022surprises} and~\citet[Theorem~4 and Remark~5.1]{misiakiewicz2024nonasymptotic}.
\end{remark}
We first present the following stability theorem regarding the deterministic backward path, whose proof is in Appendix~\ref{proof:self-consistent-path-stability}.
\begin{theorem}[Stability of the self-consistent backward path] \label{thm:self-consistent-path-stability}
    Let Assumption~\ref{assmp:features} hold and fix $z=t+\iu\eta\in\C_+\cup\R_{<0}$.
    \begin{enumerate}
        \item[(i)] (Well-posedness and monotonicity). For every positive integer $n$, the recursion~\eqref{eq:self-consistent-expected-trace-update} uniquely defines a deterministic path $\{w_k\}_{k=0}^n$ with $w_n=0$. If $\eta>0$, then
        \begin{align*}
            \Im w_0>\Im w_1>\cdots>\Im w_n=0.
        \end{align*}
        If $t<0$, then also
        \begin{align*}
            \Re w_0<\Re w_1<\cdots<\Re w_n=0.
        \end{align*}
        For $z<0$, the path remains real for all $0 \leq k \leq n$.

        \item[(ii)] (Expected-trace stability). Under Assumption~\ref{assmp:spectral-stability} at $z$ with $\kappa=1$, the expected leave-one-out traces satisfy
        \begin{align*}
            \max_{1\leq k\leq n}\brc{|\wb G_k(z,w_k)-\wb G_1(z,w_0)| \vee |\wb G_k(z,w_k)-\wb G_n(z,0)|}=O_{C_{\bx},\nu}\prn{\frac{\Lambda_n(z)^{15}}n}.
        \end{align*}

        \item[(iii)] (Path stability). Under the conditions of (ii), the backward path satisfies for $0 \leq k \leq n$,
        \begin{align*}
            |w_k-w_k\lin|=O_{C_{\bx},\nu}\prn{\Lambda_n(z)^{17} \cdot \frac{n-k}{n}},
        \end{align*}
        where $w_k\lin=-(1-k/n)z/\omega(z/n)$.
    \end{enumerate}
\end{theorem}
The above Theorem~\ref{thm:self-consistent-path-stability} controls the expected leave-one-out trace along the backward path and compares that path with the linear interpolation determined by $\omega(z/n)$. We highlight that Theorem~\ref{thm:self-consistent-path-stability} is non-random: the bounds in (ii) and (iii) are deterministic and hold simultaneously for any deterministic set of spectral parameters on which the assumptions are satisfied uniformly, in contrast to the stochastic approximation in~\cite{cheng2024dimensionfree}.

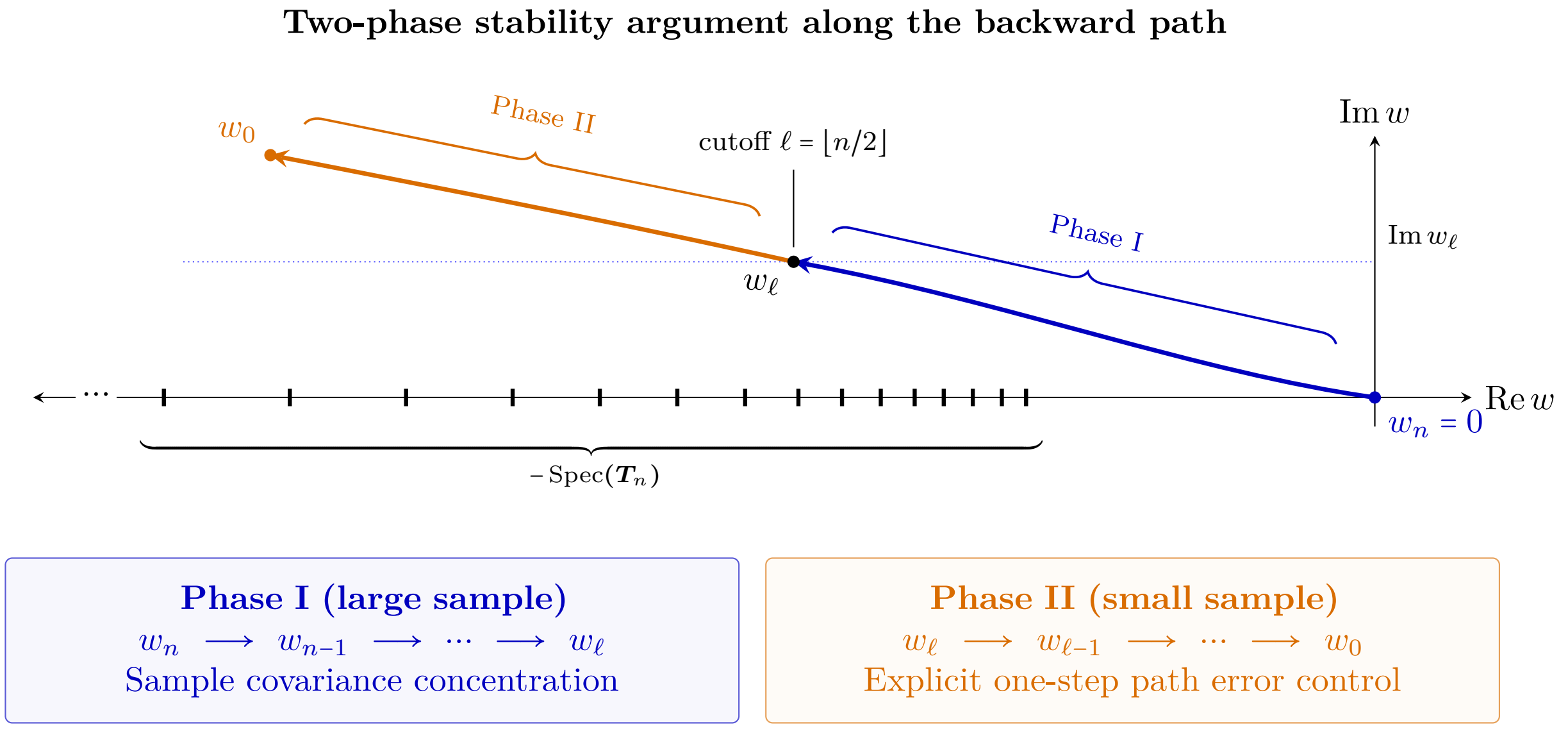
\begin{figure}[htbp!]
    \centering
    \colorlet{pathphaseone}{blue!75!black}
    \colorlet{pathphasetwo}{orange!85!black}
    \usetikzlibrary{decorations.pathreplacing}

    \begin{tikzpicture}[x=1cm,y=1cm,font=\small,>=stealth]
        \node[font=\small\bfseries] at (7.6,4.35) {Two-phase stability argument along the backward path};
        \draw[<->] (0.15,0.5)--(15,0.5) node[right] {$\Re w$};
        \draw[->] (14,0.2)--(14,3.2) node[above] {$\Im w$};
        \node[fill=white,inner sep=2pt] at (0.8,0.5) {$\cdots$};
        \foreach \x in {1.5,2.8,4,5.1,6,6.8,7.5,8.05,8.5,8.9,9.25,9.55,9.85,10.15,10.4}
        \draw[line width=1.2pt] (\x,0.41)--(\x,0.59);
        \node[anchor=north,inner sep=0pt] at (5.95,0.16)
        {$\underbrace{\hspace{9.4cm}}_{\substack{-\Spec(\bT_n)}}$};

        \draw[densely dotted,blue!65] (1.7,1.9)--(14,1.9);
        \node[above right,font=\scriptsize] at (14,1.9) {$\Im w_\ell$};
        \draw[pathphaseone,line width=1.3pt,->] (14,0.5) .. controls (12.2,0.75) and (10.15,1.5) .. (8,1.9);
        \draw[pathphasetwo,line width=1.3pt,->] (8,1.9) .. controls (6.2,2.3) and (4.4,2.65) .. (2.6,3);
        \fill[pathphaseone] (14,0.5) circle (1.8pt) node[below right] {$w_n=0$};
        \fill (8,1.9) circle (1.8pt) node[below left] {$w_\ell$};
        \fill[pathphasetwo] (2.6,3) circle (1.8pt) node[above left] {$w_0$};
        \draw[thin] (8,2.05)--(8,2.85);
        \node[above,font=\scriptsize] at (8,2.85) {cutoff $\ell=\lfloor n/2\rfloor$};
        \draw[pathphaseone,line width=0.7pt,decorate,decoration={brace,amplitude=5pt}]
        (8.4,2.2)--(13.6,1.05)
        node[midway,sloped,above=10pt,align=center,font=\scriptsize,text=pathphaseone] {Phase I};
        \draw[pathphasetwo,line width=0.7pt,decorate,decoration={brace,amplitude=5pt}]
        (2.95,3.32)--(7.65,2.37)
        node[midway,sloped,above=10pt,align=center,font=\scriptsize,text=pathphasetwo] {Phase II};

        \node[anchor=north,align=center,text width=7.15cm,minimum height=1.7cm,draw=pathphasetwo!65,fill=pathphasetwo!3,text=pathphasetwo,rounded corners=2pt,inner sep=6pt] at (11.5,-1.15)
        {\hypersetup{linkcolor=pathphasetwo}\textbf{Phase II (small sample)}\\
        $w_\ell\longrightarrow w_{\ell-1}\longrightarrow\cdots\longrightarrow w_0$\\
        Explicit one-step path error control};
        \node[anchor=north,align=center,text width=7.15cm,minimum height=1.7cm,draw=pathphaseone!65,fill=pathphaseone!3,text=pathphaseone,rounded corners=2pt,inner sep=6pt] at (3.65,-1.15)
        {\hypersetup{linkcolor=pathphaseone}\textbf{Phase I (large sample)}\\
        $w_n\longrightarrow w_{n-1}\longrightarrow\cdots\longrightarrow w_\ell$\\
        Sample covariance concentration};
    \end{tikzpicture}

    \caption{The illustration of the proof strategies along the backward path $\{w_k\}$ on $\C$ with \textcolor{pathphaseone}{Phase I (large sample)} and \textcolor{pathphasetwo}{Phase II (small sample)}. The real-axis markers show the reflected terminal ESD as a reference overlay after rescaling.}
    \label{fig:two-phase-backward-path}
\end{figure}

To provide insights from a theoretical point of view, Figure~\ref{fig:two-phase-backward-path} illustrates the two-phase argument underlying Theorem~\ref{thm:self-consistent-path-stability}. For large sample indices, sample covariance concentration controls the resolvents, while for small sample indices, we use explicit one-step path error control.

\begin{remark}[The polynomial, critical and proportional regimes]\label{rem:polynomial-critical-scales}
    Let $\zeta$ range over a fixed compact subset of $\C_+$, and take $z$ at the scales in Table~\ref{tab:spectral-stability-scales}. Proposition~\ref{prop:canonical-bulk-laws} and the discussions in Appendix~\ref{app:additional-bulk-laws} provide $\Lambda_n(z)=O((\log n)^2)$ for the power law and proportional regime, and $O((\log n)^3)$ for the critical law. Thus, in all three regimes, Theorem~\ref{thm:self-consistent-path-stability}.(ii) and (iii) control the expected traces and the backward path with only logarithmic losses.
\end{remark}
\begin{remark}[Logarithmic and local-scale refinement]
    Under the hypotheses of~\citet[Theorem~1 and Example~1]{zhivotovskiy2022dimensionfree}, the dimension-free covariance bound can potentially sharpen the real resolvent spectral parameter when the tail effective rank is comparable to the sample size. Reaching optimal local scales would require sharpening the dependence on $\Lambda_n(z)$ in both the approximation bounds and the sample-size condition to linear, up to logarithmic factors, as discussed in Appendix~\ref{app:additional-bulk-laws}. We leave these refinements for future work.
\end{remark}

\subsection{Stability of weighted single-resolvents and mixed-resolvents}
\label{sec:trace-drift-fluctuations}

The backward path $\{w_k\}_{k=0}^n$ and its stability allow us to track the random weighted resolvent traces $H_n^{\bQ}$ and $G_{n+1}^{\bQ}$ in Eqs.~\eqref{eq:relative-resolvent-trace-identity} and \eqref{eq:leave-one-out-trace-identity} to their deterministic equivalents $H_0^{\bQ}$ and $G_{1}^{\bQ}$, controlling both their bias and fluctuations with nonasymptotic error bounds. On a high level, our path construction allows us to write
\begin{align}
    H_k^{\bQ}(z,w_k)-H_{k-1}^{\bQ}(z,w_{k-1})&=\xi_k^{H,\bQ}+\beta_k^{H,\bQ} \quad \text{and} \quad G_{k+1}^{\bQ}(z,w_{k+1})-G_k^{\bQ}(z,w_k)=\xi_k^{G,\bQ}+\beta_k^{G,\bQ},
\end{align}
where the $\xi$-terms are martingale differences and the $\beta$-terms form a small accumulated drift, both controlled in $L^p$. Finally, the path stability in Theorem~\ref{thm:self-consistent-path-stability}.(iii) allows us to transfer the backward path at $w_0$ to $w_0\lin$.

Thus, the linear interpolation path makes the deterministic equivalents explicit by comparing the terminal values and the backward path gives sharp bounds. Recall Eq.~\eqref{eq:def-traces-interpolation}, at the sample terminal point, we have
\begin{align*}
    H_n^{\bQ}(z,0)=\Tr \prn{\bQ\bT_n(\bT_n-z\bI)^{-1}}
    \qquad\text{and}\qquad
    G_{n+1}^{\bQ}(z,0)=\Tr \prn{\bQ \bSigma^{\half}\bR_n(z,0)\bSigma^{\half}}.
\end{align*}
At the population terminal point, i.e.\ the linear initial point, substituting $\bT_0=0$ and $w_0\lin=-z/\omega(z/n)$ in Eq.~\eqref{eq:def-traces-interpolation}, we have
\begin{equation}\label{eq:terminal-zero-sample-traces}
    H_0^{\bQ}(z,w_0\lin)=\Tr\prn{\bQ\bSigma(\bSigma-\omega(z/n)\bI)^{-1}} \quad \text{and} \quad G_1^{\bQ}(z,w_0\lin)=\frac{1}{z}\Tr\prn{\bQ\bSigma(\bSigma/\omega(z/n)-\bI)^{-1}}.
\end{equation}
Using $|z_-|=\dist(z,\R_{\geq0})$ as in Assumption~\ref{assmp:spectral-stability} and the shorthand $H_n^{\bQ} = H_n^{\bQ}(z,0), H_0^{\bQ} = H_0^{\bQ}(z, w_0\lin)$, as well as $G_{n+1}^{\bQ}= G_{n+1}^{\bQ}(z,0), G_1^{\bQ}=G_1^{\bQ}(z, w_0\lin)$, the next theorem provides the approximation bounds at their terminal values. The proof is in Appendix~\ref{proof:terminal-weighted-traces}.
\begin{theorem}[Deterministic equivalents of single-resolvents]\label{thm:terminal-weighted-traces}
    Let $z\in\C_+\cup\R_{<0}$, $\kappa>0$ and $r$ be a nonnegative integer. Under Assumptions~\ref{assmp:features} and~\ref{assmp:spectral-stability} at $z$ with parameters $\nu,\kappa$, it holds for all $2\leq p\leq \kappa\log n$ and $\|\bQ\|<\infty$ that
    \begin{align}
        & \left\|\frac{1}{n} \cdot \frac{\de^r}{\de z^r}\prn{H_n^{\bQ}-H_0^{\bQ}}\right\|_{L^p} \vee
        \left\|\frac{\de^r}{\de z^r}\prn{G_{n+1}^{\bQ}-G_1^{\bQ}}\right\|_{L^p} =O_{C_{\bx},\nu,\kappa,r}\prn{\frac{p^2\Lambda_n(z)^{23}}{n|z_-|^r} \cdot \prn{\|\bQ\| \wedge \frac{\|\bQ\|_F}{\sqrt n}}}. \hspace{1em}\label{eq:terminal-trace-bounds}
    \end{align}
\end{theorem}
At $r=0$, Theorem~\ref{thm:terminal-weighted-traces} approximates $H_n^{\bQ}/n$ and $G_{n+1}^{\bQ}$ by $H_0^{\bQ}/n$ and $G_1^{\bQ}$ in Eq.~\eqref{eq:terminal-zero-sample-traces}. At the polynomial and centered critical bulk scales of Remark~\ref{rem:polynomial-critical-scales}, both errors are of order $n^{-1}$ up to powers of $\log n$ for $\|\bQ\|\leq1$. For $r\geq1$, differentiation includes the dependence of the population traces on $w_0\lin(z)=-z/\omega(z/n)$. For each admissible $z$ and fixed $r$, applying Markov's inequality to Eq.~\eqref{eq:terminal-trace-bounds} at $p=\kappa\log n$ gives a pointwise high probability bound for both errors with probability $1-O(n^{-D})$ for any fixed $D>0$.

For our main application, regression with a spectral filter, the prediction variance conditional on the design involves a single trace weighted by $\bSigma$. Its resolvent representation uses divided differences of $zG_{n+1}^{\bI}(z,0)$, with a derivative when the two spectral parameters coincide. Theorem~\ref{thm:terminal-weighted-traces} controls both quantities. To control the concentration of its bias term, we need traces of two weighted resolvents. We study these mixed resolvents next.

Expanding $\bT_n$ over the observations allows us to approximate the mixed trace by a divided difference of terminal single traces, with prefactors determined by the two terminal path increments. Applying Theorems~\ref{thm:self-consistent-path-stability} and~\ref{thm:terminal-weighted-traces} replaces these quantities by their linear initial values; the initial resolvent identity then yields the divided-difference factor of $\omega$. In the next theorem, we bound the trace of weighted mixed resolvents around their identified deterministic equivalents. Define the weighted resolvent $\wt{\bR}_k(z,w):=\bSigma^{\half}\bR_k(z,w)\bSigma^{\half}$. For simplicity of the statement, we write $\wt{\bR}_n=\wt{\bR}_n(z,0)$ and $\wt{\bR}_0=\wt{\bR}_0(z,w_0\lin(z))$, with $\wt{\bR}_n'=\wt{\bR}_n(z',0)$ and $\wt{\bR}_0'=\wt{\bR}_0(z',w_0\lin(z'))$ at the second parameter. We also use the shorthand $\Lambda_n = \Lambda_n(z), \Lambda_n(z')=\Lambda_n'$ and suppress the dependence on $C_{\bx}, \nu, \kappa$ in its big-Oh notation. See its proof in Appendix~\ref{proof:terminal-mixed-resolvents}.

\begin{theorem}[Deterministic equivalents of mixed-resolvents]\label{thm:terminal-mixed-resolvents}
    Fix $z,z'\in\C_+\cup\R_{<0}$, $\kappa>0$, and nonnegative integers $r,r'$. Suppose Assumption~\ref{assmp:features} holds and Assumption~\ref{assmp:spectral-stability} holds at $z,z'$ with the same $\nu,\kappa$. For every deterministic $\bQ$ with $\|\bQ\|_F<\infty$ and every $2\leq p\leq\kappa\log n$, we have
    \begin{equation}\label{eq:terminal-mixed-resolvent-bound}
        \left\|\frac{\de^{r+r'}}{\de z^r\,\de(z')^{r'}}\Tr\!\left[\bQ\!\left(\wt{\bR}_n\wt{\bR}_n'-\frac{\omega(z/n)-\omega(z'/n)}{(z-z')/n}\wt{\bR}_0\wt{\bR}_0'\right)\right]\right\|_{L^p}=O_{r,r'}\left(\frac{p^{3/2}\|\bQ\|_F \cdot (\Lambda_n\vee\Lambda_n')^{26}}{n^{5/2}|z_-|^r|z'_-|^{r'}}\right).
    \end{equation}
    When $z'=z$ the divided difference becomes $\omega'(z/n)$, and the derivatives act separately.
\end{theorem}
We point out the initial explicit factors $\wt{\bR}_0(z,w_0\lin(z))=(\omega(z/n)/z)\bSigma(\bSigma-\omega(z/n)\bI)^{-1}$. In addition, the same bound above extends to either parameter in the lower half-plane, by reflecting the corresponding path and taking adjoints (cf.~Appendix~\ref{proof:terminal-mixed-resolvents}).

    \section{Discussion}\label{sec:discussion}
Our results establish nonasymptotic deterministic equivalents for spectral laws and prediction risk in a Hilbert-space setting. As a concluding remark, we outline unresolved questions, future research directions and potential applications that merit further investigation.

\begin{description}
    \item[Optimal local laws at local scale.]
    The mesoscopic estimates in Appendix~\ref{app:additional-bulk-laws} leave a gap to the optimal local scale known in the proportional regime~\citep{knowles2017anisotropic}. As discussed there, a natural target for future refinement is linear dependence on $\Lambda_n(z)$, up to logarithmic factors, in both the error bound and the sample-size requirement. This would allow spectral heights approaching the typical bulk spacing in our canonical examples. Achieving this requires sharper control along the backward interpolation path, and most likely a completely new interpolation scheme. It is also interesting to determine how far such estimates hold under the quadratic concentration assumption alone.

    \item[Natural spikes and the transition to the bulk.]
    As the sample size naturally grows, the decaying population spectrum in a Hilbert space can produce separated leading sample eigenvalues without an added finite-rank perturbation. This produces the natural ``spiked plus noise'' random matrix model. The question of characterizing their locations, eigendirection overlaps, and eventual transition into the bulk would give a natural analogue of the spiked Wigner and Wishart theories~\citep{baik2005phase,benaychgeorges2011eigenvalues}. The self-consistent equation and weighted resolvent approximations provide a starting point for identifying the corresponding shifts and separation criteria.

    \item[Central limit theorems and uncertainty quantification.]
    Central limit theorems could provide confidence intervals for spectral quantities and prediction risk. The martingale decompositions in Section~\ref{sec:trace-drift-fluctuations} suggest a route beyond classical results for linear spectral statistics~\citep{bai2004clt}. This requires identifying and consistently estimating the fluctuation variance, and making the drift and centering errors negligible at the fluctuation scale or correcting for them. However, our current bounds do not show that these errors are negligible at the fluctuation scale.

    \item[Lower-edge behavior and universality.]
    Section~\ref{sec:applications} uses a smallest-eigenvalue bound at the correct order. A sharper characterization would locate the smallest sample eigenvalue near the lower edge of the finite-$n$ law $\mu_n$ and determine its fluctuations. At a regular positive edge, a natural question is whether Tracy-Widom limits hold after appropriate centering and rescaling, as in classical covariance models~\citep{fan2022tracywidom}. This requires finer spectral control and potentially stronger distributional assumptions, and it would enable quantification of the extent to which a negative ridge pole can approach the empirical spectrum. However, the latter consideration is likely of greater theoretical than practical interest.

    \item[PCR and nonsmooth spectral filters.]
    The theory in our paper is restricted to analytic filters. Principal component regression (PCR) uses a discontinuous spectral cutoff and therefore falls outside our analytic framework. It would be interesting to extend our deterministic risk approximations to PCR and other nonanalytic spectral filters, including filters with limited smoothness. Indeed, we can nevertheless formally substitute more general filters into Proposition~\ref{prop:filter-real-risk}. It remains open to quantify the resulting approximation errors and how these errors depend on the filter's smoothness.

    \item[Applications in learning theory.]
    Finally, our effective risk formulas could guide source-dependent regularization and stopping rules, as suggested in Section~\ref{sec:gradient-flow}. In practice, developing data-driven procedures would require estimating the relevant spectral and signal quantities, with generalized cross-validation providing a related example~\citep{misiakiewicz2024nonasymptotic}. More broadly, these tools could potentially refine analyses of random features and neural scaling laws~\citep{defilippis2024dimensionfree,bahri2024explaining}. An important and interesting future direction is to extend the theory to learned representations.
\end{description}

    \section*{Statement of AI use}
We disclose the use of AI, and if used, acknowledge it in the following two categories.

\begin{enumerate}
    \item[(i)] Primary contributions: \textbf{Yes, by C.C.}, who uses GPT-5.6 to find the covariance tensorization (bilinear form of the Efron-Stein) to control the expected traces of weighted resolvents in the proof of Lemma~\ref{lem:one-step-covariance-bounds}, which supports the proof of Lemma~\ref{lem:one-step-G-bound}.
    \item[(ii)] Supporting contributions:
    \textbf{Yes, by C.C.}, who uses GPT-5.6 and GPT-6 models for producing numerical experiments, assisting in writing, suggesting literature, and auditing and optimizing existing proofs. GPT-6 suggested the logarithmic submean inequality on $L^p$ norms for bounding the contour integral across the spectral gap in Lemma~\ref{lem:filter-contour-comparison}.
\end{enumerate}
The authors took full responsibility for the contents of the paper.

    \section*{Acknowledgment}
The author thanks Theodor Misiakiewicz for helpful discussions. This work made use of the Illinois Campus Cluster, a computing resource that is operated by the Illinois Campus Cluster Program (ICCP) in conjunction with the National Center for Supercomputing Applications (NCSA) and which is supported by funds from the University of Illinois Urbana-Champaign.

    \bibliography{bib}

@article{hastie2022surprises,
  title         = {Surprises in high-dimensional ridgeless least squares interpolation},
  author        = {Hastie, Trevor and Montanari, Andrea and Rosset, Saharon and Tibshirani, Ryan J},
  journal       = {The Annals of Statistics},
  volume        = {50},
  number        = {2},
  pages         = {949--986},
  year          = {2022}
}

@article{bartlett2021deep,
  title         = {Deep learning: a statistical viewpoint},
  author        = {Bartlett, Peter L and Montanari, Andrea and Rakhlin, Alexander},
  journal       = {Acta Numerica},
  volume        = {30},
  pages         = {87--201},
  year          = {2021},
  publisher     = {Cambridge University Press}
}

@article{rudelson2013hanson,
  author        = {Rudelson, Mark and Vershynin, Roman},
  title         = {{Hanson-Wright} inequality and {sub-Gaussian} concentration},
  journal       = {Electronic Communications in Probability},
  volume        = {18},
  number        = {82},
  pages         = {1--9},
  year          = {2013}
}

@article{cheng2024dimensionfree,
  author        = {Cheng, Chen and Montanari, Andrea},
  title         = {Dimension free ridge regression},
  journal       = {The Annals of Statistics},
  volume        = {52},
  number        = {6},
  pages         = {2879--2912},
  year          = {2024}
}

@article{misiakiewicz2024nonasymptotic,
  author        = {Misiakiewicz, Theodor and Saeed, Basil},
  title         = {A non-asymptotic theory of kernel ridge regression: deterministic equivalents, test error, and {GCV} estimator},
  journal       = {arXiv preprint arXiv:2403.08938},
  year          = {2024},
  eprint        = {2403.08938},
  archivePrefix = {arXiv},
  primaryClass  = {stat.ML}
}

@article{hachem2007deterministic,
  author        = {Hachem, Walid and Loubaton, Philippe and Najim, Jamal},
  title         = {Deterministic equivalents for certain functionals of large random matrices},
  journal       = {The Annals of Applied Probability},
  volume        = {17},
  number        = {3},
  pages         = {875--930},
  year          = {2007}
}

@article{chouard2022quantitative,
  author        = {Chouard, Cl{\'e}ment},
  title         = {Quantitative deterministic equivalent of sample covariance matrices with a general dependence structure},
  journal       = {arXiv preprint arXiv:2211.13044},
  year          = {2022},
  eprint        = {2211.13044},
  archivePrefix = {arXiv},
  primaryClass  = {math.PR}
}

@article{caponnetto2007optimal,
  author        = {Caponnetto, Andrea and De Vito, Ernesto},
  title         = {Optimal Rates for the Regularized Least-Squares Algorithm},
  journal       = {Foundations of Computational Mathematics},
  volume        = {7},
  number        = {3},
  pages         = {331--368},
  year          = {2007}
}

@inproceedings{zhang2002effective,
  author        = {Zhang, Tong},
  title         = {Effective Dimension and Generalization of Kernel Learning},
  booktitle     = {Advances in Neural Information Processing Systems},
  volume        = {15},
  pages         = {454--461},
  year          = {2002}
}

@article{blanchard2020kernel,
  author        = {Blanchard, Gilles and M{\"u}cke, Nicole},
  title         = {Kernel Regression, Minimax Rates and Effective Dimensionality: Beyond the Regular Case},
  journal       = {Analysis and Applications},
  volume        = {18},
  number        = {4},
  pages         = {683--696},
  year          = {2020}
}

@book{simon2005trace,
  author        = {Simon, Barry},
  title         = {Trace Ideals and Their Applications},
  edition       = {2},
  series        = {Mathematical Surveys and Monographs},
  volume        = {120},
  publisher     = {American Mathematical Society},
  address       = {Providence, RI},
  year          = {2005}
}

@book{schatten1960norm,
  author        = {Schatten, Robert},
  title         = {Norm Ideals of Completely Continuous Operators},
  series        = {Ergebnisse der Mathematik und ihrer Grenzgebiete},
  volume        = {27},
  publisher     = {Springer-Verlag},
  address       = {Berlin},
  year          = {1960}
}

@book{gohberg1969introduction,
  author        = {Gohberg, Israel C. and Kre{\u\i}n, Mark G.},
  title         = {Introduction to the Theory of Linear Nonselfadjoint Operators in Hilbert Space},
  translator    = {Feinstein, A.},
  series        = {Translations of Mathematical Monographs},
  volume        = {18},
  publisher     = {American Mathematical Society},
  address       = {Providence, RI},
  year          = {1969}
}

@inproceedings{defilippis2024dimensionfree,
  author        = {Defilippis, Leonardo and Loureiro, Bruno and Misiakiewicz, Theodor},
  title         = {Dimension-free deterministic equivalents and scaling laws for random feature regression},
  booktitle     = {Advances in Neural Information Processing Systems},
  volume        = {37},
  pages         = {104630--104693},
  year          = {2024}
}

@article{adamczak2015note,
  author        = {Adamczak, Rados{\l}aw},
  title         = {A note on the {Hanson-Wright} inequality for random vectors with dependencies},
  journal       = {Electronic Communications in Probability},
  volume        = {20},
  number        = {72},
  pages         = {1--13},
  year          = {2015}
}

@book{boucheron2013concentration,
  author        = {Boucheron, St{\'e}phane and Lugosi, G{\'a}bor and Massart, Pascal},
  title         = {Concentration Inequalities: A Nonasymptotic Theory of Independence},
  publisher     = {Oxford University Press},
  address       = {Oxford},
  year          = {2013}
}

@book{bakry2014analysis,
  author        = {Bakry, Dominique and Gentil, Ivan and Ledoux, Michel},
  title         = {Analysis and Geometry of {Markov} Diffusion Operators},
  series        = {Grundlehren der mathematischen Wissenschaften},
  volume        = {348},
  publisher     = {Springer},
  address       = {Cham},
  year          = {2014}
}

@article{silverstein1995analysis,
  author        = {Silverstein, Jack W. and Choi, Sang-Il},
  title         = {Analysis of the Limiting Spectral Distribution of Large Dimensional Random Matrices},
  journal       = {Journal of Multivariate Analysis},
  volume        = {54},
  number        = {2},
  pages         = {295--309},
  year          = {1995}
}

@article{bai1999exact,
  author        = {Bai, Zhidong D. and Silverstein, Jack W.},
  title         = {Exact Separation of Eigenvalues of Large Dimensional Sample Covariance Matrices},
  journal       = {The Annals of Probability},
  volume        = {27},
  number        = {3},
  pages         = {1536--1555},
  year          = {1999}
}

@article{zhivotovskiy2022dimensionfree,
  author        = {Zhivotovskiy, Nikita},
  title         = {Dimension-Free Bounds for Sums of Independent Matrices and Simple Tensors via the Variational Principle},
  journal       = {Electronic Journal of Probability},
  volume        = {29},
  number        = {13},
  pages         = {1--28},
  year          = {2024}
}

@article{yaskov2014lower,
  author        = {Yaskov, Pavel},
  title         = {Lower Bounds on the Smallest Eigenvalue of a Sample Covariance Matrix},
  journal       = {Electronic Communications in Probability},
  volume        = {19},
  number        = {83},
  pages         = {1--10},
  year          = {2014}
}

@article{boucheron2005moment,
  author        = {Boucheron, St{\'e}phane and Bousquet, Olivier and Lugosi, G{\'a}bor and Massart, Pascal},
  title         = {Moment Inequalities for Functions of Independent Random Variables},
  journal       = {The Annals of Probability},
  volume        = {33},
  number        = {2},
  pages         = {514--560},
  year          = {2005}
}

@article{kac1974r,
  author        = {Kac, Israel S. and Krein, Mark G.},
  title         = {{R}-functions---analytic functions mapping the upper halfplane into itself},
  journal       = {American Mathematical Society Translations, Series 2},
  volume        = {103},
  pages         = {1--18},
  year          = {1974}
}

@article{biane1997free,
  author        = {Biane, Philippe},
  title         = {On the free convolution with a semi-circular distribution},
  journal       = {Indiana University Mathematics Journal},
  volume        = {46},
  number        = {3},
  pages         = {705--718},
  year          = {1997}
}

@book{durrett2019probability,
  author        = {Durrett, Rick},
  title         = {Probability: Theory and Examples},
  edition       = {5},
  publisher     = {Cambridge University Press},
  year          = {2019}
}

@article{stark2026quantitative,
  author        = {Stark, Thomas and D{\"o}rnemann, Nina},
  title         = {Quantitative Resolvent Calculus for Sample Covariance Matrices},
  journal       = {arXiv preprint arXiv:2610.33087},
  year          = {2026}
}

@book{hall1980martingale,
  author        = {Hall, Peter and Heyde, C. C.},
  title         = {Martingale Limit Theory and Its Application},
  publisher     = {Academic Press},
  year          = {1980}
}

@book{bai2010spectral,
  author    = {Bai, Zhidong and Silverstein, Jack W.},
  title     = {Spectral Analysis of Large Dimensional Random Matrices},
  edition   = {2},
  publisher = {Springer},
  year      = {2010}
}

@article{erdHos2012rigidity,
  author  = {Erd{\H{o}}s, L{\'a}szl{\'o} and Yau, Horng-Tzer and Yin, Jun},
  title   = {Rigidity of eigenvalues of generalized {Wigner} matrices},
  journal = {Advances in Mathematics},
  volume  = {229},
  number  = {3},
  pages   = {1435--1515},
  year    = {2012}
}

@article{bartlett2020benign,
  author  = {Bartlett, Peter L. and Long, Philip M. and Lugosi, G{\'a}bor and Tsigler, Alexander},
  title   = {Benign overfitting in linear regression},
  journal = {Proceedings of the National Academy of Sciences},
  volume  = {117},
  number  = {48},
  pages   = {30063--30070},
  year    = {2020}
}

@article{tsigler2023benign,
  author  = {Tsigler, Alexander and Bartlett, Peter L.},
  title   = {Benign overfitting in ridge regression},
  journal = {Journal of Machine Learning Research},
  volume  = {24},
  number  = {123},
  pages   = {1--76},
  year    = {2023}
}

@inproceedings{bordelon2020spectrum,
  author    = {Bordelon, Blake and Canatar, Abdulkadir and Pehlevan, Cengiz},
  title     = {Spectrum Dependent Learning Curves in Kernel Regression and Wide Neural Networks},
  booktitle = {Proceedings of the 37th International Conference on Machine Learning},
  series    = {Proceedings of Machine Learning Research},
  volume    = {119},
  pages     = {1024--1034},
  year      = {2020}
}

@article{bahri2024explaining,
  author  = {Bahri, Yasaman and Dyer, Ethan and Kaplan, Jared and Lee, Jaehoon and Sharma, Utkarsh},
  title   = {Explaining neural scaling laws},
  journal = {Proceedings of the National Academy of Sciences},
  volume  = {121},
  number  = {27},
  pages   = {e2311878121},
  year    = {2024}
}

@article{atanasov2026scaling,
  author  = {Atanasov, Alexander and Zavatone-Veth, Jacob A. and Pehlevan, Cengiz},
  title   = {Scaling and renormalization in high-dimensional regression},
  journal = {Journal of Statistical Mechanics: Theory and Experiment},
  volume  = {2026},
  number  = {4},
  pages   = {043404},
  year    = {2026}
}

@phdthesis{cheng2025thesis,
  author = {Cheng, Chen},
  title  = {High Dimensionality in Modern Machine Learning: A Random Matrix Theory Perspective},
  school = {Stanford University},
  year   = {2025}
}

@article{logerfo2008spectral,
  author = {Lo Gerfo, Lorenzo and Rosasco, Lorenzo and Odone, Francesca and De Vito, Ernesto and Verri, Alessandro},
  title = {Spectral algorithms for supervised learning},
  journal = {Neural Computation},
  volume = {20},
  number = {7},
  pages = {1873--1897},
  year = {2008}
}

@inproceedings{ali2019continuous,
  author = {Ali, Alnur and Kolter, J. Zico and Tibshirani, Ryan J.},
  title = {A continuous-time view of early stopping for least squares regression},
  booktitle = {Proceedings of the Twenty-Second International Conference on Artificial Intelligence and Statistics},
  series = {Proceedings of Machine Learning Research},
  volume = {89},
  pages = {1370--1378},
  year = {2019}
}

@inproceedings{wu2026riskcomparisons,
  author = {Wu, Jingfeng and Bartlett, Peter L. and Kakade, Sham M. and Lee, Jason D. and Yu, Bin},
  title = {Risk comparisons in linear regression: implicit regularization dominates explicit regularization (extended abstract)},
  booktitle = {Proceedings of Thirty Ninth Conference on Learning Theory},
  series = {Proceedings of Machine Learning Research},
  volume = {336},
  pages = {6849--6851},
  year = {2026}
}

@article{kim2026monotone,
  author = {Kim, Juno and Fu, Hengyu and Bartlett, Peter and Lee, Jason D. and Wu, Jingfeng},
  title = {Principal component regression dominates all monotone spectral filters for linear regression},
  journal = {arXiv preprint arXiv:2609.39440},
  year = {2026}
}

@article{knowles2017anisotropic,
  author  = {Knowles, Antti and Yin, Jun},
  title   = {Anisotropic local laws for random matrices},
  journal = {Probability Theory and Related Fields},
  volume  = {169},
  pages   = {257--352},
  year    = {2017}
}

@article{fan2026nonseparable,
  author        = {Fan, Zhou and Ma, Renyuan and Paquette, Elliot and Wang, Zhichao},
  title         = {Anisotropic local law for non-separable sample covariance matrices},
  journal       = {arXiv preprint arXiv:2602.17960},
  year          = {2026},
  eprint        = {2602.17960},
  archivePrefix = {arXiv},
  primaryClass  = {math.PR}
}

@article{bloemendal2014isotropic,
  author  = {Bloemendal, Alex and Erd{\H{o}}s, L{\'a}szl{\'o} and Knowles, Antti and Yau, Horng-Tzer and Yin, Jun},
  title   = {Isotropic local laws for sample covariance and generalized {Wigner} matrices},
  journal = {Electronic Journal of Probability},
  volume  = {19},
  pages   = {1--53},
  year    = {2014}
}

@article{marchenko1967distribution,
  author  = {Marchenko, Vladimir A. and Pastur, Leonid A.},
  title   = {Distribution of eigenvalues for some sets of random matrices},
  journal = {Mathematics of the USSR-Sbornik},
  volume  = {1},
  number  = {4},
  pages   = {457--483},
  year    = {1967}
}

@inproceedings{ho2020denoising,
  author = {Ho, Jonathan and Jain, Ajay and Abbeel, Pieter},
  title = {Denoising diffusion probabilistic models},
  booktitle = {Advances in Neural Information Processing Systems},
  volume = {33},
  pages = {6840--6851},
  year = {2020}
}

@book{wainwright2019highdimensional,
  author = {Wainwright, Martin J.},
  title = {High-Dimensional Statistics: A Non-Asymptotic Viewpoint},
  publisher = {Cambridge University Press},
  year = {2019}
}

@article{mei2022generalization,
  author = {Mei, Song and Montanari, Andrea},
  title = {The generalization error of random features regression: Precise asymptotics and the double descent curve},
  journal = {Communications on Pure and Applied Mathematics},
  volume = {75},
  number = {4},
  pages = {667--766},
  year = {2022}
}

@article{montanari2022interpolation,
  author = {Montanari, Andrea and Zhong, Yiqiao},
  title = {The interpolation phase transition in neural networks: Memorization and generalization under lazy training},
  journal = {The Annals of Statistics},
  volume = {50},
  number = {5},
  pages = {2816--2847},
  year = {2022}
}

@article{lu2025equivalence,
  author = {Lu, Yue M. and Yau, Horng-Tzer},
  title = {An equivalence principle for the spectrum of random inner-product kernel matrices with polynomial scalings},
  journal = {The Annals of Applied Probability},
  volume = {35},
  number = {4},
  pages = {2411--2470},
  year = {2025}
}

@article{liganshilinrisk,
  author = {Li, Yicheng and Gan, Weiye and Shi, Zuoqiang and Lin, Qian},
  title = {Generalization error curves for analytic spectral algorithms under power-law decay},
  journal = {Applied and Computational Harmonic Analysis},
  volume = {87},
  pages = {101920},
  year = {2027}
}

@inproceedings{deng2009imagenet,
  author = {Deng, Jia and Dong, Wei and Socher, Richard and Li, Li-Jia and Li, Kai and Fei-Fei, Li},
  title = {{ImageNet}: A large-scale hierarchical image database},
  booktitle = {Proceedings of the IEEE Conference on Computer Vision and Pattern Recognition},
  pages = {248--255},
  year = {2009}
}

@inproceedings{schuhmann2022laion,
  author = {Schuhmann, Christoph and Beaumont, Romain and Vencu, Richard and Gordon, Cade and Wightman, Ross and Cherti, Mehdi and Coombes, Theo and Katta, Aarush and Mullis, Clayton and Wortsman, Mitchell and Schramowski, Patrick and Kundurthy, Srivatsa and Crowson, Katherine and Schmidt, Ludwig and Kaczmarczyk, Robert and Jitsev, Jenia},
  title = {{LAION-5B}: An open large-scale dataset for training next generation image-text models},
  booktitle = {Advances in Neural Information Processing Systems},
  volume = {35},
  pages = {25278--25294},
  year = {2022}
}

@inproceedings{penedo2024fineweb,
  author = {Penedo, Guilherme and Kydl{\'\i}{\v c}ek, Hynek and Ben Allal, Loubna and Lozhkov, Anton and Mitchell, Margaret and Raffel, Colin and von Werra, Leandro and Wolf, Thomas},
  title = {The {FineWeb} datasets: Decanting the web for the finest text data at scale},
  booktitle = {Advances in Neural Information Processing Systems},
  volume = {37},
  pages = {30811--30849},
  year = {2024}
}

@inproceedings{rombach2022latent,
  author = {Rombach, Robin and Blattmann, Andreas and Lorenz, Dominik and Esser, Patrick and Ommer, Bj{\"o}rn},
  title = {High-resolution image synthesis with latent diffusion models},
  booktitle = {Proceedings of the IEEE/CVF Conference on Computer Vision and Pattern Recognition},
  pages = {10684--10695},
  year = {2022}
}

@inproceedings{brown2020fewshot,
  author = {Brown, Tom B. and Mann, Benjamin and Ryder, Nick and Subbiah, Melanie and Kaplan, Jared and Dhariwal, Prafulla and others},
  title = {Language models are few-shot learners},
  booktitle = {Advances in Neural Information Processing Systems},
  volume = {33},
  pages = {1877--1901},
  year = {2020}
}

@article{grattafiori2024llama,
  author = {Grattafiori, Aaron and Dubey, Abhimanyu and Jauhri, Abhinav and Pandey, Abhinav and Kadian, Abhishek and Al-Dahle, Ahmad and others},
  title = {The {Llama 3} herd of models},
  journal = {arXiv preprint arXiv:2407.21783},
  year = {2024}
}

@book{vandervaartwellner2023weak,
  author = {van der Vaart, Aad W. and Wellner, Jon A.},
  title = {Weak Convergence and Empirical Processes: With Applications to Statistics},
  edition = {2},
  publisher = {Springer},
  year = {2023}
}

@inproceedings{paquette20244+,
  title={4+3 phases of compute-optimal neural scaling laws},
  author={Paquette, Elliot and Paquette, Courtney and Xiao, Lechao and Pennington, Jeffrey},
  booktitle={Advances in Neural Information Processing Systems},
  volume={37},
  pages={16459--16537},
  year={2024}
}

@article{liang2020just,
  author = {Liang, Tengyuan and Rakhlin, Alexander},
  title = {Just interpolate: Kernel \textquotedblleft ridgeless\textquotedblright{} regression can generalize},
  journal = {The Annals of Statistics},
  volume = {48},
  number = {3},
  pages = {1329--1347},
  year = {2020}
}

@inproceedings{cheng2022memorize,
  author = {Cheng, Chen and Duchi, John and Kuditipudi, Rohith},
  title = {Memorize to generalize: On the necessity of interpolation in high dimensional linear regression},
  booktitle = {Proceedings of the Thirty-Fifth Conference on Learning Theory},
  series = {Proceedings of Machine Learning Research},
  volume = {178},
  pages = {5528--5560},
  year = {2022}
}

@inproceedings{cheng2026memorization,
  author = {Cheng, Chen and Barber, Rina Foygel},
  title = {Is memorization helpful or harmful? Prior information sets the threshold},
  booktitle = {Proceedings of the Thirty-Ninth Conference on Learning Theory},
  series = {Proceedings of Machine Learning Research},
  volume = {336},
  pages = {1399--1433},
  year = {2026}
}

@article{bercovici1999stable,
  author = {Bercovici, Hari and Pata, Vittorino},
  title = {Stable laws and domains of attraction in free probability theory},
  journal = {Annals of Mathematics},
  year = {1999},
  volume = {149},
  number = {3},
  pages = {1023--1060},
  note = {With an appendix by Philippe Biane}
}

@article{hasebe2020properties,
  author = {Hasebe, Takahiro and Simon, Thomas and Wang, Min},
  title = {Some properties of the free stable distributions},
  journal = {Annales de l'Institut Henri Poincaré, Probabilités et Statistiques},
  year = {2020},
  volume = {56},
  number = {1},
  pages = {296--325}
}

@article{tibshirani2022divided,
  author = {Tibshirani, Ryan J.},
  title = {Divided Differences, Falling Factorials, and Discrete Splines},
  journal = {Foundations and Trends in Machine Learning},
  year = {2022},
  volume = {15},
  number = {6},
  pages = {694--846}
}

@book{aupetit1991primer,
  author = {Aupetit, Bernard},
  title = {A Primer on Spectral Theory},
  series = {Universitext},
  publisher = {Springer-Verlag},
  address = {New York},
  year = {1991}
}

@article{chen2021hanson,
  author        = {Chen, Xiaohui and Yang, Yun},
  title         = {{Hanson-Wright} inequality in {Hilbert} spaces with application to {K}-means clustering for non-{Euclidean} data},
  journal       = {Bernoulli},
  volume        = {27},
  number        = {1},
  pages         = {586--614},
  year          = {2021}
}

@article{brailovskaya2024universality,
  author = {Brailovskaya, Tatiana and van Handel, Ramon},
  title = {Universality and sharp matrix concentration inequalities},
  journal = {Geometric and Functional Analysis},
  volume = {34},
  pages = {1734--1838},
  year = {2024}
}

@inproceedings{velikanov2024generalization,
  author = {Velikanov, Maksim and Panov, Maxim and Yarotsky, Dmitry},
  title = {Generalization error of spectral algorithms},
  booktitle = {The Twelfth International Conference on Learning Representations},
  year = {2024}
}

@article{ding2023global,
  author = {Ding, Xiucai and Wang, Zhenggang},
  title = {Global and local {CLTs} for linear spectral statistics of general sample covariance matrices when the dimension is much larger than the sample size with applications},
  journal = {IEEE Transactions on Information Theory},
  volume = {71},
  number = {8},
  pages = {6248--6296},
  year = {2025}
}

@article{dobriban2018highdimensional,
  author = {Dobriban, Edgar and Wager, Stefan},
  title = {High-dimensional asymptotics of prediction: Ridge regression and classification},
  journal = {The Annals of Statistics},
  volume = {46},
  number = {1},
  pages = {247--279},
  year = {2018}
}

@book{conway1978functions,
  author = {Conway, John B.},
  title = {Functions of One Complex Variable {I}},
  edition = {2},
  series = {Graduate Texts in Mathematics},
  volume = {11},
  publisher = {Springer},
  year = {1978}
}

@book{steinwart2008support,
  author = {Steinwart, Ingo and Christmann, Andreas},
  title = {Support Vector Machines},
  series = {Information Science and Statistics},
  publisher = {Springer},
  year = {2008}
}

@article{elkaroui2010spectrum,
  author = {{El Karoui}, Noureddine},
  title = {The spectrum of kernel random matrices},
  journal = {The Annals of Statistics},
  volume = {38},
  number = {1},
  pages = {1--50},
  year = {2010}
}

@article{cheng2013spectrum,
  author = {Cheng, Xiuyuan and Singer, Amit},
  title = {The spectrum of random inner-product kernel matrices},
  journal = {Random Matrices: Theory and Applications},
  volume = {2},
  number = {4},
  pages = {1350010},
  year = {2013}
}

@article{do2013spectrum,
  author = {Do, Yen and Vu, Van},
  title = {The spectrum of random kernel matrices: Universality results for rough and varying kernels},
  journal = {Random Matrices: Theory and Applications},
  volume = {2},
  number = {3},
  pages = {1350005},
  year = {2013}
}

@article{fan2019spectral,
  author = {Fan, Zhou and Montanari, Andrea},
  title = {The spectral norm of random inner-product kernel matrices},
  journal = {Probability Theory and Related Fields},
  volume = {173},
  pages = {27--85},
  year = {2019}
}

@article{misiakiewicz2022spectrum,
  author = {Misiakiewicz, Theodor},
  title = {Spectrum of inner-product kernel matrices in the polynomial regime and multiple descent phenomenon in kernel ridge regression},
  journal = {arXiv preprint arXiv:2204.10425},
  year = {2022},
  eprint = {2204.10425},
  archivePrefix = {arXiv},
  primaryClass = {math.ST}
}

@article{dubova2026universality,
  author = {Dubova, Sofiia and Lu, Yue M. and McKenna, Benjamin and Yau, Horng-Tzer},
  title = {Universality for the global spectrum of random inner-product kernel matrices in the polynomial regime},
  journal = {The Annals of Applied Probability},
  volume = {36},
  number = {5},
  pages = {4465--4515},
  year = {2026}
}

@article{li2026riskequivalence,
  author        = {Li, Yicheng and Cheng, Yuqian and Chen, Zhuo and Lin, Qian},
  title         = {Risk Equivalence between {RKHS} Regression and Sequence Models for Lipschitz Spectral Algorithms},
  journal       = {arXiv preprint arXiv:2609.08817},
  year          = {2026}
}

@article{baik2005phase,
  author = {Baik, Jinho and Ben Arous, G{\'e}rard and P{\'e}ch{\'e}, Sandrine},
  title = {Phase transition of the largest eigenvalue for nonnull complex sample covariance matrices},
  journal = {The Annals of Probability},
  volume = {33},
  number = {5},
  pages = {1643--1697},
  year = {2005}
}

@article{benaychgeorges2011eigenvalues,
  author = {Benaych-Georges, Florent and Nadakuditi, Raj Rao},
  title = {The eigenvalues and eigenvectors of finite, low rank perturbations of large random matrices},
  journal = {Advances in Mathematics},
  volume = {227},
  number = {1},
  pages = {494--521},
  year = {2011}
}

@article{bai2004clt,
  author = {Bai, Z. D. and Silverstein, Jack W.},
  title = {{CLT} for linear spectral statistics of large-dimensional sample covariance matrices},
  journal = {The Annals of Probability},
  volume = {32},
  number = {1A},
  pages = {553--605},
  year = {2004}
}

@article{fan2022tracywidom,
  author = {Fan, Zhou and Johnstone, Iain M.},
  title = {{Tracy-Widom} at each edge of real covariance and {MANOVA} estimators},
  journal = {The Annals of Applied Probability},
  volume = {32},
  number = {4},
  pages = {2967--3003},
  year = {2022}
}

@article{wishart1928generalised,
  author = {Wishart, John},
  title = {The generalised product moment distribution in samples from a normal multivariate population},
  journal = {Biometrika},
  volume = {20A},
  number = {1-2},
  pages = {32--52},
  year = {1928}
}

@article{apostol2000calculating,
  author = {Apostol, Tom M.},
  title = {Calculating Higher Derivatives of Inverses},
  journal = {The American Mathematical Monthly},
  volume = {107},
  number = {8},
  pages = {738--741},
  year = {2000}
}
    \bibliographystyle{abbrvnat}

    \newpage
    \appendix

    \section{Technical lemmas}
\label{sec:technical-lemmas}

We record the following moment form of the Efron-Stein inequality, adapted from~\citet[Theorem~2]{boucheron2005moment}.

\begin{lemma}[Moment Efron-Stein inequality] \label{lem:moment-efron-stein} Let $\bX_1,\ldots,\bX_N$ be independent random elements, let $\bX_1',\ldots,\bX_N'$ be independent copies that are also independent of $\bX_1,\ldots,\bX_N$, and let \begin{align*}
    Z&=f(\bX_1,\ldots,\bX_N),\\ Z^{(i)}&=f(\bX_1,\ldots,\bX_{i-1},\bX_i',\bX_{i+1},\ldots,\bX_N).
\end{align*}
For every $q\geq2$ such that the right-hand side is finite and some universal constant $C > 0$,
\begin{align}
    \norm{Z-\E Z}_{L^q} \leq C\sqrt q\cdot \prn{\sum_{i=1}^N\norm{Z-Z^{(i)}}_{L^q}^2}^{1/2}.
    \label{eq:moment-efron-stein}
\end{align}
The statement holds for both $Z \in \R$ and $Z \in \C$.
\end{lemma}
\begin{proof}
    For real-valued $Z$, applying~\citet[Theorem~2]{boucheron2005moment} to $Z$ and $-Z$ gives
    \begin{align*} \norm{(Z-\E Z)_+}_{L^q}&\leq\sqrt{2\kappa_q q}\cdot \norm{V_+^{1/2}}_{L^q}, & \norm{(Z-\E Z)_-}_{L^q}&\leq\sqrt{2\kappa_q q}\cdot \norm{V_-^{1/2}}_{L^q},
    \end{align*}
    where $\kappa_q=O(1)$ and
    \begin{align*}
        V_\pm:=\sum_{i=1}^N\E_{\bX_i'}\brk{(Z-Z^{(i)})_\pm^2}.
    \end{align*}
    Set $V:=\sum_{i=1}^N\E_{\bX_i'}\brk{|Z-Z^{(i)}|^2}$. Since $V_\pm\leq V$, we have
    \begin{align*}
        \norm{Z-\E Z}_{L^q}
        &=\prn{\norm{(Z-\E Z)_+}_{L^q}^q+\norm{(Z-\E Z)_-}_{L^q}^q}^{1/q} \leq 2^{1/q}\sqrt{2\kappa_q q}\cdot \norm{V^{1/2}}_{L^q}
        =O\prn{\sqrt q\,\norm{V^{1/2}}_{L^q}},
    \end{align*}
    where the big-Oh notation subsumes the universal bounds $2^{1/q}\leq\sqrt2$ and $\sup_{q\geq2}\kappa_q<\infty$. For complex-valued $Z$, applying the real bound to $\Re Z$ and $\Im Z$ and using Minkowski's inequality gives the same estimate. It remains to bound $\norm{V^{1/2}}_{L^q}$ in terms of the replacement moments. Using $q/2\geq1$, we write
    \begin{align*}
        \norm{V^{1/2}}_{L^q}^2
        &=\norm{\sum_{i=1}^N\E_{\bX_i'}\brk{|Z-Z^{(i)}|^2}}_{L^{q/2}} \stackrel{\mathrm{(i)}}{\leq}\sum_{i=1}^N\norm{\E_{\bX_i'}\brk{|Z-Z^{(i)}|^2}}_{L^{q/2}}
        \stackrel{\mathrm{(ii)}}{\leq}\sum_{i=1}^N\norm{Z-Z^{(i)}}_{L^q}^2.
    \end{align*}
    We make use of Minkowski's inequality in $L^{q/2}$ in (i), and Jensen's inequality in (ii). Taking square roots and substituting into the preceding estimate proves Eq.~\eqref{eq:moment-efron-stein}.
\end{proof}

The following covariance tensorization inequality is the bilinear form of the Efron-Stein inequality and follows from its standard martingale proof. See, e.g., \citet[Theorem~3.1]{boucheron2013concentration}.
\begin{lemma}[Covariance tensorization]\label{lem:covariance-tensorization} Let $\bX_1,\ldots,\bX_N$ be independent random elements, let $Z=Z(\bX_1,\ldots,\bX_N)$ and $W=W(\bX_1,\ldots,\bX_N)$ be complex-valued random variables in $L^2$, and let $Z^{(i)}$ and $W^{(i)}$ be obtained by replacing $\bX_i$ by an independent copy. Then
    \begin{align*}
        \left|\E\brk{(Z-\E Z)(W-\E W)}\right|
        \leq \sum_{i=1}^N\norm{Z-Z^{(i)}}_{L^2}\norm{W-W^{(i)}}_{L^2}.
    \end{align*}
\end{lemma}

\begin{proof} Let $\mc F_i=\sigma(\bX_1,\ldots,\bX_i)$, with $\mc F_0$ trivial, and let $D_iZ=\E[Z\mid\mc F_i]-\E[Z\mid\mc F_{i-1}]$, with $D_iW$ defined analogously. Orthogonality of martingale differences gives $\E[(Z-\E Z)(W-\E W)]=\sum_{i=1}^N\E[D_iZ\,D_iW]$. If $\bX_i'$ is the independent copy used in $Z^{(i)}$, conditioning first on $\mc F_i$ and then integrating over $\bX_i'$ and the future variables shows $D_iZ=\E[Z-Z^{(i)}\mid\mc F_i]$. Hence conditional Jensen's inequality gives $\norm{D_iZ}_{L^2}\leq\norm{Z-Z^{(i)}}_{L^2}$, and similarly for $W$. The conclusion follows from Cauchy-Schwarz, applied to each martingale-difference product. For complex-valued variables, the same argument applies to the bilinear product; alternatively, apply the real-valued result to the real and imaginary parts. \end{proof}

Let $s(z)$ be defined on $z \in \C_+$. The following lemma helps us to determine when $s(z)$ is indeed a Stieltjes transform of a probability measure, often from only $s(z)$ satisfying a certain fixed-point equation (or system of equations). $s(z)$ is known to be a Herglotz-Nevanlinna function if it is analytic on $\C_+$ and $\Im s(z) \geq 0$ for all $z \in \C_+$. We have the following representation lemma~\cite{kac1974r}.
\begin{lemma}[Herglotz representation] \label{lem:herglotz-nevanlinna}
    The function $s(z)$ is a Herglotz-Nevanlinna function if and only if there are unique numbers $a \in \R, b \geq 0$ and a (positive) Borel measure $\mu$ with $\int_{\R} \frac{1}{1+t^2} \de \mu(t)< \infty$ such that
    \begin{align*}
        s(z)=a+b z+\int_{-\infty}^{+\infty} \prn{\frac{1}{t-z}-\frac{t}{1+t^2}} \mu(\de t) .
    \end{align*}
\end{lemma}

\begin{lemma}[Stability of the effective spectral parameter equation for $\omega(z)$]\label{lem:population-scalar-stability}
    Let $\bSigma\succeq0$ be trace-class. For every $z\in\C_+$,
    \begin{align}
        s(z)^{-1} = -z+\sum_{j\geq1}\frac{\sigma_j}{1+n\sigma_js(z)}
    \end{align}
    has a unique solution $s(z)=s_n(z)\in\C_+$. It is the Stieltjes transform of a probability measure $\mu_n$ supported on $[0,\infty)$, and is equal to $-(n\omega(z/n))^{-1}$. It extends analytically across $(-\infty,0)$, where it is the unique positive solution of Eq.~\eqref{eq:population-scalar-equation}. For $z,z'\in\C_+\cup\R_{<0}$,
    \begin{align*}
        |s_n(z)|\leq |z_-|^{-1}\qquad \text{and} \qquad |s_n(z')-s_n(z)|\leq\frac{|z'-z|}{|z'_-| \cdot |z_-|}.
    \end{align*}

\end{lemma}
We defer the proof of Lemma~\ref{lem:population-scalar-stability} to Appendix~\ref{proof:population-scalar-stability}. The next lemma quantifies the spectral transfer from approximation errors of Stieltjes transforms to their probability measures over an interval $E$. For a probability measure $\mu$ with Stieltjes transform $s$, the Stieltjes inversion formula~\citep[Theorem~B.8]{bai2010spectral} recovers the mass of $E=[a,b]$ when $\mu(\{a,b\})=0$, giving
\begin{align*}
    \mu(E)=\lim_{\eta\downarrow0}\frac{1}{\pi}\int_E\Im s(t+\iu\eta)\,\de t.
\end{align*}
The following elementary smoothing inequality bounds $\mu(E)$ at finite $\eta$ with explicit error terms. A related estimate for eigenvalue counts near the spectral edge follows from the local law in the proof of~\citet[Lemma~6.1]{erdHos2012rigidity}. Additionally, let the inflated and deflated intervals by $E_h^+=[a-h, b+h]$ and $E^-=[a+h, b-h]$ for $0 < h < (b-a)/2$.
\begin{lemma}\label{lem:gram-interval-smoothing}
    Let $s(z)$ be the Stieltjes transform of a probability measure $\mu$ on $\R$, and $E=[a, b]$ for $a < b$. For $E_h^\pm$ defined above, it holds for any $0 < \eta < \pi h / 2$ that
    \begin{align*}
        \frac{\frac{1}{\pi}\int_{E_h^-} \Im s(t+ \iu \eta) \de t - \frac{\eta}{\pi h} }{1 - \frac{\eta}{\pi h}} \leq \mu(E) \leq \frac{\frac{1}{\pi}\int_{E_h^+} \Im s(t+ \iu \eta) \de t}{1 - \frac{2\eta}{\pi h} }.
    \end{align*}
\end{lemma}
\begin{proof}
    By Fubini,
    \begin{align*}
        \frac{1}{\pi} \int_E \Im s(t + \iu \eta) \de t = \int_{-\infty}^{+\infty} \prn{\int_E \frac{\eta}{\pi((x-t)^2 + \eta^2)} \de x} \mu(\de t).
    \end{align*}
    We use the fact that for any $t \in E$,
    \begin{align*}
        \int_{E_h^+} \frac{\eta}{\pi((x-t)^2 + \eta^2)} \de x & \geq 1 - \frac{1}{\pi} \int_{|x| \geq h} \frac{\eta}{x^2 + \eta^2} \de x \geq 1 - \frac{2}{\pi} \arctan \frac{\eta}{h} \geq 1 - \frac{2 \eta}{\pi h},
    \end{align*}
    and
    \begin{align*}
        \frac{1}{\pi} \int_{E_h^+} \Im s(t + \iu \eta) \de t \geq \int_{-\infty}^{+\infty} \prn{ 1 - \frac{2 \eta}{\pi h}} \ind \{t \in E\} \mu(\de t) = \prn{ 1 - \frac{2 \eta}{\pi h}} \mu (E).
    \end{align*}
    Similarly for $t \not \in E$,
    \begin{align*}
        \int_{E_h^-} \frac{\eta}{\pi((x-t)^2 + \eta^2)} \de x & \leq \frac{1}{\pi} \int_{x \geq h} \frac{\eta}{x^2 + \eta^2} \de x \leq \frac{\eta}{\pi h},
    \end{align*}
    and
    \begin{align*}
        \frac{1}{\pi} \int_{E_h^-} \Im s(t + \iu \eta) \de t \leq \int_{-\infty}^{+\infty} \prn{\frac{\eta}{\pi h} \ind \{t \not \in E\} + \ind \{t \in E\}} \mu(\de t) = \frac{\eta}{\pi h} (1 - \mu(E)) + \mu(E).
    \end{align*}
    The proof is complete.
\end{proof}

In Section~\ref{sec:spectral-filter-risk}, we bound the approximation error between a spectral filter’s prediction risk and its deterministic equivalent, under the hypothesis of a spectral event $\{\what{\sigma}_n \geq \eta_n\}$, where $n\what{\sigma}_n$ is the smallest eigenvalue of the unnormalized sample Gram matrix $\bX \bX^*$. Indeed, under the  quadratic form concentration Assumption~\ref{assmp:features} of Hanson-Wright, we can already derive a sufficiently informative estimate for the power law and critical law decay regimes, as stated in Lemma~\ref{lem:smallest-eigenvalue}. We understand $d=\infty$ and $\sigma_i=0$ for $i$ beyond $\bSigma$'s rank in what follows. See its proof in Appendix~\ref{proof:smallest-eigenvalue}.
\begin{lemma}[Smallest eigenvalue bound]\label{lem:smallest-eigenvalue}
    Under Assumption~\ref{assmp:features}, there exists $\const_\tau=\const_\tau(C_{\bx})>0$ such that, for every deterministic integer $m\geq0$ and $u\geq0$, with probability at least $1-2e^{-u}$,
    \begin{align}
        n\what{\sigma}_n\geq\tau_1(m)-\const_\tau\brc{\sqrt{(n+u)\tau_2(m)}+(n+u)\sigma_{m+1}}.\label{eq:smallest-eigenvalue-tail}
    \end{align}
    In particular, if $\tau_1(m)\geq2\const_\tau\brc{\sqrt{2n\tau_2(m)}+2n\sigma_{m+1}}$, then $\P\prn{\what{\sigma}_n\geq\tau_1(m)/(2n)}\geq1-2e^{-n}$.
\end{lemma}

\subsection{Proof of Lemma~\ref{lem:population-scalar-stability}}\label{proof:population-scalar-stability}
\paragraph{Part I: Unique analytic solution as the Stieltjes transform.} We first establish existence, uniqueness, and analyticity by applying a global inversion theorem. If $\bSigma=0$, then $s=-1/z$ is the unique analytic solution. Otherwise, substituting $u=-1/s$ transforms Eq.~\eqref{eq:population-scalar-equation} into $H(u)=z$, where
\begin{align*}
    H(u):=u+\sum_j\frac{\sigma_j u}{u-n\sigma_j}
    =u+\Tr\bSigma+n\sum_j\frac{\sigma_j^2}{u-n\sigma_j},\qquad u\in\C_+.
\end{align*}
To apply~\citet[Lemma~4, p.~710]{biane1997free}, define the probability measure $\rho:=\Tr(\bSigma^2)^{-1}\sum_j\sigma_j^2\delta_{n\sigma_j}$, which is well-defined because $0<\Tr(\bSigma^2)\leq\|\bSigma\|\Tr\bSigma<\infty$. We have
\begin{align*}
    H(u)-\Tr\bSigma=u+n\Tr(\bSigma^2)\int_\R\frac{\rho(\de x)}{u-x}.
\end{align*}
This is Biane's map with probability measure $\rho$ and variance parameter $n\Tr(\bSigma^2)$, whose domain is easily checked to be $\{u\in\C_+:\Im H(u)>0\}$. Applying that lemma shows that $H-\Tr\bSigma$ maps this domain conformally onto $\C_+$. Since translation by the real number $\Tr\bSigma$ preserves $\C_+$, the same holds for $H$, which therefore has an analytic inverse.
Thus $s_n(z)=-1/H^{-1}(z)\in\C_+$ solves Eq.~\eqref{eq:population-scalar-equation}. Every solution $s\in\C_+$ corresponds to $u=-1/s$ in this domain. Hence the solution is unique and $s_n(z)=-1/H^{-1}(z)$ is analytic on $\C_+$.

We next identify the solution as a Stieltjes transform and prove the bounds. Taking imaginary parts in Eq.~\eqref{eq:population-scalar-equation} implies $(\Im z)|s_n(z)|^2\leq\Im s_n(z)\leq|s_n(z)|$, and hence $|s_n(z)|\leq(\Im z)^{-1}$, since $s_n(z)\neq0$. For $y>2n\|\bSigma\|$, we have $|n\sigma_js_n(\iu y)|\leq1/2$, and hence $\left|\sum_j\sigma_j/(1+n\sigma_js_n(\iu y))\right|\leq2\Tr\bSigma$. Multiplying Eq.~\eqref{eq:population-scalar-equation} by $s_n(\iu y)$ at $z=\iu y$ therefore implies
\begin{align*}
    \left|\iu ys_n(\iu y)+1\right|
    =\left|s_n(\iu y)\sum_j\frac{\sigma_j}{1+n\sigma_js_n(\iu y)}\right|
    \leq\frac{2\Tr\bSigma}{y} \longrightarrow 0\qquad \text{as} \qquad y\to+\infty.
\end{align*}
Thus $ys_n(\iu y)\to\iu$ as $y\to\infty$, with $n$ and $\bSigma$ fixed. Applying Lemma~\ref{lem:herglotz-nevanlinna} to the analytic map $s_n:\C_+\to\C_+$ yields a Herglotz representation. In this paragraph, let $a\in\R$, $b\geq0$, and $\mu$ denote its coefficients and measure. Taking imaginary parts at $z=\iu y$ and using the preceding limit implies
\begin{align*}
    y\Im s_n(\iu y)=by^2+\int_\R\frac{y^2}{t^2+y^2}\mu(\de t)\longrightarrow1 \qquad \text{as} \qquad y\to+\infty.
\end{align*}
Both terms on the right are nonnegative, so $b=0$. Since $y^2/(t^2+y^2)$ increases to $1$ as $y\to\infty$, applying monotone convergence then yields $\mu(\R)=1$. In particular,
\begin{align*}
    \left|\int_\R\frac{1}{t-\iu y}\mu(\de t)\right|\leq\frac1y\longrightarrow0 \qquad \text{as} \qquad y\to+\infty.
\end{align*}
We obtain by comparing the above display with $s_n(\iu y)\to0$ as $y\to\infty$ in the Herglotz representation that $a-\int_\R t(1+t^2)^{-1}\mu(\de t)=0$. This last integral is finite as $|t|/(1+t^2)\leq1/2$ and $\mu(\R)=1$. Consequently, we can establish $s_n(z)$ is a Stieltjes transform:
\begin{align*}
    s_n(z)=\int_\R\frac{1}{t-z}\mu(\de t),
\end{align*}
\paragraph{Part II: Analytic extension and quantitative bounds of $s_n(z)$.}
We next verify that $\mu$ has no negative support so that we can extend analytically to $(-\infty, 0)$ for $s_n(z)$. For a fixed $x<0$,
\begin{align*}
    1=(-x)s+\sum_j\frac{\sigma_js}{1+n\sigma_js},\qquad s>0,
\end{align*}
has a unique solution because the right-hand side is continuous and strictly increasing from $0$ to $+\infty$ as $s$ ranges over $(0,\infty)$. To apply the analytic inverse function theorem, we verify that the derivative of the right-hand side of Eq.~\eqref{eq:population-scalar-equation} is positive at this solution. Indeed,
\begin{align*}
    n\sum_j\frac{\sigma_j^2s^2}{(1+n\sigma_js)^2}\leq\sum_j\frac{\sigma_js}{1+n\sigma_js}=1+xs<1.
\end{align*}
Consequently, at this solution,
\begin{align*}
    \frac{\de}{\de s}\prn{-s^{-1}+\sum_j\frac{\sigma_j}{1+n\sigma_js}}
    &=\frac1{s^2}\prn{1-n\sum_j\frac{\sigma_j^2s^2}{(1+n\sigma_js)^2}}
    \geq-\frac{x}{s}>0.
\end{align*}
Applying the analytic inverse function theorem therefore supplies a local inverse with a positive imaginary part just above $x$. By uniqueness in $\C_+$, this inverse agrees with $s_n$. Hence, $s_n$ extends analytically and is real on $(-\infty,0)$. Applying Stieltjes inversion on compact negative intervals gives $\mu(( -\infty,0))=0$. The representation now yields the following bounds for every $z,z'\in\C_+\cup\R_{<0}$, since $\mu$ is a probability measure supported on $[0,\infty)$:
\begin{align*}
    |s_n(z)|\leq\frac1{|z_-|} \qquad \text{and} \qquad
    |s_n(z')-s_n(z)|\leq |z'-z|\int_{[0,\infty)}\frac{\mu(\de t)}{|t-z'||t-z|}\leq\frac{|z'-z|}{|z'_-|\,|z_-|}.
\end{align*}
Finally, substituting $s=-(n\omega(z/n))^{-1}$ into Eq.~\eqref{eq:population-scalar-equation} gives Eq.~\eqref{def:effective-spectral-parameter-self-consistent-equation}. Its uniqueness therefore identifies $s_n(z)=-(n\omega(z/n))^{-1}$ and completes the proof.

\subsection{Proof of Lemma~\ref{lem:smallest-eigenvalue}}\label{proof:smallest-eigenvalue}
We bound the least eigenvalue by projecting each $\bx_i$ onto the tail eigenspaces indexed by $j>m$, using the orthogonal projection $\bQ_m:=\bI-\sum_{j=1}^m\bu_j\bu_j^*$. Let $\bS\in\R^{n\times n}$ have entries $S_{ij}:=\ind\{i\ne j\}\langle\bQ_m\bx_i,\bQ_m\bx_j\rangle$. Since $\bX\bX^*\succeq\bX\bQ_m\bX^*=(\langle\bQ_m\bx_i,\bQ_m\bx_j\rangle)_{i,j\leq n}$, separating the diagonal and off-diagonal terms yields
\begin{align*}
    n\what{\sigma}_n=\lambda_{\min}(\bX\bX^*)
    \geq\min_{i\leq n}\|\bQ_m\bx_i\|^2-\|\bS\|
    \geq\tau_1(m)-\max_{i\leq n}\left|\|\bQ_m\bx_i\|^2-\tau_1(m)\right|-\|\bS\|.
\end{align*}
It suffices to control the diagonal fluctuations and $\|\bS\|$ separately. Applying Proposition~\ref{prop:hanson-wright} with $\bA=\bQ_m$ and taking a union bound yields control for the second term immediately, for $t>0$,
\begin{align*}
    \P\prn{\max_{i\leq n}\left|\|\bQ_m\bx_i\|^2-\tau_1(m)\right|>t}
    \leq2n\exp\brc{-\Omega_{C_{\bx}}\prn{\frac{t^2}{\tau_2(m)}\wedge\frac{t}{\sigma_{m+1}}}}.
\end{align*}

It remains to control $\norm{\bS}$. We will use the standard $\epsilon$-net argument, and start with the quadratic forms of $\bS$. We use the Hilbert-space Hanson-Wright inequality~\citep[Proposition~2.5]{chen2021hanson}, based on the decoupling argument in~\citep{rudelson2013hanson}. The vectors $\bQ_m\bx_i$ are independent and centered, with positive-definite, trace-class covariance $\bSigma_{>m}:=\bQ_m\bSigma\bQ_m$ on the closed span of the nonzero tail eigenvectors. Applying the rank-one case of Proposition~\ref{prop:hanson-wright} shows that, for every $\bh$ in this space, $\E e^{\langle\bh,\bQ_m\bx_i\rangle}\leq\exp\brc{O_{C_{\bx}}\prn{\bh^*\bSigma_{>m}\bh}}$, verifying the required sub-Gaussian assumption. For each unit $\bv\in\R^n$, apply the cited inequality with covariance proxy $\bSigma_{>m}$ and coefficient matrices $\pm\bv\bv^*$, whose operator and Frobenius norms equal one. Substituting $\|\bSigma_{>m}\|=\sigma_{m+1}$ and $\|\bSigma_{>m}\|_F^2=\tau_2(m)$ yields
\begin{align*}
    \P\prn{|\bv^*\bS\bv|>t}
    =\P\prn{\left|\sum_{i\ne j}v_iv_j\langle\bQ_m\bx_i,\bQ_m\bx_j\rangle\right|>t}
    \leq2\exp\brc{-\Omega_{C_{\bx}}\prn{\frac{t^2}{\tau_2(m)}\wedge\frac{t}{\sigma_{m+1}}}}.
\end{align*}

We close the proof by the standard net argument. Using a $1/4$-net with at most $9^n$ points and the union bound yields
\begin{align*}
    \P\prn{n\what{\sigma}_n<\tau_1(m)-3t}
    \leq2(n+9^n)\exp\brc{-\Omega_{C_{\bx}}\prn{\frac{t^2}{\tau_2(m)}\wedge\frac{t}{\sigma_{m+1}}}}.
\end{align*}
Taking $3t=\const_\tau\brc{\sqrt{(n+u)\tau_2(m)}+(n+u)\sigma_{m+1}}$ with $\const_\tau=\const_\tau(C_{\bx})$ sufficiently large absorbs the union-bound factor, since $\log(n+9^n)=O(n)$, and bounds the right-hand side by $2e^{-u}$, proving Eq.~\eqref{eq:smallest-eigenvalue-tail}. The final assertion follows by taking $u=n$.

    \section{Proofs for Section~\ref{sec:gram-bulk-law}: Spectral laws}\label{app:spectral-laws}

\subsection{Proof of Theorem~\ref{thm:gram-stieltjes-approx}}\label{proof:gram-bulk-law}
\paragraph{Part I: Stieltjes transform.} The proof follows from the backward path $\{w_k\}$ and relative-resolvent traces $\{H_k\}$ in Section~\ref{sec:interpolation-path}, which is made clear by Eqs.~\eqref{eq:H-and-empirical-Stieltjes-transform} and \eqref{eq:population-scalar-equation}. Indeed, the previous equations identify
\begin{align*}
    H_n(z,0)/n=1+z\what s_n(z) \qquad \text{and} \qquad H_0(z,w_0\lin)/n=1+zs_n(z).
\end{align*}
We then obtain directly from Theorem~\ref{thm:terminal-weighted-traces} with $r=0$ and $\bQ=\bI$ that
\begin{align*}
    \|\what s_n(z)-s_n(z)\|_{L^p}=\frac{\|H_n(z,0)-H_0(z,w_0\lin)\|_{L^p}}{n|z|}
    =O_{C_{\bx},\nu,\kappa}\prn{\frac{p^2\Lambda_n(z)^{23}}{n|z|}}.
\end{align*}

A slightly more refined analysis using the proof details can sharpen the exponent on $\Lambda_n(z)$, using Lemma~\ref{lem:trace-step-decomposition} and the subsequent Eq.~\eqref{eq:trace-path-H} to express $H_n(z,0)-H_0(z,w_0\lin)$ as the sum of the martingale and drift terms with indices $1\le k\le n$ and $0\le k\le n$, respectively. In the proof of Lemma~\ref{lem:trace-accumulated-drift}, Eqs.~\eqref{eq:martingale-maximal-moment} and~\eqref{eq:gram-positive-index-drift}, bound the martingale and drift sums in $L^p$ by $O(p^2\Lambda_n(z)^9)$ and $O(p^2\Lambda_n(z)^{17})$, respectively. For $k=0$, Theorem~\ref{thm:self-consistent-path-stability}.(iii) and Eqs.~\eqref{eq:path-perturbed-population-equation} and~\eqref{eq:path-spectral-perturbation-rate} in its proof imply for the initial drift $\beta_0^{H,\bI}=w_0-w_0\lin+(w_0/z)\varepsilon_n(z)=O(\Lambda_n(z)^{17})$, using $|w_0|\le n\Lambda_n(z)$ as in Part IV-i of that proof. Dividing the sum of these bounds by $n|z|$ yields the strengthened bound in Eq.~\eqref{eq:gram-stieltjes-approx}.
\paragraph{Part II: Interval mass.} We transfer the $L^p$ bound on Stieltjes transforms to their probability measure masses. Fix $E$ and $h$ as in Theorem~\ref{thm:gram-stieltjes-approx}.(ii). To compare interval masses by Lemma~\ref{lem:gram-interval-smoothing}, it suffices to control the integrated transform error $\Delta(E_h^+):=\pi^{-1}\int_{a-h}^{b+h}|\what s_n(t+\iu\eta)-s_n(t+\iu\eta)|\,\de t$. For $p=\kappa\log n$, applying Minkowski's integral inequality and then Part I at each $t+\iu\eta$, we have
\begin{align*}
    \|\Delta(E_h^+)\|_{L^p}
    &\le\frac1\pi\int_{a-h}^{b+h}\|\what s_n(t+\iu\eta)-s_n(t+\iu\eta)\|_{L^p}\,\de t
    \le O_{C_{\bx},\nu,\kappa}\prn{\frac{(\log n)^2}{n}\int_{a-h}^{b+h}\frac{\Lambda_n(t+\iu\eta)^{17}}{|t+\iu\eta|}\,\de t}.
\end{align*}
Applying Markov's inequality to $\Delta(E_h^+)$, it holds $\Delta(E_h^+)\le e^{D/\kappa}\|\Delta(E_h^+)\|_{L^p}$ with probability at least $1-n^{-D}$.

If $\eta/h>\pi/16$, Eq.~\eqref{eq:gram-interval-approx} follows from $|\what\mu_n(E)-\mu_n(E)|\le1$. Otherwise, we apply Lemma~\ref{lem:gram-interval-smoothing} to $\what\mu_n$ on $E$ with deflation/inflation width $h/2$, which implies
\begin{align*}
    \frac1\pi\int_{E_{h/2}^-}\Im\what s_n(t+\iu\eta)\,\de t-O(\eta/h)
    \le\what\mu_n(E)\le
    \frac1\pi\int_{E_{h/2}^+}\Im\what s_n(t+\iu\eta)\,\de t+O(\eta/h).
\end{align*}
The condition $\eta\le\pi h/16$ keeps both denominators in the lemma bounded away from zero, and each smoothed mass is at most one. To compare with $\mu_n$, we replace $\what{s}_n$ with $s_n$ at a cost of at most $\Delta(E_h^+)$ in either integral, since $E_{h/2}^\pm\subset E_h^+$, and use the same Poisson tail estimate as in the proof of the lemma. For $x\in E_h^-$, the two tails outside $E_{h/2}^-$ each begin at distance at least $h/2$ from $x$; and for $x\notin E_h^+$, the interval $E_{h/2}^+$ lies in one such tail. Each tail has mass at most $\pi^{-1}\int_{h/2}^\infty\eta t^{-2} \de t=2\eta/(\pi h)$. Integrating against $\mu_n$ therefore implies
\begin{align*}
    \left(1-\frac{4\eta}{\pi h}\right)\mu_n(E_h^-)
    &\le\frac1\pi\int_{E_{h/2}^-}\Im s_n(t+\iu\eta)\de t \qquad \text{and} \qquad
    \frac1\pi\int_{E_{h/2}^+}\Im s_n(t+\iu\eta)\de t \le\mu_n(E_h^+)+\frac{2\eta}{\pi h}.
\end{align*}
We therefore finish the proof by combining the above displays, concluding with $|\what\mu_n(E)-\mu_n(E)|\le\Delta(E_h^+)+O(\eta/h)+\mu_n(E_h^+\setminus E_h^-)$.

\subsection{Proof of Proposition~\ref{prop:canonical-bulk-laws}}\label{proof:canonical-bulk-laws}
We will prove the three main conclusions in the following three parts. In Part I we prove concentration of the empirical Stieltjes transform to its deterministic equivalent. In Part II we prove weak convergence of the spectral laws to $\mu_\star$ in power law and critical law, respectively. In Part III, we identify the lower edges of the limiting laws, and their density bounds.  Throughout, $z=z_n(\zeta)$ denotes the entry of Table~\ref{tab:spectral-stability-scales} for the spectrum under consideration.
\paragraph{Part I: Concentration of $\what{s}_n$ on compact subsets of $\C_+$.} We first apply Theorem~\ref{thm:gram-stieltjes-approx}.(i) at the two spectral scales in Table~\ref{tab:spectral-stability-scales} to prove Eq.~\eqref{eq:canonical-bulk-rate}. Once Assumption~\ref{assmp:spectral-stability} is verified, Eq.~\eqref{eq:gram-stieltjes-approx} then reads
\begin{align*}
    \|a_n(\what s_n(z)-s_n(z))\|_{L^p}=O_{C_{\bx},\nu,\kappa}\!\left(\frac{p^2}{n}\frac{a_n}{|z|}\Lambda_n(z)^{17}\right).
\end{align*}
It therefore suffices to verify the assumption and bound $(a_n/|z|)\Lambda_n(z)^{17}$ uniformly for $\zeta\in K$. To verify the assumption, we choose the largest admissible cutoff $m=\lfloor\rho(\lfloor n/2\rfloor-2)\rfloor=\Theta(n)$, so $\Lambda_n(z)\le\Lambda_n(z,m)$. For the power law, comparing the tail sum with its integral gives $a_n=n^{-\alpha}$ and $\tau_1(m)=\Theta(a_n)$ and $n\sigma_{m+1}=\Theta(a_n)$. For the critical law, $a_n=(1+\log n)^{-1-\beta}$ and
\begin{align*}
    b_n\sim\int_n^\infty\frac{\de x}{x(1+\log x)^{1+\beta}}=\frac{1}{\beta(1+\log n)^{\beta}},\qquad \tau_1(m)=\Theta(b_n),\qquad n\sigma_{m+1}=\Theta(a_n).
\end{align*}
Consequently $r_1(m)=\tau_1(m)/\sigma_{m+1}$ is $\Theta(n)$ and $\Theta(n \log n)$, respectively, so $\log(nr_1(m))=\Theta(\log n)$ in both cases. Thus,
\begin{align*}
    \Lambda_n(z,m)&=\frac{1}{|z_-|}\begin{cases}
    |z|+\Theta(a_n(\log n)^2),&\text{power law;}\\
    |z|+\Theta(b_n(\log n)^2),&\text{critical law.}
\end{cases}
\end{align*}
It remains to lower bound $|z_-|$ and upper bound $|z|$ at the two spectral scales. Since $K$ is compact in $\C_+$, $|z_-|\ge a_n\min_{\zeta\in K}\Im\zeta$ in both cases. The power law $z=a_n\zeta$ gives $|z|=\Theta(a_n)$, whereas the critical law $z=b_n+a_n\zeta$ implies $|z|=\Theta(b_n)$ with $b_n=\Theta(a_n\log n)$. Substituting these bounds into the preceding display gives, uniformly for $\zeta\in K$,
\begin{align*}
    \Lambda_n(z)\le\Lambda_n(z,m)=
    \begin{cases}
        O\!\left(\dfrac{a_n+a_n(\log n)^2}{a_n}\right)=O((\log n)^2),&\text{power law;} \vspace{.6em} \\
        O\!\left(\dfrac{b_n+b_n(\log n)^2}{a_n}\right)=O((\log n)^3),&\text{critical law.}
    \end{cases}
\end{align*}
Thus $\Lambda_n(z)^{16}=O((\log n)^{48})=o(n)$ in both cases, verifying the sample-size condition of Assumption~\ref{assmp:spectral-stability} holds for large enough $n$. For its distance condition, $|z_-|\ge a_n\min_{\zeta\in K}\Im\zeta$ guarantees $|z_-|\ge n^{-(1+\alpha)}$ for the power law since $a_n=n^{-\alpha}$, and $|z_-|\ge n^{-1}$ for the critical law since $a_n=(1+\log n)^{-1-\beta}$. The assumption holds uniformly on $K$ for all sufficiently large $n$ with $\nu=1+\alpha$ and $\nu=1$, respectively. Returning to the first display, the power law contributes $(\log n)^{2\times 17}=(\log n)^{34}$. For the critical law the additional factor $a_n/|z|=O((\log n)^{-1})$ gives $(\log n)^{3\times 17-1}=(\log n)^{50}$. Since these bounds are uniform on $K$, we complete the proof by taking the supremum of the $L^p$ norms.
\paragraph{Part II: Convergence of the rescaled spectral laws $\what{\mu}_n$, $\mu_n$ to $\mu_\star$.} We first show the rescaled fixed-point equations have a limit. Define $u=a_n s_n(z_n(\zeta))$ to be the Stieltjes transform at $\zeta$ of the law of $x/a_n$ or $(x-b_n)/a_n$, respectively, for $x\sim\mu_n$. To identify its limit, substitute $s_n(z_n(\zeta))=u/a_n$ and $a_n=n\sigma_n$ into Eq.~\eqref{eq:population-scalar-equation}. Since $z_n(\zeta)$ equals $a_n\zeta$ or $b_n + a_n\zeta$, dividing Eq.~\eqref{eq:population-scalar-equation} by $a_n$ yields the inverse equation $\zeta=\zeta_n(u)$, where
\begin{align*}
    \zeta_n(u)+\frac1u=\frac1n\sum_{j\ge1}\frac{1}{\sigma_n/\sigma_j+u}-\begin{cases}0,&\text{power law;}\\ b_n/a_n,&\text{critical law.}\end{cases}
\end{align*}
We show that $\zeta_n(u)\to\zeta_\star(u)$ locally uniformly on $\C_+$, where $\zeta_\star(u)$ is the spectral parameter in Eq.~\eqref{eq:canonical-bulk-inverses} at $s(\zeta)=u$.

For the power law, $\sigma_n/\sigma_j=(j/n)^{1+\alpha}$. For $u>0$ the sum is $1/u$ times the ``regularly varying spectrum'' trace sum in~\citep{cheng2024dimensionfree} at ridge parameter $\sigma_n/u$. The same Riemann sum gives, for $u\in\C_+$,
\begin{align*}
    \zeta_n(u)+\frac1u \to \int_0^\infty\frac{\de t}{t^{1+\alpha}+u}
    =\frac{\pi\,u^{-\alpha/(1+\alpha)}}{(1+\alpha)\sin(\pi/(1+\alpha))}.
\end{align*}
For the critical law, $b_n/a_n=n^{-1}\sum_{j>n}\sigma_j/\sigma_n$, and we have
\begin{align*}
    \zeta_n(u)+\frac1u
    &=\frac1n\sum_{j\le n}\frac1{\sigma_n/\sigma_j+u}
    +\frac1n\sum_{j>n}\left(\frac1{\sigma_n/\sigma_j+u}-\frac{\sigma_j}{\sigma_n}\right) \to \int_0^1\frac{\de t}{t+u}+\int_1^\infty\left(\frac1{t+u}-\frac1t\right)\de t=-\log u.
\end{align*}
Indeed, $\sigma_n/\sigma_{\lceil nt\rceil}\to t$ for $t>0$ and this ratio is at least $t$ for $t>1$. As the power law summand at $j=\lceil nt\rceil$ is $O((1+t^{1+\alpha})^{-1})$, while the critical law summand is bounded for $0<t\leq1$ and its centered tail is $O(t^{-2})$ for $t>1$, locally uniform convergence follows on compact subsets of $\C_+$ by dominated convergence.

We next show the convergence of the laws and identify $\mu_\star$. To this end, it suffices to show that every vague subsequential limit of the rescaled laws has mass one and that all such limits agree. Fix one such limit $\mu$, and let $s$ be its Stieltjes transform. For $u\in\C_+$ with $\zeta=\zeta_\star(u)\in\C_+$, we have $\zeta_n(u)\in\C_+$ for all large $n$, and thus applying Lemma~\ref{lem:population-scalar-stability} at $z_n(\zeta_n(u))$ implies $a_ns_n(z_n(\zeta_n(u)))=u$. Invoking the resolvent identity at $\zeta$ and $\zeta_n(u)$, we have
\begin{align*}
    |a_ns_n(z_n(\zeta))-u|\leq\frac{|\zeta-\zeta_n(u)|}{\Im\zeta \cdot \Im\zeta_n(u)}\to0.
\end{align*}
Along the chosen subsequence, vague convergence implies $a_ns_n(z_n(\zeta))\to s(\zeta)$ since $x\mapsto(x-\zeta)^{-1}$ is continuous and vanishes at infinity. Hence $s(\zeta_\star(u))=u$, proving uniqueness whenever a solution exists. To verify that $\mu$ is a probability measure, take $u=\iu/y$ and $\zeta=\zeta_\star(u)=\iu y+o(y)$ as $y\to\infty$. Taking imaginary parts as in the proof of Lemma~\ref{lem:population-scalar-stability} yields
\begin{align*}
    \frac{\Im\zeta}{y}=(\Im\zeta)\Im s(\zeta)=\int\frac{(\Im\zeta)^2}{|x-\zeta|^2}\,\mu(\de x)\to\mu(\R).
\end{align*}
The convergence follows by dominated convergence as $y \to \infty$. Since $\Im\zeta/y\to1$, we obtain $\mu(\R)=1$. Thus $s(\zeta)\in\C_+$ for every $\zeta\in\C_+$. Passing to the limit along the subsequence in $\zeta_n(a_ns_n(z_n(\zeta)))=\zeta$, it holds $\zeta_\star(s(\zeta))=\zeta$ by locally uniform convergence of $\zeta_n$. By the uniqueness of $s$ established above, all vague subsequential limits therefore have the same transform $s$, proving weak convergence to $\mu = \mu_\star$. The resolvent identity makes the rescaled transforms equicontinuous on each compact $K\subset\C_+$, so their convergence is locally uniform.

Finally, setting $p=2$ in Eq.~\eqref{eq:canonical-bulk-rate} yields a summable sequence of mean-square errors for each fixed $\zeta$. By invoking Chebyshev's inequality together with the Borel-Cantelli lemma, we deduce almost sure convergence of the empirical Stieltjes transforms to $s$ on a countable dense subset of $\C_+$. This is sufficient to conclude almost sure weak convergence via the Stieltjes continuity theorem~\citep[Theorem~B.9]{bai2010spectral}.

\paragraph{Part III: Lower edges and density bounds.} To obtain the edge formula and the density bound, we reparameterize by $u=s(\zeta)=re^{\iu\theta}$, $0<\theta<\pi$. In the power law case, we additionally set $\gamma=1/(1+\alpha)$ and $c=\pi\gamma/\sin(\pi\gamma)$. Then, take imaginary parts for both cases in Eq.~\eqref{eq:canonical-bulk-inverses}, and since $\Im\zeta>0$, we have
\begin{align*}
    0<\Im\zeta=\frac{\sin\theta}{r}-\begin{cases}c r^{\gamma-1}\sin((1-\gamma)\theta),&\text{power law;} \vspace{.6em} \\ \theta,&\text{critical law.}\end{cases}
    \qquad r<\begin{cases}\prn{\dfrac{\sin\theta}{c\sin((1-\gamma)\theta)}}^{1/\gamma},&\text{power law;} \vspace{.6em} \\ \dfrac{\sin\theta}{\theta},&\text{critical law.}\end{cases}
\end{align*}
Using $\sin((1-\gamma)\theta)\geq(1-\gamma)\sin\theta$ and $\sin\theta\leq\theta$, we have $\Im s(\zeta)\leq r\leq (c(1-\gamma))^{-1/\gamma}$ for the power law decay and $\Im s(\zeta) \leq 1$ for the critical law decay. By Stieltjes inversion on any interval $(a,b)$,
\begin{align*}
    \mu_\star((a,b))\leq\lim_{\eta\downarrow0}\frac1\pi\int_a^b\Im s(t+\iu\eta)\de t \leq \frac{b-a}{\pi} \cdot \begin{cases}
    (c(1-\gamma))^{-1/\gamma},&\text{power law;}\\ 1,&\text{critical law.}\end{cases}
\end{align*}
which proves the density bound on $\|\rho_\star\|_\infty$ for $\mu_\star$ .

Next we find $\kappa_\star$ by differentiating Eq.~\eqref{eq:canonical-bulk-inverses} on $u>0$. The edge equation (cf.~\citep[Section~4]{silverstein1995analysis}) is
\begin{align}\label{eq:canonical-edge-equation}
    0=u^2\zeta_\star'(u)=\begin{cases}
    1-c(1-\gamma)u^\gamma,&\text{power law;}\\
    1-u,&\text{critical law.}
\end{cases}
\end{align}
On the increasing branch from $u=0$, the inverse agrees with $s$ by Part II, so there is no mass below $\max_{u>0}\zeta_\star(u)$. Evaluating this maximum at $u=(c(1-\gamma))^{-1/\gamma}$ and $u=1$, respectively, we obtain the two edges $\kappa_\star = (c(1-\gamma))^{1/\gamma}/\alpha$ and $\kappa_\star = -1$ in Eq.~\eqref{eq:canonical-power-edge-density}.

It only remains to show that the support contains every $x>\kappa_\star$. Indeed, we will derive the exact densities below using $(r, \theta)$ parameterization. Taking equality in the bound for $r$ sets $\Im\zeta_\star(re^{\iu\theta})=0$ in the first display of Part III. Thus $x=\zeta_\star(re^{\iu\theta})$ is real and continuous in $\theta$. As $\theta\downarrow0$, $r$ tends to the positive critical point computed above (use $\sin\theta\sim\theta$), so $x\to\kappa_\star$. As $\theta\uparrow\pi$, $r\to0$ and $x\sim1/r\to\infty$. The intermediate value theorem therefore gives such a $\theta$ for every $x>\kappa_\star$. Since $\zeta_\star'(re^{\iu\theta})\ne0$ by the derivative formula, its local inverse at $x$ takes values in $\C_+$ near $x$ and agrees with $s$ there for $\Im\zeta>0$ by Part II. Hence $s(x+\iu0)=re^{\iu\theta}$. Stieltjes inversion yields $\rho_\star(x)=r\sin\theta/\pi>0$, proving that the support is $[\kappa_\star,\infty)$.

\subsection{Mesoscopic scales and the generalized Marchenko-Pastur law}\label{app:additional-bulk-laws}
We provide two additional secondary implications of Theorem~\ref{thm:gram-stieltjes-approx} below.
\paragraph{Extension I: Mesoscopic scale.} To set the stage, we write $h=\Im\zeta=\Im z/a_n$ for the spectral height in units of the bulk width $a_n$. The results presented in Proposition~\ref{prop:canonical-bulk-laws} treat the macroscopic scale $h=\Theta(1)$. At a fixed interior bulk energy, typical finer mesoscopic results concern heights satisfy $h\to0$ and $nh\to\infty$. In fact, our Theorem~\ref{thm:gram-stieltjes-approx} can already give such a refinement at $h=n^{-\delta}$ with $0<\delta<1/17$. To our knowledge, such a dimension-free Stieltjes transform estimate was previously unavailable for the infinite-rank power and critical laws under Assumption~\ref{assmp:features}. Related dimension-free estimates for ridge regression concern prediction risk at nonnegative real regularization parameters~\citep{cheng2024dimensionfree}.

We briefly discuss what is known to be optimal in the classical proportional regime, where $d/n$ is bounded above and away from zero, the population eigenvalues lie in a fixed compact interval, and the normalized Gram matrix $\bX\bX^*/n$ has bulk width of order one. At a regular bulk energy, the averaged local law gives a Stieltjes transform error of order $1/(nh)$ in the sense of stochastic domination down to $h=n^{-1+\varepsilon}$ for any fixed $\varepsilon>0$~\citep[Theorem~3.6]{knowles2017anisotropic}. This scale is optimal up to $n^\varepsilon$, since the mean bulk spacing is of order $1/n$. For identity covariance, the aspect ratio may vary polynomially~\citep{bloemendal2014isotropic}. In the proportional regime, the averaged local law also holds for dependent coordinates under quadratic-form concentration~\citep[Theorem~2.5]{fan2026nonseparable}. It is reasonable to expect such the similar height refinement can hold in our setup; however, neither result covers an infinite-rank covariance with decaying eigenvalues and a shrinking bulk width as in Proposition~\ref{prop:canonical-bulk-laws}.

In Theorem~\ref{thm:gram-stieltjes-approx}.(i), $\Lambda_n(z)=\Omega(1/h)$ and $(a_n/|z|)\Lambda_n(z)=O((\log n)^2/h)$, so one power in the error bound would give $O(p^2(\log n)^2/(nh))$. This provides an explicit target to refine, i.e.\ replacing the sample-size condition $n=\Omega(\Lambda_n(z)^{16})$ by $n=\Omega(\Lambda_n(z))$, which would admit $h=n^{-1+\varepsilon}$ for every fixed $\varepsilon>0$. Along the interpolation path, this calls for (a potentially completely new strategy for) stability in Theorem~\ref{thm:self-consistent-path-stability} at those heights and an unweighted trace error in Theorem~\ref{thm:terminal-weighted-traces} linear in $\Lambda_n(z)$ up to logarithmic factors. We leave the optimal local law open and give the mesoscopic bound below.
\begin{proposition}[Mesoscopic bulk scales]\label{prop:mesoscopic-bulk}
    Use the setup and scaling $z=z_n(\zeta)$ of Proposition~\ref{prop:canonical-bulk-laws}, with $\Re\zeta$ fixed in the interior of the limiting bulk and $h=\Im\zeta\in(0,1]$. If $z$ satisfies Assumption~\ref{assmp:spectral-stability}, then for all large $n$ and $2\leq p\leq\kappa\log n$,
    \begin{align}
        \|a_n(\what s_n(z)-s_n(z))\|_{L^p}=\frac{p^2}{nh^{17}}\begin{cases} O_{C_{\bx},\nu,\kappa,\Re\zeta,\alpha}((\log n)^{34}),&\text{power law;}\\ O_{C_{\bx},\nu,\kappa,\Re\zeta,\beta}((\log n)^{50}),&\text{critical law.}\end{cases}\label{eq:bulk-mesoscopic-estimate}
    \end{align}
\end{proposition}
This implies mesoscopic refinement, as for $h=n^{-\delta}$ with fixed $0<\delta<1/17$, the assumption holds for all large $n$ with $\nu>\alpha+\delta$ or $\nu>\delta$, respectively, and the bound is $o(1)$ uniformly in $p$.
\paragraph{Extension II: Proportional regime.} In the main paper, our examples focus on the less well-understood settings (power law and critical law decay) involving infinitely many population eigenvalues. Still, it is useful to illustrate how our results apply to the standard proportional framework, which follows as a more straightforward application of the same theorem. Suppose $d/n$ stays within a fixed compact subset of $(0,+\infty)$, and let $\pi_n=d^{-1}\sum_{j=1}^d\delta_{\sigma_j}$ be the population spectral law. At $z=n\zeta$ with $\zeta$ in a fixed compact set $K\subset\C_+$, simply choosing the cutoff $m=1$ verifies Assumption~\ref{assmp:spectral-stability}. The finite-$n$ comparison law for $\bX\bX^*/n$ is the generalized Marchenko-Pastur law on the Gram scale determined by $d/n$ and $\pi_n$. By Eq.~\eqref{eq:population-scalar-equation}, its Stieltjes transform $u=ns_n(n\zeta)$ satisfies
\begin{align}
    \zeta=-\frac1u+\frac dn\int\frac{x}{1+xu}\pi_n(\de x).\label{eq:proportional-finite-n-inverse}
\end{align}
These results are standardly known. At fixed $\zeta\in\C_+$, sub-Gaussian concentration for all $1$-Lipschitz functions of the full data matrix in Frobenius norm gives an almost sure $O(\sqrt{\log n}/n)$ transform error for the same finite-$n$ law~\citep[Corollary~2.5]{chouard2022quantitative}. To see how our result applies, under Assumption~\ref{assmp:features}, which allows independent sub-Gaussian coordinates or convex concentration of each standardized observation, Proposition~\ref{prop:proportional-gram-law} applies Theorem~\ref{thm:gram-stieltjes-approx} directly to recover this finite-$n$ law at the $1/n$ scale up to logarithmic factors, with an explicit $L^p$ bound on $K$. If $d/n\to\phi$ and $\pi_n\Rightarrow\pi$, the limiting law is the classical Marchenko-Pastur law for identity covariance~\cite{marchenko1967distribution} and the generalized Marchenko-Pastur law for population spectrum $\pi$~\cite{silverstein1995analysis}. However, a quantitative rate to this limit would also depend on how fast $d/n$ and $\pi_n$ converge, which is not covered by the statement of the result itself in what follows.

\begin{proposition}[Generalized Marchenko-Pastur law]\label{prop:proportional-gram-law}
    In the proportional setup above, suppose Assumption~\ref{assmp:features} holds and fix $\kappa>0$. For all large $n$ and $2\leq p\leq\kappa\log n$,
    \begin{align*}
        \sup_{\zeta\in K}\|n\what s_n(n\zeta)-ns_n(n\zeta)\|_{L^p}=O_{C_{\bx},\kappa,K}\prn{\frac{p^2(\log n)^{34}}{n}}.
    \end{align*}
    If $d/n\to\phi\in(0,\infty)$ and $\pi_n$ converges weakly to $\pi$, the empirical law of $\bX\bX^*/n$ converges almost surely to the generalized Marchenko-Pastur law. Its Stieltjes transform uniquely satisfies in $\C_+$ that
    \begin{align}
        \zeta=-\frac1s+\phi\int\frac{x}{1+xs}\pi(\de x).\label{eq:companion-mp-limit}
    \end{align}
\end{proposition}

\subsubsection{Proof of Proposition~\ref{prop:mesoscopic-bulk}}
The proof is a simple application of Theorem~\ref{thm:gram-stieltjes-approx}.(i), using the cutoff and the bound on $\Lambda_n(z,m)$ from Part I of the proof of Proposition~\ref{prop:canonical-bulk-laws}. In the power law case, $\Re\zeta>\kappa_\star>0$ by Eq.~\eqref{eq:canonical-power-edge-density}, so $|z|=a_n|\zeta|=\Theta(a_n)$. In the critical law case, $z=b_n+a_n\zeta$ and $a_n/b_n\to0$ give $|z|=\Theta(b_n)$, uniformly for $0<h\leq1$. In both cases $\Re z>0$ for large $n$, so $|z_-|=a_n h$. The same bound on $\Lambda_n(z,m)$ therefore acquires a factor $h^{-1}$. Applying Theorem~\ref{thm:gram-stieltjes-approx}.(i) gives Eq.~\eqref{eq:bulk-mesoscopic-estimate}, with the factor $h^{-17}$ and the same logarithmic powers as Eq.~\eqref{eq:canonical-bulk-rate}.

\subsubsection{Proof of Proposition~\ref{prop:proportional-gram-law}}

We apply Theorem~\ref{thm:gram-stieltjes-approx}.(i) at $z=n\zeta$ by taking $m=1$, which is admissible for all large $n$ since $d/n$ is bounded away from zero. The normalization $\|\bSigma\|=1$ gives $\sigma_2\leq1$, $\tau_1(1)\leq d$ and $r_1(1)\leq d$. Uniformly for $\zeta\in K$, $|z|=O_K(n)$ and $|z_-|\geq n\min_{\zeta\in K}\Im\zeta$. Hence, when $\sigma_2>0$,
\begin{align*}
    \Lambda_n(n\zeta)\leq\Lambda_n(n\zeta,1)\leq\frac{n\sup_{\zeta\in K}|\zeta|+n+d\log n\log(nd)}{n\min_{\zeta\in K}\Im\zeta}=O_K((\log n)^2).
\end{align*}
Here and in the statement, the implicit constant may depend on the fixed range of $d/n$. If $\sigma_2=0$, the defining convention gives $\Lambda_n(n\zeta,1)=|z|/|z_-|=O_K(1)$. Thus $\Lambda_n(n\zeta)^{16}=o(n)$ and $|z_-|\geq n^{-1}$ for all large $n$, verifying Assumption~\ref{assmp:spectral-stability} with $\nu=1$. The empirical transform of $\bX\bX^*/n$ is $n\what s_n(n\zeta)$. Multiplying Eq.~\eqref{eq:gram-stieltjes-approx} in Theorem~\ref{thm:gram-stieltjes-approx}.(i) by $n$ gives the stated bound since $|\zeta|$ is bounded away from zero on $K$.

To identify the limit, we apply the inverse comparison in Part II of the proof of Proposition~\ref{prop:canonical-bulk-laws}. For $u=ns_n(n\zeta)$, the finite-$n$ inverse equation is Eq.~\eqref{eq:proportional-finite-n-inverse}. Its right-hand side converges locally uniformly in $u\in\C_+$ to the right-hand side of Eq.~\eqref{eq:companion-mp-limit}, since $d/n\to\phi$, $\pi_n$ converges weakly to $\pi$. At $u=\iu/y$ the limiting right-hand side is $\iu y+O(1)$ as $y\to\infty$, since $|x/(1+x\iu/y)|\leq1$ on $[0,1]$. These are precisely the two inputs to the inverse comparison and mass-one calculation in Part II. Applying them with $z_n(\zeta)=n\zeta$ shows that the deterministic laws converge weakly to the probability law whose transform is the unique upper-half-plane solution of Eq.~\eqref{eq:companion-mp-limit}. At $p=2$ the displayed rate in Proposition~\ref{prop:proportional-gram-law} gives summable mean-square errors for every fixed $\zeta$. The same countable-set argument as in Part II therefore gives almost sure weak convergence of the empirical laws.

For $\pi=\delta_1$, Eq.~\eqref{eq:companion-mp-limit} reduces to $\zeta s^2+(\zeta+1-\phi)s+1=0$, the Marchenko-Pastur equation for the $n\times n$ Gram matrix~\citep{silverstein1995analysis}.

    \section{Proofs for Section~\ref{sec:spectral-filter-risk}: Risk of spectral filters}\label{app:spectral-filters}

\subsection{Proof of Lemma~\ref{lem:conditional-risk-resolvents}}\label{proof:conditional-risk-resolvents}
Starting from the conditional variance in Eq.~\eqref{eq:conditional-filter-risks}, we have
\begin{align*}
    \sampleVariance[\upphi]
    &=\frac{\tau^2}{n}\Tr\prn{\bSigma\what{\bSigma}\upphi(\what{\bSigma})^2}\\
    &\stackrel{\mathrm{(i)}}{=}\tau^2\iint_{\Gamma^2}\Tr\prn{\bSigma\bT_n(\bT_n-z\bI)^{-1}(\bT_n-z'\bI)^{-1}}\,\alpha_{\upphi}(\de z)\alpha_{\upphi}(\de z')\\
    &\stackrel{\mathrm{(ii)}}{=}\tau^2\iint_{\Gamma^2}\frac{z\Tr\prn{\bSigma(\bT_n-z\bI)^{-1}}-z'\Tr\prn{\bSigma(\bT_n-z'\bI)^{-1}}}{z-z'}\,\alpha_{\upphi}(\de z)\alpha_{\upphi}(\de z'),
\end{align*}
where in (i), we substitute Eq.~\eqref{eq:analytic-filter-representation} for the squared filter factors with a double integral and use $\what{\bSigma}=\bT_n/n$, and in (ii) we apply the resolvent identity at $z\ne z'$ and use $\bT_n(\bT_n-z\bI)^{-1}=\bI+z(\bT_n-z\bI)^{-1}$. At $z'=z$ the quotient equals the derivative of $z\Tr(\bSigma(\bT_n-z\bI)^{-1})$, which can be verified directly at step (ii). Similarly, for the conditional squared bias in Eq.~\eqref{eq:conditional-filter-risks}, we have
\begin{align*}
    \sampleBias[\upphi] = \iint_{\Gamma^2}zz'\bbeta^*(\bT_n-z\bI)^{-1}\bSigma(\bT_n-z'\bI)^{-1}\bbeta \alpha_{\upphi}(\de z)\alpha_{\upphi}(\de z'),
\end{align*}
where we substitute $\uppsi(\what{\bSigma})=-\int_\Gamma z(\bT_n-z\bI)^{-1} \alpha_{\upphi}(\de z)$ from Eq.~\eqref{eq:analytic-filter-representation} on both sides of $\bSigma$. Since $0\in\Omega$, Cauchy's formula at zero covers $\ker\what{\bSigma}$ as well.

We implicitly used the dominated convergence theorem to interchange the double integral and the trace, since the resolvent kernels are bounded on $\Gamma^2$ and $\tvnorms{\alpha_{\upphi}}^2\Tr\bSigma<\infty$.

\subsection{Proof of Theorem~\ref{thm:analytic-filter-risk}}\label{proof:analytic-filter-risk}
Our goal is to replace the two integrands in Lemma~\ref{lem:conditional-risk-resolvents} with their deterministic equivalents in Theorems~\ref{thm:terminal-weighted-traces} and~\ref{thm:terminal-mixed-resolvents}, respectively. We then integrate the errors over $\Gamma^2$. Accordingly, we advise readers to refer to Section~\ref{sec:interpolation-path} first, as the present proof relies on the path-interpolation estimates established there for the weighted resolvent trace and the mixed-resolvent trace. Throughout, we adopt the notation introduced in that section. We begin by computing the differences $\sampleVariance[\upphi] - \effectiveVariance[\upphi]$ and $\sampleBias[\upphi]-\effectiveBias[\upphi]$ in a form to which we can apply Theorems~\ref{thm:terminal-weighted-traces} and~\ref{thm:terminal-mixed-resolvents} of Section~\ref{sec:interpolation-path}. We first work on the event $\{\what{\sigma}_n\ge\eta_n\}$.

For the variance terms, we will take $\bQ=\bI$ in Theorem~\ref{thm:terminal-weighted-traces}. Recall that $G_{n+1}(z)=G_{n+1}(z,0)=\Tr(\bSigma(\bT_n-z\bI)^{-1})$ and $G_1(z)=G_1(z,w_0\lin(z))$ at the two endpoints. Combining Eqs.~\eqref{def:effective-spectral-parameter-self-consistent-equation} and~\eqref{eq:terminal-zero-sample-traces}, we have
\begin{align*}
    zG_1(z)=\sum_j\frac{\sigma_j\omega(z/n)}{\sigma_j-\omega(z/n)}=\wt{\omega}(z)-z.
\end{align*}
Therefore, denote by $\chi^G(z):=z\big(G_{n+1}(z,0)-G_1(z,w_0\lin(z))\big)$ and its divided difference by $\chi^G[z,z']$, we can write
\begin{align*}
    \sampleVariance[\upphi]-\effectiveVariance[\upphi]&=\tau^2\iint_{\Gamma^2}\chi^G[z,z']\,\alpha_{\upphi}(\de z)\alpha_{\upphi}(\de z').
\end{align*}

For bias terms, we take $\bQ=\btheta\btheta^*$ in Theorem~\ref{thm:terminal-mixed-resolvents}. Since $\bbeta=\bSigma^{1/2}\btheta$, multiplying the sample trace by $zz'$ recovers $zz'\bbeta^*(\bT_n-z\bI)^{-1}\bSigma(\bT_n-z'\bI)^{-1}\bbeta$. As $\omega(z/n)(\omega(z/n)\bI-\bSigma)^{-1}=(\bI-n\bSigma/\wt{\omega}(z))^{-1}$ at each of the two arguments $\wt{\bR}_0 \wt{\bR}_0'$, we recover the integrand of $\effectiveBias[\upphi]$ in Eq.~\eqref{eq:analytic-filter-risk}. Denoting the trace difference in Eq.~\eqref{eq:terminal-mixed-resolvent-bound} at $r=r'=0$ by (cf.~the opening of the proof of Theorem~\ref{thm:terminal-mixed-resolvents} in Appendix~\ref{proof:terminal-mixed-resolvents} for the same notation used)
\begin{align*}
    \chi^{\bQ}(z,z')=\Tr \prn{\bQ \left(\wt{\bR}_n\wt{\bR}_n'-\wt{\omega}[z,z']\wt{\bR}_0\wt{\bR}_0'\right)},
\end{align*}
With $\bQ=\btheta\btheta^*$, subtracting Eq.~\eqref{eq:analytic-filter-risk} from Lemma~\ref{lem:conditional-risk-resolvents} yields
\begin{align*}
    \sampleBias[\upphi]-\effectiveBias[\upphi]&=\iint_{\Gamma^2}zz'\chi^{\bQ}(z,z')\,\alpha_{\upphi}(\de z)\alpha_{\upphi}(\de z').
\end{align*}

It remains to bound these errors for each pair on $\Gamma^2$. As sketched in the paragraph before Theorem~\ref{thm:analytic-filter-risk}, when the contour crosses the gap $(0,n\eta_n)$, Theorems~\ref{thm:terminal-weighted-traces} and~\ref{thm:terminal-mixed-resolvents} do not apply directly. The following lemma extends the approximation errors to the entire contour. We defer its proof to Appendix~\ref{proof:filter-contour-comparison}.
\begin{lemma}\label{lem:filter-contour-comparison}
    Under the assumptions of Theorem~\ref{thm:analytic-filter-risk}, for $z,z'\in\Gamma$ and $\bQ=\btheta\btheta^*$,
    \begin{subequations}
        \begin{align}
            \left\|\ind\{\what{\sigma}_n\ge\eta_n\}\chi^G[z,z']\right\|_{L^p}
            &=O\left(\frac{p^2\Delta_\Gamma^{23}}n \cdot \prn{\frac{|z|\vee n\eta_n}{\dist(z,I_n)}+\frac{|z'|\vee n\eta_n}{\dist(z',I_n)}}\right),\label{eq:filter-projected-single}\\
            \left\|\ind\{\what{\sigma}_n\ge\eta_n\}zz'\chi^{\bQ}(z,z')\right\|_{L^p}
            &=O\left(\frac{p^{3/2}\Delta_\Gamma^{26}\|\btheta\|^2}{n^{5/2}} \cdot (|z|\vee n\eta_n) \cdot (|z'|\vee n\eta_n)\right), \label{eq:filter-projected-mixed}
        \end{align}
    \end{subequations}
    where $\chi^G[z,z]=(\chi^G)'(z)$.
\end{lemma}
Applying Minkowski's integral inequality to the two contour identities and invoking Lemma~\ref{lem:filter-contour-comparison}, Eqs.~\eqref{eq:filter-projected-single} and~\eqref{eq:filter-projected-mixed}, yields both bounds in Theorem~\ref{thm:analytic-filter-risk}. For the variance, the two summands have equal integrals by symmetry, and integration in the other variable gives $\tvnorms{\alpha_{\upphi}}$. For the squared bias, the double integral factors as $\big(\int_\Gamma (|z|\vee n\eta_n)\,|\alpha_{\upphi}|(\de z)\big)^2$.

\subsubsection{Proof of Lemma~\ref{lem:filter-contour-comparison}}\label{proof:filter-contour-comparison}
Our proof consists of three parts. To simplify the notation, we use $\mc{A}_n = \{\what{\sigma}_n \geq \eta_n\}$ in the proof.
\begin{enumerate}
    \item In Part I, for variance, we reduce to a bound on $\chi^G$ in a local disk centered at each $z \in \Gamma$.
    \item In Part II, we apply Theorems~\ref{thm:terminal-weighted-traces} and~\ref{thm:terminal-mixed-resolvents} where spectral stability holds, and use a crude bound at smaller spectral heights.
    \item In Part III, we upgrade the crude bounds by applying submean inequality, averaging the logarithm of each normalized error over the disks to recover the desired bounds.
\end{enumerate}

\paragraph{Part I: Reduction for $\chi^G$.} The key observation is that we can reduce the divided difference into contour integrals on the boundary of two disks. For $z,z'\in\Gamma$, applying Cauchy's formula on the union of the local disks centered at $z,z'$ with radii $\dist(z,I_n)/4,\dist(z',I_n)/4$, we immediately have
\begin{align*}
    \chi^G[z,z']=\frac{1}{2\pi\iu}\oint_{\Gamma_{z \mid z'} \cup \Gamma_{z' \mid z}}\frac{\chi^G(\zeta)}{(\zeta-z)(\zeta-z')} \de\zeta,
\end{align*}
where $\Gamma_{z \mid z'}$ and $\Gamma_{z' \mid z}$ are the two arcs on the boundary, centered at $z$ and $z'$ respectively. The integral is positively oriented on their union and $\chi^G$ is zero outside $\mc{A}_n$. This naturally recovers the derivative convention, as at $z=z'$, this is Cauchy's derivative formula. For $\zeta \in \Gamma_{z \mid z'}$, $|\zeta-z|=\dist(z,I_n)/4$ and $|\zeta-z'|\ge\dist(z',I_n)/4$. Hence
\begin{align*}
    \tfrac34\dist(z,I_n)\le\dist(\zeta,I_n)\le |\zeta-z'|+\dist(z',I_n)\le5|\zeta-z'|.
\end{align*}
Since $\Gamma_{z \mid z'}$ lies on the circle centered at $z$ with radius $\dist(z,I_n)/4$, its length is at most the circumference $2\pi \cdot \dist(z,I_n)/4 = O(\dist(z, I_n))$. We can apply Minkowski's inequality and the preceding denominator bound to conclude that
\begin{align*}
    \left\|\frac{1}{2\pi\iu}\int_{\Gamma_{z \mid z'}}\frac{\chi^G(\zeta)}{(\zeta-z)(\zeta-z')}\,\de\zeta\right\|_{L^p}
    =O\!\left(\frac{1}{\dist(z,I_n)}\sup_{|\zeta-z|\le\dist(z,I_n)/4}\|\chi^G(\zeta)\|_{L^p}\right).
\end{align*}
Similarly, we have the integral bound over $\Gamma_{z' \mid z}$. Combine those two bounds and compare with Eq.~\eqref{eq:filter-projected-single}, it is then clear that it suffices to prove the following local disk bound for the variance:
\begin{align}
    \sup_{|z'-z|\le\dist(z,I_n)/4}\big\|\ind_{\mc{A}_n} \chi^G(z')\big\|_{L^p}
    =O\!\left(\frac{p^2\Delta_\Gamma^{23}}{n}(|z|\vee n\eta_n)\right),\qquad z\in\Gamma.\label{eq:filter-single-on-disk}
\end{align}

\paragraph{Part II: Separate bounds inside and outside the horizontal strip.} We bound $\chi^G(z')$ on the local disk in Eq.~\eqref{eq:filter-single-on-disk} and $zz'\chi^{\bQ}(z,z')$ for $z,z'\in\Gamma$. Theorems~\ref{thm:terminal-weighted-traces} and~\ref{thm:terminal-mixed-resolvents} apply directly outside the horizontal strip
\begin{align*}
    \{\zeta\in\C:|\Im\zeta|<h_n\dist(z,I_n)\}
\end{align*}
associated with each $z\in\Gamma$. We will select the relative height $h_n>0$ later.

Our first step is to apply Theorems~\ref{thm:terminal-weighted-traces} and~\ref{thm:terminal-mixed-resolvents} outside the strip, by verifying Assumption~\ref{assmp:spectral-stability} uniformly there. Let $D_z:=\{\zeta\in\C:|\zeta-z|\le\dist(z,I_n)/4\}$. For variance, we take $z'\in D_z$ and $\zeta\in D_{z'}$. For squared bias we take $(\zeta,\zeta')\in D_z\times D_{z'}$ with $z,z'\in\Gamma$. Since distance to $I_n$ is $1$-Lipschitz, $\dist(z',I_n)\asymp\dist(z,I_n)$ when $z'\in D_z$. Thus in either case $|\zeta-z|=O(\dist(z,I_n))$ and $|\zeta|=O(|z|\vee n\eta_n)$. Outside the strip, $|\zeta_-|\ge|\Im\zeta|\ge h_n\dist(z,I_n)$.

To verify Assumption~\ref{assmp:spectral-stability}, we make use of the fact that $|\zeta_-|\Lambda_n(\zeta)-|\zeta|$ is the minimized spectral tail and hence takes the same value at $\zeta$ and at $-\dist(z,I_n)$. Since $r\mapsto\Lambda_n(-r)$ is nonincreasing for $r>0$, we have $\Lambda_n(-\dist(z,I_n))\le\Delta_\Gamma$ and $|z|\le\Delta_\Gamma\dist(z,I_n)$ from Theorem~\ref{thm:analytic-filter-risk}. For any $|\zeta-z|\le3\dist(z,I_n)/4$ and $\Im\zeta\ne0$, it holds
\begin{equation}\label{eq:filter-shifted-stability-controls}
    \begin{aligned}
        |\zeta_-|\Lambda_n(\zeta)
        &=|\zeta|+\dist(z,I_n)\cdot \prn{\Lambda_n(-\dist(z,I_n))-1} \le 2\Delta_\Gamma\dist(z,I_n).
    \end{aligned}
\end{equation}
To see that we can use Eq.~\eqref{eq:filter-shifted-stability-controls} for $\zeta \in D_{z'}$ with $z'\in D_z$. By $\dist(z',I_n)\le\dist(z,I_n)+|z'-z|\le5\dist(z,I_n)/4$, it holds for $\zeta\in D_{z'}$,
\begin{align*}
    |\zeta-z|\le |\zeta-z'|+|z'-z|\le\frac{\dist(z',I_n)}4+\frac{\dist(z,I_n)}4\le\left(\frac5{16}+\frac14\right)\dist(z,I_n)<\frac34\dist(z,I_n).
\end{align*}
Therefore, for the variance term and $\zeta \in D_z \cup D_{z'}$ outside the strip, Eq.~\eqref{eq:filter-shifted-stability-controls} implies $\Lambda_n(\zeta)\le 2\Delta_\Gamma/h_n$. We can choose $h_n=C\Delta_\Gamma n^{-1/16}$ and $\nu=C_\Gamma+1$, which yields
\begin{align*}
    \Lambda_n(\zeta)^{16}\le (2/C)^{16}n\qquad\text{and}\qquad |\zeta_-|\ge h_n\dist(\Gamma,I_n)\ge C\Delta_\Gamma n^{-1/16-C_\Gamma}\ge n^{-\nu}.
\end{align*}
Choosing $C$ sufficiently large as a function of $C_{\bx},\kappa,C_\Gamma$ ensures that Assumption~\ref{assmp:spectral-stability} holds uniformly for all sufficiently large $n$. Moreover, the requirement on $C_\Gamma$ in Theorem~\ref{thm:analytic-filter-risk} implies that $h_n=O(1/\log n)$.

With Assumption~\ref{assmp:spectral-stability} verified, we apply Theorem~\ref{thm:terminal-weighted-traces} with $\bQ=\bI$ and, if $|\Im\zeta'|\ge h_n\dist(z',I_n)$ as well, Theorem~\ref{thm:terminal-mixed-resolvents} with $\bQ=\btheta\btheta^*$ to bound $\chi^G$ and $\chi^{\bQ}$. At derivative order zero, using reflection for lower-half-plane arguments (for opposite half-planes, see the extension after Theorem~\ref{thm:terminal-mixed-resolvents}), we have results outside the strip,
\begin{equation}\label{eq:filter-direct-off-strip}
    \begin{aligned}
        \big\|\ind_{\mc{A}_n}\chi^G(\zeta)\big\|_{L^p}
        &=O\!\left(\frac{p^2\Delta_\Gamma^{23}}{n}(|z|\vee n\eta_n)\brc{1\vee\frac{\dist(z,I_n)}{|\Im\zeta|}}^{23}\right),\\
        \big\|\ind_{\mc{A}_n}\zeta\zeta'\chi^{\bQ}(\zeta,\zeta')\big\|_{L^p}
        &=O\!\left(\frac{p^{3/2}\Delta_\Gamma^{26}\|\btheta\|^2}{n^{5/2}}(|z|\vee n\eta_n)(|z'|\vee n\eta_n)\brc{1\vee\frac{\dist(z,I_n)}{|\Im\zeta|}\vee\frac{\dist(z',I_n)}{|\Im\zeta'|}}^{26}\right).
    \end{aligned}
\end{equation}
When either $\zeta$ or $\zeta'$ lies inside its strip, we will instead use the polynomial bounds in the following lemma, proved in Appendix~\ref{proof:filter-crude-on-disks}.

\begin{lemma}\label{lem:filter-crude-on-disks}
    Under the setup of Lemma~\ref{lem:filter-contour-comparison}, when $\mc{A}_n$ holds, the following bounds hold for some fixed $K$ and all sufficiently large $n$, uniformly in the contour points.
    \begin{subequations}
        \begin{align}
            \sup_{z'\in D_z}\sup_{\zeta\in D_{z'}}|\chi^G(\zeta)|&\le n^K\,\frac{|z|\vee n\eta_n}{n},\label{eq:filter-crude-single-on-disks}\\
            \sup_{(\zeta,\zeta')\in D_z\times D_{z'}}|\zeta\zeta'\chi^{\bQ}(\zeta,\zeta')|&\le n^K\,\frac{\|\btheta\|^2}{n^{5/2}}(|z|\vee n\eta_n)(|z'|\vee n\eta_n).\label{eq:filter-crude-mixed-on-disks}
        \end{align}
    \end{subequations}
\end{lemma}
Combining the crude and fine-grained estimates by the logarithmic submean inequality, we can close the proof in the next part.

\paragraph{Part III: Combining bounds by the submean inequality.}
For $z\in\Gamma$ and $z'\in D_z$, on $\mc{A}_n$ the scaled trace error is analytic samplewise on $D_{z'}$ by spectral separation and $\wt{\omega}=-1/s_n$, and the uniform bound in Lemma~\ref{lem:filter-crude-on-disks} makes $\ind_{\mc{A}_n}\chi^G$ analytic there in $L^p$ by Cauchy's formula. Applying the submean inequality to its log-norm~\citep[Lemma~3.4.6]{aupetit1991primer}, with $\log^+x:=\log(1\vee x)$, yields
\begin{align*}
    \log^+\frac{n\|\ind_{\mc{A}_n}\chi^G(z')\|_{L^p}}{p^2\Delta_\Gamma^{23}(|z|\vee n\eta_n)}
    \le \frac1{|D_{z'}|}\int_{D_{z'}}\log^+\frac{n\|\ind_{\mc{A}_n}\chi^G(\zeta)\|_{L^p}}{p^2\Delta_\Gamma^{23}(|z|\vee n\eta_n)} \de A(\zeta).
\end{align*}
Here $\de A$ denotes area measure, and the denominator depends only on the fixed contour point $z$.

It suffices to bound the average on the right by $O(1)$. By splitting at the strip boundary, i.e.\ $|\Im\zeta|=h_n\dist(z,I_n)$, we apply Lemma~\ref{lem:filter-crude-on-disks} inside the strip and Eq.~\eqref{eq:filter-direct-off-strip} outside to obtain
\begin{align*}
    &\frac1{|D_{z'}|}\int_{D_{z'}}\log^+\frac{n\|\ind_{\mc{A}_n}\chi^G(\zeta)\|_{L^p}}{p^2\Delta_\Gamma^{23}(|z|\vee n\eta_n)} \de A(\zeta)\le \underbrace{\vphantom{O\!\left(1+23\int_{h_n}^1\log(1/y)\,\de y\right)}O(h_n K\log n)}_{\text{inside the strip}}
    +\underbrace{O\!\left(1+23\int_{h_n}^1\log(1/y)\,\de y\right)}_{\text{outside the strip}}=O(1).
\end{align*}
Inside the strip, Lemma~\ref{lem:filter-crude-on-disks} bounds the ratio inside $\log^+$ by $n^K/(p^2\Delta_\Gamma^{23})\le n^K$, since the factor $(|z|\vee n\eta_n)/n$ cancels and $p,\Delta_\Gamma\ge1$. Outside the strip, the full prefactor in Eq.~\eqref{eq:filter-direct-off-strip} cancels, leaving at most $\log^+ C+23\log^+(1/y)$ with $y=|\Im\zeta|/\dist(z,I_n)$. The within strip term thus contributes $O(h_nK\log n)$, and the remaining integral is bounded by $\int_0^1\log(1/y) \de y=1$. Since $h_n\log n=O(1)$, both contributions are bounded. Exponentiating the submean bound proves Eq.~\eqref{eq:filter-single-on-disk} uniformly for $z'\in D_z$, and Part I closes the variance bound of Lemma~\ref{lem:filter-contour-comparison}.

For the squared bias, we repeat the same argument by applying the submean inequality to $\ind_{\mc{A}_n}\zeta\zeta'\chi^{\bQ}(\zeta,\zeta')$ successively on $D_z$ and $D_{z'}$. We can therefore have Eq.~\eqref{eq:filter-projected-mixed}.

\subsubsection{Proof of Lemma~\ref{lem:filter-crude-on-disks}}\label{proof:filter-crude-on-disks}

\paragraph{Part I: Proof of Eq.~\eqref{eq:filter-crude-single-on-disks}.} On $\mc{A}_n$, we bound $|\chi^G(\zeta)|$ for $z'\in D_z$ and $\zeta\in D_{z'}$. We work on the larger disk $|\zeta-z|\le 3\dist(z,I_n)/4$, which contains every such $D_{z'}$ by the disk inclusion following Eq.~\eqref{eq:filter-shifted-stability-controls}. Part II will use the bound for $|\wt{\omega}(\zeta)-\zeta|$ on these same disks to estimate a divided difference, so we provide bounds on them together. Taking $\bQ=\bI$ in Eq.~\eqref{eq:terminal-zero-sample-traces} and using $\wt{\omega}=-1/s_n$, we can write both endpoints as weighted traces.
\begin{align*}
    \zeta G_{n+1}(\zeta)&=\Tr\prn{\bSigma\,\zeta(\bT_n-\zeta\bI)^{-1}},\\
    \zeta G_1(\zeta)&=\Tr\prn{\bSigma(\bSigma/\omega(\zeta/n)-\bI)^{-1}}
    =-\Tr\prn{\bSigma(\bI+n s_n(\zeta)\bSigma)^{-1}}=\wt{\omega}(\zeta)-\zeta.
\end{align*}
It thus suffices to bound the two operator norms. For the sample endpoint, the spectrum of $\bT_n$ lies in $\{0\}\cup I_n$ on $\mc{A}_n$, with $\zeta/(0-\zeta)=-1$ by continuation. For the population endpoint, diagonalize $\bSigma$. If $\Re\zeta\le n\eta_n$, then $\Re s_n(\zeta)\ge0$ since $\supp\mu_n\subset I_n$, so $|1+n\sigma_j s_n(\zeta)|\ge1$. Otherwise $\dist(\zeta,I_n)=|\Im\zeta|$, and $\Im(\zeta s_n(\zeta))=\Im\zeta\int_{I_n}x|x-\zeta|^{-2}\mu_n(\de x)$ implies $|\zeta|\cdot |1+n\sigma_j s_n(\zeta)|\ge|\Im (\zeta(1+n\sigma_j s_n(\zeta)))|\ge|\Im\zeta|$. Hence
\begin{align*}
    \|\zeta(\bT_n-\zeta\bI)^{-1}\|\vee\|(\bI+n s_n(\zeta)\bSigma)^{-1}\|
    \le 1\vee\frac{|\zeta|}{\dist(\zeta,I_n)}=O(\Delta_\Gamma),
\end{align*}
where $\dist(\zeta,I_n)\ge\dist(z,I_n)/4$ and $|z|\le\Delta_\Gamma\dist(z,I_n)$. Applying $|\Tr(\bSigma\bA)|\le\Tr(\bSigma)\|\bA\|$ to each endpoint, we have
\begin{align*}
    |\chi^G(\zeta)|&\le |\zeta G_{n+1}(\zeta)|+|\zeta G_1(\zeta)| \le\Tr(\bSigma)\big(\|\zeta(\bT_n-\zeta\bI)^{-1}\|+\|(\bI+n s_n(\zeta)\bSigma)^{-1}\|\big)
    =O(\Delta_\Gamma\Tr\bSigma).
\end{align*}
Note that the same bound holds for $|\wt{\omega}(\zeta)-\zeta|$.

We are now ready to close the proof of Eq.~\eqref{eq:filter-crude-single-on-disks}. Taking a minimizing cutoff $m$ in $\Lambda_n(-\dist(\Gamma,I_n))$ and using $\|\bSigma\|=1$ and Eq.~\eqref{eq:path-spectral-envelope}, we have $\Tr\bSigma\le m+\tau_1(m)\le n+\dist(\Gamma,I_n)\Delta_\Gamma$. Since $\dist(\Gamma,I_n)\le |z-n\eta_n|\le2(|z|\vee n\eta_n)$,
\begin{align*}
    \frac{n|\chi^G(\zeta)|}{|z|\vee n\eta_n}
    =O\left(n\Delta_\Gamma \cdot \brc{\frac{n}{\dist(\Gamma,I_n)}+\Delta_\Gamma}\right).
\end{align*}
As Theorem~\ref{thm:analytic-filter-risk} assumes $\dist(\Gamma,I_n)\ge n^{-C_\Gamma}$ and $\Delta_\Gamma=O(n^{1/16}/\log n)$, we conclude that the right-hand side is at most $n^K$ for some fixed $K$. We complete the proof for the first bound.

\paragraph{Part II: Proof of Eq.~\eqref{eq:filter-crude-mixed-on-disks}.} For $(\zeta,\zeta')\in D_z\times D_{z'}$, applying the two resolvent bounds from Part I to the mixed-trace identity in Appendix~\ref{proof:analytic-filter-risk}, with $\bbeta=\bSigma^{1/2}\btheta$ and $\|\bSigma\|=1$, yields
\begin{align*}
    |\zeta\zeta'\chi^{\bQ}(\zeta,\zeta')|
    \le O\!\left(\|\btheta\|^2\Delta_\Gamma^2\right)\big(1+|\wt{\omega}[\zeta,\zeta']|\big).
\end{align*}
It remains to bound the divided difference. Since $\wt{\omega}[\zeta,\zeta']-1$ is the divided difference of $\wt{\omega}(\zeta)-\zeta$, we employ the same strategy on $\Gamma_{z\mid z'}$ and $\Gamma_{z'\mid z}$ from Part I of the proof of Lemma~\ref{lem:filter-contour-comparison}, with $(z,z')$ there replaced by $(\zeta,\zeta')$. Each arc lies in the corresponding larger disk from Part I, and $\dist(\zeta,I_n)\ge 3\dist(z,I_n)/4$, with the same comparison for $\zeta'$. Inserting the bound from Part I on each arc gives
\begin{align*}
    |\wt{\omega}[\zeta,\zeta']-1|
    &=O\left(\Delta_\Gamma\Tr\bSigma\cdot\brc{\frac1{\dist(z,I_n)}+\frac1{\dist(z',I_n)}}\right)
    =O\left(\Delta_\Gamma\cdot\brc{\frac{n}{\dist(\Gamma,I_n)}+\Delta_\Gamma}\right).
\end{align*}
Finally, $(|z|\vee n\eta_n)(|z'|\vee n\eta_n)\ge\dist(\Gamma,I_n)^2/4$. Thus multiplying the mixed-trace bound by $n^{5/2}/((|z|\vee n\eta_n)(|z'|\vee n\eta_n))$ yields at most $n^K\|\btheta\|^2$ by the same polynomial bounds as in Part I. We prove Eq.~\eqref{eq:filter-crude-mixed-on-disks} and thus the lemma.

\subsection{Proof of Proposition~\ref{prop:filter-real-risk}}\label{proof:filter-real-risk}
We shall begin with the Stieltjes representation in Eq.~\eqref{eq:filter-directional-measures}, with associated probability measure $\mu_{\omega, j}$ and weighted measure $\mu_\omega$. After establishing the formulas in Eq.~\eqref{eq:filter-real-risk} from Eq.~\eqref{eq:analytic-filter-risk} in Parts I and II, we come back to verify that $s_{\omega, j}$ are indeed Stieltjes transforms of probability measures and that $\mu_\omega$ is finite.

\paragraph{Part I: Proof of Eq.~\eqref{eq:filter-real-variance}.} Starting from Eq.~\eqref{eq:analytic-filter-risk-variance}, we have
\begin{align*}
    \effectiveVariance[\upphi]=\tau^2\iint_{\Gamma^2}\big(\wt{\omega}[z,z']-1\big)\,\alpha_{\upphi}(\de z)\alpha_{\upphi}(\de z').
\end{align*}
Next, we use the fixed-point equation for $\omega$ and Eq.~\eqref{eq:filter-directional-measures} at $z/n$ to write
\begin{align*}
    \wt{\omega}(z)-z+\Tr\bSigma
    &=\sum_{j}\sigma_j\left(1+\frac{z}{n}s_{\omega,j}(z/n)\right)
    =\sum_{j}\sigma_j\int\frac{nx}{nx-z}\,\mu_{\omega,j}(\de x)
    =n^2\int\frac{\mu_\omega(\de x)}{nx-z}.
\end{align*}
Taking the difference of the above identity between $z$ and $z'$ and dividing by $z-z'$, we have
\begin{equation}\label{eq:omega-divided-difference-measure}
    \wt{\omega}[z,z']-1
    =\frac{(\wt{\omega}(z)-z)-(\wt{\omega}(z')-z')}{z-z'}
    =n^2\int\frac{\mu_\omega(\de x)}{(nx-z)(nx-z')}.
\end{equation}
As we will show later in Part III, $\supp\mu_\omega\subseteq[\eta_n,\infty)$, so we have $nx\in I_n$ on its support. Thus, by substituting into the variance and applying Eq.~\eqref{eq:analytic-filter-representation} in $z$ and $z'$, we recover Eq.~\eqref{eq:filter-real-variance}:
\begin{align*}
    \effectiveVariance[\upphi]
    &=\tau^2\int\left(\int_\Gamma\frac{n\,\alpha_{\upphi}(\de z)}{nx-z}\right)
    \left(\int_\Gamma\frac{n\,\alpha_{\upphi}(\de z')}{nx-z'}\right)\mu_\omega(\de x)=\tau^2\int\upphi(x)^2\,\mu_\omega(\de x).
\end{align*}
\paragraph{Part II: Proof of Eq.~\eqref{eq:filter-real-bias}.} By Eq.~\eqref{eq:analytic-filter-risk-bias} and writing the two population resolvents in the eigenbasis of $\bSigma$, we have
\begin{align*}
    \effectiveBias[\upphi]=\sum_j\sigma_j|\langle\bu_j,\bbeta\rangle|^2\iint_{\Gamma^2}
    \frac{\wt{\omega}[z,z']\,\alpha_{\upphi}(\de z)\alpha_{\upphi}(\de z')}
    {(1-n\sigma_j/\wt{\omega}(z))(1-n\sigma_j/\wt{\omega}(z'))}.
\end{align*}
To identify the two terms in Eq.~\eqref{eq:filter-real-bias}, we apply Eq.~\eqref{eq:omega-divided-difference-measure} to $\wt{\omega}[z,z']=1+(\wt{\omega}[z,z']-1)$. Therefore it follows that
\begin{align*}
    \effectiveBias[\upphi]=\sum_j\sigma_j|\langle\bu_j,\bbeta\rangle|^2 \Biggl\{\Bigl(\underbrace{\int_\Gamma\frac{\alpha_{\upphi}(\de z)}{1-n\sigma_j/\wt{\omega}(z)}}_{\mathrm{(I)}}\Bigr)^2
    +\int\Bigl(\underbrace{\int_\Gamma\frac{n}{nx'-z}\,\frac{\alpha_{\upphi}(\de z)}{1-n\sigma_j/\wt{\omega}(z)}}_{\mathrm{(II)}}\Bigr)^2\mu_\omega(\de x')\Biggr\}.
\end{align*}
Indeed, the $z$ and $z'$ integrals coincide in each summand, with $x'$ fixed in the second. By Eq.~\eqref{eq:filter-directional-measures} at $z/n$, the common resolvent factor is
\begin{align*}
    \left(1-\frac{n\sigma_j}{\wt{\omega}(z)}\right)^{-1}=-\frac{z}{n}s_{\omega,j}(z/n)=\int\frac{z}{z-nx}\,\mu_{\omega,j}(\de x).
\end{align*}
Substituting this representation and applying Eq.~\eqref{eq:analytic-filter-representation}, we have
\begin{align*}
    \mathrm{(I)}&=\int\!\left[\int_\Gamma\frac{z}{z-nx}\,\alpha_{\upphi}(\de z)\right]\mu_{\omega,j}(\de x)
    =\int\uppsi(x)\,\mu_{\omega,j}(\de x),\\
    \mathrm{(II)}&=\int\!\left[\frac{1}{x'-x}\int_\Gamma\left(\frac{z}{z-nx}-\frac{z}{z-nx'}\right)\alpha_{\upphi}(\de z)\right]\mu_{\omega,j}(\de x)
    =\int\frac{\uppsi(x)-\uppsi(x')}{x'-x}\,\mu_{\omega,j}(\de x).
\end{align*}
In the second row we use partial fractions for $x\ne x'$ and the continuous extension $-\uppsi'(x')$ at $x=x'>0$. Substituting (I) and (II) into the expansion of $\effectiveBias[\upphi]$ proves Eq.~\eqref{eq:filter-real-bias}.

We implicitly applied Fubini's theorem in Parts I and II. We can apply because the contour bounds in Theorem~\ref{thm:analytic-filter-risk} make both kernels uniformly bounded for fixed $n$, the measures have finite total variation, and $\sum_j\sigma_j|\langle\bu_j,\bbeta\rangle|^2<\infty$.

\paragraph{Part III: Probability measures in Eq.~\eqref{eq:filter-directional-measures}.} Lastly, we check $\mu_{\omega, j}$ are indeed probability measures and the representation is valid. From Eq.~\eqref{eq:filter-directional-measures}, $s_{\omega,j}(z)^{-1}=-z(1+n\sigma_j s_n(nz))$  for $\sigma_j>0$. Using the Stieltjes representation of $s_n$, we obtain that
\begin{align*}
    \Im\prn{z(1+n\sigma_j s_n(nz))}
    &=\Im z\left(1+n\sigma_j\int\frac{x}{|x-nz|^2}\,\mu_n(\de x)\right)>0,\qquad z\in\C_+,\\
    (-\iu y)s_{\omega,j}(\iu y)&=\frac{1}{1+n\sigma_j s_n(n\iu y)}\to1,\qquad y\to\infty,
\end{align*}
where the second equation makes use of $|s_n(n\iu y)|\le1/(ny)$.  By the first equation, $z(1+n\sigma_j s_n(nz))$ never vanishes on $\C_+$. Its negative reciprocal $s_{\omega,j}$ is therefore analytic there, with
\begin{align*}
    \Im s_{\omega,j}(z)=\frac{\Im\prn{z(1+n\sigma_j s_n(nz))}}{|z(1+n\sigma_j s_n(nz))|^2}>0.
\end{align*}
Thus $s_{\omega,j}$ is Herglotz-Nevanlinna, and the second equation gives $y\Im s_{\omega,j}(\iu y)\to1$ and $s_{\omega,j}(\iu y)\to0$ as $y\to\infty$. Applying Lemma~\ref{lem:herglotz-nevanlinna} with these limits, as in the proof of Lemma~\ref{lem:population-scalar-stability}, yields a unique positive measure $\mu_{\omega,j}$ satisfying
\begin{align*}
    s_{\omega,j}(z)=\int_{\R}\frac{\mu_{\omega,j}(\de x)}{x-z}
    \qquad\text{and}\qquad
    \mu_{\omega,j}(\R)=\lim_{y\to\infty}y\Im s_{\omega,j}(\iu y)=1.
\end{align*}
By Assumption~\ref{assmp:analytic-spectral-filter}, $s_n(nz)>0$ for real $z<\eta_n$, so $s_{\omega,j}$ is real analytic there except possibly at zero. Stieltjes inversion gives $\supp\mu_{\omega,j}\subseteq\{0\}\cup[\eta_n,\infty)$. Applying monotone convergence to $\lambda x/(x+\lambda) \to x$ and using Eq.~\eqref{eq:filter-directional-measures} at $-\lambda$ implies
\begin{align*}
    \int x\,\mu_{\omega,j}(\de x)=\lim_{\lambda\to\infty}\lambda\big(1-\lambda s_{\omega,j}(-\lambda)\big)=\sigma_j,
\end{align*}
since $n\lambda s_n(-n\lambda)\to1$. Thus $\mu_\omega(\R)=n^{-1}\sum_j\sigma_j^2=\Tr(\bSigma^2)/n<\infty$, and the factor $x$ in its definition removes any mass at zero, concluding the proof with $\supp\mu_\omega\subseteq[\eta_n,\infty)$.

    \section{Proofs for Applications}\label{app:proofs-for-applications}
Throughout the application proofs, we leave dependence on the decay exponent $\alpha$ or $\beta$ implicit for simplicity.

We first introduce the lemma that provides the hard edge guarantee on the event $\{\what{\sigma}_n \geq \eta_n\}$, which serves as the common input for ridge-type filters and gradient flow. The lemma also gives quantitative bounds that helps us check assumptions and conditions to apply Theorem~\ref{thm:analytic-filter-risk}. We shall take $\eta_n = \Theta_{C_{\bx}}(\tau_1(n)/n)$ next for both propositions. See its proof in Appendix~\ref{proof:application-spectral-bounds}.
\begin{lemma}\label{lem:application-spectral-bounds}
    Under Assumption~\ref{assmp:features}, take either spectrum in Eq.~\eqref{eq:spectrum-decay-examples}. We can take $\eta_n$ above such that for all $n=\Omega_{C_{\bx}}(1)$, it holds $\P\prn{\what\sigma_n\geq\eta_n} \geq1-2e^{-n}$ and $\supp\mu_n\subseteq[n\eta_n,\infty)$, and the following consequences hold under both laws
    \begin{align*}
        \Lambda_n(-n\eta_n/4)&=O_{C_{\bx}}((\log n)^2)\qquad \text{and} \qquad -\log(1\wedge n\eta_n)=O_{C_{\bx}}(\log n).
    \end{align*}
\end{lemma}

\subsection{Proof of Proposition~\ref{prop:ridge-type-risk}}\label{proof:ridge-type-risk}

\paragraph{Part I: Bounding approximation errors in (i).}
In the proof, we write $\lambda_+:=\lambda\vee0$. We will first apply Theorem~\ref{thm:analytic-filter-risk} to obtain the error control before computing their explicit formulas and specializing to $k=1$ in the next part. To this end, Lemma~\ref{lem:application-spectral-bounds} supplies the event probability and the inclusion $\supp\mu_n\subseteq I_n$. It remains to check the contour conditions and bound the two contour integrals. We use the two-circle construction (cf.~Figure~\ref{fig:ridge-filter-contours} for illustration) with the explicit choice
\begin{align*}
    \Omega=\brc{|z|<n\eta_n/4}  \cup \brc{|z+n\lambda_+|>n(\lambda_++3\eta_n/4)}=:D_{\mathrm{in}}\cup D_{\mathrm{out}}.
\end{align*}
We point out that the orientation convention is that boundaries $\partial D_{\mathrm{in}}$ and $\partial D_{\mathrm{out}}$ are counterclockwise and clockwise, respectively. Define the residual $\uppsi_{k,\lambda}(z/n)$ to equal $1$ in $D_{\mathrm{in}}$ and $(\lambda/(z/n+\lambda))^k$ on $D_{\mathrm{out}}$. The pole $-n\lambda$ lies strictly inside the excluded outer disk. Thus the corresponding filter is holomorphic on both components, agrees with $z\mapsto\upphi_{k,\lambda}(z/n)$ on $I_n$, and has residual tending to zero at infinity, verifying Assumption~\ref{assmp:analytic-spectral-filter}. Both circles lie in $\Re z\leq3n\eta_n/4$. Consequently,
\begin{align*}
    \dist(\Gamma,I_n)=n\eta_n/4\qquad\text{and}\qquad
    \sup_{z\in\Gamma}\frac{|z|\vee n\eta_n}{\dist(z,I_n)}\leq5.
\end{align*}
Indeed, $\dist(z,I_n)=|z-n\eta_n|$ on $\Gamma$ and $|z|\leq|z-n\eta_n|+n\eta_n$. Applying Lemma~\ref{lem:application-spectral-bounds} yields
\begin{align*}
    \Delta_\Gamma\leq\Lambda_n(-n\eta_n/4)\vee5=O_{C_{\bx}}((\log n)^2)\qquad\text{and}\qquad
    -\log(1\wedge\dist(\Gamma,I_n))=O_{C_{\bx}}(\log n).
\end{align*}
Thus the two numerical conditions of Theorem~\ref{thm:analytic-filter-risk} hold for all sufficiently large $n$.  It remains to evaluate the three quantities entering the error bounds of Theorem~\ref{thm:analytic-filter-risk}:
\begin{align*}
    \tvnorms{\alpha_{\upphi_{k,\lambda}}}, &  &
    \int_\Gamma\frac{|z|\vee n\eta_n}{\dist(z,I_n)}|\alpha_{\upphi_{k,\lambda}}|(\de z),  & &
    \int_\Gamma(|z|\vee n\eta_n)|\alpha_{\upphi_{k,\lambda}}|(\de z).
\end{align*}

First, we control the term $\tvnorms{\alpha_{\upphi_{k,\lambda}}}$. On $\partial D_{\mathrm{in}}$, we have $|\alpha_{\upphi_{k,\lambda}}|(\de z)=|\de z|/(2\pi|z|)$, so this circle contributes exactly $1$. On $\partial D_{\mathrm{out}}$, by a change of variable $z/n=(\lambda_++3\eta_n/4)e^{\iu\theta}-\lambda_+$, if $-\eta_n/2\leq\lambda\leq0$, then $|(3\eta_n/4)e^{\iu\theta}+\lambda|\geq3\eta_n/4-|\lambda|\geq\eta_n/4$, so
\begin{align*}
    |\alpha_{\upphi_{k,\lambda}}|(\de z)
    =\left|\frac{\lambda}{(3\eta_n/4)e^{\iu\theta}+\lambda}\right|^k\frac{\de\theta}{2\pi}
    \leq2^k\frac{\de\theta}{2\pi}.
\end{align*}
If $\lambda>0$, the density on $\partial D_{\mathrm{out}}$ is instead
\begin{align*}
    |\alpha_{\upphi_{k,\lambda}}|(\de z)
    =\left(\frac{\lambda}{\lambda+3\eta_n/4}\right)^k
    \frac{\lambda+3\eta_n/4}{|(\lambda+3\eta_n/4)e^{\iu\theta}-\lambda|}\frac{\de\theta}{2\pi}.
\end{align*}
Using $|(\lambda+3\eta_n/4)e^{\iu\theta}-\lambda|^2=(3\eta_n/4)^2+4\lambda(\lambda+3\eta_n/4)\sin^2(\theta/2)$ and $\sin(\theta/2)\geq\theta/\pi$ on $[0,\pi]$, we obtain
\begin{align*}
    \int_{\partial D_{\mathrm{out}}}|\alpha_{\upphi_{k,\lambda}}|(\de z)
    &\stackrel{\mathrm{(i)}}{\leq}\frac{\lambda}{\pi}\int_0^\pi\frac{\de\theta}{\sqrt{(3\eta_n/4)^2+4\lambda^2\theta^2/\pi^2}}
    \leq\frac12\log\prn{1+\frac{16\lambda}{3\eta_n}},
\end{align*}
where we use in (i) that $(\lambda/(\lambda+3\eta_n/4))^k \cdot (\lambda+3\eta_n/4)\leq\lambda$. Combining the two cases, and adding the contribution from $\partial D_{\mathrm{in}}$, we have
\begin{align}
    \tvnorms{\alpha_{\upphi_{k,\lambda}}}=O_k\prn{1+\log(1+\lambda_+/\eta_n)}=O_{C_{\bx},k}(\log n),\label{eq:ridge-contour-total-variation}
\end{align}
where the last step uses $\lambda=O(1)$ and $-\log\eta_n=O_{C_{\bx}}(\log n)$. The second quantity then follows immediately, applying the geometric bound already proved on $\Gamma$,
\begin{align}
    \int_\Gamma\frac{|z|\vee n\eta_n}{\dist(z,I_n)}\,|\alpha_{\upphi_{k,\lambda}}|(\de z)
    \leq5\tvnorms{\alpha_{\upphi_{k,\lambda}}}=O_{C_{\bx},k}(\log n).\label{eq:ridge-contour-separation-integral}
\end{align}

Thirdly, we control the integral weighted by $|z|\vee n\eta_n$. On $\partial D_{\mathrm{in}}$, this weight equals $n\eta_n$. On $\partial D_{\mathrm{out}}$, the preceding density bounds imply
\begin{align*}
    \int_{\partial D_{\mathrm{out}}}|z|\,|\alpha_{\upphi_{k,\lambda}}|(\de z)
    &\leq\begin{cases}
    2^k\cdot3n\eta_n/4,&-\eta_n/2\leq\lambda\leq0,\\
    n(\lambda+3\eta_n/4)\left(\dfrac{\lambda}{\lambda+3\eta_n/4}\right)^k,&\lambda>0,
\end{cases}
\end{align*}
which are both bounded by $O_k(n(\lambda+\eta_n))$. Combining with the inner contribution and using $|z|\vee n\eta_n\leq|z|+n\eta_n$ on the outer circle, we conclude that
\begin{align}
    \int_\Gamma(|z|\vee n\eta_n)|\alpha_{\upphi_{k,\lambda}}|(\de z)
    &\leq n\eta_n\tvnorms{\alpha_{\upphi_{k,\lambda}}}+\int_{\partial D_{\mathrm{out}}}|z|\,|\alpha_{\upphi_{k,\lambda}}|(\de z)\nonumber\\
    &=O_k\prn{n(\lambda+\eta_n)+n\eta_n\log(1+\lambda_+/\eta_n)}=O_k(n(\lambda+\eta_n)),\label{eq:ridge-contour-weighted-integral}
\end{align}
where the last step uses $\log(1+x)\leq x$ and $\lambda_+\leq\lambda+\eta_n$. Substituting Eqs.~\eqref{eq:ridge-contour-total-variation}, \eqref{eq:ridge-contour-separation-integral} and~\eqref{eq:ridge-contour-weighted-integral} into Theorem~\ref{thm:analytic-filter-risk} there yields
\begin{align*}
    \|\ind\{\what\sigma_n\geq\eta_n\}(\sampleVariance[\upphi_{k,\lambda}]-\effectiveVariance[\upphi_{k,\lambda}])\|_{L^p}
    &=O_{C_{\bx},\kappa,k}\prn{\frac{\tau^2p^2\Delta_\Gamma^{23}(\log n)^2}{n}},\\
    \|\ind\{\what\sigma_n\geq\eta_n\}(\sampleBias[\upphi_{k,\lambda}]-\effectiveBias[\upphi_{k,\lambda}])\|_{L^p}
    &=O_{C_{\bx},\kappa,k}\prn{\frac{p^{3/2}\Delta_\Gamma^{26}(\lambda+\eta_n)^2\norm{\btheta}^2}{\sqrt n}}.
\end{align*}
Here $\Delta_\Gamma=O_{C_{\bx}}((\log n)^2)$ and $p\leq\kappa\log n$ contribute only powers of $\log n$. This proves the two approximation bounds in (i).

\paragraph{Part II: Evaluating the effective formulas in (i) and (ii).} It remains to concretely compute the effective risk $\effectiveRisk{}$, and the Cauchy's formula is helpful in this case. We first reduce both contour integrals to the pole at $-\lambda$. For $F$ analytic on a neighborhood of the closed outer disk in units $\zeta=z/n$, deform the inner integral of $F(\zeta)/\zeta$ to its outer boundary. Reversing the exterior contribution gives
\begin{align}
    \int_\Gamma F(z/n)\,\alpha_{\upphi_{k,\lambda}}(\de z)
    &=-\frac1{2\pi\iu}\oint_{\partial D_{\mathrm{out}}}
    F(z/n)\frac{1-(\lambda/(z/n+\lambda))^k}{z}\,\de z=\sum_{a=0}^{k-1}\frac{\lambda^a}{a!}F^{(a)}(-\lambda).\label{eq:ridge-filter-residue}
\end{align}
Here $\partial D_{\mathrm{out}}$ is clockwise, and after substituting $\zeta=z/n$, the last equality follows by expanding the kernel as $\sum_{a=0}^{k-1}\lambda^a/(\zeta+\lambda)^{a+1}$ and applying Cauchy's derivative formula. At $\lambda=0$ it is simply $1/\zeta$, giving $F(0)$ for every $k$. By Eq.~\eqref{eq:analytic-filter-risk}, for variance and squared bias we first apply Eq.~\eqref{eq:ridge-filter-residue} in $\zeta=z/n$, with $\zeta'=z'/n$ fixed, taking
\begin{align*}
    F(\zeta)=\begin{cases}
    \tau^2\prn{\omega[\zeta,\zeta']-1},&\text{for computing }\effectiveVariance,\\[.4em]
    \omega[\zeta,\zeta']\,\bbeta^*\bSigma\prn{\bI-\dfrac{\bSigma}{\omega(\zeta)}}^{-1}\prn{\bI-\dfrac{\bSigma}{\omega(\zeta')}}^{-1}\bbeta,&\text{for computing }\effectiveBias.
\end{cases}
\end{align*}
Substituting the second kernel into Eq.~\eqref{eq:ridge-filter-residue} in both variables directly recovers the squared-bias formula in (i). For variance, applying the same identity yields
\begin{align*}
    \effectiveVariance[\upphi_{k,\lambda}]
    &\stackrel{\mathrm{(i)}}{=}\tau^2\brc{\sum_{a,b=0}^{k-1}\frac{\lambda^{a+b}\omega^{(a+b+1)}(-\lambda)}{(a+b+1)!}-1}= \tau^2\brc{\sum_{m=1}^{2k-1}\frac{(k-|k-m|)\lambda^{m-1}\omega^{(m)}(-\lambda)}{m!}-1}.
\end{align*}
In (i), the Taylor expansion of $\omega$ at $-\lambda$ gives the coefficient $\omega^{(a+b+1)}(-\lambda)/(a+b+1)!$ of $(\zeta+\lambda)^a(\zeta'+\lambda)^b$ in $\omega[\zeta,\zeta']$; the constant $-1$ contributes only when $a=b=0$. In the last equality, we group the terms with $a+b=m-1$, of which there are $k-|k-m|$ for $0\leq a,b\leq k-1$. This proves the variance formula in (i).

We move on to specialization for (ii), using $\lambda_\star=-\omega(-\lambda)=(ns_n(-n\lambda))^{-1}>0$. Rearranging the fixed-point equation in Eq.~\eqref{eq:ridge-effective-regularization} and differentiating in $\lambda_\star = \lambda_\star(\lambda)$ implies
\begin{align*}
    \lambda=\lambda_\star\cdot \prn{1-\frac1n\Tr(\bSigma(\bSigma+\lambda_\star\bI)^{-1})}
    \qquad\text{and}\qquad
    \omega'(-\lambda)=\frac1{1-\mathsf{J}_n(\lambda)}>0.
\end{align*}
Substituting $k=1$ into the formulas in (i) yields both effective risks in (ii). Finally, we have from $\bbeta=\bSigma^{1/2}\btheta$ that
\begin{align*}
    \frac{\tau^2}{\effectiveVariance[\upphi_{1,\lambda}]}=\frac{1-\mathsf{J}_n(\lambda)}{\mathsf{J}_n(\lambda)}\qquad\text{and}\qquad
    \frac{(\lambda+\eta_n)^2\|\btheta\|^2}{\effectiveBias[\upphi_{1,\lambda}]}=\frac{(\lambda+\eta_n)^2(1-\mathsf{J}_n(\lambda))\|\btheta\|^2}{\lambda_\star^2\|\bSigma(\bSigma+\lambda_\star\bI)^{-1}\btheta\|^2}.
\end{align*}
We close the proof with the multiplicative bounds.

\subsection{Proof of Proposition~\ref{prop:gradient-flow-risk}}\label{proof:gradient-flow-risk}
\paragraph{Part I: Bounding approximation errors.}
We repeat largely the same proof in Appendix~\ref{proof:ridge-type-risk}. Applying Theorem~\ref{thm:analytic-filter-risk} and Lemma~\ref{lem:application-spectral-bounds}, we compute the contour integrals on the positively oriented boundary $\Gamma=\partial\Omega$ of $\Omega=\{z:|\arg(z+n/t)|<\pi/4\}$, as in Figure~\ref{fig:analytic-filter-support}. The filter $\upphi_t(z/n)$ is entire, and its residual $e^{-tz/n}$ decays exponentially along the two rays, verifying Assumption~\ref{assmp:analytic-spectral-filter}. Here
\begin{align*}
    \dist(\Gamma,I_n)=\frac{n(\eta_n+t^{-1})}{\sqrt2}
    \qquad\text{and}\qquad
    \sup_{z\in\Gamma}\frac{|z|\vee n\eta_n}{\dist(z,I_n)}\leq3.
\end{align*}
For the second bound, use $|z|\leq|z-n\eta_n|+n\eta_n\leq(1+\sqrt2)\dist(z,I_n)\leq3\dist(z,I_n)$ when $\Re z\leq n\eta_n$, where $\sqrt2=1/\sin(\pi/4)$ comes from the angle of the rays, and $\dist(z,I_n)=|\Im z|=\Re z+n/t$ otherwise. Since $\dist(\Gamma,I_n)\geq n\eta_n/4$, applying Lemma~\ref{lem:application-spectral-bounds} and the monotonicity of $\Lambda_n(-a)$ in Eq.~\eqref{eq:path-spectral-envelope} yields $\Delta_\Gamma=O_{C_{\bx}}((\log n)^2)$ and $-\log(1\wedge\dist(\Gamma,I_n))=O_{C_{\bx}}(\log n)$. Thus the conditions hold for all sufficiently large $n$, independently of $t$. As in Appendix~\ref{proof:ridge-type-risk}, the remaining quantities are
\begin{align*}
    \tvnorms{\alpha_{\upphi_t}}, & &
    \int_\Gamma\frac{|z|\vee n\eta_n}{\dist(z,I_n)}|\alpha_{\upphi_t}|(\de z), & &
    \int_\Gamma(|z|\vee n\eta_n)|\alpha_{\upphi_t}|(\de z).
\end{align*}
We bound them one-by-one next. For the first, parametrize each ray by $z=(n/t)(-1+r e^{\pm\iu\pi/4})$, $r\geq0$, so that
\begin{align*}
    |\alpha_{\upphi_t}|(\de z)=\frac{e^{1-r/\sqrt2}}{2\pi|-1+r e^{\pm\iu\pi/4}|}\,\de r
    \qquad\text{and}\qquad
    |-1+r e^{\pm\iu\pi/4}|^2=(r-1/\sqrt2)^2+1/2.
\end{align*}
Integrating over both rays yields
\begin{align}
    \tvnorms{\alpha_{\upphi_t}}\leq\frac{\sqrt2}{\pi}\int_0^\infty e^{1-r/\sqrt2}\,\de r=O(1).\label{eq:gradient-flow-contour-total-variation}
\end{align}
For the second quantity, applying the geometric bound above yields
\begin{align}
    \int_\Gamma\frac{|z|\vee n\eta_n}{\dist(z,I_n)}|\alpha_{\upphi_t}|(\de z)
    \leq3\tvnorms{\alpha_{\upphi_t}}=O(1).\label{eq:gradient-flow-contour-separation-integral}
\end{align}
For the third quantity, it follows from $|z|\vee n\eta_n\leq|z|+n\eta_n$ and the cancellation of $|z|$ that
\begin{align}
    \int_\Gamma(|z|\vee n\eta_n)|\alpha_{\upphi_t}|(\de z)
    &\leq\frac{n}{\pi t}\int_0^\infty e^{1-r/\sqrt2}\,\de r+n\eta_n\tvnorms{\alpha_{\upphi_t}}
    =O\prn{n(\eta_n+t^{-1})}.\label{eq:gradient-flow-contour-weighted-integral}
\end{align}
Substituting Eqs.~\eqref{eq:gradient-flow-contour-total-variation}, \eqref{eq:gradient-flow-contour-separation-integral} and~\eqref{eq:gradient-flow-contour-weighted-integral} into Theorem~\ref{thm:analytic-filter-risk} proves the two approximation bounds. As in Part I of Appendix~\ref{proof:ridge-type-risk}, the factors $p^2\Delta_\Gamma^{23}$ and $p^{3/2}\Delta_\Gamma^{26}$ are absorbed into $\wt O_{C_{\bx},\kappa}$, with implied constants independent of $t$.

\paragraph{Part II: Evaluating the effective formulas.}
To obtain the effective-risk formulas from Proposition~\ref{prop:filter-real-risk}, write the two filter expressions there as
\begin{align*}
    x\upphi_t(x)^2=\int_0^t(e^{-rx}-e^{-(r+t)x})\de r\qquad\text{and}\qquad
    \frac{\uppsi_t(x)-\uppsi_t(x')}{x'-x}=\int_0^t e^{-rx-(t-r)x'}\de r.
\end{align*}
The identities extend continuously to zero and to $x=x'$. Integrating against $\mu_{\omega,j}$ and applying Fubini-Tonelli by nonnegativity yields
\begin{align*}
    \int_{[0,\infty)}x\upphi_t(x)^2\mu_{\omega,j}(\de x)&=\int_0^t\bigl(\mathcal L_{\omega,j}(r)-\mathcal L_{\omega,j}(r+t)\bigr)\de r,\\
    \int_{[0,\infty)}\frac{\uppsi_t(x)-\uppsi_t(x')}{x'-x}\mu_{\omega,j}(\de x)&=\int_0^t e^{-(t-r)x'}\mathcal L_{\omega,j}(r)\de r.
\end{align*}
Substituting the first identity into Eq.~\eqref{eq:filter-real-variance} proves the variance formula. Substituting the second into Eq.~\eqref{eq:filter-real-bias}, with $\int\uppsi_t(x)\mu_{\omega,j}(\de x)=\mathcal L_{\omega,j}(t)$, proves the squared-bias formula.

\subsection{Proof of Lemma~\ref{lem:application-spectral-bounds}}\label{proof:application-spectral-bounds}
\paragraph{Part I: $\eta_n$ for the hard edge event.}
By Lemma~\ref{lem:smallest-eigenvalue}, to obtain $\P\prn{\what\sigma_n\geq\eta_n}\geq1-2e^{-n}$, it suffices to take $\eta_n:=\tau_1(m)/(2n)$ for some $m$, provided that
\begin{align*}
    \tau_1(m)=\Omega_{C_{\bx}}\prn{\sqrt{2n\tau_2(m)}+2n\sigma_{m+1}}.
\end{align*}
We verify this condition for both laws when $n=\Omega_{C_{\bx}}(1)$. We take $m=\lceil\const n\rceil$ and $m=n$ in the power law and critical law, respectively. Here $\const \geq1$ depends only on $C_{\bx}$ and $\alpha$. By standard elementary estimates to the tails of the relevant series,
\begin{align*}
    \frac{\sqrt{2n\tau_2(m)}+2n\sigma_{m+1}}{\tau_1(m)}
    =\begin{cases}O_\alpha(\const^{-1/2}+\const^{-1}),&\text{power law,}\\ O_\beta((1+\log n)^{-1}),&\text{critical law.}\end{cases}
\end{align*}
We meet the condition by choosing $\const$ large enough and then taking $n$ large enough. Thus the above choices of $\eta_n = \tau_1(m)/(2n)$ suffice and satisfy $n\eta_n=\Theta_{C_{\bx}}(\tau_1(n))$.

For the deterministic support, we apply the edge characterization in Eq.~\eqref{eq:canonical-edge-equation} and the paragraph following it in Part III of the proof of Proposition~\ref{prop:canonical-bulk-laws} to Eq.~\eqref{eq:population-scalar-equation} to obtain
\begin{align*}
    n\eta_n&\leq\tau_1(m)-2\sqrt{n\tau_2(m)}\leq\inf\supp\mu_n =\max_{s>0}\brc{-\frac1s+\sum_{j\geq1}\frac{\sigma_j}{1+n\sigma_j s}}\leq\tau_1(n)=\Theta_{C_{\bx}}(n\eta_n).
\end{align*}
We explain the bounds above. For the lower bound, we preserve the tail terms, use $(1+x)^{-1}\geq1-x$, and evaluate at $s=1/\sqrt{n\tau_2(m)}$. For the upper bound, it is directly obtained by bounding the first $n$ summands by $1/(ns)$ and the remaining sum by $\tau_1(n)$.

\paragraph{Part II: Preparing bounds for invoking Theorem~\ref{thm:analytic-filter-risk}.}
In the previous Part I, we choose $\eta_n$ explicitly. For the spectral stability parameter, we use the admissible cutoff $m=\lfloor\rho n/4\rfloor$ in Eq.~\eqref{eq:path-spectral-envelope} (note that this $m$ is different from the $m$ we used for the eigenvalue lower bound lemma). From the tail bounds, $n\sigma_{m+1}+\tau_1(m)=O_{C_{\bx}}(\tau_1(n))$ and $\log(nr_1(m))=O_{C_{\bx}}(\log n)$, so
\begin{align*}
    \Lambda_n(-n\eta_n/4)&\leq1+\frac{4\brc{n\sigma_{m+1}+\tau_1(m)\log n\log(nr_1(m))}}{n\eta_n}=O_{C_{\bx}}((\log n)^2).
\end{align*}
For the remaining logarithmic bound, the scale of $\eta_n$ established in Part I implies $-\log(1\wedge n\eta_n)=O_{C_{\bx}}(\log n)$, completing the proof.

    \section{Proofs for the deterministic self-consistent interpolation} \label{proof:interpolation-path}

In the proofs, we will repeatedly use the elementary identity $\Im(\bM^{-1})=-\bM^{-*}(\Im\bM)\bM^{-1}$ for any bounded invertible operator $\bM$ on (the complexification of) $\mc H$. If, in addition, $\Im\bM$ is bounded and invertible, then $\Im(\bM^{-1})^{-1}=-\bM(\Im\bM)^{-1}\bM^*$.

\subsection{Proof of Theorem~\ref{thm:self-consistent-path-stability}} \label{proof:self-consistent-path-stability}
We prove the theorem for $z=t+\iu\eta\in\C_+\cup\R_{<0}$ in four parts.
\begin{enumerate}
    \item In Part I, we prove (i) and collect the weighted resolvent bounds needed for stability.
    \item In Part II, we bound the one-step change $|\wb{G}_{k+1}(z,w_{k+1})-\wb{G}_k(z,w_k)|$ for every $1\leq k\leq n-1$ by canceling its leading drift through leave-one-out identities.
    \item In Part III, we prove (ii), the pathwise expected-trace stability for $\{\wb{G}_k(z,w_k)\}$ by controlling the resolvent norms separately for large and small sample sizes, then sum the one-step errors.
    \item In Part IV, we prove (iii) for the backward path stability for $\{w_k\}$, by converting expected-trace stability into nearly constant increments and comparing the initial endpoint with the effective spectral parameter equation.
\end{enumerate}

\paragraph{Part I: Proof of (i) and resolvent bounds.} We first invoke the following lemma to establish well-posedness and monotonicity of the imaginary part and, when $t<0$, the real part of the backward path. We then prepare the weighted resolvent estimates needed for stability in the following parts. The lemma collects the resolvent and expected-trace bounds from which we derive these path properties. See its proof in Appendix~\ref{proof:resolvent-bounds-strict-positivity}.
\begin{lemma}\label{lem:resolvent-bounds-strict-positivity}
    For any fixed $z=t+\iu\eta\in\C_+ \cup \R_{<0}$, the backward path $\{w_k\}_{k=0}^n$ is uniquely defined, and the following bounds hold.
    \begin{subequations}
        \begin{enumerate}
            \item[(i)] (Resolvent bounds). For $j\in[n]$ and $(z,w)\in\mc D$ (and in particular, $z \in \C_+$), the resolvent exists and $\norm{\bR_{j-1}(z,w)}\leq\eta^{-1}$. In addition, allowing $\eta=0$, along the path, we have for every $j\in[n]$, $0\leq k\leq n$ and $\bx\in\mc H$,
            \begin{align}
                \norm{\bR_{j-1}(z,w_k)}\leq\frac1{|z_-|}\qquad \text{and} \qquad
                |1+\bx^*\bR_{j-1}(z,w_k)\bx|^{-1}\leq\frac{|z|}{|z_-|}.
                \label{eq:helper-path-resolvent-quadratic-bounds}
            \end{align}
            \item[(ii)] (Expected-trace bounds). For $j\in[n]$ and $(z,w)\in\mc D$, the imaginary parts of the leave-one-out trace and expected-trace are positive, i.e.\ $\Im G_j(z,w)>0$ and $\Im\wb G_j(z,w)>0$. In addition, allowing $\eta=0$, along the path, we have for every $j\in[n]$ and $0\leq k\leq n$,
            \begin{align}
                |1+\wb G_j(z,w_k)|\geq\frac{|z_-|}{|z|}.
                \label{eq:helper-denominator-separation}
            \end{align}
            If $t<0$, then also $\Re\wb G_j(z,w_k)\geq0$.
            \item[(iii)] Along the backward path, we also have for $0 \leq k \leq n$
            \begin{align}
                \Im(w_k\wb{z})\geq0 \qquad \text{and} \qquad \Im(w_0\wb{z})\geq \frac{\eta|w_0|^2}{n}.
                \label{eq:helper-path-cone-bound}
            \end{align}
        \end{enumerate}
    \end{subequations}
\end{lemma}
We are now equipped to prove the monotonicity results of the backward path. By Lemma~\ref{lem:resolvent-bounds-strict-positivity}, the backward path is well-defined. For $\eta>0$, applying the recursion and Lemma~\ref{lem:resolvent-bounds-strict-positivity}.(ii) yields
\begin{align*}
    \Im w_{k-1}-\Im w_k =\frac{\Im\wb G_k(z,w_k)}{|1+\wb G_k(z,w_k)|^2}>0.
\end{align*}For $t<0$, applying the same recursion and using $\Re\wb G_k(z,w_k)\geq0$ from Lemma~\ref{lem:resolvent-bounds-strict-positivity}.(ii), and we have similarly
\begin{align*}
    \Re w_k-\Re w_{k-1} =\frac{1+\Re\wb G_k(z,w_k)}{|1+\wb G_k(z,w_k)|^2}>0.
\end{align*}
Combined with $w_n=0$, this shows that the real parts are decreasing, namely $\Re w_0<\Re w_1<\cdots<\Re w_n=0$. For real $z<0$, the backward induction directly implies $w_k\in\R$ for every $k$, and the recursion determines $w_0<w_1<\cdots<w_n=0$. Together with $w_n=0$, these inequalities prove Theorem~\ref{thm:self-consistent-path-stability}.(i).

The remainder of Part I supplies additional bounds for the weighted resolvents $\wt{\bR}_k(z,w):=\bSigma^{1/2}\bR_k(z,w)\bSigma^{1/2}$. For large sample sizes, we compare them with the real resolvent $\bA_k(\zeta):=\bSigma^{1/2}(\bT_k+\zeta\bI)^{-1}\bSigma^{1/2}$ at $\zeta=|z_-|$. For small sample sizes, we use the coefficient $-\Re(e^{\iu\theta}w_j)$ of $\bSigma$ in the rotated inverse resolvent. The following lemma provides both estimates, whose proof is in Appendix~\ref{proof:weighted-resolvent-complex-to-real}.
\begin{lemma}[Weighted resolvent bounds]\label{lem:weighted-resolvent-complex-to-real}
    For any fixed $z=t+\iu\eta\in\C_+\cup\R_{<0}$ and every $0\leq k,j\leq n$, let $\wt{\bR}_k=\wt{\bR}_k(z,w_j)$. Along the backward path, the real resolvents provide the first type of bounds:
    \begin{subequations}
        \begin{align}
            \|\wt{\bR}_k\|\leq\frac{3|z|}{|z_-|} \cdot \|\bA_k(|z_-|)\|,\qquad
            \|\wt{\bR}_k\|_F\leq\frac{3|z|}{|z_-|} \cdot \|\bA_k(|z_-|)\|_F,\qquad
            |\Tr\wt{\bR}_k|\leq\frac{3|z|}{|z_-|} \cdot \Tr\bA_k(|z_-|).
            \label{eq:path-real-comparisons}
        \end{align}
        Set $\theta=\arctan(1+2(t\vee0)/\eta)$ for $\eta>0$ and $\theta=\pi/4$ for real $z<0$. For $j<n$, $-\Re(e^{\iu\theta}w_j)>0$ and, for every integer $1\le m\le d$, the path parameter itself provides the second type of bounds:
        \begin{align}
            \|\wt{\bR}_k\|\le\frac{1}{-\Re(e^{\iu\theta}w_j)}\qquad\text{and}\qquad
            \|\wt{\bR}_k\|_F^2\le\frac{m}{\Re(e^{\iu\theta}w_j)^2}+\frac{4\tau_2(m)}{|z_-|^2}.
            \label{eq:path-buffer-norms}
        \end{align}
    \end{subequations}
\end{lemma}
To use Lemma~\ref{lem:weighted-resolvent-complex-to-real}, Eq.~\eqref{eq:path-real-comparisons}, we need moment bounds for $\bA_k(\zeta)$. We obtain them by separating the first $m$ population eigenvectors from the tail, following the truncation and intrinsic-dimension matrix Bernstein argument of \citet[Lemma~7.2]{cheng2024dimensionfree}. The tail block is controlled by the following real resolvent spectral parameter, expressed using the tail quantities in Eq.~\eqref{eq:tail-parameters-trace}: for $k\geq2$ and $1\leq m<k$,
\begin{align*}
    \const_{\bA}(k,m):=k\sigma_{m+1}+\sqrt{k\sigma_{m+1}\tau_1(m)\log(k)\log\prn{k r_1(m)}}+\tau_1(m)\log(k)\log\prn{k r_1(m)},
\end{align*}
with $\const_{\bA}(k,m)=0$ when $\sigma_{m+1}=0$. The next lemma converts this tail estimate into the required resolvent moments. Its proof is in Appendix~\ref{proof:weighted-resolvent-norm-control}.
\begin{lemma}[$L^p$-type norm bounds for $\bA_k(\zeta)$]\label{lem:weighted-resolvent-norm-control}
    Let Assumption~\ref{assmp:features} hold. With $\rho=\rho(C_{\bx})$ from Assumption~\ref{assmp:spectral-stability}, for every fixed $\kappa,D>0$, $k=\Omega_{C_{\bx},\kappa,D}(1)$, $1\leq m\leq\rho k\wedge d$, $\zeta>0$, and $1\leq p\leq\kappa\log k$, we have
    \begin{subequations}\label{eq:real-resolvent-uniform-moments}
        \begin{align}
            \bigl\|\|\bA_k(\zeta)\|\bigr\|_{L^p(\bT_k)}&\leq O_{C_{\bx}, \kappa, D} \prn{\frac1k\left(1+\frac{\const_{\bA}(k,m)}{\zeta}\right)+\frac{k^{-D}}{\zeta}},\label{eq:real-resolvent-operator-moment}\\
            \bigl\|\|\bA_k(\zeta)\|_F^2\bigr\|_{L^p(\bT_k)}&\leq O_{C_{\bx}, \kappa, D} \prn{\frac{m}{k^2}\left(1+\frac{\const_{\bA}(k,m)}{\zeta}\right)^2+\frac{\tau_2(m)}{\zeta^2}+\frac{\|\bSigma\|_F^2k^{-D}}{\zeta^2}},\label{eq:real-resolvent-frobenius-moment}\\
            \|\Tr\bA_k(\zeta)\|_{L^p(\bT_k)}&\leq O_{C_{\bx}, \kappa, D} \prn{\frac{m}{k}\left(1+\frac{\const_{\bA}(k,m)}{\zeta}\right)+\frac{\tau_1(m)}{\zeta}+\frac{\Tr(\bSigma)k^{-D}}{\zeta}}.\label{eq:real-resolvent-trace-moment}
        \end{align}
    \end{subequations}
\end{lemma}
This wraps up Part I, in which we have proved Theorem~\ref{thm:self-consistent-path-stability}.(i), and Lemmas~\ref{lem:weighted-resolvent-complex-to-real} and~\ref{lem:weighted-resolvent-norm-control} provide the weighted resolvent bounds needed in the next steps. Part II first expresses the one-step error in terms of these norm bounds before we evaluate them in Part III.
\paragraph{Part II: One-step stability by leave-one-out.} We next bound $|\wb G_{k+1}(z,w_{k+1})-\wb G_k(z,w_k)|$ in terms of the operator and Frobenius norms of the weighted resolvents. The leave-one-out calculation uses sample sizes $k-2,k-1,k$ and the two adjacent path parameters. Accordingly, for $k\in[n-1]$, let $\mc I_k:=\{(k-2)\vee0,k-1,k\}\times\{w_k,w_{k+1}\}$ and define
\begin{align*}
    \gamma_k:=\max_{(j,\zeta)\in\mc I_k}\bigl\|\|\wt{\bR}_j(z,\zeta)\|\bigr\|_{L^4(\bT_j)}\qquad\text{and}\qquad\varphi_k:=\max_{(j,\zeta)\in\mc I_k}\bigl\|\|\wt{\bR}_j(z,\zeta)\|_F\bigr\|_{L^4(\bT_j)}.
\end{align*}
The next lemma bounds the one-step error
\begin{align}
    \sconst_{\wb G}(k)=\left|\wb G_{k+1}(z,w_{k+1})-\wb G_k(z,w_k)\right|
\end{align}
in terms of $\gamma_k$ and $\varphi_k$. Its proof, including the leading-order cancellation, is in Appendix~\ref{proof:one-step-G-bound}.
\begin{lemma}\label{lem:one-step-G-bound}
    Let Assumption~\ref{assmp:features} hold. For a fixed $z=t+\iu\eta\in\C_+\cup\R_{<0}$ and all $1\leq k\leq n-1$,
    \begin{subequations}
        \begin{align}
            \sconst_{\wb G}(k) & =  O_{C_{\bx}}\prn{|z|\gamma_k^2/|z_-| + |z|^2(k\gamma_k^3+\gamma_k\varphi_k^2)/|z_-|^2 + |z|^3\prn{\varphi_k^2(\sqrt{k}\gamma_k+\varphi_k)^2+\gamma_k(\sqrt{k}\gamma_k+\varphi_k)^3}/|z_-|^3 }. \label{eq:one-step-explicit-denominator-delta-free}
        \end{align}
        If in addition, for some $\Delta>0$ satisfying $\Delta=\Omega(\sqrt{k}\gamma_k+\varphi_k)$ and $|1+\wb G_{k+1}(z,w_{k+1})|\geq\Delta$, we have
        \begin{align}
            \sconst_{\wb G}(k) =O_{C_{\bx},\Delta}\prn{(\gamma_k+\varphi_k^2)^2+k(\gamma_k+\varphi_k^2)^3}.
            \label{eq:one-step-explicit-denominator-delta}
        \end{align}
    \end{subequations}

\end{lemma}
Heuristically, we can view the ``typical scaling'' above as $\gamma_k=\Theta(\varphi_k^2)$ and $\varphi_k^2=\wt{O}(n^{-1})$, where $\wt O$ suppresses polynomial factors in the path stability parameter $\Lambda_n(z)$. Since $k\leq n$ and $|z|/|z_-|\leq\Lambda_n(z)$, substituting this scaling into Lemma~\ref{lem:one-step-G-bound}, Eq.~\eqref{eq:one-step-explicit-denominator-delta-free}, gives $\sconst_{\wb G}(k)=\wt O_{C_{\bx},\nu}(n^{-2})$, and hence an accumulated error of order $\wt O_{C_{\bx},\nu}(n^{-1})$ over the path, which is the desired conclusion of Theorem~\ref{thm:self-consistent-path-stability}.(ii). We will justify this scaling in Part III uniformly in $k$, using two different approaches in the large and small sample size regimes.

\paragraph{Part III: Proof of (ii), expected-trace stability.}
To prove Theorem~\ref{thm:self-consistent-path-stability}.(ii), we need to compare every $\wb G_k(z,w_k)$ with $\wb G_1(z,w_0)$. The norm estimates from Lemma~\ref{lem:weighted-resolvent-norm-control} above are sufficient to control the moments in Lemma~\ref{lem:one-step-G-bound} when the sample size is on the order of $n$, but they do not apply when the sample size is close to zero. Accordingly, we further divide this part into two steps.
\begin{enumerate}
    \item In Part III-i, we first aggregate the large-sample error terms, and then leverage the resulting trace stability to obtain a uniform lower bound on the rotated coefficient $-\Re(e^{\iu\theta}w_k)$ for all $k\le\ell$.
    \item In Part III-ii, we use this lower bound for the small sample indices, where the upper bounds in Eq.~\eqref{eq:path-buffer-norms} of Lemma~\ref{lem:weighted-resolvent-complex-to-real} provide the remaining norm bounds required for telescoping.
\end{enumerate}
The proof of (ii) then follows by combining the steps above, and accounting separately for the final change from $w_1$ to $w_0$ in $\wb G_1$. To illustrate this strategy, Figure~\ref{fig:two-phase-backward-path} summarizes the two sample-size regimes. Throughout, fix $\nu>0$ as in Assumption~\ref{assmp:spectral-stability} and choose $m$ attaining the minimum in Eq.~\eqref{eq:path-spectral-envelope}. The same $m$ works for every $z$, since $|z|$ and $|z_-|$ are independent of the cutoff. We write
\begin{align}
    \delta_k = \delta_k(z):=\frac{1}{1+\wb G_k(z,w_k)} = w_k-w_{k-1} \label{eq:def-delta-path-parameter}
\end{align}
for the backward increments when $1 \leq k \leq n$. We split the path at $\ell:=\lfloor n/2\rfloor$.

\subparagraph{Part III-i: Large sample sizes and a lower bound for the path parameter.}
We begin by controlling the segment from $n$ to $\ell$, since the cumulative increments there will in turn govern the rest of the segment. To do this, we first need uniform bounds on $\gamma_k$ and $\varphi_k^2$ for $\ell\leq k\leq n-1$, which then allows us to sum the one-step errors from Lemma~\ref{lem:one-step-G-bound}.

For $\ell\leq k\leq n-1$, every index $j\in\{k-2,k-1,k\}$ defining $\gamma_k,\varphi_k$ satisfies $j\geq \ell-2 = \lfloor n/2\rfloor-2 \asymp n$. Thus the cutoff chosen above satisfies $m\leq\rho(C_{\bx})(\ell-2)\leq\rho(C_{\bx})j$ for every such $j$. For fixed $D>0$, apply Lemma~\ref{lem:weighted-resolvent-norm-control}, Eq.~\eqref{eq:real-resolvent-operator-moment}, with $p=4$ for $\gamma_k$, and Eq.~\eqref{eq:real-resolvent-frobenius-moment} with $p=2$ for $\varphi_k^2$. Using $\const_{\bA}(j,m)\leq\const_{\bA}(n,m)$ and then Eq.~\eqref{eq:path-real-comparisons} in Lemma~\ref{lem:weighted-resolvent-complex-to-real} implies

\begin{align*}
    \gamma_k&=O_{C_{\bx},D}\prn{\frac{|z|}{|z_-|}\cdot \brc{\frac1n\prn{1+\frac{\const_{\bA}(n,m)}{|z_-|}}+\frac{n^{-D}}{|z_-|}}},\\
    \varphi_k^2&=O_{C_{\bx},D}\prn{\frac{|z|^2}{|z_-|^2}\cdot \brc{\frac{m}{n^2}\prn{1+\frac{\const_{\bA}(n,m)}{|z_-|}}^2+\frac{\tau_2(m)}{|z_-|^2}+\frac{\|\bSigma\|_F^2n^{-D}}{|z_-|^2}}}.
\end{align*}
To express these bounds in terms of $\Lambda_n(z)$, we need to control $\const_{\bA}(n,m)$ and $\tau_2(m)$ by $\Lambda_n(z)$. If $\sigma_{m+1}=0$, then $\const_{\bA}(n,m)=\tau_2(m)=0$ and both comparisons below follow from $\Lambda_n(z)\ge1$. Otherwise $r_1(m)\ge1$ and $\log(n)\log(nr_1(m))\ge1$ for $n\ge3$. The middle term $\sqrt{n\sigma_{m+1}\tau_1(m)\log(n)\log\prn{nr_1(m)}}$ in $\const_{\bA}(n,m)$ is the geometric mean of its first and last terms. We then use $\tau_2(m)\leq\sigma_{m+1}\tau_1(m)$ and $2\sqrt{ab}\leq a+b$ in the definitions of $\const_{\bA}(n,m)$ and $\Lambda_n(z)$ to conclude that
\begin{align}
    \frac{n\tau_2(m)}{|z_-|^2}\leq\Lambda_n(z)^2\qquad \text{and} \qquad
    1+\frac{\const_{\bA}(n,m)}{|z_-|}\leq\frac32\Lambda_n(z).
    \label{eq:path-stability-parameter-controls}
\end{align}
We then choose $D$ so that the last remainder terms in $\gamma_k$ and $\varphi_k^2$ preserve the $n^{-1}$ scale of the norm estimates inside the brackets. Since $\|\bSigma\|=1$, $\|\bSigma\|_F^2\leq\Tr\bSigma\leq m+\tau_1(m)\leq n+\Lambda_n(z)|z_-|$. Using $|z_-|\geq n^{-\nu}$, it follows that
\begin{align*}
    \frac{n^{-D}}{|z_-|}\leq n^{\nu-D}\qquad \text{and} \qquad
    \frac{\|\bSigma\|_F^2n^{-D}}{|z_-|^2}\leq n^{1+2\nu-D}+\Lambda_n(z)n^{\nu-D}.
\end{align*}
Taking $D=2\nu+2$ makes $n^{1+2\nu-D}=n^{-1}$, so both remainders are $O(\Lambda_n(z)/n)$. Since all implicit constants from this point on depend only on $C_{\bx}$ and $\nu$, we will suppress them from the big-Oh notations for simplicity. Consequently, by collecting all the preceding bounds, for all $\ell\leq k\leq n-1$, we have
\begin{align}
    \gamma_k=O\prn{\frac{\Lambda_n(z)^2}{n}}\qquad \text{and} \qquad
    \varphi_k^2=O\prn{\frac{\Lambda_n(z)^4}{n}}.
    \label{eq:path-large-moments}
\end{align}
Since $k\leq n$, these two bounds give $\sqrt{k}\gamma_k+\varphi_k=O(\Lambda_n(z)^2/\sqrt n)$. Substituting the preceding bounds into Eq.~\eqref{eq:one-step-explicit-denominator-delta-free} in Lemma~\ref{lem:one-step-G-bound} implies the following without the additional denominator stability assumption on $\Delta$:
\begin{align}
    \sconst_{\wb G}(k) & = O\prn{\frac{|z|\gamma_k^2}{|z_-|} + \frac{|z|^2(k\gamma_k^3+\gamma_k\varphi_k^2)}{|z_-|^2} + \frac{|z|^3\prn{\varphi_k^2(\sqrt{k}\gamma_k+\varphi_k)^2+\gamma_k(\sqrt{k}\gamma_k+\varphi_k)^3}}{|z_-|^3} } \nonumber \\
    & =  O\prn{\frac{\Lambda_n(z)^5}{n^2}+ \frac{\Lambda_n(z)^8}{n^2} + \frac{\Lambda_n(z)^{11}}{n^2}} =O \prn{\frac{\Lambda_n(z)^{11}}{n^2}}
    \label{eq:path-large-one-step-rate}
\end{align}
Summing $\wb G_{k+1}(z,w_{k+1})-\wb G_k(z,w_k)$ from $k=j$ to $n-1$ and applying Eq.~\eqref{eq:path-large-one-step-rate} to the $n-j$ summands, with $w_n=0$, yields, for $\ell\leq j\leq n$,
\begin{align}
    |\wb G_j(z,w_j)-\wb G_n(z,0)|\leq\sum_{k=j}^{n-1}\sconst_{\wb G}(k)=O\prn{\frac{\Lambda_n(z)^{11}}n}.
    \label{eq:path-large-increment-stability}
\end{align}
For the small-sample segment, Lemma~\ref{lem:weighted-resolvent-complex-to-real}, Eq.~\eqref{eq:path-buffer-norms}, reduces the norm estimates to a lower bound on $-\Re(e^{\iu\theta}w_k)$ for $k\le\ell$. Since $w_k=-\sum_{j=k+1}^n\delta_j$, it suffices to show that every rotated increment has positive real part and that the last $n-\ell$ increments contribute at least order $n/\Lambda_n(z)^3$ in total.

Fix $j\in[n]$ and use the rotation in the proof of Lemma~\ref{lem:weighted-resolvent-complex-to-real}, Eq.~\eqref{eq:weighted-common-rotation}, at $\bU=e^{\iu\theta}\bR_{j-1}(z,w_j)^{-1}$. Since $\Re\bU\succ0$, we have $\Re(e^{-\iu\theta}\bR_{j-1}(z,w_j))=\bU^{-*}(\Re\bU)\bU^{-1}\succeq0$. Taking the $\bSigma$-weighted trace and expectation gives $\Re(e^{-\iu\theta}\wb G_j(z,w_j))\ge0$. Using $\delta_j=(1+\wb G_j(z,w_j))^{-1}$ therefore yields
\begin{align*}
    \Re(e^{\iu\theta}\delta_j)=\frac{\cos\theta+\Re(e^{-\iu\theta}\wb G_j(z,w_j))}{|1+\wb G_j(z,w_j)|^2}\ge\cos\theta\,|\delta_j|^2>0.
\end{align*}
Since we need a lower bound on $-\Re(e^{\iu\theta}w_k)$ for $k\leq\ell$, we retain the last $n-\ell$ terms in its sum. The preceding inequality reduces this to bounding $|\delta_j|$ from below, or equivalently $|\delta_j^{-1}|$ from above, for $j>\ell$. By Eq.~\eqref{eq:def-delta-path-parameter}, $\delta_j^{-1}-\delta_n^{-1}=\wb G_j(z,w_j)-\wb G_n(z,0)$. Applying the triangle inequality and Eq.~\eqref{eq:path-large-increment-stability} yields $|\delta_j^{-1}|\leq|\delta_n^{-1}|+O(\Lambda_n(z)^{11}/n)$, so it suffices to bound the terminal reciprocal.

Since $w_n=0$, we can compare the terminal trace directly with the real resolvent at $|z|$. This choice cancels the factor $|z|$ in the comparison against the denominators $\zeta=|z|$ in the real trace bound, leaving a single factor $|z_-|^{-1}$. For $x\geq0$, we have $(x+|z|)/|x-z|\leq1+2|z|/|z_-|\leq3|z|/|z_-|$. Applying this inequality in the eigenbasis of $\bT_{n-1}$ and then applying Lemma~\ref{lem:weighted-resolvent-norm-control}, Eq.~\eqref{eq:real-resolvent-trace-moment}, at $\zeta=|z|$, sample size $n-1$, $p=1$ and $D=2\nu+2$ yields
\begin{align}
    |\delta_n^{-1}|=|1+\wb G_n(z,0)|
    &\stackrel{\mathrm{(i)}}{\leq}1+\frac{3|z|}{|z_-|}\E\Tr\bA_{n-1}(|z|)\nonumber\\
    &\stackrel{\mathrm{(ii)}}{\leq}O\prn{\frac{|z|+\const_{\bA}(n,m)+\tau_1(m)+n^{-D}\Tr\bSigma}{|z_-|}}
    =O(\Lambda_n(z)).\label{eq:path-terminal-reciprocal-bound}
\end{align}
In (i), the trace assigns nonnegative weights from $\bSigma$ to the eigenvalues of $\bT_{n-1}$. In (ii), multiplying Lemma~\ref{lem:weighted-resolvent-norm-control}, Eq.~\eqref{eq:real-resolvent-trace-moment}, by $|z|$ yields
\begin{align*}
    |z|\E\Tr\bA_{n-1}(|z|)
    =O\prn{\frac{m}{n-1}\prn{|z|+\const_{\bA}(n-1,m)}+\tau_1(m)+(n-1)^{-D}\Tr\bSigma}.
\end{align*}
We use $m\leq n-1$, $\const_{\bA}(n-1,m)\leq\const_{\bA}(n,m)$. The last bound follows from Eq.~\eqref{eq:path-stability-parameter-controls}, $\Tr\bSigma\leq m+\tau_1(m)\leq n+\Lambda_n(z)|z_-|$ and $|z_-|\geq n^{-\nu}$, as in the choice of $D$ after Eq.~\eqref{eq:path-stability-parameter-controls}. Using $n=\Omega(\Lambda_n(z)^{16})$ and $\Lambda_n(z)\ge1$ in Eqs.~\eqref{eq:path-large-increment-stability} and~\eqref{eq:path-terminal-reciprocal-bound}, we have
\begin{align*}
    \max_{\ell<j\le n}|\delta_j^{-1}|\le |\delta_n^{-1}|+\max_{\ell<j\le n}|\wb G_j(z,w_j)-\wb G_n(z,0)|=O\prn{\Lambda_n(z)+\Lambda_n(z)^{11}/n}=O(\Lambda_n(z)).
\end{align*}
For every $k\le\ell$, positivity allows us to retain just the last $n-\ell$ terms. The preceding bound gives $\min_{\ell<j\le n}|\delta_j|=\Omega(1/\Lambda_n(z))$. Using $\cos\theta\ge |z_-|/(3|z|)\ge 1/(3\Lambda_n(z))$, as in Part II of the proof of Lemma~\ref{lem:weighted-resolvent-complex-to-real}, and $n-\ell\ge n/2$ yields, for $0\le k\le\ell$,
\begin{align}
    -\Re(e^{\iu\theta}w_k)&=\sum_{j=k+1}^n\Re(e^{\iu\theta}\delta_j)\ge\frac{n-\ell}{3\Lambda_n(z)}\min_{\ell<j\le n}|\delta_j|^2=\Omega\prn{\frac{n}{\Lambda_n(z)^3}}.\label{eq:path-imaginary-buffer-rate}
\end{align}

\subparagraph{Part III-ii: Small sample sizes and stability of the expected traces.}
Our goal is now to bound every remaining one-step error and then sum along the entire path. For $1\leq k<\ell$, the weighted resolvents defining $\gamma_k$ and $\varphi_k$ use path indices $j\in\{k,k+1\}\subseteq\{0,\ldots,\ell\}$.
Thus Eq.~\eqref{eq:path-imaginary-buffer-rate} applies to both. Applying Eq.~\eqref{eq:path-buffer-norms} in Lemma~\ref{lem:weighted-resolvent-complex-to-real}, we have the following estimates as
\begin{align*}
    \gamma_k&=O\prn{\frac{\Lambda_n(z)^{3}}{n}} \qquad \text{and} \qquad
    \varphi_k^2=O\prn{\frac{m\Lambda_n(z)^{6}}{n^2}+\frac{\tau_2(m)}{|z_-|^2}}=O\prn{\frac{\Lambda_n(z)^{6}}n},
\end{align*}
where we use $m\le n$ and Eq.~\eqref{eq:path-stability-parameter-controls}.
In particular, this implies $\sqrt{k}\gamma_k+\varphi_k=O(\Lambda_n(z)^3/\sqrt n)$, and these bounds for $\gamma_k$, $\varphi_k$ and $\sqrt{k}\gamma_k+\varphi_k$ combined with Eq.~\eqref{eq:one-step-explicit-denominator-delta-free} in Lemma~\ref{lem:one-step-G-bound} imply that
\begin{align*}
    \sconst_{\wb G}(k) & =  O\prn{\frac{\Lambda_n(z)^{7}}{n^2}+ \frac{\Lambda_n(z)^{11}}{n^2} + \frac{\Lambda_n(z)^{15}}{n^2}} =O \prn{\frac{\Lambda_n(z)^{15}}{n^2}}
\end{align*}
Combining this estimate with the large-sample bound $\sconst_{\wb G}(k)=O(\Lambda_n(z)^{11}/n^2)$ from Part III-i, Eq.~\eqref{eq:path-large-one-step-rate}, and using $\Lambda_n(z)\geq1$, we obtain $\sconst_{\wb G}(k)=O(\Lambda_n(z)^{15}/n^2)$ for every $1\leq k\leq n-1$. Summing the increments between any two indices proves
\begin{align}
    \max_{1\leq i,j\leq n}|\wb G_i(z,w_i)-\wb G_j(z,w_j)|\leq\sum_{k=1}^{n-1}\sconst_{\wb G}(k)=O\prn{\frac{\Lambda_n(z)^{15}}n}.
    \label{eq:path-G-telescope-rate}
\end{align}
To establish Theorem~\ref{thm:self-consistent-path-stability}.(ii), which also compares $\wb G_k(z,w_k)$ with $\wb G_1(z,w_0)$, we also need to bound the offset $\wb G_1(z,w_1)-\wb G_1(z,w_0)$. For this last change of parameter, we use the same path lower bound at the two indices $0$ and $1$. Applying Lemma~\ref{lem:weighted-resolvent-complex-to-real}, Eq.~\eqref{eq:path-buffer-norms}, with sample index $0$ and using Eqs.~\eqref{eq:path-stability-parameter-controls} and~\eqref{eq:path-imaginary-buffer-rate} yields, for $i\in\{0,1\}$,
\begin{align*}
    \|\wt{\bR}_0(z,w_i)\|_F^2\leq\frac{m}{\Re(e^{\iu\theta}w_i)^2}+\frac{4\tau_2(m)}{|z_-|^2}
    =O\prn{\frac{\Lambda_n(z)^{6}}n}.
\end{align*}
Here $|\delta_1|\leq\Lambda_n(z)$, and therefore by trace Cauchy-Schwarz,
\begin{align*}
    |\wb G_1(z,w_1)-\wb G_1(z,w_0)|
    &\leq|\delta_1|\,\|\wt{\bR}_0(z,w_1)\|_F\|\wt{\bR}_0(z,w_0)\|_F
    =O\prn{\frac{\Lambda_n(z)^{7}}n}.
\end{align*}
Combining the bound for $|\wb G_1(z,w_1)-\wb G_1(z,w_0)|$ with Eq.~\eqref{eq:path-G-telescope-rate} yields

\begin{align}
    \max_{1\leq k\leq n}|\wb G_k(z,w_k)-\wb G_1(z,w_0)|=O\prn{\frac{\Lambda_n(z)^{15}}n}.
    \label{eq:path-G-population-endpoint}
\end{align}
Combining with taking $j=n$ in Eq.~\eqref{eq:path-G-telescope-rate} completes the proof of (ii).

\paragraph{Part IV: Proof of (iii), path stability.} To prove (iii), we use Eq.~\eqref{eq:path-G-population-endpoint} to compare backward increments with $(1+\wb G_1(z,w_0))^{-1}$. Summing yields an approximate equation for $w_0$, and Lemma~\ref{lem:population-scalar-stability} then compares $w_0$ and $w_0\lin$. The proof concludes by summing the increment differences.

\subparagraph{Part IV-i: Stable increments and the endpoint equation.} Let $\delta_0:=(1+\wb G_1(z,w_0))^{-1}$. Applying Lemma~\ref{lem:resolvent-bounds-strict-positivity}.(ii), Eq.~\eqref{eq:helper-denominator-separation}, we have $\max_{0\leq j\leq n}|\delta_j|\leq\Lambda_n(z)$. Using the reciprocal identity and Eq.~\eqref{eq:path-G-population-endpoint}, we obtain, for $1\leq j\leq n$,
\begin{align}
    \delta_0-\delta_j&=\delta_0\delta_j\prn{\wb G_j(z,w_j)-\wb G_1(z,w_0)} \qquad \text{and} \qquad
    |\delta_j-\delta_0|=O\prn{\frac{\Lambda_n(z)^{17}}n}.
    \label{eq:path-increment-common-value}
\end{align}
Here $|\delta_0\delta_j|\leq\Lambda_n(z)^2$, while Eq.~\eqref{eq:path-G-population-endpoint} bounds the trace difference by $O(\Lambda_n(z)^{15}/n)$ uniformly in $j$. It then follows from summing these differences, with $\big|\sum_{j=1}^n(\delta_0-\delta_j)\big|\leq n\max_{1\leq j\leq n}|\delta_0-\delta_j|$, and using $w_0=-\sum_{j=1}^n\delta_j$ from Eq.~\eqref{eq:def-delta-path-parameter} and $w_n=0$ that
\begin{align}
    w_0+n\delta_0=\sum_{j=1}^n(\delta_0-\delta_j)=O(\Lambda_n(z)^{17}).
    \label{eq:path-population-residual-rate}
\end{align}
We next convert this residual into a perturbation of the effective spectral parameter equation, which will allow us to compare $w_0$ with $w_0\lin$.

\subparagraph{Part IV-ii: Comparison with the population endpoint.} To compare $w_0$ with $w_0\lin=nzs_n(z)$, we will identify $w_0/(nz)$ with $s_n$ using an approximate fixed-point argument and apply Lemma~\ref{lem:population-scalar-stability}. Indeed, $w_0\lin/(nz)$ satisfies the exact fixed-point equation in Eq.~\eqref{eq:population-scalar-equation},
\begin{align*}
    z=-\prn{\frac{w_0\lin}{nz}}^{-1}+\sum_{j\geq1}\frac{\sigma_j}{1+\sigma_jw_0\lin/z}.
\end{align*}
This boils down to controlling the perturbed equation, which we will identify next.

By Eq.~\eqref{eq:path-imaginary-buffer-rate}, we have $w_0\neq0$. For $z\in\C_+$, applying Lemma~\ref{lem:resolvent-bounds-strict-positivity}.(iii) gives $\Im(w_0/(nz))\geq\eta|w_0/(nz)|^2>0$.
Since $\bT_0=0$, the trace $\wb G_1(z,w_0)$ is deterministic and satisfies $\wb G_1(z,w_0)=-z^{-1}\sum_{j\geq1}\sigma_j/(1+\sigma_jw_0/z)$. Thus $\delta_0^{-1}=1+\wb G_1(z,w_0)$ implies $\sum_{j\geq1}\sigma_j/(1+\sigma_jw_0/z)=z(1-\delta_0^{-1})$. Subtracting $nz/w_0$ from both sides of the preceding equation and using $z(1-\delta_0^{-1})-nz/w_0=z-z(w_0+n\delta_0)/(w_0\delta_0)$ yields the perturbed equation
\begin{align}
    z+\varepsilon_n(z)&=-\prn{\frac{w_0}{nz}}^{-1}+\sum_{j\geq1}\frac{\sigma_j}{1+\sigma_jw_0/z}\qquad \text{where} \qquad
    \varepsilon_n(z):=-\frac{z}{w_0\delta_0}(w_0+n\delta_0).
    \label{eq:path-perturbed-population-equation}
\end{align}
The residual $\varepsilon_n(z) = -z/(w_0 \delta_0) \cdot \sum_{j=1}^n (\delta_0-\delta_j)$ by Eq.~\eqref{eq:path-population-residual-rate}, which can be controlled by trace differences already bounded in Part III. From Eq.~\eqref{eq:path-increment-common-value},
\begin{align*}
    \varepsilon_n(z)=-\frac{z}{w_0}\sum_{j=1}^n\delta_j\prn{\wb G_j(z,w_j)-\wb G_1(z,w_0)}.
\end{align*}
To bound the ratio of $\sum_j|\delta_j|$ to $|w_0|$, compare each increment with $\delta_n$ before summing. By Lemma~\ref{lem:resolvent-bounds-strict-positivity}.(ii), Eq.~\eqref{eq:helper-denominator-separation}, we have $|\delta_j|\leq\Lambda_n(z)$, and Eq.~\eqref{eq:path-G-telescope-rate} bounds the trace difference with the terminal index by $O(\Lambda_n(z)^{15}/n)$. Thus
\begin{align*}
    \max_{1\le j\le n}\left|\frac{\delta_j}{\delta_n}-1\right|
    &=\max_{1\le j\le n}\left|\delta_j\prn{\wb G_n(z,0)-\wb G_j(z,w_j)}\right|=O\prn{\frac{\Lambda_n(z)^{16}}n}\le\frac12.
\end{align*}
The last inequality follows by choosing the implicit constant in Assumption~\ref{assmp:spectral-stability}. Since $w_0=-\sum_{j=1}^n\delta_j$, we have $|w_0+n\delta_n|\le n|\delta_n|/2$, so $|w_0|\ge n|\delta_n|/2$ and $\sum_{j=1}^n|\delta_j|\le3n|\delta_n|/2$. Hence $\sum_{j=1}^n|\delta_j|/|w_0|\leq3$. Applying this bound and Eq.~\eqref{eq:path-G-population-endpoint} to the preceding identity yields
\begin{align}
    |\varepsilon_n(z)|&\leq\frac{|z|}{|w_0|}\sum_{j=1}^n|\delta_j|\max_{1\leq j\leq n}|\wb G_j(z,w_j)-\wb G_1(z,w_0)| =O\prn{\frac{|z|\Lambda_n(z)^{15}}n} =O\prn{\frac{|z_-|\Lambda_n(z)^{16}}n}. \hspace{1em}
    \label{eq:path-spectral-perturbation-rate}
\end{align}
By Eq.~\eqref{eq:path-spectral-perturbation-rate}, choosing the implicit constant in $n=\Omega(\Lambda_n(z)^{16})$ sufficiently large ensures $|\varepsilon_n(z)|\leq|z_-|/4$. For $z\in\C_+$, since $w_0/(nz)\in\C_+$ and $\bSigma\succeq0$ is trace-class, assuming $z+\varepsilon_n(z)\in\C_+$ and applying uniqueness in Lemma~\ref{lem:population-scalar-stability} to Eq.~\eqref{eq:path-perturbed-population-equation}, we have $w_0/(nz)=s_n(z+\varepsilon_n(z))$ and further
\begin{align*}
    |(z+\varepsilon_n(z))_-|=\dist(z+\varepsilon_n(z),\R_{\geq0})\geq|z_-|-|\varepsilon_n(z)|\geq\frac34|z_-|.
\end{align*}
The fixed-point stability bound in the same lemma at $z$ and $z+\varepsilon_n(z)$ implies then
\begin{align*}
    \left|\frac{w_0}{nz}-s_n(z)\right|
    &=|s_n(z+\varepsilon_n(z))-s_n(z)|
    \leq\frac{|\varepsilon_n(z)|}{|(z+\varepsilon_n(z))_-|\,|z_-|}
    \leq\frac{4|\varepsilon_n(z)|}{3|z_-|^2}.
\end{align*}
By $w_0\lin=nzs_n(z)$, Eq.~\eqref{eq:path-spectral-perturbation-rate}, and $|z|/|z_-|\leq\Lambda_n(z)$, we can control the terminal path difference by
\begin{align}
    |w_0-w_0\lin|&=n|z|\left|\frac{w_0}{nz}-s_n(z)\right|
    =O\prn{\frac{|z|}{|z_-|}\Lambda_n(z)^{16}}
    =O(\Lambda_n(z)^{17}).
    \label{eq:path-population-comparison-rate}
\end{align}
It thus remains only to verify $z+\varepsilon_n(z)\in\C_+$, which can be shown by the following short proof by contradiction. Suppose that $\Im(z+\varepsilon_n(z))\leq0$. Since $\eta\leq|\varepsilon_n(z)|<|z_-|$, we must have $t<0$. Moreover, $\Re(z+\varepsilon_n(z))<0$, since otherwise $|\varepsilon_n(z)|^2\geq t^2+\eta^2=|z_-|^2$, yielding a contradiction. By the path signs from Part I, $\Re(w_0/(nz))=(t\Re w_0+\eta\Im w_0)/(n|z|^2)>0$. Together with $\Im(w_0/(nz))>0$, these signs imply
\begin{align*}
    \Im((z+\varepsilon_n(z))w_0/(nz)) = \Im(z+\varepsilon_n(z)) \Re(w_0/(nz)) + \Re(z+\varepsilon_n(z)) \Im(w_0/(nz))<0.
\end{align*} However, multiplying Eq.~\eqref{eq:path-perturbed-population-equation} by $w_0/(nz)$ and taking imaginary parts implies
\begin{align*}
    \Im\prn{(z+\varepsilon_n(z))w_0/(nz)}
    =\Im\prn{w_0/(nz)} \cdot \sum_{j\geq1}\frac{\sigma_j}{|1+\sigma_jw_0/z|^2}\geq0,
\end{align*}
yielding a contradiction. Thus $z+\varepsilon_n(z)\in\C_+$, which proves Eq.~\eqref{eq:path-population-comparison-rate} for $z\in\C_+$.

For $z=t<0$, we have from Part I $w_0<0$ and $\delta_0>0$. Hence $w_0/(nt)>0$ and $\varepsilon_n(t)\in\R$. Since $|\varepsilon_n(t)|\leq|t|/4$, we also have $t+\varepsilon_n(t)<0$. Applying the uniqueness of the positive solution in Lemma~\ref{lem:population-scalar-stability} to Eq.~\eqref{eq:path-perturbed-population-equation} identifies $w_0/(nt)=s_n(t+\varepsilon_n(t))$, and subsequently, we obtain from its difference bound that
\begin{align*}
    |w_0-w_0\lin|
    &=n|t|\,|s_n(t+\varepsilon_n(t))-s_n(t)|
    \leq\frac{n|\varepsilon_n(t)|}{|t+\varepsilon_n(t)|}
    \leq\frac{4n|\varepsilon_n(t)|}{3|t|}
    =O\prn{\Lambda_n(t)^{16}},
\end{align*}
where we use Eq.~\eqref{eq:path-spectral-perturbation-rate} and $|z_-|=|t|$. Since $\Lambda_n(t)\geq1$, Eq.~\eqref{eq:path-population-comparison-rate} also follows for $z=t<0$.

We also transfer the terminal reciprocal bound in Eq.~\eqref{eq:path-terminal-reciprocal-bound} to the constant increment $\delta\lin(z):=-w_0\lin/n=-zs_n(z)$ of the linear path. For both $z\in\C_+$ and $z<0$, we have just identified $w_0/(nz)=s_n(z+\varepsilon_n(z))$. Rearranging the definition of $\varepsilon_n(z)$ in Eq.~\eqref{eq:path-perturbed-population-equation} gives $-n/w_0=\delta_0^{-1}+\varepsilon_n(z)/z$. Substitution yields the exact identity
\begin{equation}
    \frac1{\delta\lin(z)}=\frac{s_n(z+\varepsilon_n(z))}{s_n(z)}\left(\frac1{\delta_0}+\frac{\varepsilon_n(z)}{z}\right).\label{eq:path-linear-reciprocal-transfer}
\end{equation}
To bound the multiplier, write $s_n(z)=\int(x-z)^{-1}\mu(\de x)$ as in Lemma~\ref{lem:population-scalar-stability}. Since $|\varepsilon_n(z)|\leq|z_-|/4$ by Eq.~\eqref{eq:path-spectral-perturbation-rate}, subtracting the two Stieltjes transforms yields
\begin{align*}
    \left|\frac{s_n(z+\varepsilon_n(z))}{s_n(z)}-1\right|
    \leq\frac{4|\varepsilon_n(z)|}{3|s_n(z)|}\int\frac{\mu(\de x)}{|x-z|^2}
    \leq\frac{4|\varepsilon_n(z)|}{3|z_-|}\leq\frac13.
\end{align*}
For the second inequality, $\Im s_n(z)=\eta\int |x-z|^{-2}\mu(\de x)$ suffices when $t\geq0$. When $t<0$, also use $\Re s_n(z)\geq(-t)\int |x-z|^{-2}\mu(\de x)$ and $|z_-|=\sqrt{t^2+\eta^2}$. Finally, Eq.~\eqref{eq:path-G-population-endpoint} gives $|\delta_0^{-1}-\delta_n^{-1}|=O(\Lambda_n(z)^{15}/n)=O(1)$. Substituting into Eq.~\eqref{eq:path-linear-reciprocal-transfer} and using $|\varepsilon_n(z)/z|\leq1/4$ and Eq.~\eqref{eq:path-terminal-reciprocal-bound} yields $|\delta\lin(z)|^{-1}=O(\Lambda_n(z))$.

\subparagraph{Part IV-iii: Comparison of the paths.}
We are now ready to finish the proof of Theorem~\ref{thm:self-consistent-path-stability}.(iii) for path comparisons when $1 \leq k \leq n$. To compare the path parameters, we first compare each increment $\delta_j$ with $-w_0\lin/n$, i.e.\ the constant backward increment of $w_k\lin$. By Eqs.~\eqref{eq:path-increment-common-value}, \eqref{eq:path-population-residual-rate}, and~\eqref{eq:path-population-comparison-rate}, and $\Lambda_n(z)\geq1$, we have for $1 \leq j \leq n$,
\begin{align*}
    \left|\delta_j+\frac{w_0\lin}{n}\right|
    &\leq|\delta_j-\delta_0|+\frac{|w_0+n\delta_0|+|w_0-w_0\lin|}{n}
    =O\prn{\frac{\Lambda_n(z)^{17}}n}.
\end{align*}
Since $w_k=-\sum_{j=k+1}^n\delta_j$ and $w_k\lin=(n-k)w_0\lin/n$, summing from $j=k+1$ to $n$ yields
\begin{align*}
    |w_k-w_k\lin|
    =\left|\sum_{j=k+1}^n\prn{\delta_j+\frac{w_0\lin}{n}}\right|
    \leq\sum_{j=k+1}^n\left|\delta_j+\frac{w_0\lin}{n}\right|
    =O\prn{\Lambda_n(z)^{17} \cdot \frac{n-k}{n}},
\end{align*}
for every $0 \leq k \leq n$. The last display proves (iii) for $z\in\C_+\cup\R_{<0}$. With Parts I and III, the proof of Theorem~\ref{thm:self-consistent-path-stability} is complete.

\subsection{Proof of Theorem~\ref{thm:terminal-weighted-traces}}\label{proof:terminal-weighted-traces}
We track the random variables $H_k^{\bQ}$ and $G_k^{\bQ}$ along the backward path to their deterministic equivalents. We will bound the two fluctuations below, including the additional correction from $w_0$ to $w_0\lin$.
\paragraph{Part I: The exact path decomposition.} We first identify this decomposition in a one-step update. Each trace increment combines a rank-one sample update with a change in the parameter. We first isolate the sample update by fixing $w$, since the conditional centering of this update will produce the martingale difference. For $\bR=\bR_{k-1}(z,w)$, the Sherman-Morrison formula identifies the corresponding trace updates as
\begin{align}
    V_k^{H,\bQ}(w):=\frac{z\bx_k^*\bR\bQ\bR\bx_k}{1+\bx_k^*\bR\bx_k}\qquad \text{and} \qquad
    V_k^{G,\bQ}(w):=\frac{\bx_k^*\bR\bSigma^{1/2}\bQ\bSigma^{1/2}\bR\bx_k}{1+\bx_k^*\bR\bx_k}.
    \label{eq:trace-rank-one-ratios}
\end{align}
Because the backward path is deterministic, the sample updates can be analyzed in the forward filtration $\{\mc F_k\}$ where $\bx_k$ is independent of $\mc{F}_{k-1}$, and all path parameters are fixed. Writing $\delta_k=w_k-w_{k-1}=(1+\wb G_k(z,w_k))^{-1}$, we obtain the following decomposition. Its proof is in Appendix~\ref{proof:trace-step-decomposition}.

\begin{lemma}[One-step trace decomposition]\label{lem:trace-step-decomposition}
    Under Assumption~\ref{assmp:features}, along the backward path $\{w_k\}$,
    \begin{subequations}\label{eq:trace-step-decomposition}
        \begin{align}
            H_k^{\bQ}(z,w_k)-H_{k-1}^{\bQ}(z,w_{k-1})&=\xi_k^{H,\bQ}+\beta_k^{H,\bQ}, &&1\leq k\leq n,\label{eq:trace-step-H}\\
            G_{k+1}^{\bQ}(z,w_{k+1})-G_k^{\bQ}(z,w_k)&=\xi_k^{G,\bQ}+\beta_k^{G,\bQ}, &&1\leq k<n,\label{eq:trace-step-G}
        \end{align}
    \end{subequations}
    where
    \begin{align*}
        \xi_k^{H,\bQ}&=\E_{k-1}V_k^{H,\bQ}(w_k)-V_k^{H,\bQ}(w_k) \qquad \text{and} \qquad
        \xi_k^{G,\bQ}=\E_{k-1}V_k^{G,\bQ}(w_{k+1})-V_k^{G,\bQ}(w_{k+1}),
    \end{align*}
    and
    \begin{align}
        \beta_k^{H,\bQ}&=z\delta_k\Tr\brk{\bQ\bR_{k-1}(z,w_k)\bSigma\bR_{k-1}(z,w_{k-1})}-\E_{k-1}V_k^{H,\bQ}(w_k),\nonumber\\
        \beta_k^{G,\bQ}&=\delta_{k+1}\Tr\brk{\bQ\bSigma^{1/2}\bR_{k-1}(z,w_{k+1})\bSigma\bR_{k-1}(z,w_k)\bSigma^{1/2}}-\E_{k-1}V_k^{G,\bQ}(w_{k+1}),
        \label{eq:trace-predictable-drifts}
    \end{align}
    Both the $\{\xi_k^{H,\bQ}\}$ and $\{\xi_k^{G,\bQ}\}$ sequences are martingale differences for $\{\mc F_k\}$, and both $\{\beta_k^{H,\bQ}\}$ and $\{\beta_k^{G,\bQ}\}$ sequences (drift terms) are $\mc F_{k-1}$-measurable.
\end{lemma}
We can now telescope from the backward initial point $w_0$, and further connect to the effective spectral parameter initial point (i.e., linear interpolation) $w_0\lin$ by including the additional corrections $\beta_0^{H, \bQ} := H_0^{\bQ}(z, w_0) - H_0^{\bQ}(z, w_0\lin)$ and $\beta_0^{G, \bQ} := G_1^{\bQ}(z, w_1) - G_1^{\bQ}(z, w_0\lin)$, which implies
\begin{subequations}\label{eq:trace-path-decomposition}
    \begin{align}
        H_k^{\bQ}(z,w_k)-H_0^{\bQ}(z,w_0\lin)&=\sum_{j=1}^k\xi_j^{H,\bQ}+ \sum_{j=0}^k \beta_j^{H,\bQ}, &&0\leq k\leq n,\label{eq:trace-path-H}\\
        G_k^{\bQ}(z,w_k)-G_1^{\bQ}(z,w_0\lin)&=\sum_{j=1}^{k-1}\xi_j^{G,\bQ}+ \sum_{j=0}^{k-1} \beta_j^{G,\bQ}, &&1\leq k\leq n.\label{eq:trace-path-G}
    \end{align}
\end{subequations}
\paragraph{Part II: Controlling the two error terms.} The backward path $\{w_k\}_{k=0}^n$ allows us to apply Theorem~\ref{thm:self-consistent-path-stability} to bound the drift and martingale contributions in the exact decompositions above. The two-phase argument in its proof additionally supplies resolvent moments for all sample indices, and the endpoint comparison in Theorem~\ref{thm:self-consistent-path-stability}.(iii) controls the correction from $w_0$ to $w_0\lin$.
\begin{lemma}\label{lem:trace-accumulated-drift}
    Under the same setup of Theorem~\ref{thm:terminal-weighted-traces}, we control the errors in Eq.~\eqref{eq:trace-path-decomposition}.
    \begin{subequations}
        \begin{enumerate}
            \item[(i)]The martingale partial sums obey
            \begin{align}
                \left\|\max_{0\leq k\leq n}\frac1n\left|\sum_{j=1}^k\xi_j^{H,\bQ}\right|\right\|_{L^p}
                \vee\left\|\max_{1\leq k\leq n}\left|\sum_{j=1}^{k-1}\xi_j^{G,\bQ}\right|\right\|_{L^p}
                &=O_{C_{\bx},\nu,\kappa}\prn{\frac{p^2\Lambda_n(z)^{11}}n \cdot \prn{\|\bQ\| \wedge \frac{\|\bQ\|_F}{\sqrt n}}}.
                \label{eq:trace-martingale}
            \end{align}
            \item[(ii)] The accumulated drifts, including the initial corrections, obey
            \begin{align}
                \left\|\max_{0\leq k\leq n}\frac1n\left|\sum_{j=0}^k\beta_j^{H,\bQ}\right|\right\|_{L^p}
                \vee\left\|\max_{1\leq k\leq n}\left|\sum_{j=0}^{k-1}\beta_j^{G,\bQ}\right|\right\|_{L^p}
                &=O_{C_{\bx},\nu,\kappa}\prn{\frac{p^2\Lambda_n(z)^{23}}n \cdot \prn{\|\bQ\| \wedge \frac{\|\bQ\|_F}{\sqrt n}}}.
                \label{eq:trace-drift-bounds}
            \end{align}
        \end{enumerate}
    \end{subequations}
\end{lemma}
The proof is in Appendix~\ref{proof:trace-accumulated-drift}. We are now ready to closes the trace bounds in the next part, and we will obtain the additional derivative estimates by Cauchy's formula.

\paragraph{Part III: Terminal traces and derivatives.}
Taking $k=n$ in Eqs.~\eqref{eq:trace-path-H} and~\eqref{eq:trace-path-G}, applying the triangle inequality and invoking Lemma~\ref{lem:trace-accumulated-drift} implies
\begin{align*}
    \frac{\|H_n^{\bQ}-H_0^{\bQ}\|_{L^p}}n\vee\|G_n^{\bQ}(z,0)-G_1^{\bQ}\|_{L^p}
    =O\prn{\frac{p^2\Lambda_n(z)^{23}}n\min\brc{\|\bQ\|,\frac{\|\bQ\|_F}{\sqrt n}}}.
\end{align*}
For the last $G$ update, by Sherman-Morrison at $w_n=0$, it holds $G_{n+1}^{\bQ}(z,0)-G_n^{\bQ}(z,0)=-V_n^{G,\bQ}(0)$. We additionally invoke Eq.~\eqref{eq:trace-ratio-moments} from Part III of the proof of Lemma~\ref{lem:trace-centered-increments} at $(s,j)=(n-1,n)$ to control
\begin{align*}
    \|G_{n+1}^{\bQ}(z,0)-G_n^{\bQ}(z,0)\|_{L^p}=\|V_n^{G,\bQ}(0)\|_{L^p}=O\prn{\frac{\Lambda_n(z)^{7}}n\min\brc{\|\bQ\|,\frac{\|\bQ\|_F}{\sqrt n}}}.
\end{align*}
Adding this bound to the preceding estimate for $G_n^{\bQ}(z,0)-G_1^{\bQ}$ proves Theorem~\ref{thm:terminal-weighted-traces} for $r=0$.

Both terminal errors are analytic on $\C\setminus\R_{\geq0}$ by the sample resolvent formulas and Eq.~\eqref{eq:terminal-zero-sample-traces}. We can then transfer the analytic function value bounds to their derivatives when $r \geq 1$ by Cauchy's formula on the circle $|\zeta-z|=|z_-|/4$, as we have $3|z_-|/4\leq\dist(\zeta,\R_{\geq0})\leq5|z_-|/4$ and $|\zeta|\asymp|z|$. With the same minimizing cutoff $m$, Eq.~\eqref{eq:path-spectral-envelope} of Assumption~\ref{assmp:spectral-stability} guarantees $\Lambda_n(\zeta)=O(\Lambda_n(z))$ for $\Im\zeta\geq0$, and $3|z_-|/4\geq n^{-(\nu+1)}$ for sufficiently large $n$. Thus the $r=0$ bound applies uniformly on the circle, by conjugation with $\bQ^*$ when $\Im\zeta<0$. Applying Cauchy's formula and Minkowski's inequality for $r\geq1$ gives
\begin{align*}
    \frac1n\left\|\frac{\de^r}{\de z^r}(H_n^{\bQ}-H_0^{\bQ})\right\|_{L^p}
    &\leq\frac{r!}{2\pi n}\int_{|\zeta-z|=|z_-|/4}\frac{\|H_n^{\bQ}(\zeta,0)-H_0^{\bQ}(\zeta,w_0\lin(\zeta))\|_{L^p}}{|\zeta-z|^{r+1}}\,|\de \zeta|\\
    &\leq r!\left(\frac4{|z_-|}\right)^r\sup_{|\zeta-z|=|z_-|/4}\frac{\|H_n^{\bQ}(\zeta,0)-H_0^{\bQ}(\zeta,w_0\lin(\zeta))\|_{L^p}}{n}.
\end{align*}
The same argument applies to $G_{n+1}^{\bQ}(\zeta,0)-G_1^{\bQ}(\zeta,w_0\lin(\zeta))$ with the same multiplicative factor $r!(4/|z_-|)^r$. Substituting the $r=0$ bounds on the circle completes the proof.

\begin{remark}[Derivative estimates]
    It is easy to control the derivatives by Cauchy's formula, as we allow the spectral parameter to take complex values from the $r=0$ bound in Theorem~\ref{thm:terminal-weighted-traces}. For comparison, to estimate a first derivative from values on the real line, \citet{cheng2024dimensionfree} expand the function at equally spaced points and invert the resulting Vandermonde matrix. To bound the higher derivatives, they apply Fa\`a di Bruno's formula to the composition of the deterministic trace with the implicit spectral parameter. Here Cauchy's formula gives all derivative orders directly from the $r=0$ bound.
\end{remark}

\subsection{Proof of Theorem~\ref{thm:terminal-mixed-resolvents}}\label{proof:complex-mixed}\label{proof:terminal-mixed-resolvents}
Let $\chi^{\bQ}(z,z')$ denote the mixed-resolvent trace difference on the left-hand side of Eq.~\eqref{eq:terminal-mixed-resolvent-bound} at $r=r'=0$. Throughout the proof, we work by taking $\|\bQ\|_F=\sqrt{n}$ without loss of generality, by linearity in $\bQ$. We begin by establishing concentration of $\Tr(\bQ \wt{\bR}_n \wt{\bR}_n')$ around its deterministic equivalent at zeroth order, i.e., for $r=r'=0$, and subsequently extend the result to arbitrary derivative orders using Cauchy’s integral formula. We extend our shorthand notations for $\wt{\bR}_{n-1} = \wt{\bR}_{n-1}(z,0), \wt{\bR}_{n-1}' = \wt{\bR}_{n-1}(z',0)$ and $\delta_k = \delta_k(z), \delta_k' = \delta_k(z')$ (cf.~Eq.~\eqref{eq:def-delta-path-parameter}), and additionally define $\delta\lin(z)=-w_0\lin(z)/n=z/(n\omega(z/n))$ and similarly for $\delta\lin(z')$. We are now ready to first explain the heuristic and outline of the proof below.
\begin{align*}
    \Tr(\bQ\wt{\bR}_n\wt{\bR}_n') & \approx \E\Tr(\bQ\wt{\bR}_{n-1}\wt{\bR}_{n-1}') \approx \frac{z\E G_{n+1}^{\bQ}(z,0)-z'\E G_{n+1}^{\bQ}(z',0)}{n\delta_n \delta_n' (z-z')} &  \textrm{(Part I: LOO concentration)}\\
    & \approx \frac{1}{\delta\lin(z) \delta\lin(z')}\frac{zG_1^{\bQ}(z,w_0\lin(z))-z'G_1^{\bQ}(z',w_0\lin(z'))}{n(z-z')} &  \textrm{(Part II: Stability by path)} \\
    & = \frac{\omega(z/n)-\omega(z'/n)}{(z-z')/n}\Tr(\bQ\wt{\bR}_0\wt{\bR}_0'). & \textrm{(Exact resolvent identity)}
\end{align*}
In the first line of approximation, we will use leave-one-out on the terminal sample covariance matrix $\bT_n$ and concentration by Hanson-Wright, whose details will be in Part I of the proof; in Part II, we establish the second line of approximation, where path and single resolvent stability allow us to rewrite the previous step in its deterministic equivalent form at the zero sample point; finally, in Part III (as well as the reduction within the steps of Part I), we make use of the exact mixed-resolvent identity expressed by single resolvent divided differences to close the proof. Concretely, we will prove by adding the errors in the first two lines of the roadmap that
\begin{align*}
    \|\chi^{\bQ}(z,z')\|_{L^p}\leq
    \underbrace{O\prn{\frac{p^{3/2}(\Lambda_n\vee\Lambda_n')^{13}+(\Lambda_n\vee\Lambda_n')^{16}}{n^2}}}_{\text{Part I}}
    +\underbrace{O\prn{\frac{(\Lambda_n\vee\Lambda_n')^{26}}{n^2}}}_{\text{Part II}}.
\end{align*}
For the third line, we have the exact identity formulas below, reducing our remaining work to establishing the above inequality.
\begin{lemma}\label{lem:mixed-endpoint-identities}
    Under the setup of Theorem~\ref{thm:terminal-mixed-resolvents}, the sample and zero-sample endpoints satisfy
    \begin{subequations}
        \begin{align}
            \frac{zG_{n+1}^{\bQ}(z,0)-z'G_{n+1}^{\bQ}(z',0)}{z-z'}&=\Tr\prn{\bQ\bSigma^{1/2}\bR_n(z,0)\bT_n\bR_n(z',0)\bSigma^{1/2}},\label{eq:mixed-terminal-T-insertion}\\
            \frac{zG_1^{\bQ}(z,w_0\lin(z))-z'G_1^{\bQ}(z',w_0\lin(z'))}{z-z'}
            &=\frac{zz'\big(\omega(z/n)-\omega(z'/n)\big)}{(z-z')\omega(z/n)\omega(z'/n)}\Tr(\bQ\wt{\bR}_0\wt{\bR}_0').\label{eq:mixed-initial-identity}
        \end{align}
    \end{subequations}
    At $z=z'$ the divided differences are evaluated by differentiation.
\end{lemma}
The proof is in Appendix~\ref{proof:mixed-endpoint-identities}. Note that $\delta\lin(z)=z/(n\omega(z/n))$, and multiplying Eq.~\eqref{eq:mixed-initial-identity} of Lemma~\ref{lem:mixed-endpoint-identities} by $1/(n\delta\lin(z)\delta\lin(z'))=n\omega(z/n)\omega(z'/n)/(zz')$ leaves the prefactor $(\omega(z/n)-\omega(z'/n))/((z-z')/n)$. This yields the final equality in the roadmap. We now focus on its first line of approximation.

\paragraph{Part I: Leave-one-out concentration at the terminal sample.} As promised, we study the errors in the two leave-one-out approximations outlined at the opening of the proof through
\begin{align*}
    \chi_1&:=\Tr(\bQ\wt{\bR}_n\wt{\bR}_n')-\E\Tr(\bQ\wt{\bR}_{n-1}\wt{\bR}_{n-1}'),\\*
    \chi_2&:=\E\Tr(\bQ\wt{\bR}_{n-1}\wt{\bR}_{n-1}')-\dfrac{z\E G_{n+1}^{\bQ}(z,0)-z'\E G_{n+1}^{\bQ}(z',0)}{n\delta_n\delta_n'(z-z')}.
\end{align*}
The first residual error $\chi_1$ compares the full mixed trace with the expected trace after removing one observation. For the second residual $\chi_2$, taking expectation in Eq.~\eqref{eq:mixed-terminal-T-insertion} of Lemma~\ref{lem:mixed-endpoint-identities}, and expanding $\bT_n=\sum_{i=1}^n\bx_i\bx_i^*$ gives $n$ equal expected contributions. After removing $\bx_n$, the row numerator has conditional mean $\Tr(\bQ\wt{\bR}_{n-1}\wt{\bR}_{n-1}')$ given $\mc F_{n-1}$. Thus $\chi_2$ measures the error of replacing the two row reciprocals by $\delta_n$ and $\delta_n'$ in expectation. We explicitly bound the terms below, and defer these detailed proofs to Appendix~\ref{proof:terminal-mixed-replacement}.
\begin{lemma}\label{lem:terminal-mixed-replacement}
    Under the setup within the proof, we have
    \begin{align*}
        \|\chi_1\|_{L^p}=O\prn{\frac{p^{3/2}(\Lambda_n\vee\Lambda_n')^{13}}{n^2}}\qquad\text{and}\qquad
        |\chi_2|=O\prn{\frac{(\Lambda_n\vee\Lambda_n')^{16}}{n^2}}.
    \end{align*}
\end{lemma}

\paragraph{Part II: Stability of the increments and single traces.} We turn to the second line of the roadmap by first replacing the terminal increments by their linear values while keeping the terminal traces fixed. We then replace the expected terminal traces by their zero-sample values with the linear increments fixed. The corresponding residuals are
\begin{align*}
    \chi_3&:=\left(\frac1{\delta_n\delta_n'}-\frac1{\delta\lin(z)\delta\lin(z')}\right) \cdot \frac{z\E G_{n+1}^{\bQ}(z,0)-z'\E G_{n+1}^{\bQ}(z',0)}{n(z-z')},\\*
    \chi_4&:=\frac{z\prn{\E G_{n+1}^{\bQ}(z,0)-G_1^{\bQ}(z,w_0\lin(z))}-z'\prn{\E G_{n+1}^{\bQ}(z',0)-G_1^{\bQ}(z',w_0\lin(z'))}}{n\delta\lin(z)\delta\lin(z')(z-z')}.
\end{align*}
The next lemma controls the terms above, whose proof details are in Appendix~\ref{proof:mixed-single-trace-comparison}.
\begin{lemma}\label{lem:mixed-single-trace-comparison}
    Under the setup within the proof, we have
    \begin{align*}
        |\chi_3|=O\prn{\frac{(\Lambda_n\vee\Lambda_n')^{24}}{n^2}}\qquad\text{and}\qquad |\chi_4|=O\prn{\frac{(\Lambda_n\vee\Lambda_n')^{26}}{n^2}}.
    \end{align*}
\end{lemma}
Identifying $\chi^{\bQ}(z,z')=\chi_1+\chi_2+\chi_3+\chi_4$, we combine Lemmas~\ref{lem:terminal-mixed-replacement} and~\ref{lem:mixed-single-trace-comparison} with Minkowski's inequality and restore the factor $\|\bQ\|_F/\sqrt n$ by linearity to obtain Eq.~\eqref{eq:terminal-mixed-resolvent-bound} for $r=r'=0$,
\begin{equation}
    \|\chi^{\bQ}(z,z')\|_{L^p}\leq\|\chi_1\|_{L^p}+|\chi_2|+|\chi_3|+|\chi_4|=O\prn{\frac{p^{3/2}\|\bQ\|_F(\Lambda_n\vee\Lambda_n')^{26}}{n^{5/2}}}.\label{eq:mixed-trace-order-zero}
\end{equation}

\paragraph{Part III: Derivatives in the two spectral parameters.}
Let $\Gamma_z$ and $\Gamma_{z'}$ be the positively oriented circles centered at $z$ and $z'$ with radii $|z_-|/4$ and $|z'_-|/4$, respectively. We can then apply the same circle arguments from Part III of the proof of Theorem~\ref{thm:terminal-weighted-traces} in each variable, which enables us to use Eq.~\eqref{eq:mixed-trace-order-zero} on $\Gamma_z\times\Gamma_{z'}$ (taking adjoints for reflected resolvents when needed). The trace difference is analytic in each variable, with the divided difference extended by continuity on the diagonal. Applying Cauchy's formula on the two circles and then Minkowski's inequality, we have
\begin{align*}
    \left\|\frac{\de^{r+r'}}{\de z^r\,\de(z')^{r'}}\chi^{\bQ}(z,z')\right\|_{L^p}
    &=\frac{r!\cdot r'!}{(2\pi)^2}\left\|\oint_{\Gamma_z}\oint_{\Gamma_{z'}}
    \frac{\chi^{\bQ}(\zeta,\zeta')}{(\zeta-z)^{r+1}(\zeta'-z')^{r'+1}}\,\de\zeta'\,\de\zeta\right\|_{L^p}\\
    &\leq r!\cdot r'!\left(\frac{4}{|z_-|}\right)^r\left(\frac{4}{|z'_-|}\right)^{r'}\sup_{\Gamma_z\times\Gamma_{z'}}\|\chi^{\bQ}(\zeta,\zeta')\|_{L^p}.
\end{align*}
Substituting the result for $r=r'=0$ into the preceding bound completes the proof for any fixed order of mixed derivatives.

\subsection{Auxiliary proofs for Theorem~\ref{thm:self-consistent-path-stability}} \label{proof:auxiliary-theorem1}

\subsubsection{Proof of Lemma~\ref{lem:resolvent-bounds-strict-positivity}} \label{proof:resolvent-bounds-strict-positivity}
We will first prove the path-dependent bounds in (i) and (ii) at fixed parameters satisfying $\Im(w\wb{z})\geq0$ (or $w\leq0$ when $z<0$) in the first two parts, and then verify these conditions along the backward path (as well as its well-posedness) in the proof of (iii). Recall that $z= t+ \iu \eta$, $w = u+ \iu v$ and therefore $\Im(w \wb{z}) = tv -\eta u$.

\paragraph{Part I: Resolvent bounds in (i).}
We first prove the existence and the bound $\|\bR_{j-1}(z,w)\|\leq\eta^{-1}$ at a general fixed $(z,w)\in\mc D$, where $\eta>0$. Fixing $j\in[n]$ and using the shorthand $\bM_{j-1}=\bM_{j-1}(z,w):=\bT_{j-1}-w\bSigma-z\bI$ and $\bR_{j-1}=\bR_{j-1}(z,w)$ whenever the inverse exists, we have for $\bh\in\mc H$,
\begin{align}
    -\Im\langle\bh,\bM_{j-1}\bh\rangle_{\mc H}
    =\eta\|\bh\|_{\mc H}^2+v\|\bh\|_{\bSigma}^2
    \geq\eta\|\bh\|_{\mc H}^2.
    \label{eq:upper-half-plane-coercivity}
\end{align}
It follows that $\|\bM_{j-1}\bh\|_{\mc H}\geq\eta\|\bh\|_{\mc H}$. Hence, $\bM_{j-1}$ is bounded below by $\eta$ and, in particular, is injective and has closed range. Applying the same argument to its adjoint implies $\bM_{j-1}^*$ has a trivial kernel, and by Hilbert space duality, $\bM_{j-1}$ has a dense range that equals $\mc H$. Thus the inverse exists and $\|\bR_{j-1}\|\leq\eta^{-1}$. For the path-dependent results, we shall assume $\Im(w\wb{z})\geq0$ and prove the two bounds in (i) at this fixed parameter $w$, as is mentioned in the proof outline.

We begin with proving the resolvent norm bound $\norm{\bR_{j-1}} \leq |z_-|^{-1}$. When $t \geq 0$, this is implied by the previous general bound $\|\bR_{j-1}\|\leq\eta^{-1}$ as $z_-= \iu \eta$. When $t<0$, this condition implies $u\leq tv/\eta \leq0$. Consequently, applying Cauchy-Schwarz for every $\bh\in\mc H$,
\begin{align*}
    \|\bM_{j-1}\bh\|_{\mc H}^2\|\bh\|_{\mc H}^2
    &\geq|\langle\bh,\bM_{j-1}\bh\rangle_{\mc H}|^2 = \prn{\Re\langle\bh,\bM_{j-1}\bh\rangle_{\mc H}}^2 + \prn{\Im\langle\bh,\bM_{j-1}\bh\rangle_{\mc H}}^2 \geq(t^2+\eta^2)\|\bh\|_{\mc H}^4.
\end{align*}
In the last step, we additionally make use of the parameter condition $u \leq 0$ for the real part, i.e., $\Re\langle \bh, \bM_{j-1} \bh\rangle_{\mc{H}}=\langle \bh, \bT_{j-1} \bh\rangle_{\mc{H}}-u\|\bSigma^{1 / 2} \bh\|_{\mc{H}}^2-t\|\bh\|_{\mc{H}}^2 \geq(-t)\|\bh\|_{\mc{H}}^2$. Hence, $\|\bM_{j-1}\bh\|_{\mc H}\geq |z_-|\|\bh\|_{\mc H}$, and using the inverse already established concludes with $\|\bR_{j-1}\|\leq|z|^{-1}=|z_-|^{-1}$.

Next, we establish the quadratic-form bound by exploiting the positive semidefinite property of $\Im(z\bR_{j-1})$. Since $\bM_{j-1}\bR_{j-1}=\bI$ and $\bR_{j-1}^*\bM_{j-1}^*=\bI$, we can write
\begin{align}
    \Im(z\bR_{j-1})
    & = \bR_{j-1}^* \cdot \frac{z\bM_{j-1}^*-\wb{z}\bM_{j-1}}{2\iu} \cdot \bR_{j-1} =\bR_{j-1}^*\prn{\eta\bT_{j-1}+\Im(w\wb{z})\bSigma}\bR_{j-1}.
    \label{eq:helper-rotated-resolvent-identity}
\end{align}
Since $\bT_{j-1}\succeq0$ and $\Im(w\wb{z})\geq0$, the right-hand side is positive semidefinite. Hence, for every $\bx\in\mc H$,
\begin{align*}
    \Im\prn{z(1+\bx^*\bR_{j-1}\bx)}
    =\eta+\bx^*\Im(z\bR_{j-1})\bx\geq\eta,
\end{align*}
which implies $|1+\bx^*\bR_{j-1}\bx|\geq\eta/|z|$. Similarly, we only need to prove the case when $t<0$ as $\eta = |z_-|$ when $t \geq 0$. We again use $u\leq tv/\eta \leq 0$ and additionally the inverse real-part identity:
\begin{align}
    \Re\bR_{j-1}
    =\bR_{j-1}^*\prn{\bT_{j-1}-u\bSigma-t\bI}\bR_{j-1}\succeq0.
    \label{eq:helper-path-resolvent-real-part}
\end{align}
It then follows that $|1+\bx^*\bR_{j-1}\bx|\geq\Re(1+\bx^*\bR_{j-1}\bx)\geq1 = |z_-|/|z|$. Thus $|1+\bx^*\bR_{j-1}\bx|^{-1}\leq |z|/|z_-|$ in both cases.

We finally cover the case of $z<0$ and real $w\leq0$. Now the operator is self-adjoint and
\begin{align}
    \bM_{j-1}(z,w)=\bT_{j-1}-w\bSigma+|z|\bI\succeq|z|\bI,
    \qquad 0\preceq\bR_{j-1}(z,w)\preceq|z|^{-1}\bI.
    \label{eq:helper-negative-real-resolvent}
\end{align}
Thus the inverse exists and $1+\bx^*\bR_{j-1}(z,w)\bx\geq1$. Since $|z_-|=|z|$, this proves both bounds in (i) at these fixed real parameters as well. It remains to verify the required conditions on $w$ along the path in Part III, i.e.\ $\Im(w_k\wb{z}) \geq 0$ when $\eta > 0$ and $w_k \leq 0$ when $t < 0, \eta = 0$.
\paragraph{Part II: Expected-trace bounds in (ii).}
We first prove the strict positivity results at a general $(z,w)\in\mc D$. By the inverse imaginary-part identity, we have
\begin{align}
    \Im\bR_{j-1}
    =\bR_{j-1}^* \cdot \frac{\bM_{j-1}^*-\bM_{j-1}}{2\iu} \cdot \bR_{j-1}
    =\bR_{j-1}^*(\eta\bI+v\bSigma)\bR_{j-1}\succeq0.
    \label{eq:upper-half-plane-resolvent-imaginary-part}
\end{align}
Recalling the definition of $G_j$ in Eq.~\eqref{eq:def-traces-interpolation}, we obtain from $v \geq 0$ that
\begin{align*}
    \Im G_j(z,w)
    & = \Tr\prn{\bSigma^{\half}\bR_{j-1}^*(\eta\bI+v\bSigma)\bR_{j-1}\bSigma^{\half}}
    \geq \eta\|\bR_{j-1}\bSigma^{\half}\|_F^2
    \geq \frac{\eta\Tr(\bSigma)}{\|\bM_{j-1}\|^2}>0,
\end{align*}
where the last lower bound follows from $\bM_{j-1}\bR_{j-1} = \bI $ and $\Tr(\bSigma)^\half = \|\bSigma^{\half}\|_F\leq\|\bM_{j-1}\|\,\|\bR_{j-1}\bSigma^{\half}\|_F$. Applying the norm bound from Part I gives $|G_j(z,w)|=\left|\Tr\prn{\bSigma\bR_{j-1}(z,w)}\right|\leq\Tr(\bSigma)\norm{\bR_{j-1}(z,w)}\leq\Tr(\bSigma)/\eta$, so we can take expectation by boundedness, and $\Im\wb G_j(z,w)=\E[\Im G_j(z,w)]>0$.

We then use the operator identities from Part I to obtain the denominator and real-part bounds under $\Im(w\wb{z})\geq0$. Taking the weighted trace and expectation in Eq.~\eqref{eq:helper-rotated-resolvent-identity} implies
\begin{align}
    \Im\prn{z(1+\wb G_j(z,w))}
    =\eta+\E\Tr\prn{\bSigma^{\half}\Im(z\bR_{j-1})\bSigma^{\half}}\geq\eta.
    \label{eq:helper-pointwise-trace-signs}
\end{align}
We obtain the denominator bound when $t \geq 0$ since $|z| \cdot |1+\wb G_j(z,w)| \geq\eta=|z_-|$. For $t<0$, taking the weighted trace and expectation in Eq.~\eqref{eq:helper-path-resolvent-real-part} instead yields
\begin{align*}
    \Re\wb G_j(z,w)=\E\Tr\prn{\bSigma^{\half}\Re\bR_{j-1}\bSigma^{\half}}\geq0
    \qquad \text{and} \qquad |1+\wb G_j(z,w)|\geq1+\Re\wb G_j(z,w)\geq1=\frac{|z_-|}{|z|}.
\end{align*}
This also proves the real-part positivity in (ii).

Lastly, for the case of $z<0$ and real $w\leq0$, we obtain from taking the weighted trace and expectation in Eq.~\eqref{eq:helper-negative-real-resolvent} that
\begin{align}
    0\leq\wb G_j(z,w)\leq\frac{\Tr(\bSigma)}{|z|},\qquad
    |1+\wb G_j(z,w)|=1+\wb G_j(z,w)\geq1=\frac{|z_-|}{|z|}.
    \label{eq:helper-negative-real-trace}
\end{align}
Thus (ii) holds at each fixed parameter satisfying the same path parameter conditions of Part I. We now verify these conditions along the backward path.

\paragraph{Part III: Path and endpoint bounds in (iii).}
For $z\in\C_+$, we prove by backward induction that $(z,w_k)\in\mc D$ and $\Im(w_k\wb{z})\geq0$, which also justifies applying Parts I and II at every path parameter. The exact analysis also yields a quantitative lower bound for the latter increment, whose sum will give the second endpoint estimate conclusion in (iii).

We start the induction backward at $w_n=0$. Suppose the two conditions hold at some $1\leq k\leq n$. From the strict positivity proved in Part II, namely $\Im\wb G_k(z,w_k)>0$, we have that the denominator in the backward recursion is nonzero and the imaginary part increases strictly
\begin{align*}
    \Im w_{k-1}
    =\Im w_k+\frac{\Im\wb G_k(z,w_k)}{|1+\wb G_k(z,w_k)|^2}>0.
\end{align*}
Thus $(z,w_{k-1})\in\mc D$. For the second condition, it holds from the same recursion that
\begin{align}
    \Im(w_{k-1}\wb{z}) - \Im(w_k\wb{z})
    &= \frac{\Im\prn{z(1+\wb G_k(z,w_k))}}{|1+\wb G_k(z,w_k)|^2} \geq \frac{\eta}{|1+\wb G_k(z,w_k)|^2} > 0
    \label{eq:path-cone-update}
\end{align}
The last inequality applies Eq.~\eqref{eq:helper-pointwise-trace-signs} with $j=k$ and $w=w_k$, and we close the positivity results in Lemma~\ref{lem:resolvent-bounds-strict-positivity}.(iii) by induction. Substituting $w=w_k$ into Parts I and II we also close the loop and prove Lemma~\ref{lem:resolvent-bounds-strict-positivity}.(i) and (ii).

It only remains to establish the endpoint bound. To this end, we first sum Eq.~\eqref{eq:path-cone-update} from $w_n=0$,
\begin{align*}
    \Im(w_0\wb{z})\geq\eta\sum_{k=1}^n\frac1{|1+\wb G_k(z,w_k)|^2}.
\end{align*}
Telescoping the backward recursion and then applying Cauchy-Schwarz implies
\begin{align*}
    |w_0|^2
    & = \left|\sum_{k=1}^n\frac1{1+\wb G_k(z,w_k)}\right|^2
    \leq n\sum_{k=1}^n\frac1{|1+\wb G_k(z,w_k)|^2}
    \leq \frac n\eta\Im(w_0\wb{z}),
\end{align*}
where the last inequality uses the preceding display. Therefore, rearranging proves the second bound in Lemma~\ref{lem:resolvent-bounds-strict-positivity}.(iii), Eq.~\eqref{eq:helper-path-cone-bound}.

Finally, for $z<0$, it suffices to verify that the recursion stays on the real half-line $w\leq0$. Starting from $w_n=0$, suppose $w_k\leq0$ is real. Applying Eq.~\eqref{eq:helper-negative-real-trace} with $j=k$ gives $\wb G_k(z,w_k)\geq0$, so the next value is real and satisfies
\begin{align*}
    w_{k-1}=w_k-\frac1{1+\wb G_k(z,w_k)}<w_k\leq0.
\end{align*}
This closes the real-parameter induction. Finally, the backward recursion uniquely defines $\{w_k\}_{k=0}^n$, and all resolvents and expected traces appearing in it are well-defined, and the proof is completed.

\subsubsection{Proof of Lemma~\ref{lem:weighted-resolvent-complex-to-real}}
\label{proof:weighted-resolvent-complex-to-real}
We prove both types of bounds using the same rotation argument of the inverse resolvent. In Part I, we bound the norms and trace of the weighted resolvent in terms of the real part of the rotated inverse. In Part II, keeping the $\bT_k$ term explicit provides the comparison with $\bA_k(|z_-|)$. In Part III, using $\bSigma$ term implies the bounds in terms of $-\Re(e^{\iu\theta}w_j)$.

Recall that $z=t+\iu\eta$, $w=u+\iu v$, and $\wt{\bR}_k=\wt{\bR}_k(z,w)$. We work first at $(z,w)\in\mc D$ with $\Im(w\wb{z})=tv-\eta u\geq0$, or at real $z<0$ and $w\leq0$. The path properties established in Part I of Appendix~\ref{proof:self-consistent-path-stability} verify these conditions at every $w=w_j$.

\paragraph{Part I: The common rotation and comparison.}
For $t_+ = t \vee 0$, we set the rotation angle $\theta=\arctan(1+2t_+/\eta)$ when $\eta>0$, and $\theta=\pi/4$ when $z<0$ is real. For $\bU=e^{\iu\theta}(\bT_k-w\bSigma-z\bI)$, we will prove the two lower bounds
\begin{align}
    \Re\bU\succeq\cos\theta \cdot (\bT_k+|z_-|\bI)\qquad\text{and}\qquad\Re\bU\succeq\tfrac12|z_-|\bI-\Re(e^{\iu\theta}w)\bSigma.\label{eq:weighted-common-rotation}
\end{align}
First, for $t\geq0$, we have $\eta>0$ and $u\leq tv/\eta$, so
\begin{align*}
    \Re\bU
    &=\cos\theta \cdot \prn{\bT_k+((1+2t/\eta)v-u)\bSigma+(\eta+t)\bI} \succeq\cos\theta \cdot \bT_k+ \frac{\eta+t}{\sqrt{\eta^2+(\eta+2t)^2}}(v\bSigma+\eta\bI).
\end{align*}
The coefficient in the lower bound is at least $1/2$, because $\eta^2+(\eta+2t)^2\leq4(\eta+t)^2$. In the first equality above, the coefficient of $\bSigma$ is $\cos\theta \cdot ((1+2t/\eta)v-u)=-\Re(e^{\iu\theta}w)$, and the coefficient of $\bI$ is $(\eta+t)\cos\theta=\eta(\eta+t)/\sqrt{\eta^2+(\eta+2t)^2}\geq\eta/2$. This proves the second bound in Eq.~\eqref{eq:weighted-common-rotation}, since $|z_-|=\eta$. The first equality also proves the first bound in Eq.~\eqref{eq:weighted-common-rotation}, as its $\bSigma$ coefficient is nonnegative and $\eta+t\geq|z_-|$. Second, for $t<0$ and $\eta>0$, the same condition and $v\geq0$ imply $u\leq tv/\eta\leq0$. Third, for real $z<0$, we assumed $u=w\leq0$ and $v=0$. In both cases $\theta=\pi/4$, so
\begin{align*}
    \Re\bU=\frac{\bT_k+(v-u)\bSigma+(\eta-t)\bI}{\sqrt2} \succeq\frac{\bT_k+(v-u)\bSigma+|z_-|\bI}{\sqrt2},
\end{align*}
where we use $\eta-t\geq|z_-|$ and $v-u\geq0$. Here $(v-u)/\sqrt2=-\Re(e^{\iu\theta}w)$ and $|z_-|/\sqrt2\geq |z_-|/2$. Both lower bounds in Eq.~\eqref{eq:weighted-common-rotation} still follow.

Next we derive the common norm comparisons used in Parts II and III, by setting $\bV:=\Re\bU\succ0$. The inverse of $\bU$ exists by Part I of the proof of Lemma~\ref{lem:resolvent-bounds-strict-positivity}. We use $\bU+\bU^*=2\bV$, $\bV^* = \bV$ and then take the congruence by $\bU^{-1}\bSigma^{1/2}$ to obtain that
\begin{align}
    \bU^*\bV^{-1}\bU-\bV&=(\bU-\bV)^*\bV^{-1}(\bU-\bV)\succeq0 \quad \text{and} \quad \bSigma^{1/2}\bU^{-*}\bV\bU^{-1}\bSigma^{1/2} \preceq\bSigma^{1/2}\bV^{-1}\bSigma^{1/2}.
    \label{eq:weighted-comparison-order}
\end{align}
Identify that $\wt{\bR}_k=e^{\iu\theta}(\bSigma^{1/2}\bV^{-1/2})(\bV^{1/2}\bU^{-1}\bSigma^{1/2})$, and we apply sub-multiplicativity for the operator norm and trace Cauchy-Schwarz for the Frobenius norm and trace to obtain
\begin{subequations}
    \begin{align}
        \|\wt{\bR}_k\|&\leq\|\bSigma^{1/2}\bV^{-1/2}\|\,\|\bV^{1/2}\bU^{-1}\bSigma^{1/2}\|
        \leq\|\bSigma^{1/2}\bV^{-1}\bSigma^{1/2}\|,\label{eq:weighted-comparison-operator}\\
        \|\wt{\bR}_k\|_F^2
        &=\Tr\prn{(\bV^{-1/2}\bSigma\bV^{-1/2})(\bV^{1/2}\bU^{-1}\bSigma\bU^{-*}\bV^{1/2})}\nonumber\\
        &\leq\|\bSigma^{1/2}\bV^{-1}\bSigma^{1/2}\|_F\,
        \|\bSigma^{1/2}\bU^{-*}\bV\bU^{-1}\bSigma^{1/2}\|_F
        \leq\|\bSigma^{1/2}\bV^{-1}\bSigma^{1/2}\|_F^2,\label{eq:weighted-comparison-frobenius}\\
        |\Tr\wt{\bR}_k|&\leq\|\bSigma^{1/2}\bV^{-1/2}\|_F\,\|\bV^{1/2}\bU^{-1}\bSigma^{1/2}\|_F
        \leq\Tr(\bSigma^{1/2}\bV^{-1}\bSigma^{1/2}).\label{eq:weighted-comparison-trace}
    \end{align}
\end{subequations}
Each final inequality above follows from Eq.~\eqref{eq:weighted-comparison-order}. In particular, the Frobenius norm is monotone for positive operators since $0\preceq\bC\preceq\bD$ implies $\Tr(\bD^2)-\Tr(\bC^2)=\Tr((\bD-\bC)(\bD+\bC))\geq0$.

\paragraph{Part II: Comparison with the real weighted resolvent.}
To prove Eq.~\eqref{eq:path-real-comparisons}, it suffices to bound $\bSigma^{1/2}\bV^{-1}\bSigma^{1/2}$ by $3|z|\bA_k(|z_-|)/|z_-|$. We use the first lower bound in Eq.~\eqref{eq:weighted-common-rotation}. Taking inverses reverses the positive-operator order, and taking the congruence by $\bSigma^{1/2}$ yields
\begin{align*}
    \bSigma^{1/2}\bV^{-1}\bSigma^{1/2}
    \preceq\frac1{\cos\theta}\bSigma^{1/2}(\bT_k+|z_-|\bI)^{-1}\bSigma^{1/2}
    =\frac1{\cos\theta}\bA_k(|z_-|).
\end{align*}
By monotonicity in p.s.d.\ cone for Eqs.~\eqref{eq:weighted-comparison-operator} to \eqref{eq:weighted-comparison-trace}, we establish all three bounds with factor $(\cos\theta)^{-1}$. Therefore it only remains to bound this factor. If $t<0$, then $(\cos\theta)^{-1}=\sqrt2$ and $|z|=|z_-|$. If $t\geq0$, then $|z_-|=\eta$ and
\begin{align*}
    (\cos\theta)^{-1}
    =\frac{\sqrt{\eta^2+(\eta+2t)^2}}\eta
    \leq\frac{\sqrt{4\eta^2+6t^2}}\eta
    \leq\frac{\sqrt6\,|z|}\eta
    \leq\frac{3|z|}{|z_-|}.
\end{align*}
The first inequality uses $4\eta t\leq2\eta^2+2t^2$. Taking $w=w_j$ now proves Eq.~\eqref{eq:path-real-comparisons}.

\paragraph{Part III: Bounds using the path parameter.}
Finally, we prove Eq.~\eqref{eq:path-buffer-norms}. Write $w_j=u_j+\iu v_j$. For $j<n$, the path monotonicity established in Part I of Appendix~\ref{proof:self-consistent-path-stability} gives $v_j>0$ when $z\in\C_+$, and $u_j=w_j<0$, $v_j=0$ when $z<0$ is real. For $t\geq0$, using $u_j\leq tv_j/\eta$ in the first calculation in Part I yields $-\Re(e^{\iu\theta}w_j)\geq\cos\theta(1+t/\eta)v_j>0$. For $t<0$, the same coefficient is $(v_j-u_j)/\sqrt2>0$, since $u_j\leq0$ and either $v_j>0$ or $w_j=u_j<0$ is real. We may therefore first derive the bounds for any $w$ with $-\Re(e^{\iu\theta}w)>0$ and then substitute $w=w_j$.

The second lower bound in Eq.~\eqref{eq:weighted-common-rotation} compares $\bV$ with an operator diagonal in the population eigenbasis. Taking inverses and the congruence by $\bSigma^{1/2}$ gives
\begin{align*}
    \bSigma^{1/2}\bV^{-1}\bSigma^{1/2}
    \preceq2\bSigma^{1/2} \cdot \prn{|z_-|\bI-2\Re(e^{\iu\theta}w)\bSigma}^{-1}\bSigma^{1/2}.
\end{align*}
From Eqs.~\eqref{eq:weighted-comparison-operator} and~\eqref{eq:weighted-comparison-frobenius}, monotonicity in the p.s.d.\ cone, and splitting the eigenvalue sum at $m$, we have
\begin{align*}
    \|\wt{\bR}_k\|&\leq2\sup_i\frac{\sigma_i}{|z_-|-2\Re(e^{\iu\theta}w)\sigma_i}
    \leq \frac1{-\Re(e^{\iu\theta}w)},\\
    \|\wt{\bR}_k\|_F^2&\leq4\sum_i\frac{\sigma_i^2}{(|z_-|-2\Re(e^{\iu\theta}w)\sigma_i)^2}\leq \sum_{i\leq m}\frac1{\Re(e^{\iu\theta}w)^2}
    +4\sum_{i>m}\frac{\sigma_i^2}{|z_-|^2}= \frac{m}{\Re(e^{\iu\theta}w)^2}+\frac{4\tau_2(m)}{|z_-|^2}.
\end{align*}
Taking $w=w_j=u_j+\iu v_j$ proves Eq.~\eqref{eq:path-buffer-norms} and the proof is complete.

\subsubsection{Proof of Lemma~\ref{lem:weighted-resolvent-norm-control}}
\label{proof:weighted-resolvent-norm-control}

We need bounds for the operator norm, squared Frobenius norm, and trace of $\bA_k(\zeta)$, uniformly over $1\leq p\leq\kappa\log k$. In Part I, we reduce these quantities to the norm of the lower population block on an event where the top Gram matrix is well-conditioned. In Part II, we integrate a tail bound for the lower block. In Part III, we combine these estimates and bound the exceptional-event contribution by the deterministic inequality $0\preceq\bA_k(\zeta)\preceq\bSigma/\zeta$.

\paragraph{Part I: Reduction to the lower block.} The following deterministic inequality separates the top and lower population eigenspaces without requiring independence between them.

\begin{lemma}[{\citep[Lemma~C.1]{cheng2024dimensionfree}}]
    \label{lem:cheng-montanari-truncation}
    Let $\bU,\bW:\mc H\to\R^k$ be bounded operators. For every $\zeta>0$,
    \begin{align*}
        \zeta\bI+(\bU+\bW)^*(\bU+\bW)
        \succeq\frac{\zeta}{2}\bI+\prn{1+\frac{2\norm{\bW^*\bW}}{\zeta}}^{-1}\bU^*\bU.
    \end{align*}
\end{lemma}

Let $\bP_m$ be the orthogonal projection onto the first $m$ population eigenvectors, let $\bQ_m=\bI-\bP_m$, and write $\bSigma_{\leq m}=\bP_m\bSigma\bP_m$ and $\bSigma_{>m}=\bQ_m\bSigma\bQ_m$. Recall that $\bX_k:\mc H\to\R^k$ and $\bT_k=\bX_k^*\bX_k$. Decompose into lower and upper blocks,
\begin{align*}
    \bX_k=\bU_{k,m}+\bW_{k,m},\qquad \bU_{k,m}:=\bX_k\bP_m,\qquad \bW_{k,m}:=\bX_k\bQ_m.
\end{align*}
Letting $\bZ_{k,m}$ be the $k\times m$ matrix whose $i$-th row consists of the first $m$ coordinates of $\bz_i$ (we identify $\R^m$ with $\bP_m\mc H$ through its ordered population eigenbasis), then
\begin{align*}
    \bU_{k,m}^*\bU_{k,m}=\bSigma_{\leq m}^{\half}\bZ_{k,m}^*\bZ_{k,m}\bSigma_{\leq m}^{\half}.
\end{align*}
Applying Lemma~\ref{lem:cheng-montanari-truncation} with $\bU=\bU_{k,m}$ and $\bW=\bW_{k,m}$ gives the deterministic reduction
\begin{align}
    \bA_k(\zeta)
    \preceq\bSigma^{\half}\prn{\frac{\zeta}{2}\bI+\prn{1+\frac{2\norm{\bW_{k,m}^*\bW_{k,m}}}{\zeta}}^{-1}\bSigma_{\leq m}^{\half}\bZ_{k,m}^*\bZ_{k,m}\bSigma_{\leq m}^{\half}}^{-1}\bSigma^{\half}.
    \label{eq:norm-bound-truncation-reduction}
\end{align}
We therefore need a lower bound for $\bZ_{k,m}^*\bZ_{k,m}$ and an upper bound for $\|\bW_{k,m}^*\bW_{k,m}\|$, as in \citet[Proof of Lemma~7.2]{cheng2024dimensionfree}. We first establish the former. Under Assumption~\ref{assmp:features}, we can repeat the fourth-moment verification in~\citet[Appendix~C.2, Part~II]{cheng2024dimensionfree} with sample size $k$ and cutoff $m$, which yields isotropic rows of $\bZ_{k,m}$ whose unit linear forms have fourth moments bounded by $\const^2$, for $\const=O_{C_{\bx}}(1)$. Then we invoke~\citet[Theorem~2.2]{yaskov2014lower} with dimension $m$, sample size $k$, and lower-tail parameter $\sqrt{k}/(4\const)$. Taking $\const=\const(C_{\bx})\ge1$ uniformly over the upper block coordinates and setting $\rho(C_{\bx})=(16\const)^{-2}$ gives $4\const\sqrt{\rho(C_{\bx})}=1/4$. With this choice, we obtain the following bounds for $\mc E(k,m):=\{\bZ_{k,m}^*\bZ_{k,m}\succeq(k/2)\bI_m\}$ when $m\leq\rho k$:
\begin{align}
    \P\prn{\mc E(k,m)^c}
    &\leq\exp\{-\Omega_{C_{\bx}}(k)\},&
    \bZ_{k,m}^*\bZ_{k,m}
    &\succeq\frac{k}{2}\bI_m,&
    \norm{(\bZ_{k,m}^*\bZ_{k,m})^{-1}}
    &\leq\frac{2}{k}.
    \label{eq:norm-bound-head-event}
\end{align}
Here the last two inequalities hold on $\mc E(k,m)$.
We next make explicit how this top Gram matrix bound controls the two resolvent norms. On $\mc E(k,m)$, substituting $\bZ_{k,m}^*\bZ_{k,m}\succeq(k/2)\bI_m$ into Eq.~\eqref{eq:norm-bound-truncation-reduction} gives
\begin{align*}
    \bA_k(\zeta)
    &\preceq\bSigma^{1/2}\prn{\frac\zeta2\bI+\frac{k}{2}\prn{1+\frac{2\|\bW_{k,m}^*\bW_{k,m}\|}\zeta}^{-1}\bSigma_{\leq m}}^{-1}\bSigma^{1/2} \preceq\frac2k\prn{1+\frac{2\|\bW_{k,m}^*\bW_{k,m}\|}\zeta}\bP_m+\frac2\zeta\bSigma_{>m}.
\end{align*}
For the second inequality, diagonalizing in the population eigenbasis gives, on the first $m$ coordinates,
\begin{align*}
    \frac{\sigma_i}{\frac{\zeta}{2}+\frac{k}{2}(1+2\|\bW_{k,m}^*\bW_{k,m}\|/\zeta)^{-1}\sigma_i}
    \leq\frac{2}{k}\left(1+\frac{2\|\bW_{k,m}^*\bW_{k,m}\|}{\zeta}\right).
\end{align*}
On the complement, the eigenvalues of the upper bound are $2\sigma_i/\zeta$. The two terms in the upper bound of $\bA_k(\zeta)$ above act on orthogonal subspaces, and the first has rank at most $m$. Taking spectral and squared Frobenius norms respectively therefore gives on the event $\mc E(k,m)$
\begin{align}
    \|\bA_k(\zeta)\|&\leq\frac2k\prn{1+\frac{2\|\bW_{k,m}^*\bW_{k,m}\|}\zeta}+\frac{2\sigma_{m+1}}\zeta,\notag\\
    \|\bA_k(\zeta)\|_F^2&\leq\frac{4m}{k^2}\prn{1+\frac{2\|\bW_{k,m}^*\bW_{k,m}\|}\zeta}^2+\frac{4\tau_2(m)}{\zeta^2}.
    \label{eq:norm-bound-headwise-estimates}
\end{align}
These are the estimates whose $L^p$ norms we need to control. The random tail norm enters the spectral norm bound linearly and the squared Frobenius norm bound quadratically, so an $L^{2p}$ bound on $\|\bW_{k,m}^*\bW_{k,m}\|$ will control both contributions on $\mc E(k,m)$. Part II derives a uniform bound on this moment for all $1 \leq p \leq \kappa \log k$. Part III subsequently applies this estimate to obtain the stated upper bounds on the corresponding norms and the associated trace bound, and then controls the remaining contribution using the high probability event $\mc{E}(k,m)$ in Eq.~\eqref{eq:norm-bound-head-event}.

\paragraph{Part II: Moments of the lower block.}
By Eq.~\eqref{eq:norm-bound-headwise-estimates}, it suffices to prove $\|\|\bW_{k,m}^*\bW_{k,m}\|\|_{L^p}=O_{C_{\bx},\kappa}(\const_{\bA}(k,m))$ for $1\leq p\leq2\kappa\log k$. Replacing $p$ by $2p$ will then supply the moment required for the squared Frobenius norm. We first derive a tail bound with an explicit tail parameter $a\geq1$, and then integrate it. The lower block is $\bW_{k,m}^*\bW_{k,m}=\sum_{i=1}^k(\bQ_m\bx_i)(\bQ_m\bx_i)^*$, with covariance parameters
\begin{align*}
    \E[(\bQ_m\bx_i)(\bQ_m\bx_i)^*]=\bSigma_{>m},\qquad
    \|\bSigma_{>m}\|=\sigma_{m+1},\qquad
    \Tr\bSigma_{>m}=\tau_1(m),\qquad
    \|\bSigma_{>m}\|_F^2=\tau_2(m).
\end{align*}
If $\sigma_{m+1}=0$, the lower block vanishes almost surely and its moment bound is immediate. Otherwise, its covariance effective rank is $r_1(m)=\tau_1(m)/\sigma_{m+1}\geq1$. Recall the deterministic threshold
\begin{align*}
    \const_{\bA}(k,m):=k\sigma_{m+1}+\sqrt{k\sigma_{m+1}\tau_1(m)\log (k)\log(kr_1(m))}+\tau_1(m)\log (k)\log(kr_1(m)).
\end{align*}
We use the truncation argument of \citet[Lemma~C.2]{cheng2024dimensionfree}, keeping the truncation level dependent on $a$ so that the resulting probability bound can be integrated at the required moment orders. Applying Proposition~\ref{prop:hanson-wright}, Eq.~\eqref{eq:hanson-wright}, with $\bA=\bQ_m$, we have for $u>0$,
\begin{align*}
    \P\left(\|\bQ_m\bx_i\|^2>\tau_1(m)+u\right)
    \leq2\exp\left\{-\Omega_{C_{\bx}}\left(\frac{u^2}{\tau_2(m)}\wedge\frac{u}{\sigma_{m+1}}\right)\right\}.
\end{align*}
Since $\sqrt{\tau_2(m)}\leq\tau_1(m)$ and $\sigma_{m+1}\leq\tau_1(m)$, we can choose a deterministic truncation level satisfying
\begin{align*}
    L_a=\const\tau_1(m)(\log k+a),\qquad \P(\|\bQ_m\bx_i\|^2>L_a)\leq2e^{-a}/k,
\end{align*}
where $\const$ is sufficiently large and depends only on $C_{\bx}$. Define the independent truncated summands, for $i\in[k]$, by
\begin{align*}
    \bY_{i,a}:=(\bQ_m\bx_i)(\bQ_m\bx_i)^*\ind\{\|\bQ_m\bx_i\|^2\leq L_a\}.
\end{align*}
Their positivity and rank-one form imply
\begin{align*}
    0\preceq\bY_{i,a}\preceq L_a\bI,\qquad
    \E\bY_{i,a}\preceq\bSigma_{>m},\qquad
    (\bY_{i,a})^2=\|\bQ_m\bx_i\|^2\bY_{i,a}\preceq L_a\bY_{i,a}.
\end{align*}
Thus, the centered summands have a norm of at most $L_a$, and their variance satisfies
\begin{align*}
    \sum_{i=1}^k\E(\bY_{i,a}-\E\bY_{i,a})^2
    &=\sum_{i=1}^k\left[\E(\bY_{i,a})^2-(\E\bY_{i,a})^2\right]
    \preceq kL_a\bSigma_{>m}.
\end{align*}
The deterministic variance upper bound $kL_a\bSigma_{>m}$ has norm $kL_a\sigma_{m+1}$ and intrinsic dimension $r_1(m)$. We can then apply the intrinsic-dimension Bernstein inequality in~\citet[Proof of Lemma~C.2]{cheng2024dimensionfree} with summand bound $L_a$ and this variance upper bound for $s\geq\sqrt{kL_a\sigma_{m+1}}+L_a/3$,
\begin{align*}
    \P\prn{\left\|\sum_{i=1}^k(\bY_{i,a}-\E\bY_{i,a})\right\|>s}
    \leq8r_1(m) \cdot \exp\left\{-\frac{s^2}{2(kL_a\sigma_{m+1}+L_as/3)}\right\}.
\end{align*}
Taking $s$ to be a sufficiently large constant times the expression below, and using $\log r_1(m)+a\geq1$, bounds the failure probability by $O(e^{-a})$, which implies
\begin{align*}
    \left\|\sum_{i=1}^k(\bY_{i,a}-\E\bY_{i,a})\right\|
    =O\left(\sqrt{kL_a\sigma_{m+1}(\log r_1(m)+a)}+L_a(\log r_1(m)+a)\right).
\end{align*}
By the union bound, with probability at least $1-2e^{-a}$ we have $\bW_{k,m}^*\bW_{k,m}=\sum_i\bY_{i,a}$, while $\|\sum_{i=1}^k\E\bY_{i,a}\|\leq k\sigma_{m+1}$. Adding the mean and the centered fluctuation consequently proves, with probability $1-O(e^{-a})$,
\begin{align}
    & \|\bW_{k,m}^*\bW_{k,m}\| \\
    &=O_{C_{\bx}}\prn{k\sigma_{m+1}+\sqrt{k\sigma_{m+1}\tau_1(m)(\log k+a)(\log r_1(m)+a)} +\tau_1(m)(\log k+a)(\log r_1(m)+a)}.
    \hspace{1em} \label{eq:norm-bound-lower-block-tail}
\end{align}
To express this threshold in terms of $\const_{\bA}(k,m)$, observe that
\begin{align*}
    (\log k+a)(\log r_1(m)+a)\leq\left(1+\frac a{\log k}\right)^2\log (k) \log(kr_1(m)).
\end{align*}
Thus Eq.~\eqref{eq:norm-bound-lower-block-tail} bounds $\|\bW_{k,m}^*\bW_{k,m}\|$ by $O_{C_{\bx}}((1+a/\log k)^2\const_{\bA}(k,m))$ with failure probability $O(e^{-a})$ for every $a\geq1$. We are now ready to combine this tail bound with its failure probability and bound the $L^p$ norm by integrating the tail. For every $p\geq1$,
\begin{align}
    \norm{\|\bW_{k,m}^*\bW_{k,m}\|}_{L^p}
    &=\brc{p\int_0^\infty t^{p-1}\P\left(\|\bW_{k,m}^*\bW_{k,m}\|>t\right) \de t } ^{1/p}\nonumber\\*
    &\stackrel{\mathrm{(i)}}{\leq}O_{C_{\bx}}(\const_{\bA}(k,m))\brc{\left(1+\frac1{\log k}\right)^{2p}+\frac{2p}{\log k}\int_1^\infty\left(1+\frac a{\log k}\right)^{2p-1}e^{-a}\de a}^{1/p}\nonumber\\
    &\stackrel{\mathrm{(ii)}}{\leq}O_{C_{\bx}}(\const_{\bA}(k,m))\brc{\int_0^\infty\left(1+\frac a{\log k}\right)^{2p}e^{-a}\de a}^{1/p}\nonumber\\
    &\stackrel{\mathrm{(iii)}}{\leq}O_{C_{\bx}}(\const_{\bA}(k,m))\left(1+\frac{\Gamma(2p+1)^{1/(2p)}}{\log k}\right)^2.
    \label{eq:norm-bound-lower-block-moment}
\end{align}
In (i), we split at the threshold corresponding to $a=1$, use the bound $1$ below it, and substitute $t=\const\const_{\bA}(k,m)(1+a/\log k)^2$ above it. The fixed constants in the tail probability are absorbed using $p\geq1$. For (ii), integration by parts gives
\begin{align*}
    \frac{2p}{\log k}\int_1^\infty\left(1+\frac a{\log k}\right)^{2p-1}e^{-a}\de a
    =-e^{-1}\left(1+\frac1{\log k}\right)^{2p}+\int_1^\infty\left(1+\frac a{\log k}\right)^{2p}e^{-a}\de a.
\end{align*}
The remaining boundary term is absorbed because $\int_1^\infty(1+a/\log k)^{2p}e^{-a}\de a\geq e^{-1}(1+1/\log k)^{2p}$ and $e^{1/p}\leq e$. In (iii), we apply Minkowski in $L^{2p}([0,\infty),e^{-a}\de a)$ and then square the resulting bound. Finally, using $\Gamma(2p+1)^{1/(2p)}=O(p)$ in Eq.~\eqref{eq:norm-bound-lower-block-moment} proves
\begin{align*}
    \bigl\|\|\bW_{k,m}^*\bW_{k,m}\|\bigr\|_{L^p}=O_{C_{\bx},\kappa}(\const_{\bA}(k,m)),\qquad 1\leq p\leq2\kappa\log k.
\end{align*}
This is the moment bound required in Part I.
\begin{remark}
    Under the stronger effective-rank condition $k=\Omega_{C_{\bx}}(r_1(m))$, applying~\citep[Theorem~1 and Example~1]{zhivotovskiy2022dimensionfree} instead gives $\|\bW_{k,m}^*\bW_{k,m}\|=O_{C_{\bx}}(k\sigma_{m+1})$ with failure probability $\exp\{-\Omega_{C_{\bx}}(k)\}$. Thus one can remove the logarithmic factors in the tail estimate in this regime. However, we do not pursue this refinement, since it improves the downstream bounds only by logarithmic factors.
\end{remark}
\paragraph{Part III: Completing the resolvent estimates.}
We now apply the moment bound from Part II to the resolvent estimates on $\mc E(k,m)$ from Part I, and then control the contribution of $\mc E(k,m)^c$. Throughout this part, $1\leq p\leq\kappa\log k$. For the first contribution, we apply the Minkowski inequality and the $L^{2p}$ bound from Part II to get
\begin{align*}
    \norm{1+\frac{2\|\bW_{k,m}^*\bW_{k,m}\|}{\zeta}}_{L^{2p}}
    \leq1+\frac2\zeta\norm{\|\bW_{k,m}^*\bW_{k,m}\|}_{L^{2p}}
    =O_{C_{\bx},\kappa}\prn{1+\frac{\const_{\bA}(k,m)}\zeta}.
\end{align*}
Thus it controls both terms involving the tail block in Eq.~\eqref{eq:norm-bound-headwise-estimates}. To obtain the trace estimate as well, taking the trace in the operator comparison preceding that equation gives, on $\mc E(k,m)$,
\begin{align*}
    \Tr\bA_k(\zeta)\leq\frac{2m}{k}\prn{1+\frac{2\|\bW_{k,m}^*\bW_{k,m}\|}\zeta}+\frac{2\tau_1(m)}\zeta.
\end{align*}
We again apply Minkowski's inequality, using $\sigma_{m+1}\leq\const_{\bA}(k,m)/k$ to obtain
\begin{align*}
    \prn{\int_{\mc E(k,m)}\|\bA_k(\zeta)\|^p\de\P}^{1/p}
    &=O_{C_{\bx},\kappa}\prn{\frac1k\prn{1+\frac{\const_{\bA}(k,m)}\zeta}},\\
    \prn{\int_{\mc E(k,m)}\|\bA_k(\zeta)\|_F^{2p}\de\P}^{1/p}
    &=O_{C_{\bx},\kappa}\prn{\frac m{k^2}\prn{1+\frac{\const_{\bA}(k,m)}\zeta}^2+\frac{\tau_2(m)}{\zeta^2}},\\
    \prn{\int_{\mc E(k,m)}(\Tr\bA_k(\zeta))^p\de\P}^{1/p}
    &=O_{C_{\bx},\kappa}\prn{\frac m k\prn{1+\frac{\const_{\bA}(k,m)}\zeta}+\frac{\tau_1(m)}\zeta}.
\end{align*}

It remains to control the same three quantities on $\mc E(k,m)^c$. Here we control deterministically by $0\preceq\bA_k(\zeta)\preceq\bSigma/\zeta$. Applying Eq.~\eqref{eq:norm-bound-head-event} and using $p\leq\kappa\log k$ implies
\begin{align*}
    \P(\mc E(k,m)^c)^{1/p}
    \leq\exp\{-\Omega_{C_{\bx}}(k/p)\}
    \leq\exp\{-\Omega_{C_{\bx}}(k/(\kappa\log k))\}\leq k^{-D}.
\end{align*}
Since $k/(\log k)^2\to\infty$ as $k\to\infty$, the last inequality holds for $k=\Omega_{C_{\bx},\kappa,D}(1)$, uniformly over the stated range of $p$. Using $\|\bSigma\|=1$, we therefore obtain
\begin{align*}
    \prn{\int_{\mc E(k,m)^c}\|\bA_k(\zeta)\|^p\de\P}^{1/p}
    &\leq\frac{\P(\mc E(k,m)^c)^{1/p}}\zeta\leq\frac{k^{-D}}\zeta,\\
    \prn{\int_{\mc E(k,m)^c}\|\bA_k(\zeta)\|_F^{2p}\de\P}^{1/p}
    &\leq\frac{\|\bSigma\|_F^2\P(\mc E(k,m)^c)^{1/p}}{\zeta^2}\leq\frac{\|\bSigma\|_F^2k^{-D}}{\zeta^2},\\
    \prn{\int_{\mc E(k,m)^c}(\Tr\bA_k(\zeta))^p\de\P}^{1/p}
    &\leq\frac{\Tr(\bSigma)\P(\mc E(k,m)^c)^{1/p}}\zeta\leq\frac{\Tr(\bSigma)k^{-D}}\zeta.
\end{align*}
Finally, splitting each quantity over $\mc E(k,m)$ and $\mc E(k,m)^c$ and applying the Minkowski inequality in $L^p$ one last time proves all stated results in Lemma~\ref{lem:weighted-resolvent-norm-control}, Eqs.~\eqref{eq:real-resolvent-uniform-moments} uniformly in $p$.

\subsubsection{Proof of Lemma~\ref{lem:one-step-G-bound}}\label{proof:one-step-G-bound}
To bound the one-step difference, we compare both endpoints with their common Euler approximation. Their leading terms then cancel exactly, leaving two remainders to estimate. Recall that $\delta_{k+1}=(1+\wb G_{k+1}(z,w_{k+1}))^{-1}$. For real $z<0$, the path parameters are nonpositive, so $\bT_j-w_k\bSigma-z\bI\succeq|z|\bI$ and both the quadratic-form and expected-trace denominators are at least one. Thus the resolvent identities and moment estimates below apply directly in this case as well. We insert the Euler approximation as follows:
\begin{align}
    \wb G_{k+1}(z,w_{k+1})
    \quad\longleftrightarrow\quad
    \wb G_k(z,w_{k+1})-\delta_{k+1}\frac{\partial}{\partial w}\wb G_k(z,w_{k+1})
    \quad\longleftrightarrow\quad
    \wb G_k(z,w_k).
    \label{eq:one-step-three-points}
\end{align}
The first comparison adds one sample while keeping $w=w_{k+1}$ fixed, whose leading contribution is
\begin{align*}
    \wb G_{k+1}(z,w_{k+1})-\wb G_k(z,w_{k+1})
    =-\delta_{k+1}\frac{\partial}{\partial w}\wb G_k(z,w_{k+1})+\text{higher-order fluctuation}.
\end{align*}
The second comparison moves $w_{k+1}$ to $w_k=w_{k+1}-\delta_{k+1}$ while keeping the sample size fixed, so
\begin{align*}
    \wb G_k(z,w_k)-\wb G_k(z,w_{k+1})
    =-\delta_{k+1}\frac{\partial}{\partial w}\wb G_k(z,w_{k+1})+\text{higher-order resolvent remainder}.
\end{align*}
The two increments have the same leading term, which cancels when we subtract them. Only the two higher-order remainders remain. The identity $\delta_{k+1}=(1+\wb G_{k+1}(z,w_{k+1}))^{-1}$ is exactly what matches the two leading terms. Viewed in this way, the comparisons in~\eqref{eq:one-step-three-points} are the expected-trace analogue of the cancellation strategy in~\citet[Proof of Theorem~6, Part~II]{cheng2024dimensionfree}. Their strategy requires inserting an auxiliary resolvent between consecutive iterates and splitting the increment into a rank-one update and a coefficient update, and their update rule explicitly cancels these conditional leading terms.

In the comparisons in Eq.~\eqref{eq:one-step-three-points}, subtracting the common Euler approximation leaves a nonlinear rank-one fluctuation in the first comparison and a resolvent remainder in the second. We now derive their exact formulas and bound them separately.

\paragraph{Part I: Explicit higher-order terms $\chi_1$ and $\chi_2$.}  We first write the two endpoints and their common Euler approximation as expected traces:
\begin{align*}
    \wb G_{k+1}(z,w_{k+1})&=\E\Tr\prn{\bSigma^{\half}\bR_k(z,w_{k+1})\bSigma^{\half}},\\
    \wb G_k(z,w_{k+1})-\delta_{k+1}\frac{\partial}{\partial w}\wb G_k(z,w_{k+1})
    &=\E\Tr\prn{\bSigma^{\half}\bR_{k-1}(z,w_{k+1})\bSigma^{\half}} -\delta_{k+1}\E\Tr\prn{\bSigma^{\half}\bR_{k-1}(z,w_{k+1})\bSigma^{\half}}^2,\\
    \wb G_k(z,w_k)&=\E\Tr\prn{\bSigma^{\half}\bR_{k-1}(z,w_k)\bSigma^{\half}}.
\end{align*}
Using the weighted resolvent notation $\wt{\bR}_j(z,w)$ from Lemma~\ref{lem:weighted-resolvent-complex-to-real}, applying the Sherman-Morrison identity, and then using the isotropy of $\bz_k$ conditional on $\mc F_{k-1}$ gives
\begin{align*}
    \wb G_{k+1}(z,w_{k+1})
    &=\E\Tr\prn{\wt{\bR}_{k-1}(z,w_{k+1})}
    -\E\brk{\frac{\bz_k^*\prn{\wt{\bR}_{k-1}(z,w_{k+1})}^2\bz_k}{1+\bz_k^*\wt{\bR}_{k-1}(z,w_{k+1})\bz_k}} \\
    & = \E\Tr\prn{\wt{\bR}_{k-1}(z,w_{k+1})}-\delta_{k+1}\E\Tr\prn{\wt{\bR}_{k-1}(z,w_{k+1})}^2 \\
    & \hspace{8em} + \underbrace{\prn{\frac{\E \brk{\bz_k^*\prn{\wt{\bR}_{k-1}(z,w_{k+1})}^2\bz_k}}{1+\wb G_{k+1}(z,w_{k+1})} -\E\brk{\frac{\bz_k^*\prn{\wt{\bR}_{k-1}(z,w_{k+1})}^2\bz_k}{1+\bz_k^*\wt{\bR}_{k-1}(z,w_{k+1})\bz_k}}}}_{=:\chi_1}
\end{align*}
It is equally instructive to expand the third equation around the same trace at $w_{k+1}$. Applying the resolvent identity twice gives
\begin{align*}
    \wb G_k(z,w_k)
    &=\E\Tr\prn{\wt{\bR}_{k-1}(z,w_{k+1})}- \delta_{k+1}\E\Tr\prn{\wt{\bR}_{k-1}(z,w_{k+1})}^2\\
    & \hspace{12em} +\underbrace{\delta_{k+1}^2\E\Tr\brk{
    \wt{\bR}_{k-1}(z,w_k)
    \prn{\wt{\bR}_{k-1}(z,w_{k+1})}^2}}_{=:\chi_2}.
\end{align*}
Thus the second equation by Euler approximation exactly recover the common first-order approximation to both endpoints in the two displays above, and $\chi_1, \chi_2$ are the remaining two higher-order terms.  From the resolvent identity, we can write
\begin{align*}
    \wt{\bR}_{k-1}(z,w_{k+1})-\wt{\bR}_{k-1}(z,w_k)
    &=\delta_{k+1}\wt{\bR}_{k-1}(z,w_{k+1})\wt{\bR}_{k-1}(z,w_k)\\
    &=\delta_{k+1}\wt{\bR}_{k-1}(z,w_k)\wt{\bR}_{k-1}(z,w_{k+1}).
\end{align*}
In this notation, the comparisons in~\eqref{eq:one-step-three-points} combine into the exact decomposition
\begin{align}
    \wb G_{k+1}(z,w_{k+1})-\wb G_k(z,w_k)
    &=\brc{\wb G_{k+1}(z,w_{k+1})
    -\prn{\wb G_k(z,w_{k+1})-\delta_{k+1}\frac{\partial}{\partial w}\wb G_k(z,w_{k+1})}}\nonumber\\
    &\quad+\brc{\prn{\wb G_k(z,w_{k+1})-\delta_{k+1}\frac{\partial}{\partial w}\wb G_k(z,w_{k+1})}
    -\wb G_k(z,w_k)}\nonumber\\
    &=\chi_1-\chi_2.
    \label{eq:one-step-euler-decomposition}
\end{align}
The two remainders also take the following equivalent forms,
\begin{align*}
    \chi_1
    &=
    \frac{\E\brk{\bz_k^*\wt{\bR}_{k-1}(z,w_{k+1})^2\bz_k}}{1+\wb G_{k+1}(z,w_{k+1})}
    -\E\brk{\frac{\bz_k^*\wt{\bR}_{k-1}(z,w_{k+1})^2\bz_k}{1+\bz_k^*\wt{\bR}_{k-1}(z,w_{k+1})\bz_k}},\\*
    \chi_2&=
    \delta_{k+1}^2\E\Tr(\wt{\bR}_{k-1}(z,w_k)\wt{\bR}_{k-1}(z,w_{k+1})^2).
\end{align*}
The parameter error term $\chi_2$ is immediate to control. From its explicit formula, applying the trace Cauchy-Schwarz and then the H\"older's inequality gives
\begin{align}
    |\chi_2|&\leq|\delta_{k+1}|^2\E\brk{\norm{\wt{\bR}_{k-1}(z,w_k)}_F\norm{\wt{\bR}_{k-1}(z,w_{k+1})}\norm{\wt{\bR}_{k-1}(z,w_{k+1})}_F} =O\prn{|\delta_{k+1}|^2\gamma_k\varphi_k^2}.
    \label{eq:one-step-euler-remainder}
\end{align}
In contrast, the analysis of $\chi_1$ requires a more fine-grained approach. It is a nonlinear rank-one fluctuation whose numerator and denominator depend on the same fresh sample. The subsequent parts analyze $\chi_1$ in detail, and we first present a useful decomposition of $\chi_1$.

\paragraph{Part II: A bias-covariance decomposition for $\chi_1$.}
We start from the expression for $\chi_1$ obtained in Part I. In the remainder of this proof and the proof of Lemma~\ref{lem:one-step-covariance-bounds}, we suppress the fixed arguments $(z,w_{k+1})$ from $\wb G_{k+1}$ and the weighted resolvents $\wt{\bR}_j(z,w_{k+1})$. Set
\begin{align*}
    U_k:=\bz_k^*\wt{\bR}_{k-1}\bz_k,
    \quad \text{and} \quad
    V_k:=\frac{\bz_k^*\wt{\bR}_{k-1}^2\bz_k}{1+\bz_k^*\wt{\bR}_{k-1}\bz_k}.
\end{align*}
Using $\delta_{k+1}(1+\wb G_{k+1})=1$ gives
\begin{align*}
    \chi_1 = \delta_{k+1} \E \brk{\bz_k^*\wt{\bR}_{k-1}^2\bz_k \prn{1- \frac{1 + \wb{G}_{k+1}}{1+\bz_k^*\wt{\bR}_{k-1}\bz_k}}}=\delta_{k+1}\E\brk{(U_k-\wb G_{k+1})V_k}.
\end{align*}
We can simplify the expectation of the first factor $U_k - \wb{G}_{k+1}$ using $\E[\bz_k\bz_k^*]=\bI$ and then applying Sherman-Morrison on $\wt{\bR}_{k-1}$:
\begin{align}
    \E U_k-\wb G_{k+1}
    =\E\Tr(\wt{\bR}_{k-1})-\wb G_{k+1}
    =\E V_k.
    \label{eq:one-step-mean-identity}
\end{align}
Subtracting the product of the means therefore yields the exact bias-covariance decomposition
\begin{align}
    \chi_1
    =\delta_{k+1}\prn{(\E V_k)^2+\Cov(U_k,V_k)}.
    \label{eq:one-step-bias-covariance}
\end{align}
The bias term is relatively easy to estimate. Applying the rank-one identity gives a direct bound for the bias. The weighted resolvent $\wt{\bR}_k = \wt{\bR}_k(z,w_{k+1})$ is obtained from $\wt{\bR}_{k-1}$ by adding the $k$-th sample. Applying Sherman-Morrison shows $V_k = \Tr(\wt{\bR}_{k-1} - \wt{\bR}_k)$, and
\begin{align}
    |V_k|
    =\left|\Tr(\wt{\bR}_{k-1}-\wt{\bR}_k)\right|
    \leq \norm{\wt{\bR}_{k-1}}+\norm{\wt{\bR}_k}.
    \label{eq:one-step-robust-ratio}
\end{align}
The last inequality is a consequence of $\wt{\bR}_{k-1}-\wt{\bR}_k$ having rank at most one. Both resolvents are indexed by $\mc I_k$, so $\norm{V_k}_{L^4}=O_{C_{\bx}}(\gamma_k)$. Consequently,
\begin{align}
    \left|\delta_{k+1}(\E V_k)^2\right|=O_{C_{\bx}}\prn{|\delta_{k+1}|\gamma_k^2}.
    \label{eq:one-step-bias-bound}
\end{align}
Thus it remains only to control the covariance term, which we analyze in the next part.

\paragraph{Part III: Bounding the covariance term $\delta_{k+1}\Cov(U_k,V_k)$.} We first derive a decomposition around the stable denominator $1+\wb{G}_{k+1}$ to replace the random denominator $1+U_k$. The decomposition allows us to control the covariance between $U_k$ and $V_k$ by $U_k$ with smaller and low-rank terms. We first derive the exact denominator expansion that organizes the proof. Since $1+U_k=1+\wb G_{k+1}+(U_k-\wb G_{k+1})$, we can explicitly expand
\begin{align*}
    \frac{1}{1+U_k}
    &=\frac{1}{1+\wb G_{k+1}}
    -\frac{U_k-\wb G_{k+1}}{(1+\wb G_{k+1})(1+U_k)} =\frac{1}{1+\wb G_{k+1}}
    -\frac{U_k-\wb G_{k+1}}{(1+\wb G_{k+1})^2}
    +\frac{(U_k-\wb G_{k+1})^2}{(1+\wb G_{k+1})^2(1+U_k)}.
\end{align*}
Multiplying this identity by $\bz_k^*\wt{\bR}_{k-1}^2\bz_k$ on both sides, by the definition of $V_k$ we can write
\begin{align}
    V_k
    &=\frac{\bz_k^*\wt{\bR}_{k-1}^2\bz_k}{1+U_k}=\frac{\bz_k^*\wt{\bR}_{k-1}^2\bz_k}{1+\wb G_{k+1}}
    -\frac{\bz_k^*\wt{\bR}_{k-1}^2\bz_k (U_k-\wb G_{k+1})}{\prn{1+\wb G_{k+1}}^2}
    +\frac{V_k(U_k-\wb G_{k+1})^2}{\prn{1+\wb G_{k+1}}^2}.
    \label{eq:one-step-deterministic-denominator-expansion}
\end{align}
Taking covariance with $U_k$ in Eq.~\eqref{eq:one-step-deterministic-denominator-expansion} and by its linearity in the second argument, we derive
\begin{align}
    &\Cov(U_k,V_k)=\frac{\Cov\prn{U_k,\bz_k^*\wt{\bR}_{k-1}^2\bz_k}}{1+\wb G_{k+1}}
    -\frac{\Cov\prn{U_k,\bz_k^*\wt{\bR}_{k-1}^2\bz_k (U_k-\wb G_{k+1})}}{\prn{1+\wb G_{k+1}}^2}
    +\frac{\Cov\prn{U_k,V_k(U_k-\wb G_{k+1})^2}}{\prn{1+\wb G_{k+1}}^2}\nonumber\\
    &=\delta_{k+1}\underbrace{\Cov\prn{U_k,\bz_k^*\wt{\bR}_{k-1}^2\bz_k}}_{=:\psi_1}
    -\delta_{k+1}^2\underbrace{\Cov\prn{U_k,\bz_k^*\wt{\bR}_{k-1}^2\bz_k(U_k-\wb G_{k+1})}}_{=:\psi_2}
    +\delta_{k+1}^2\underbrace{\Cov\prn{U_k,V_k(U_k-\wb G_{k+1})^2}}_{=:\psi_3}.
    \label{eq:one-step-three-covariances} \hspace{2em}
\end{align}
This decomposition replaces the random denominator $1+U_k$ with powers of the deterministic denominator $1+\wb G_{k+1}$. Thus $\Cov(U_k,V_k)=\delta_{k+1}\psi_1-\delta_{k+1}^2\psi_2+\delta_{k+1}^2\psi_3$. We bound the three covariances separately from the deterministic factors involving $\delta_{k+1}$, which we collect in the following lemma. We defer its proof to Appendix~\ref{proof:one-step-covariance-bounds}.
\begin{lemma}\label{lem:one-step-covariance-bounds}
    Under the assumptions of Lemma~\ref{lem:one-step-G-bound}, let $\psi_1,\psi_2,\psi_3$ be the covariances defined in Eq.~\eqref{eq:one-step-three-covariances}. Then
    \begin{align*}
        |\psi_1|&=O_{C_{\bx}}\prn{k\gamma_k^3+\gamma_k\varphi_k^2},&|\psi_2|&=O_{C_{\bx}}\prn{\varphi_k^2(\sqrt{k}\gamma_k+\varphi_k)^2}, & |\psi_3|&=O_{C_{\bx}}\prn{\gamma_k(\sqrt{k}\gamma_k+\varphi_k)^3}.
    \end{align*}
\end{lemma}
We are now ready to combine all bounds in the final step.

\paragraph{Part IV: Combining the bounds.} Applying Lemma~\ref{lem:one-step-covariance-bounds} to the decomposition in Eq.~\eqref{eq:one-step-three-covariances} and using the triangle inequality gives

\begin{align}
    \left|\delta_{k+1}\Cov(U_k,V_k)\right|
    &=O_{C_{\bx}}\Bigl(|\delta_{k+1}|^2 \cdot (k\gamma_k^3+\gamma_k\varphi_k^2) +|\delta_{k+1}|^3\cdot \prn{\varphi_k^2(\sqrt{k}\gamma_k+\varphi_k)^2+\gamma_k(\sqrt{k}\gamma_k+\varphi_k)^3}\Bigr). \hspace{2em}
    \label{eq:one-step-covariance-bound}
\end{align}

Combining Eq.~\eqref{eq:one-step-covariance-bound} with the decomposition in Eqs.~\eqref{eq:one-step-euler-decomposition} and~\eqref{eq:one-step-bias-covariance}, the Euler-remainder bound in Eq.~\eqref{eq:one-step-euler-remainder}, and the bias bound in Eq.~\eqref{eq:one-step-bias-bound}, then applying $|\delta_{k+1}|\leq |z|/|z_-|$ from Lemma~\ref{lem:resolvent-bounds-strict-positivity}, Eq.~\eqref{eq:helper-denominator-separation}, proves Lemma~\ref{lem:one-step-G-bound}, Eq.~\eqref{eq:one-step-explicit-denominator-delta-free}. For its second bound, Eq.~\eqref{eq:one-step-explicit-denominator-delta}, we instead use $|\delta_{k+1}|\leq\Delta^{-1}$ and

\begin{align*}
    \gamma_k^2+\gamma_k\varphi_k^2+k\gamma_k^3+\varphi_k^2(\sqrt{k}\gamma_k+\varphi_k)^2 &=O\prn{(\gamma_k+\varphi_k^2)^2+k(\gamma_k+\varphi_k^2)^3},\\
    \gamma_k(\sqrt{k}\gamma_k+\varphi_k)^3 & \leq2(\sqrt{k}\gamma_k+\varphi_k)(k\gamma_k^3+\gamma_k\varphi_k^2).
\end{align*}
The condition $\Delta=\Omega(\sqrt{k}\gamma_k+\varphi_k)$ subsumes the last factor into the big-Oh notation $O_{C_{\bx},\Delta}$. Combining these estimates once again with Eqs.~\eqref{eq:one-step-euler-decomposition}, \eqref{eq:one-step-euler-remainder}, \eqref{eq:one-step-bias-covariance}, and~\eqref{eq:one-step-bias-bound} proves Eq.~\eqref{eq:one-step-explicit-denominator-delta} under its additional assumptions. The proof is complete.

\subsubsection{Proof of Lemma~\ref{lem:one-step-covariance-bounds}}\label{proof:one-step-covariance-bounds}
We use the same notations in Appendix~\ref{proof:one-step-G-bound} (i.e. proof of Lemma~\ref{lem:one-step-G-bound}). The covariances $\psi_2$ and $\psi_3$ are more straightforward to control, while $\psi_1$ contains two sources of randomness---the fresh sample $\bz_k$ and the leave-one-out resolvent $\wt{\bR}_{k-1}$. We will condition on the first $k-1$ samples to separate the two contributions. We control the first source by quadratic-form concentration, and the second source by covariance tensorization (cf. Lemma~\ref{lem:covariance-tensorization}).

\paragraph{Part I: Bounding $\psi_1$.}
Since $\wt{\bR}_{k-1}$ depends only on the first $k-1$ samples, conditional isotropy provides $\E[U_k\mid\bT_{k-1}]=\Tr(\wt{\bR}_{k-1})$ and $\E[\bz_k^*\wt{\bR}_{k-1}^2\bz_k\mid\bT_{k-1}]=\Tr(\wt{\bR}_{k-1}^2)$. Subtracting these conditional means in the two quadratic forms gives
\begin{align}
    \psi_1&=\E\brk{(U_k-\Tr\wt{\bR}_{k-1})\prn{\bz_k^*\wt{\bR}_{k-1}^2\bz_k-\Tr(\wt{\bR}_{k-1}^2)}}+\Cov\prn{\Tr(\wt{\bR}_{k-1}),\Tr(\wt{\bR}_{k-1}^2)}.
    \label{eq:one-step-polynomial-covariance-decomposition}
\end{align}
The first term measures the fresh-sample fluctuation, whereas the second measures the fluctuation of the two leave-one-out traces.

For the first term, applying the conditional Hanson-Wright in Proposition~\ref{prop:hanson-wright} and integrating over tail yields
\begin{align}
    \norm{U_k-\Tr\wt{\bR}_{k-1}}_{L^2(\bz_k\mid\bT_{k-1})}&=O_{C_{\bx}}\prn{\norm{\wt{\bR}_{k-1}}_F},\nonumber\\
    \norm{\bz_k^*\wt{\bR}_{k-1}^2\bz_k-\Tr(\wt{\bR}_{k-1}^2)}_{L^2(\bz_k\mid\bT_{k-1})}&=O_{C_{\bx}}\prn{\norm{\wt{\bR}_{k-1}}\norm{\wt{\bR}_{k-1}}_F},
    \label{eq:one-step-conditional-quadratic-moments}
\end{align}
where we use $\norm{\wt{\bR}_{k-1}^2}_F\leq\norm{\wt{\bR}_{k-1}}\norm{\wt{\bR}_{k-1}}_F$ and $\norm{\wt{\bR}_{k-1}}\leq\norm{\wt{\bR}_{k-1}}_F$. For complex operators, the same estimates follow by repeating the same analyses to their Hermitian real and imaginary parts. The H\"older's inequality with exponents $(2,4,4)$ therefore implies
\begin{align*}
    \left|\E\brk{(U_k-\Tr\wt{\bR}_{k-1})\prn{\bz_k^*\wt{\bR}_{k-1}^2\bz_k-\Tr(\wt{\bR}_{k-1}^2)}}\right|
    &=O_{C_{\bx}}\prn{\E\brk{\norm{\wt{\bR}_{k-1}}\norm{\wt{\bR}_{k-1}}_F^2}}=O_{C_{\bx}}(\gamma_k\varphi_k^2).
\end{align*}
This completes the control of the first term in Eq.~\eqref{eq:one-step-polynomial-covariance-decomposition}. To control the second covariance of traces term, we estimate how each trace changes when one of the first $k-1$ observations is replaced. For $1\leq i\leq k-1$, let $\wt{\bR}_{k-1}^{(i)}$ denote $\wt{\bR}_{k-1}$ with $\bz_i$ replaced by an independent copy. The two sample covariance operators differ by an operator of rank at most two; applying the resolvent identity therefore shows that $\wt{\bR}_{k-1}-\wt{\bR}_{k-1}^{(i)}$ also has rank at most two. Consequently,
\begin{align*}
    |\Tr(\wt{\bR}_{k-1})-\Tr(\wt{\bR}_{k-1}^{(i)})|&\leq2\norm{\wt{\bR}_{k-1}-\wt{\bR}_{k-1}^{(i)}}\leq2\prn{\norm{\wt{\bR}_{k-1}}+\norm{\wt{\bR}_{k-1}^{(i)}}}.
\end{align*}
For the trace of the squared resolvent, using
$\wt{\bR}_{k-1}^2-(\wt{\bR}_{k-1}^{(i)})^2=\wt{\bR}_{k-1}(\wt{\bR}_{k-1}-\wt{\bR}_{k-1}^{(i)})+(\wt{\bR}_{k-1}-\wt{\bR}_{k-1}^{(i)})\wt{\bR}_{k-1}^{(i)}$
shows that this difference has rank at most four and operator norm at most $(\|\wt{\bR}_{k-1}\|+\|\wt{\bR}_{k-1}^{(i)}\|) \cdot \|\wt{\bR}_{k-1}-\wt{\bR}_{k-1}^{(i)}\|$. Since $\wt{\bR}_{k-1}^{(i)}$ has the same distribution as $\wt{\bR}_{k-1}$, taking moments and applying H\"older's inequality with exponents $(4,4)$ to the product yields
\begin{align}
    \norm{\Tr(\wt{\bR}_{k-1})-\Tr(\wt{\bR}_{k-1}^{(i)})}_{L^4}&=O_{C_{\bx}}(\gamma_k) \quad \text{and} \quad
    \norm{\Tr(\wt{\bR}_{k-1}^2)-\Tr((\wt{\bR}_{k-1}^{(i)})^2)}_{L^2}=O_{C_{\bx}}(\gamma_k^2).
    \label{eq:one-step-trace-replacement-moments}
\end{align}
Both traces are functions of the independent samples $\bz_1,\ldots,\bz_{k-1}$ and belong to $L^2$. We can thus apply Efron-Stein via covariance tensorization in Lemma~\ref{lem:covariance-tensorization} with these two traces. Using Eq.~\eqref{eq:one-step-trace-replacement-moments}, we have
\begin{align*}
    \left|\Cov\prn{\Tr(\wt{\bR}_{k-1}),\Tr(\wt{\bR}_{k-1}^2)}\right|
    &\leq\sum_{i=1}^{k-1}\norm{\Tr(\wt{\bR}_{k-1})-\Tr(\wt{\bR}_{k-1}^{(i)})}_{L^2}\norm{\Tr(\wt{\bR}_{k-1}^2)-\Tr((\wt{\bR}_{k-1}^{(i)})^2)}_{L^2}\\*
    &=O_{C_{\bx}}(k\gamma_k^3).
\end{align*}
Combining the two contributions in Eq.~\eqref{eq:one-step-polynomial-covariance-decomposition} completes the proof for
\begin{align}
    |\psi_1|=O_{C_{\bx}}\prn{k\gamma_k^3+\gamma_k\varphi_k^2}.
    \label{eq:one-step-polynomial-covariance-bound}
\end{align}

\paragraph{Part II: Moment bounds for the remaining covariances.}
The covariances $\psi_2$ and $\psi_3$ contain moments of $U_k-\wb G_{k+1}$, which we aim to control first. Applying Eq.~\eqref{eq:one-step-mean-identity} gives the decomposition
\begin{align*}
    U_k-\wb G_{k+1}=(U_k-\Tr\wt{\bR}_{k-1})+(\Tr\wt{\bR}_{k-1}-\E\Tr\wt{\bR}_{k-1})+\E V_k.
\end{align*}
The first summand is the fresh quadratic-form fluctuation. Applying Hanson-Wright in Proposition~\ref{prop:hanson-wright} conditionally with moment order $4$ and then averaging over $\bT_{k-1}$ gives $\norm{U_k-\Tr\wt{\bR}_{k-1}}_{L^4}=O_{C_{\bx}}(\varphi_k)$. For the second summand, we use the replacement estimate already established in the previous part (cf.~Eq.~\eqref{eq:one-step-trace-replacement-moments}). By Efron-Stein in Lemma~\ref{lem:moment-efron-stein} with $q=4$, the conditional Jensen's inequality and Minkowski's inequality for the square-function bound, we can bound
\begin{align*}
    \norm{\Tr\wt{\bR}_{k-1}-\E\Tr\wt{\bR}_{k-1}}_{L^4}
    &=O\prn{\prn{\sum_{i=1}^{k-1}\norm{\Tr(\wt{\bR}_{k-1})-\Tr(\wt{\bR}_{k-1}^{(i)})}_{L^4}^2}^{1/2}}
    =O_{C_{\bx}}(\sqrt{k}\gamma_k).
\end{align*}
The bound $|\E V_k|=O_{C_{\bx}}(\gamma_k)$ was already proved in Part II of Lemma~\ref{lem:one-step-G-bound}, using Eq.~\eqref{eq:one-step-robust-ratio}. Combining these estimates by the triangle inequality yields
\begin{align}
    \norm{U_k-\wb G_{k+1}}_{L^4}=O_{C_{\bx}}\prn{\sqrt{k}\gamma_k+\varphi_k}.
    \label{eq:one-step-denominator-fluctuation}
\end{align}

The remaining components in $\psi_2$ and $\psi_3$ include $\bz_k^*\wt{\bR}_{k-1}^2\bz_k$ and $V_k$. Using the conditional estimate in Eq.~\eqref{eq:one-step-conditional-quadratic-moments} with $|\Tr(\wt{\bR}_{k-1}^2)|\leq\norm{\wt{\bR}_{k-1}}_F^2$ and $\norm{\wt{\bR}_{k-1}}\leq\norm{\wt{\bR}_{k-1}}_F$, then averaging over $\bT_{k-1}$, gives the first bound below. The second follows from the rank-one bound~\eqref{eq:one-step-robust-ratio}, as already used above:
\begin{align}
    \norm{\bz_k^*\wt{\bR}_{k-1}^2\bz_k}_{L^2}=O_{C_{\bx}}(\varphi_k^2)\qquad \text{and} \qquad
    \norm{V_k}_{L^4}=O_{C_{\bx}}(\gamma_k).
    \label{eq:one-step-numerator-ratio-moments}
\end{align}

\paragraph{Part III: Bounding $\psi_2$ and $\psi_3$.}
Recall we can write
\begin{align*}
    \psi_2&=\E\brk{(U_k-\E U_k)(\bz_k^*\wt{\bR}_{k-1}^2\bz_k)(U_k-\wb G_{k+1})} \qquad \text{and} \qquad \psi_3=\E\brk{(U_k-\E U_k)V_k(U_k-\wb G_{k+1})^2}.
\end{align*}
Although $U_k-\wb G_{k+1}$ need not be centered, applying Jensen's inequality gives $|\E U_k-\wb G_{k+1}|\leq\norm{U_k-\wb G_{k+1}}_{L^4}$ and thus $\norm{U_k-\E U_k}_{L^4}\leq2\norm{U_k-\wb G_{k+1}}_{L^4}$. Invoking H\"older's inequality with exponents $(4,2,4)$ to $\psi_2$ and using Eqs.~\eqref{eq:one-step-denominator-fluctuation} and~\eqref{eq:one-step-numerator-ratio-moments} implies
\begin{align}
    |\psi_2|&\leq2\norm{\bz_k^*\wt{\bR}_{k-1}^2\bz_k}_{L^2} \cdot \norm{U_k-\wb G_{k+1}}_{L^4}^2
    =O_{C_{\bx}}\prn{\varphi_k^2(\sqrt{k}\gamma_k+\varphi_k)^2}.
    \label{eq:one-step-first-higher-order-covariance}
\end{align}
For $\psi_3$, applying H\"older with exponents $(4,4,4,4)$ yields
\begin{align}
    |\psi_3|&\leq2\norm{V_k}_{L^4} \cdot \norm{U_k-\wb G_{k+1}}_{L^4}^3
    =O_{C_{\bx}}\prn{\gamma_k(\sqrt{k}\gamma_k+\varphi_k)^3}.
    \label{eq:one-step-second-higher-order-covariance}
\end{align}
This proves all three bounds and we complete the proof.

\subsection{Auxiliary proofs for Theorem~\ref{thm:terminal-weighted-traces}}\label{proof:auxiliary-theorem2}

\subsubsection{Proof of Lemma~\ref{lem:trace-step-decomposition}}\label{proof:trace-step-decomposition}
The proof is a pure algebraic verification. We compute directly along the backward path, first for the resolvents
\begin{align*}
    & \bR_{k}(z,w_k) - \bR_{k-1}(z,w_{k-1}) = \prn{\bR_{k}(z,w_k) - \bR_{k-1}(z,w_k)} + \prn{\bR_{k-1}(z,w_k) - \bR_{k-1}(z,w_{k-1})} \\
    & = -\frac{\bR_{k-1}(z,w_k)\bx_k\bx_k^*\bR_{k-1}(z,w_k)}{1+\bx_k^*\bR_{k-1}(z,w_k)\bx_k} + \delta_k\bR_{k-1}(z,w_k)\bSigma\bR_{k-1}(z,w_{k-1}),
\end{align*}
where we use the Sherman-Morrison formula and the path update $\delta_k=w_k-w_{k-1}$. Next, we can recall the definitions of $H_k^{\bQ}$ and $G_k^{\bQ}$ in Eq.~\eqref{eq:def-traces-interpolation}. Using $(\bT_j-w\bSigma)\bR_j(z,w)=\bI+z\bR_j(z,w)$, we can make the one-step update explicit using the calculation above. For $1\leq k\leq n$, the resolvent update implies
\begin{align*}
    & H_k^{\bQ}(z,w_k)-H_{k-1}^{\bQ}(z,w_{k-1}) =z\Tr\prn{\bQ\prn{\bR_k(z,w_k)-\bR_{k-1}(z,w_{k-1})}} \\
    & =-V_k^{H,\bQ}(w_k)+z\delta_k\Tr\prn{\bQ\bR_{k-1}(z,w_k)\bSigma\bR_{k-1}(z,w_{k-1})}=\xi_k^{H,\bQ}+\beta_k^{H,\bQ}.
\end{align*}
For $G_k^{\bQ}$, the parameter pair is $(w_{k+1},w_k)$. Applying the same resolvent calculation with this pair and taking the trace weighted by $\bSigma^{1/2}\bQ\bSigma^{1/2}$, we have for $1\leq k<n$,
\begin{align*}
    &G_{k+1}^{\bQ}(z,w_{k+1})-G_k^{\bQ}(z,w_k)
    =\Tr\prn{\bSigma^{1/2}\bQ\bSigma^{1/2}\prn{\bR_k(z,w_{k+1})-\bR_{k-1}(z,w_k)}}\\
    &=-V_k^{G,\bQ}(w_{k+1})+\delta_{k+1}\Tr\prn{\bQ\bSigma^{1/2}\bR_{k-1}(z,w_{k+1})\bSigma\bR_{k-1}(z,w_k)\bSigma^{1/2}} =\xi_k^{G,\bQ}+\beta_k^{G,\bQ}.
\end{align*}
The deterministic path makes the parameter-change terms $\mc F_{k-1}$-measurable. Thus both drift term sequences are predictable, and conditional centering gives $\E_{k-1}\xi_k^{H,\bQ}=\E_{k-1}\xi_k^{G,\bQ}=0$. We therefore complete the proof.

\subsubsection{Proof of Lemma~\ref{lem:trace-accumulated-drift}}\label{proof:trace-accumulated-drift}
We first bound the martingale partial sums in (i). For (ii), we separate each predictable increment into its full expectation and its centered fluctuation. Indeed, taking full expectation recovers the Euler decomposition in Lemma~\ref{lem:one-step-G-bound}, which makes the role of our backward path explicit. It suffices to work under either $\|\bQ\|\leq1$ or $\|\bQ\|_F\leq\sqrt n$, as rescaling restores the factor $\|\bQ\| \wedge \|\bQ\|_F/\sqrt n$.
\paragraph{Part I: Bounding the martingale partial sums.} To control the sums of martingale differences $\xi_k^{H, \bQ}$ and $\xi_k^{G, \bQ}$, it suffices to bound their $L^p$ norms, as the martingale inequality below will show. By Lemma~\ref{lem:trace-step-decomposition}, each $\xi$-term is the negative of a conditionally centered $V$-term. The next lemma supplies a general bound for $V$-term fluctuations, which will also be used later in the proof for the drift $\beta$-terms analysis. We defer its proof to Appendix~\ref{proof:trace-centered-increments}.
\begin{lemma}\label{lem:trace-centered-increments}
    Under the same setup as Theorem~\ref{thm:terminal-weighted-traces}, let $0\leq j\leq n$ and $(j-3)\vee0\leq s\leq j$. Form the ratios $V_{s+1}^{H,\bQ}(w_j)$ and $V_{s+1}^{G,\bQ}(w_j)$ of Lemma~\ref{lem:trace-step-decomposition}, Eq.~\eqref{eq:trace-rank-one-ratios}, using $\bR_s(z,w_j)$ and a fresh observation $\bx_{s+1}$ independent of $\mc F_s$. Then, for $2\leq p\leq\kappa\log n$,
    \begin{align}
        \max\brc{\frac{\|V_{s+1}^{H,\bQ}(w_j)-\E_s V_{s+1}^{H,\bQ}(w_j)\|_{L^p}}n,\|V_{s+1}^{G,\bQ}(w_j)-\E_s V_{s+1}^{G,\bQ}(w_j)\|_{L^p}} & \nonumber \\
        & \hspace{-10em} =O_{C_{\bx},\nu,\kappa}\prn{\frac{p\Lambda_n(z)^{11}}{n^{3/2}} \cdot \prn{\|\bQ\| \wedge \frac{\|\bQ\|_F}{\sqrt n}}}.
        \label{eq:trace-ratio-fluctuations}
    \end{align}
\end{lemma}
We immediately recover the $\xi$-terms in the martingale partial sums by taking $(s,j)=(k-1,k)$ for $H$ and $(s,j)=(k-1,k+1)$ for $G$ in the above Lemma~\ref{lem:trace-centered-increments}. For $H$ we use the separate power $9$ in Part II of the proof of that lemma. Therefore, we obtain
\begin{align}
    \frac1n\left\|\max_{0\leq j\leq n}\left|\sum_{k=1}^j\xi_k^{H,\bQ}\right|\right\|_{L^p}
    &\leq \frac{O(p)}n\prn{\sum_{k=1}^n\|\xi_k^{H,\bQ}\|_{L^p}^2}^{1/2}
    =O\prn{\frac{p^2\Lambda_n(z)^{9}}n},\nonumber\\*
    \left\|\max_{0\leq j<n}\left|\sum_{k=1}^j\xi_k^{G,\bQ}\right|\right\|_{L^p}
    &\leq O(p)\prn{\sum_{k=1}^{n-1}\|\xi_k^{G,\bQ}\|_{L^p}^2}^{1/2}
    =O\prn{\frac{p^2\Lambda_n(z)^{11}}n}.
    \label{eq:martingale-maximal-moment}
\end{align}
In both inequalities above, we apply Doob's maximal inequality~\citep[Theorem~4.4.4]{durrett2019probability} and Burkholder's inequality~\citep[Theorem~2.10]{hall1980martingale} to the real and imaginary parts, and then apply the Minkowski inequality in $L^{p/2}$ to reduce to sum of $L^p$ norms of $\xi$-terms. We suppress dependence on $C_{\bx},\nu,\kappa$ in the big-Oh notation to simplify presentation, and we shall keep this convention throughout the remainder of this proof. This completes the proof of Lemma~\ref{lem:trace-accumulated-drift}.(i).

\paragraph{Part II: Bounding the accumulated drift.}
We next derive an upper bound for the accumulated $\beta$-terms. This entails establishing bounds on both the full expectation and the centered fluctuation of each predictable increment, in addition to controlling the initial correction term. We start with the exact decomposition
\begin{align}
    \beta_k^{H,\bQ}=\E\beta_k^{H,\bQ}+(\beta_k^{H,\bQ}-\E\beta_k^{H,\bQ})
    \qquad \text{and} \qquad
    \beta_k^{G,\bQ}=\E\beta_k^{G,\bQ}+(\beta_k^{G,\bQ}-\E\beta_k^{G,\bQ}).
    \label{eq:trace-drift-mean-fluctuation}
\end{align}
We point out how the drift terms connect to our proof for the path stability of $\{w_k\}_{k=0}^n$. Taking full expectation in Lemma~\ref{lem:trace-step-decomposition}, Eq.~\eqref{eq:trace-step-G}, with $\bQ=\bI$ recovers Eq.~\eqref{eq:one-step-euler-decomposition} in the proof of Lemma~\ref{lem:one-step-G-bound}:
\begin{align*}
    \E\beta_k^{G,\bI}
    =\wb G_{k+1}(z,w_{k+1})-\wb G_k(z,w_k)=\chi_1-\chi_2,
\end{align*}
since $\E\xi_k^{G,\bI}=0$. This gives rise to our proof strategy by repeating the Euler calculation with the test $\bQ$ under conditional expectation $\E_{k-1}$, so that the conditional trace of the weighted resolvent, $\Tr\wt{\bR}_{k-1}=\E_{k-1}[\bz_k^*\wt{\bR}_{k-1}\bz_k]$, remains explicit. Below, we shall suppress $(z,w_k)$ in the $\beta_k^{H, \bQ}$ terms and $(z,w_{k+1})$ in the $\beta_k^{G, \bQ}$ terms, and write $U_k=\bz_k^*\wt{\bR}_{k-1}\bz_k$ at the corresponding evaluation point. Applying the Euler identity from Part I of the proof of Lemma~\ref{lem:one-step-G-bound}, Eq.~\eqref{eq:one-step-euler-decomposition}, to Eq.~\eqref{eq:trace-predictable-drifts}, we can write
\begin{align}
    \beta_k^{H,\bQ}&=\chi_1^{H,\bQ}-\chi_2^{H,\bQ} \qquad \text{and} \qquad \beta_k^{G,\bQ}=\chi_1^{G,\bQ}-\chi_2^{G,\bQ},\label{eq:trace-three-drift-terms}
\end{align}
where conditional isotropy and $\delta_j(1+\wb G_j)=1$ identify
\begin{align*}
    \chi_1^{H,\bQ}&=\delta_k\E_{k-1}[(U_k-\wb G_k)V_k^{H,\bQ}],
    &
    \chi_1^{G,\bQ} &=\delta_{k+1}\E_{k-1}[(U_k-\wb G_{k+1})V_k^{G,\bQ}],\\
    \chi_2^{H,\bQ}&=z\delta_k^2\Tr[\bQ\bR_{k-1}\bSigma\bR_{k-1}(z,w_{k-1})\bSigma\bR_{k-1}],  & \chi_2^{G,\bQ}&=\delta_{k+1}^2\Tr[\bQ\wt{\bR}_{k-1}\wt{\bR}_{k-1}(z,w_k)\wt{\bR}_{k-1}].
\end{align*}
We use $\bR_{k-1}(z,w_{k-1})=\bR_{k-1}-\delta_k\bR_{k-1}(z,w_{k-1})\bSigma\bR_{k-1}$ inside the trace for the $H$-terms. The $G$-terms calculation is the weighted version of the same identity. We first bound the remainders $\chi_2^{H,\bQ}$ and $\chi_2^{G,\bQ}$ directly. For $\chi_2^{G,\bQ}$, trace Cauchy-Schwarz and the operator-Frobenius product inequality give
\begin{align*}
    \left|\Tr\prn{\bQ\wt{\bR}_{k-1}\wt{\bR}_{k-1}(z,w_k)\wt{\bR}_{k-1}}\right|
    &\leq\|\wt{\bR}_{k-1}(z,w_k)\|\prn{\|\wt{\bR}_{k-1}\|_F^2\vee\sqrt n\|\wt{\bR}_{k-1}\|\|\wt{\bR}_{k-1}\|_F}.
\end{align*}
The two alternatives correspond to $\|\bQ\|\leq1$ and $\|\bQ\|_F\leq\sqrt n$, respectively. For $\chi_2^{H,\bQ}$, we can factor the trace through $\bR_{k-1}\bSigma^{1/2}$, $\wt{\bR}_{k-1}(z,w_{k-1})$, and $\bSigma^{1/2}\bR_{k-1}$. Applying the same product inequalities and Eq.~\eqref{eq:trace-H-factor-bounds} from Part I-ii of the proof of Lemma~\ref{lem:trace-centered-increments} gives
\begin{align*}
    \frac{|\chi_2^{H,\bQ}|}n
    &\leq\frac{2|z||\delta_k|^2}{n|z_-|}\|\wt{\bR}_{k-1}(z,w_{k-1})\|\prn{\Tr|\wt{\bR}_{k-1}|\vee\sqrt{n\|\wt{\bR}_{k-1}\|\Tr|\wt{\bR}_{k-1}|}}.
\end{align*}
Using $|\delta_j|\leq\Lambda_n(z)$ and $|z|/|z_-|\leq\Lambda_n(z)$, then applying H\"older and the resolvent moments in Eq.~\eqref{eq:trace-growing-resolvent-moments} from Part II of the proof of Lemma~\ref{lem:trace-centered-increments}, we obtain
\begin{align}
    \|\chi_2^{G,\bQ}\|_{L^p}
    &=O\prn{\Lambda_n(z)^2\frac{\Lambda_n(z)^{3}}n\frac{\Lambda_n(z)^{6}}n}
    =O\prn{\frac{\Lambda_n(z)^{11}}{n^2}},\nonumber\\
    \frac{\|\chi_2^{H,\bQ}\|_{L^p}}n
    &=O\prn{\frac{\Lambda_n(z)^3}n\frac{\Lambda_n(z)^{3}}n\Lambda_n(z)^{3}}
    =O\prn{\frac{\Lambda_n(z)^{9}}{n^2}},
    \label{eq:trace-euler-drift-bounds}
\end{align}
using the H\"older inequality with exponents $(3,3,3)$ in the first line, and in the second line $(2,2)$ first, followed by $(2,2)$ for the square-root product. The weighted resolvent moments used above have order at most $3p\leq4\kappa\log n$, and thus inequalities in Eq.~\eqref{eq:trace-growing-resolvent-moments} apply.

It remains to bound $\chi_1^{H, \bQ}$ and $\chi_1^{G,\bQ}$. Both expansions contain the conditional trace fluctuation $\Tr\wt{\bR}_{k-1}-\E\Tr\wt{\bR}_{k-1}$. Using $\E_{k-1}U_k=\Tr\wt{\bR}_{k-1}$ and the mean identity in Part II of the proof of Lemma~\ref{lem:one-step-G-bound}, Eq.~\eqref{eq:one-step-mean-identity}, we have for $\chi_1^{G,\bQ}$,
\begin{align}
    & \chi_1^{G,\bQ}
    =\delta_{k+1}\E_{k-1}[(U_k-\E U_k)V_k^{G,\bQ}]
    +\delta_{k+1} \cdot \E V_k^{G,\bI} \cdot \E_{k-1}V_k^{G,\bQ}\nonumber\\
    &=\delta_{k+1}(\Tr\wt{\bR}_{k-1}-\E\Tr\wt{\bR}_{k-1})\E_{k-1}V_k^{G,\bQ}
    +\delta_{k+1}\Cov(U_k,V_k^{G,\bQ}\mid\mc F_{k-1}) +\delta_{k+1}\cdot \E V_k^{G,\bI} \cdot \E_{k-1}V_k^{G,\bQ}. \hspace{1em}
    \label{eq:trace-conditional-drift-components}
\end{align}
For the second line, we split $U_k-\E U_k=(U_k-\E_{k-1}U_k)+(\Tr\wt{\bR}_{k-1}-\E\Tr\wt{\bR}_{k-1})$ and use $\E_{k-1}(U_k-\E_{k-1}U_k)=0$. For $\chi_1^{H,\bQ}$, we have similarly
\begin{align}
    & \chi_1^{H,\bQ}
    =\delta_k\E_{k-1}[(U_k-\E U_k)V_k^{H,\bQ}] \\ & =\delta_k(\Tr\wt{\bR}_{k-1}-\E\Tr\wt{\bR}_{k-1})\E_{k-1}V_k^{H,\bQ}
    +\delta_k\Cov(U_k,V_k^{H,\bQ}\mid\mc F_{k-1}).
    \label{eq:trace-conditional-drift-components-H}
\end{align}
We note that only the $\chi_1^{G, \bQ}$ expansion has the additional mean correction. The following lemma bounds $\chi_1$-terms and the initial correction, making use of the decomposition in the above displays. We postpone its proof to Appendix~\ref{proof:trace-drift-components}.

\begin{lemma}\label{lem:trace-drift-components}
    Under the same setup as Theorem~\ref{thm:terminal-weighted-traces}, for $2\leq p\leq\kappa\log n$,
    \begin{subequations}
        \begin{align}
            \frac{\|\chi_1^{H,\bQ}\|_{L^p}}n\vee\|\chi_1^{G,\bQ}\|_{L^p}
            &=O_{C_{\bx},\nu,\kappa}\prn{\frac{p^2\Lambda_n(z)^{19}}{n^2} \cdot \prn{\|\bQ\| \wedge \frac{\|\bQ\|_F}{\sqrt n}}},\label{eq:trace-chi-one-bounds}\\
            \frac{|\beta_0^{H,\bQ}|}n\vee|\beta_0^{G,\bQ}|
            &=O_{C_{\bx},\nu,\kappa}\prn{\frac{\Lambda_n(z)^{23}}n \cdot \prn{\|\bQ\| \wedge \frac{\|\bQ\|_F}{\sqrt n}}}.
            \label{eq:trace-initial-corrections}
        \end{align}
    \end{subequations}
\end{lemma}
We can now combine the two parts of each drift. Applying Lemma~\ref{lem:trace-drift-components} above to the decomposition in Eq.~\eqref{eq:trace-three-drift-terms}, we obtain
\begin{align*}
    \frac{\|\beta_k^{H,\bQ}\|_{L^p}}n\vee\|\beta_k^{G,\bQ}\|_{L^p}
    &\leq\prn{\frac{\|\chi_1^{H,\bQ}\|_{L^p}}n\vee\|\chi_1^{G,\bQ}\|_{L^p}}+\prn{\frac{\|\chi_2^{H,\bQ}\|_{L^p}}n\vee\|\chi_2^{G,\bQ}\|_{L^p}} =O\prn{\frac{p^2\Lambda_n(z)^{19}}{n^2}}.
\end{align*}
For the scalar estimate in Theorem~\ref{thm:gram-stieltjes-approx} we preserve a separate $H$ bound here. For $\chi_2^{H,\bQ}$, the second line of Eq.~\eqref{eq:trace-euler-drift-bounds} gives the exponent $9$. For $\chi_1^{H,\bQ}$, the last display in Part II of the proof of Lemma~\ref{lem:trace-drift-components} bounds the two terms of Eq.~\eqref{eq:trace-conditional-drift-components-H} separately. The conditional trace term has exponent $17$ and the conditional covariance has power $13$. Summing from $k=1$ and applying Minkowski's inequality to the $n$ one-step drifts yields a sharpened bound
\begin{align}
    \frac1n\left\|\sum_{k=1}^n|\beta_k^{H,\bQ}|\right\|_{L^p} =O\prn{\frac{p^2\Lambda_n(z)^{17}}n}.\label{eq:gram-positive-index-drift}
\end{align}
Applying Minkowski's inequality to the $n$ increments for $H$ and the $n-1$ increments for $G$, we have
\begin{align}
    \max\brc{\frac1n\left\|\sum_{k=1}^n|\beta_k^{H,\bQ}|\right\|_{L^p},
    \left\|\sum_{k=1}^{n-1}|\beta_k^{G,\bQ}|\right\|_{L^p}}
    =O\prn{\frac{p^2\Lambda_n(z)^{19}}{n}}.
    \label{eq:trace-summed-predictable-drifts}
\end{align}
The proof is complete by adding the initial corrections bounded in Lemma~\ref{lem:trace-drift-components}, Eq.~\eqref{eq:trace-initial-corrections}.

\subsubsection{Proof of Lemma~\ref{lem:trace-centered-increments}}\label{proof:trace-centered-increments}
Within the proof, we assume either $\|\bQ\|\leq1$ or $\|\bQ\|_F\leq\sqrt n$ as linearity in $\bQ$ allows rescaling.
\paragraph{Part I: Reduction to numerator and denominator moments.}
For any fixed $s, j$ under the setup of Lemma~\ref{lem:trace-centered-increments}, let $\wt{\bR}_s(z, w_j) = \bSigma^{\half} \bR_s(z, w_j) \bSigma^{\half}$, and we shall suppress the dependence on $(z,w_j)$ for notational simplicity. Starting from the left-hand side of Eq.~\eqref{eq:trace-ratio-fluctuations}, we apply conditional Jensen's inequality by freezing the denominator first, and then the triangle inequality to obtain
\begin{align*}
    \max\left\{\frac{\|V_{s+1}^{H,\bQ}-\E_sV_{s+1}^{H,\bQ}\|_{L^p}}n,\ \|V_{s+1}^{G,\bQ}-\E_sV_{s+1}^{G,\bQ}\|_{L^p}\right\} &\\
    &\hspace{-12em} \leq2\max\left\{\frac1n\left\|V_{s+1}^{H,\bQ}-\frac{z\Tr(\bQ\bR_s\bSigma\bR_s)}{1+\Tr\wt{\bR}_s}\right\|_{L^p},\ \left\|V_{s+1}^{G,\bQ}-\frac{\Tr(\bQ\wt{\bR}_s^2)}{1+\Tr\wt{\bR}_s}\right\|_{L^p}\right\}.
\end{align*}
It then boils down to bounding the right hand side in the above display. To see what estimates we need, for the $G$ term, we expand
\begin{align}
    V_{s+1}^{G,\bQ}-\frac{\Tr(\bQ\wt{\bR}_s^2)}{1+\Tr\wt{\bR}_s}
    &=\frac{\bz_{s+1}^*\wt{\bR}_s\bQ\wt{\bR}_s\bz_{s+1}-\Tr(\bQ\wt{\bR}_s^2)}
    {1+\bz_{s+1}^*\wt{\bR}_s\bz_{s+1}} -\frac{\Tr(\bQ\wt{\bR}_s^2)(\bz_{s+1}^*\wt{\bR}_s\bz_{s+1}-\Tr\wt{\bR}_s)}
    {(1+\Tr\wt{\bR}_s)(1+\bz_{s+1}^*\wt{\bR}_s\bz_{s+1})}.
    \label{eq:trace-ratio-centered-identity}
\end{align}
We bound the denominators above first. The quadratic-form can be controlled by Lemma~\ref{lem:resolvent-bounds-strict-positivity}.(i). For the trace denominator, note that $\Tr\wt{\bR}_s=G_{s+1}(z,w_j)$. Repeating the argument for $\wb G_{s+1}$ in Part II of that lemma's proof before taking expectation gives $\Im\prn{z(1+\Tr\wt{\bR}_s)}\geq\eta$, as in Eq.~\eqref{eq:helper-pointwise-trace-signs}, and $\Re\Tr\wt{\bR}_s\geq0$ when $t<0$. Thus both denominators satisfy
\begin{align*}
    |1+\bz_{s+1}^*\wt{\bR}_s\bz_{s+1}|^{-1}
    \vee|1+\Tr\wt{\bR}_s|^{-1}\leq |z|/|z_-|\leq\Lambda_n(z).
\end{align*}
Applying the denominator bounds and Cauchy-Schwarz to Eq.~\eqref{eq:trace-ratio-centered-identity} yields
\begin{subequations}\label{eq:trace-centered-ratio-reduction}
    \begin{align}
        \|V_{s+1}^{G,\bQ}-\E_sV_{s+1}^{G,\bQ}\|_{L^p}
        \leq2\Lambda_n(z)\|\bz_{s+1}^*\wt{\bR}_s\bQ\wt{\bR}_s\bz_{s+1}-\Tr(\bQ\wt{\bR}_s^2)\|_{L^p} &\nonumber\\
        & \hspace{-14em} +2\Lambda_n(z)^2\|\Tr(\bQ\wt{\bR}_s^2)\|_{L^{2p}}
        \|\bz_{s+1}^*\wt{\bR}_s\bz_{s+1}-\Tr\wt{\bR}_s\|_{L^{2p}}.
        \label{eq:trace-centered-ratio-G}
    \end{align}
    Similarly, for $V_{s+1}^{H,\bQ}$, we have
    \begin{align}
        \frac{\|V_{s+1}^{H,\bQ}-\E_sV_{s+1}^{H,\bQ}\|_{L^p}}n
        &\leq\frac{2\Lambda_n(z)}n\|z\bz_{s+1}^*\bSigma^{1/2}\bR_s\bQ\bR_s\bSigma^{1/2}\bz_{s+1}-z\Tr(\bQ\bR_s\bSigma\bR_s)\|_{L^p}\nonumber\\
        &\quad+\frac{2\Lambda_n(z)^2}n\|z\Tr(\bQ\bR_s\bSigma\bR_s)\|_{L^{2p}}
        \|\bz_{s+1}^*\wt{\bR}_s\bz_{s+1}-\Tr\wt{\bR}_s\|_{L^{2p}}.
        \label{eq:trace-centered-ratio-H}
    \end{align}
\end{subequations}
Thus it suffices to control the numerator moments on the right-hand sides above.
\subparagraph{Part I-i: Reduction of $G$-terms numerator moments.} For $2\leq q\leq2\kappa\log n$, applying conditional Hanson-Wright and integrating the tail, as in Eq.~\eqref{eq:one-step-conditional-quadratic-moments}, we have
\begin{align}
    \|\bz_{s+1}^*\wt{\bR}_s\bQ\wt{\bR}_s\bz_{s+1}-\Tr(\bQ\wt{\bR}_s^2)\|_{L^q} \leq & O_{C_{\bx}}(q)\bigl\|\|\wt{\bR}_s\bQ\wt{\bR}_s\|_F\bigr\|_{L^q} \nonumber \\
    \stackrel{\mathrm{(i)}}{\leq} & O_{C_{\bx}}(q) \cdot \prn{\bigl\|\|\wt{\bR}_s\|\|\wt{\bR}_s\|_F\bigr\|_{L^q} \vee \bigl\|\sqrt n\,\|\wt{\bR}_s\|^2\bigr\|_{L^q}},\nonumber\\
    \|\bz_{s+1}^*\wt{\bR}_s\bz_{s+1}-\Tr\wt{\bR}_s\|_{L^q}
    \leq & O_{C_{\bx}}(q) \cdot \bigl\|\|\wt{\bR}_s\|_F\bigr\|_{L^q},
    \label{eq:trace-conditional-hanson-wright}
\end{align}
where we invoke the operator-Frobenius product inequality in (i), using either $\|\bQ\|\leq1$ or $\|\bQ\|_F\leq\sqrt n$. For the conditional trace in Eq.~\eqref{eq:trace-centered-ratio-G}, we apply the trace Cauchy-Schwarz and the same product inequality to obtain
\begin{align*}
    \|\Tr(\bQ\wt{\bR}_s^2)\|_{L^q}
    \leq \bigl\|\|\wt{\bR}_s\|_F^2\bigr\|_{L^q} \vee \bigl\|\sqrt n\,\|\wt{\bR}_s\|\|\wt{\bR}_s\|_F\bigr\|_{L^q}.
\end{align*}
\subparagraph{Part I-ii: Reduction of $H$-terms numerator moments.}
For the $H$-terms numerator, we will express the bound through the same weighted resolvent $\wt{\bR}_s$. We define the operator modulus by its semidefinite square-root $|\bA| := (\bA^* \bA)^{1/2}$ when $\bA$ is Hilbert-Schmidt. Taking congruences by $\bR_s$ and $\bR_s^*$ in the lower bound $\Re(e^{\iu\theta}\bR_s^{-1})\succeq|z_-|\bI/2$ from Part I of the proof of Lemma~\ref{lem:weighted-resolvent-complex-to-real}, Eq.~\eqref{eq:weighted-common-rotation} and multiplying by $\bSigma^{1/2}$ and taking norms or traces imply
\begin{align}
    \|\bR_s\bSigma^{1/2}\|^2\vee\|\bSigma^{1/2}\bR_s\|^2
    &\leq \frac{2\|\wt{\bR}_s\|}{|z_-|} \qquad \text{and} \qquad
    \|\bR_s\bSigma^{1/2}\|_F^2\vee\|\bSigma^{1/2}\bR_s\|_F^2
    \leq \frac{2\Tr|\wt{\bR}_s|}{|z_-|}.
    \label{eq:trace-H-factor-bounds}
\end{align}
Therefore, using $|z|/|z_-|\leq\Lambda_n(z)$, we have
\begin{align*}
    \|z\bSigma^{1/2}\bR_s\bQ\bR_s\bSigma^{1/2}\|_F
    &\leq2\Lambda_n(z) \cdot \prn{\sqrt{\|\wt{\bR}_s\|\Tr|\wt{\bR}_s|} \vee (\sqrt{n}\|\wt{\bR}_s\|)},\\
    |z\Tr(\bQ\bR_s\bSigma\bR_s)|
    &\leq2\Lambda_n(z) \cdot \prn{\Tr|\wt{\bR}_s| \vee \sqrt{n\|\wt{\bR}_s\|\Tr|\wt{\bR}_s|}}.
\end{align*}
The first alternative in each maximum applies when $\|\bQ\|\leq1$, and the second when $\|\bQ\|_F\leq\sqrt n$.
Applying the same conditional Hanson-Wright calculation for $G$-terms, it holds
\begin{align*}
    \|z\bz_{s+1}^*\bSigma^{1/2}\bR_s\bQ\bR_s\bSigma^{1/2}\bz_{s+1}-z\Tr(\bQ\bR_s\bSigma\bR_s)\|_{L^q} &\\
    &\hspace{-10em} \leq O_{C_{\bx}}(q\Lambda_n(z)) \cdot \prn{\Bigl\|\sqrt{\|\wt{\bR}_s\|\Tr|\wt{\bR}_s|} \Bigr\|_{L^q} \vee \Bigl\| \sqrt{n} \|\wt{\bR}_s\| \Bigr\|_{L^q}}.
\end{align*}
Combining the reduction to weighted resolvent moments above and substituting into Eqs.~\eqref{eq:trace-centered-ratio-G} and~\eqref{eq:trace-centered-ratio-H}, we have
\begin{align}
    \|V_{s+1}^{G,\bQ}-\E_sV_{s+1}^{G,\bQ}\|_{L^p}
    &\leq O_{C_{\bx}}(p) \cdot \biggl[\Lambda_n(z)\bigl\|\|\wt{\bR}_s\|\|\wt{\bR}_s\|_F+\sqrt n\|\wt{\bR}_s\|^2\bigr\|_{L^p}\nonumber\\
    &\qquad+\Lambda_n(z)^2\bigl\|\|\wt{\bR}_s\|_F^2+\sqrt n\|\wt{\bR}_s\|\|\wt{\bR}_s\|_F\bigr\|_{L^{2p}}\bigl\|\|\wt{\bR}_s\|_F\bigr\|_{L^{2p}}\biggr],\nonumber\\
    \frac{\|V_{s+1}^{H,\bQ}-\E_sV_{s+1}^{H,\bQ}\|_{L^p}}n
    &\leq\frac{O_{C_{\bx}}(p)}n \cdot \biggl[\Lambda_n(z)^2\bigl\|\sqrt{\|\wt{\bR}_s\|\Tr|\wt{\bR}_s|}+\sqrt n\|\wt{\bR}_s\|\bigr\|_{L^p}\nonumber\\
    &\qquad+\Lambda_n(z)^3\bigl\|\Tr|\wt{\bR}_s|+\sqrt{n\|\wt{\bR}_s\|\Tr|\wt{\bR}_s|}\bigr\|_{L^{2p}}\bigl\|\|\wt{\bR}_s\|_F\bigr\|_{L^{2p}}\biggr].
    \label{eq:trace-centered-resolvent-reduction}
\end{align}
It remains to control the operator, Frobenius, and trace-norm moments of $\wt{\bR}_s$.

\paragraph{Part II: Bounding the resolvent moments.} Throughout this part, we suppress dependence on $C_{\bx},\nu,\kappa$ in big-Oh notations. To close Eq.~\eqref{eq:trace-centered-resolvent-reduction}, we prove, for $2\leq q\leq4\kappa\log n$,
\begin{align}
    \bigl\|\|\wt{\bR}_s\|\bigr\|_{L^q}
    &=O\prn{\frac{\Lambda_n(z)^{3}}{n}},&
    \bigl\|\|\wt{\bR}_s\|_F\bigr\|_{L^q}
    &=O\prn{\frac{\Lambda_n(z)^{3}}{\sqrt n}},&
    \bigl\|\Tr|\wt{\bR}_s|\bigr\|_{L^q}
    &=O\prn{\Lambda_n(z)^{3}}. \hspace{1.5em}
    \label{eq:trace-growing-resolvent-moments}
\end{align}
The first two bounds follow from our existing arguments. For $j>\ell:=\lfloor n/2\rfloor$, combining Lemma~\ref{lem:weighted-resolvent-complex-to-real}, Eq.~\eqref{eq:path-real-comparisons}, with Lemma~\ref{lem:weighted-resolvent-norm-control}, Eqs.~\eqref{eq:real-resolvent-operator-moment} and~\eqref{eq:real-resolvent-frobenius-moment}, at $\zeta=|z_-|$ and moment orders $q,q/2$ gives the (strengthened) first two bounds, with $\Lambda_n(z)^{3}$ replaced by $\Lambda_n(z)^2$. To see how these lemmas apply, note that $s\geq\ell-2\asymp n$, so the minimizing cutoff in Assumption~\ref{assmp:spectral-stability} satisfies $m\leq\rho(C_{\bx})(\ell-2)\leq\rho(C_{\bx}) s$ and, for sufficiently large $n$, $q\leq4\kappa\log n\leq8\kappa\log s$. We therefore apply Lemma~\ref{lem:weighted-resolvent-norm-control} with $\kappa$ replaced by $8\kappa$ for both orders $q$ and $q/2$. To evaluate the resulting bounds, use Eq.~\eqref{eq:path-stability-parameter-controls} from Part III-i of the proof of Theorem~\ref{thm:self-consistent-path-stability} to control the spectral terms by $\Lambda_n(z)$, and take $D=2\nu+2$ to absorb the exceptional-event terms, as verified immediately after that equation. Lastly, for $j\leq\ell$, combining Lemma~\ref{lem:weighted-resolvent-complex-to-real}, Eq.~\eqref{eq:path-buffer-norms}, with the path lower bound in Eq.~\eqref{eq:path-imaginary-buffer-rate} from that same part gives the first two bounds pointwise.

For the third estimate, we use the rotation from Part I of the proof of Lemma~\ref{lem:weighted-resolvent-complex-to-real}, with $\bU=e^{\iu\theta}\bR_s^{-1}$ and $\bV=\Re\bU$ and obtain that
\begin{align*}
    & \Tr|\wt{\bR}_s|
    =\Tr\left|(\bSigma^{1/2}\bV^{-1/2})(\bV^{1/2}\bU^{-1}\bSigma^{1/2})\right|\\
    &\stackrel{\mathrm{(i)}}{\leq}\|\bSigma^{1/2}\bV^{-1/2}\|_F\,\|\bV^{1/2}\bU^{-1}\bSigma^{1/2}\|_F =\sqrt{\Tr(\bSigma^{1/2}\bV^{-1}\bSigma^{1/2})}\,
    \sqrt{\Tr(\bSigma^{1/2}\bU^{-*}\bV\bU^{-1}\bSigma^{1/2})}\\
    &\stackrel{\mathrm{(ii)}}{\leq}\Tr(\bSigma^{1/2}\bV^{-1}\bSigma^{1/2}).
\end{align*}
Here in (i), we apply the trace-norm Cauchy-Schwarz, i.e.\ $\Tr|\bA\bB|\leq\|\bA\|_F\|\bB\|_F$, to the two displayed factors. In (ii), we take the trace in the positive-operator inequality from Eq.~\eqref{eq:weighted-comparison-order}. We similarly discuss the large sample and small sample regimes separately. For $j>\ell$, the first lower bound in Eq.~\eqref{eq:weighted-common-rotation} in Part I of the proof of Lemma~\ref{lem:weighted-resolvent-complex-to-real} implies
\begin{align*}
    \bigl\|\Tr|\wt{\bR}_s|\bigr\|_{L^q}
    &\leq\frac{3|z|}{|z_-|}\|\Tr\bA_s(|z_-|)\|_{L^q}
    =O(\Lambda_n(z)^2).
\end{align*}
Using the same parameter choices as for the first two estimates, we apply Lemma~\ref{lem:weighted-resolvent-norm-control}, Eq.~\eqref{eq:real-resolvent-trace-moment}, at order $q$. The bounds $\tau_1(m)/|z_-|\leq\Lambda_n(z)$ and $\Tr\bSigma\leq n+\tau_1(m)$ control the last two terms, giving the displayed $O(\Lambda_n(z)^2)$ bound. For $j\leq\ell$, we use the second lower bound in Eq.~\eqref{eq:weighted-common-rotation} from Part I of the proof of Lemma~\ref{lem:weighted-resolvent-complex-to-real} instead to estimate a pointwise bound:
\begin{align*}
    \Tr|\wt{\bR}_s|
    &\leq2\sum_i\frac{\sigma_i}{|z_-|-2\Re(e^{\iu\theta}w_j)\sigma_i} \leq\frac{m}{-\Re(e^{\iu\theta}w_j)}+\frac{2\tau_1(m)}{|z_-|}
    =O(\Lambda_n(z)^{3}).
\end{align*}
The last step uses $m\leq n$ and the path lower bound in Eq.~\eqref{eq:path-imaginary-buffer-rate} from Part III-i of the proof of Theorem~\ref{thm:self-consistent-path-stability}. We complete the proof of the third estimate in Eq.~\eqref{eq:trace-growing-resolvent-moments} by taking its $L^q$ norm. Therefore, substituting all weighted resolvent bounds in Eq.~\eqref{eq:trace-centered-resolvent-reduction} and using Cauchy-Schwarz, we conclude that
\begin{align*}
    \|V_{s+1}^{G,\bQ}-\E_sV_{s+1}^{G,\bQ}\|_{L^p}
    &=O\prn{\frac{p\Lambda_n(z)^{7}}{n^{3/2}}+\frac{p\Lambda_n(z)^{11}}{n^{3/2}}},\\
    \frac{\|V_{s+1}^{H,\bQ}-\E_sV_{s+1}^{H,\bQ}\|_{L^p}}n
    &=O\prn{\frac{p\Lambda_n(z)^{5}}{n^{3/2}}+\frac{p\Lambda_n(z)^{9}}{n^{3/2}}}.
\end{align*}
All resolvent moments above have order at most $4p\leq4\kappa\log n$. Since $\Lambda_n(z)\geq1$, we complete the proof by rescaling $\bQ$.

\paragraph{Part III: Additional uncentered bounds.}
For the later estimation of $\beta$-terms, we record a useful (and cruder) consequence for $2\leq q\leq2\kappa\log n$. Substituting the weighted resolvent bounds in Eq.~\eqref{eq:trace-growing-resolvent-moments} into the numerator inequalities from Part I and applying H\"older, we have
\begin{align}
    \|\Tr(\bQ\wt{\bR}_s^2)\|_{L^q}
    &=O(\Lambda_n(z)^{6}/n),&
    \bigl\|\|\wt{\bR}_s\bQ\wt{\bR}_s\|_F\bigr\|_{L^q}
    &=O(\Lambda_n(z)^{6}/n^{3/2}),\nonumber\\
    \|z\Tr(\bQ\bR_s\bSigma\bR_s)\|_{L^q}
    &=O(\Lambda_n(z)^{4}),&
    \bigl\|\|z\bSigma^{1/2}\bR_s\bQ\bR_s\bSigma^{1/2}\|_F\bigr\|_{L^q}
    &=O(\Lambda_n(z)^{4}/\sqrt n).
    \label{eq:trace-numerator-moments}
\end{align}
To bound the uncentered ratios, apply the reciprocal bound from Part I, then conditional Hanson-Wright and Eq.~\eqref{eq:trace-numerator-moments}:
\begin{align*}
    \|V_{s+1}^{G,\bQ}\|_{L^q}
    &\leq\Lambda_n(z)\prn{\|\Tr(\bQ\wt{\bR}_s^2)\|_{L^q}+O_{C_{\bx}}(q)\bigl\|\|\wt{\bR}_s\bQ\wt{\bR}_s\|_F\bigr\|_{L^q}} =O\prn{\frac{\Lambda_n(z)^{7}}n},\\
    \frac{\|V_{s+1}^{H,\bQ}\|_{L^q}}n
    &\leq\frac{\Lambda_n(z)}n\prn{\|z\Tr(\bQ\bR_s\bSigma\bR_s)\|_{L^q}+O_{C_{\bx}}(q)\bigl\|\|z\bSigma^{1/2}\bR_s\bQ\bR_s\bSigma^{1/2}\|_F\bigr\|_{L^q}}\\
    &=O\prn{\frac{\Lambda_n(z)^{5}}n}.
\end{align*}
where we use $q\leq2\kappa\log n$. We thus obtain from $\Lambda_n(z)\geq1$ that
\begin{align}
    \frac{\|V_{s+1}^{H,\bQ}\|_{L^q}}{n} \vee \|V_{s+1}^{G,\bQ}\|_{L^q}
    =O\prn{\frac{\Lambda_n(z)^{7}}{n}},
    \label{eq:trace-ratio-moments}
\end{align}
Here, as in the preceding parts, either $\|\bQ\|\leq1$ or $\|\bQ\|_F\leq\sqrt n$.

\subsubsection{Proof of Lemma~\ref{lem:trace-drift-components}}\label{proof:trace-drift-components}
Our goal is to control each term in the expansion of $\chi_1$-terms in Eqs.~\eqref{eq:trace-conditional-drift-components} and~\eqref{eq:trace-conditional-drift-components-H}, and additionally provide the zero-sample endpoint errors. Throughout the proof, we shall assume either $\|\bQ\|\leq1$ or $\|\bQ\|_F\leq\sqrt n$, and suppress dependence on $C_{\bx}, \nu, \kappa$ in big-Oh notations. We first bound the conditional trace fluctuation appearing in the first terms of Eqs.~\eqref{eq:trace-conditional-drift-components} and~\eqref{eq:trace-conditional-drift-components-H}.
\paragraph{Part I: The conditional trace fluctuation.} To upper-bound the leading term in each expansion in Part II of the proof of Lemma~\ref{lem:trace-accumulated-drift}, namely in Eqs.~\eqref{eq:trace-conditional-drift-components} and~\eqref{eq:trace-conditional-drift-components-H}, we first require a concentration result for $\Tr\wt{\bR}_{k-1}(z,w_j)$ for $j=k,k+1$. We would like to use Efron-Stein in Lemma~\ref{lem:moment-efron-stein}, and therefore for $k\geq2$, we replace observation $i<k$ and condition on the common remaining sample. This is equivalent to invoking Lemma~\ref{lem:trace-centered-increments} with $V_{s+1}^{G,\bI}$. Applying Sherman-Morrison with common resolvent $(\bT_{k-1}-\bx_i\bx_i^*-w_j\bSigma-z\bI)^{-1}$ allows us to use Lemma~\ref{lem:trace-centered-increments} with $\kappa$ replaced by $6\kappa$, Eq.~\eqref{eq:trace-ratio-fluctuations}, with $s=k-2$ and $j=k$ or $k+1$ since $j-s\leq3$, and implies
\begin{align*}
    \|\Tr(\wt{\bR}_{k-1}-\wt{\bR}_{k-1}^{(i)})\|_{L^q}
    =O(q\Lambda_n(z)^{11}/n^{3/2}).
\end{align*}
We can then use Efron-Stein (cf.~Lemma~\ref{lem:moment-efron-stein}) to the first $k-1$ observations and conclude that for $2\leq q\leq6\kappa\log n$,
\begin{align}
    \|\Tr\wt{\bR}_{k-1}-\E\Tr\wt{\bR}_{k-1}\|_{L^q}
    &\leq O(\sqrt q) \times \prn{\sum_{i=1}^{k-1}\|\Tr(\wt{\bR}_{k-1}-\wt{\bR}_{k-1}^{(i)})\|_{L^q}^2}^{1/2}
    =O\prn{\frac{q^{3/2}\Lambda_n(z)^{11}}{n}}.
    \label{eq:trace-centered-old-samples}
\end{align}
We point out this generalizes the $L^4$ bound already obtained in Part II of the proof of Lemma~\ref{lem:one-step-covariance-bounds}, and for $k=1$ the trace is deterministic.

\paragraph{Part II: Bounding $\chi_1^{H,\bQ}$ and $\chi_1^{G,\bQ}$.}
We begin by bounding the first terms in Eqs.~\eqref{eq:trace-conditional-drift-components} and~\eqref{eq:trace-conditional-drift-components-H}. By Cauchy-Schwarz and the conditional Jensen's inequality, it follows that
\begin{align*}
    \frac{\|\delta_k(\Tr\wt{\bR}_{k-1}-\E\Tr\wt{\bR}_{k-1})\E_{k-1}V_k^{H,\bQ}\|_{L^p}}n&\leq\Lambda_n(z)\|\Tr\wt{\bR}_{k-1}-\E\Tr\wt{\bR}_{k-1}\|_{L^{2p}}\frac{\|V_k^{H,\bQ}\|_{L^{2p}}}n, \\
    \|\delta_{k+1}(\Tr\wt{\bR}_{k-1}-\E\Tr\wt{\bR}_{k-1})\E_{k-1}V_k^{G,\bQ}\|_{L^p}&\leq\Lambda_n(z)\|\Tr\wt{\bR}_{k-1}-\E\Tr\wt{\bR}_{k-1}\|_{L^{2p}}\|V_k^{G,\bQ}\|_{L^{2p}}.
\end{align*}
Here the resolvents are evaluated at $(z,w_k)$ and $(z,w_{k+1})$, respectively, and we use $|\delta_k|\vee|\delta_{k+1}|\leq\Lambda_n(z)$. Applying the conditional trace bound in Eq.~\eqref{eq:trace-centered-old-samples} from Part I and Eq.~\eqref{eq:trace-ratio-moments} from Part III of the proof of Lemma~\ref{lem:trace-centered-increments}, both at order $2p$, bounds each right-hand side by
\begin{align*}
    O\prn{\Lambda_n(z)\frac{p^{3/2}\Lambda_n(z)^{11}}n\frac{\Lambda_n(z)^{7}}n}
    =O\prn{\frac{p^{3/2}\Lambda_n(z)^{19}}{n^2}}.
\end{align*}
For $H$, the uncentered ratio estimate in Part III of the proof of Lemma~\ref{lem:trace-centered-increments} gives $\|V_k^{H,\bQ}\|_{L^{2p}}/n=O(\Lambda_n(z)^5/n)$. Thus the first $H$ term has order $O(p^{3/2}\Lambda_n(z)^{1+11+5}/n^2)$. For the second terms in the same two expansions, recall $U_k=\bz_k^*\wt{\bR}_{k-1}\bz_k$ and we have similarly
\begin{align*}
    \frac{\|\delta_k\Cov(U_k,V_k^{H,\bQ}\mid\mc F_{k-1})\|_{L^p}}n&\leq\Lambda_n(z)\|U_k-\E_{k-1}U_k\|_{L^{2p}}\frac{\|V_k^{H,\bQ}-\E_{k-1}V_k^{H,\bQ}\|_{L^{2p}}}n, \\
    \|\delta_{k+1}\Cov(U_k,V_k^{G,\bQ}\mid\mc F_{k-1})\|_{L^p}&\leq\Lambda_n(z)\|U_k-\E_{k-1}U_k\|_{L^{2p}}\|V_k^{G,\bQ}-\E_{k-1}V_k^{G,\bQ}\|_{L^{2p}}.
\end{align*}
Here $\E_{k-1}U_k = \Tr \wt{\bR}_{k-1}$ exactly, and we can therefore apply Eq.~\eqref{eq:trace-conditional-hanson-wright} from Part I of the proof of Lemma~\ref{lem:trace-centered-increments} to bound the first factors following $\Lambda_n(z)$ above directly. By Eq.~\eqref{eq:trace-growing-resolvent-moments}, we have
\begin{align*}
    \|U_k-\E_{k-1}U_k\|_{L^{2p}} \leq O(p) \cdot \norm{\norm{\wt{\bR}_{k-1}}_F}_{L^{2p}} = O \prn{\frac{p\Lambda_n(z)^{3}}{\sqrt n}}.
\end{align*}

For the last factor, apply Lemma~\ref{lem:trace-centered-increments} with $\kappa$ replaced by $2\kappa$ and moment order $2p$, taking $(s,j)=(k-1,k+1)$ for $G$ and $(s,j)=(k-1,k)$ for $H$. This yields the desired control, so we can bound each right-hand side by
\begin{align*}
    O\prn{\Lambda_n(z)\frac{p\Lambda_n(z)^{3}}{\sqrt n}\frac{p\Lambda_n(z)^{11}}{n^{3/2}}}
    =O\prn{\frac{p^2\Lambda_n(z)^{15}}{n^2}}.
\end{align*}
For $H$, the centered ratio estimate in Part II of the proof of Lemma~\ref{lem:trace-centered-increments} has power $9$, giving $O(p^2\Lambda_n(z)^{1+3+9}/n^2)$ for the second $H$ term.
It only remains to bound the third term in Eq.~\eqref{eq:trace-conditional-drift-components}. Applying Eq.~\eqref{eq:one-step-robust-ratio} from Part II of the proof of Lemma~\ref{lem:one-step-G-bound}, followed by Eq.~\eqref{eq:trace-growing-resolvent-moments} from Part II of the proof of Lemma~\ref{lem:trace-centered-increments}, we obtain
\begin{align*}
    |\E V_k^{G,\bI}|\leq\E|V_k^{G,\bI}|
    \leq\E\|\wt{\bR}_{k-1}(z,w_{k+1})\|+\E\|\wt{\bR}_k(z,w_{k+1})\|
    =O\prn{\frac{\Lambda_n(z)^{3}}n}.
\end{align*}
Further by the conditional Jensen's inequality and Eq.~\eqref{eq:trace-ratio-moments} from Part III of the proof of Lemma~\ref{lem:trace-centered-increments}, we bound the last term by
\begin{align*}
    \|\delta_{k+1} \cdot \E V_k^{G,\bI} \cdot \E_{k-1}V_k^{G,\bQ}\|_{L^p}
    &\leq\Lambda_n(z)|\E V_k^{G,\bI}|\|V_k^{G,\bQ}\|_{L^p}
    =O\prn{\frac{\Lambda_n(z)^{11}}{n^2}}.
\end{align*}
Collecting all the bounds from the preceding displays and substituting them into Eqs.~\eqref{eq:trace-conditional-drift-components} and~\eqref{eq:trace-conditional-drift-components-H} enables us to complete the proof of Eq.~\eqref{eq:trace-chi-one-bounds} by
\begin{align*}
    \|\chi_1^{G,\bQ}\|_{L^p}
    &=O\prn{\frac{p^{3/2}\Lambda_n(z)^{19}+p^2\Lambda_n(z)^{15}+\Lambda_n(z)^{11}}{n^2}},\\
    \frac{\|\chi_1^{H,\bQ}\|_{L^p}}n
    &=O\prn{\frac{p^{3/2}\Lambda_n(z)^{17}+p^2\Lambda_n(z)^{13}}{n^2}}.
\end{align*}
Since $p\geq2$ and $\Lambda_n(z)\geq1$, both bounds are $O(p^2\Lambda_n(z)^{19}/n^2)$.

\paragraph{Part III: The initial corrections.} Finally, we control the endpoint corrections $\beta_0^{H, \bQ} = H_0^{\bQ}(z, w_0) - H_0^{\bQ}(z, w_0\lin)$ and $\beta_0^{G, \bQ} = G_1^{\bQ}(z, w_1) - G_1^{\bQ}(z, w_0\lin)$. We apply the resolvent identity and identify
\begin{align*}
    \beta_0^{H,\bQ}&=(w_0-w_0\lin)\cdot z\Tr\prn{\bQ\bR_0(z,w_0)\bSigma\bR_0(z,w_0\lin)},\\
    \beta_0^{G,\bQ}&=(w_1-w_0\lin)\cdot \Tr\prn{\bQ\wt{\bR}_0(z,w_1)\wt{\bR}_0(z,w_0\lin)}.
\end{align*}
We bound the first path differences factors by applying Eq.~\eqref{eq:path-population-comparison-rate} from the proof of Theorem~\ref{thm:self-consistent-path-stability}, i.e., $|w_0-w_0\lin|=O(\Lambda_n(z)^{17})$, and $|\delta_1|\leq\Lambda_n(z)$ implying $|w_1-w_0\lin|\leq|w_0-w_0\lin|+|\delta_1|=O(\Lambda_n(z)^{17})$.

To bound the next two tracial terms, we recall Eq.~\eqref{eq:trace-numerator-moments} from Part III of the proof of Lemma~\ref{lem:trace-centered-increments} with argument $(z, w_0)$ as follows
\begin{align*}
    \|z\Tr(\bQ\bR_0\bSigma\bR_0)\|_{L^p}
    &=O(\Lambda_n(z)^{4}) \qquad \text{and} \qquad
    \|\Tr(\bQ\wt{\bR}_0^2)\|_{L^p} =O(\Lambda_n(z)^{6}/n).
\end{align*}
It thus remains to replace $\bR_0(z,w_0)$ with $\bR_0(z,w_0\lin)$ and $\wt{\bR}_0(z,w_0)^2$ with $\wt{\bR}_0(z,w_1) \wt{\bR}_0(z,w_0\lin)$ in the above display. As in the zero-sample endpoint, $\bT_0=0$, the resolvents in both traces are diagonal in the eigenbasis of $\bSigma$. On the eigenspace of $\sigma_i$, replacing one resolvent at $w_0$ by the one at $w_0\lin$ multiplies its entry by $(z+w_0\sigma_i)/(z+w_0\lin\sigma_i)$. For the $G$ trace we also replace the other factor at $w_0$ by the one at $w_1$, which contributes $(z+w_0\sigma_i)/(z+w_1\sigma_i)$. It suffices to bound both ratios uniformly in $i$.

We start from the second ratio, using $w_1-w_0=\delta_1$ and applying Eq.~\eqref{eq:trace-growing-resolvent-moments} from the proof of Lemma~\ref{lem:trace-centered-increments} at $(s,j)=(0,1)$. Since $|\delta_1|\leq\Lambda_n(z)$, we obtain
\begin{align*}
    \left|\frac{z+w_0\sigma_i}{z+w_1\sigma_i}\right|=\left|1-\frac{\delta_1\sigma_i}{z+w_1\sigma_i}\right|\leq1+|\delta_1|\,\|\wt{\bR}_0(z,w_1)\|=1+O\prn{\frac{\Lambda_n(z)^{4}}{n}}=O(1).
\end{align*}
For the first ratio, recall the fixed-point equations comparison with Eq.~\eqref{eq:path-perturbed-population-equation} from Part IV-ii of the proof of Theorem~\ref{thm:self-consistent-path-stability}, i.e., $s_n(z)=w_0\lin/(nz)$ and $s_n(z+\varepsilon_n(z)) = w_0/(nz)$ with $|\varepsilon_n(z)|\leq|z_-|/4$ by Eq.~\eqref{eq:path-spectral-perturbation-rate}. Writing $s_n(z)=\int(x-z)^{-1}\mu(\de x)$ as in Lemma~\ref{lem:population-scalar-stability}, we obtain
\begin{align*}
    \left|\frac{z+w_0\sigma_i}{z+w_0\lin\sigma_i}-1\right|&=\frac{n\sigma_i\left|s_n(z+\varepsilon_n(z))-s_n(z)\right|}{|1+n\sigma_i s_n(z)|}\\
    &\stackrel{\mathrm{(i)}}{\leq}\frac{4|\varepsilon_n(z)|}{3}\frac{n\sigma_i\int |x-z|^{-2}\mu(\de x)}{|1+n\sigma_i s_n(z)|}\stackrel{\mathrm{(ii)}}{\leq}\frac{4|\varepsilon_n(z)|}{3|z_-|}\leq \frac{1}{3}.
\end{align*}
In (i) we subtract the two Stieltjes transforms and use $|x-z-\varepsilon_n(z)|\geq3|x-z|/4$. For (ii), taking the imaginary part of the Stieltjes representation yields $\Im s_n(z)=\eta\int |x-z|^{-2}\mu(\de x)$. If $t\geq0$, then $|z_-|=\eta$, so $|1+n\sigma_i s_n(z)|\geq n\sigma_i\Im s_n(z)$ and we have the desired bound. If $t<0$, then $x-t\geq-t$ on the support of $\mu$, so $\Re s_n(z)=\int (x-t)|x-z|^{-2}\mu(\de x)\geq(-t)\int |x-z|^{-2}\mu(\de x)$. Combining the two parts and using $|z_-|=\sqrt{t^2+\eta^2}$ implies (ii) as
\begin{align*}
    |1+n\sigma_i s_n(z)|\geq n\sigma_i|s_n(z)|\geq n\sigma_i|z_-|\int |x-z|^{-2}\mu(\de x).
\end{align*}
Multiplying $\bQ$ by the bounded diagonal multipliers preserves both test norm bounds up to a constant. Taken together, the above displays complete the proof as
\begin{align*}
    \frac{|\beta_0^{H,\bQ}|}{n}&=O\prn{|w_0-w_0\lin| \cdot \sup_{i \geq 1} \left|\frac{z+w_0\sigma_i}{z+w_0\lin\sigma_i}\right| \cdot \frac{\Lambda_n(z)^{4}}{n}}=O\prn{\frac{\Lambda_n(z)^{21}}{n}},\\
    |\beta_0^{G,\bQ}|&=O\prn{|w_1-w_0\lin| \cdot \sup_{i \geq 1} \left|\frac{z+w_0\sigma_i}{z+w_1\sigma_i}\right| \cdot  \sup_{i \geq 1} \left|\frac{z+w_0\sigma_i}{z+w_0\lin\sigma_i}\right| \cdot \frac{\Lambda_n(z)^{6}}{n}}=O\prn{\frac{\Lambda_n(z)^{23}}{n}}.
\end{align*}

\subsection{Auxiliary proofs for Theorem~\ref{thm:terminal-mixed-resolvents}}

\subsubsection{Proof of Lemma~\ref{lem:mixed-endpoint-identities}}\label{proof:mixed-endpoint-identities}
At the full sample endpoint, $z\bR_n(z,0)=\bT_n\bR_n(z,0)-\bI$. Therefore, subtraction of the same identity at $z'$ and applying the resolvent identity implies the following:
\begin{align*}
    z\bR_n(z,0)-z'\bR_n(z',0)=(z-z')\bR_n(z,0)\bT_n\bR_n(z',0).
\end{align*}
We recover Eq.~\eqref{eq:mixed-terminal-T-insertion} by taking the weighted trace.

At the zero sample endpoint, we make use of Eq.~\eqref{eq:terminal-zero-sample-traces} to express $zG_1^{\bQ}(z,w_0\lin(z))$ as the trace against $\bQ$ of a diagonal operator in the eigenbasis of $\bSigma$. On the eigenspace for $\sigma_i$, it holds
\begin{align*}
    \frac{\omega(z/n)\sigma_i}{\sigma_i-\omega(z/n)}-\frac{\omega(z'/n)\sigma_i}{\sigma_i-\omega(z'/n)}
    =\frac{(\omega(z/n)-\omega(z'/n))\sigma_i^2}{(\sigma_i-\omega(z/n))(\sigma_i-\omega(z'/n))}.
\end{align*}
Dividing the display by $z-z'$, using $\wt{\bR}_0=(\omega(z/n)/z)\bSigma(\bSigma-\omega(z/n)\bI)^{-1}$ and that of $\wt{\bR}_0'$ at both spectral parameters and taking the trace against $\bQ$ proves Eq.~\eqref{eq:mixed-initial-identity}. It follows by taking $z \to z'$ to see that both endpoint identities extend to $z=z'$ by differentiation.

\subsubsection{Proof of Lemma~\ref{lem:terminal-mixed-replacement}}\label{proof:terminal-mixed-replacement}
\paragraph{Part I: Bounding the term $\chi_1$.} We write $\chi_1$ into the following two parts:
\begin{align*}
    \chi_1
    &=\underbrace{\E\Tr\prn{\bQ\big(\wt{\bR}_n\wt{\bR}_n'-\wt{\bR}_{n-1}\wt{\bR}_{n-1}'\big)}}_{:=\xi_1} + \underbrace{\Tr(\bQ\wt{\bR}_n\wt{\bR}_n')-\E\Tr(\bQ\wt{\bR}_n\wt{\bR}_n')}_{:=\xi_2},
\end{align*}
where $\xi_1=\E \chi_1$. Therefore, $\xi_2$ can be viewed as its fluctuation, and $\xi_1$ as its residual bias. We will bound those two terms separately. The residual bias $\xi_1$ can be bounded by a leave-one-out argument of $\bx_n$, and the fluctuation term $\xi_2$ can be bounded by reusing the analysis of $\xi_1$ applying the Efron-Stein inequality. We start with $\xi_1$ next.
\subparagraph{Part I-i: Bounding $\chi_1$'s residual bias $\xi_1$.} Following the definition of $U_k$ at the start of Part II of the proof of Lemma~\ref{lem:one-step-G-bound}, we write
\begin{align*}
    U_n:=\bz_n^*\wt{\bR}_{n-1}\bz_n\qquad\text{and}\qquad U_n':=\bz_n^*\wt{\bR}_{n-1}'\bz_n
\end{align*}
for the two terminal parameters $(z,0)$ and $(z',0)$, respectively. By applying the Sherman-Morrison formula to both resolvents, we can write
\begin{align}
    & \Tr \prn{\bQ\big(\wt{\bR}_n\wt{\bR}_n'-\wt{\bR}_{n-1}\wt{\bR}_{n-1}'\big)} \nonumber \\
    & =-\frac{\bz_n^*\wt{\bR}_{n-1}\wt{\bR}_{n-1}'\bQ\wt{\bR}_{n-1}\bz_n}{1+U_n}-\frac{\bz_n^*\wt{\bR}_{n-1}'\bQ\wt{\bR}_{n-1}\wt{\bR}_{n-1}'\bz_n}{1+U_n'}  +\frac{(\bz_n^*\wt{\bR}_{n-1}\wt{\bR}_{n-1}'\bz_n)(\bz_n^*\wt{\bR}_{n-1}'\bQ\wt{\bR}_{n-1}\bz_n)}{(1+U_n)(1+U_n')}. \hspace{2em}\label{eq:mixed-terminal-restore}
\end{align}
To bound $\xi_1=\E\chi_1$, it suffices to bound the $L^p$ norm of the right-hand side of Eq.~\eqref{eq:mixed-terminal-restore}, since the absolute value of its expectation is at most its $L^p$ norm by H\"older's inequality. To bound each term, we start by bounding the reciprocals of the denominators. By Lemma~\ref{lem:resolvent-bounds-strict-positivity}.(i) with $j=k=n$, $w_n=0$ and $\bx=\bx_n$, since $U_n=\bx_n^*\bR_{n-1}(z,0)\bx_n$ and the same identity holds at $z'$, we obtain
\begin{align*}
    |1+U_n|^{-1}\leq |z|/|z_-|\leq\Lambda_n\qquad\text{and}\qquad |1+U_n'|^{-1}\leq |z'|/|z'_-|\leq\Lambda_n'.
\end{align*}
For the numerators, they are the same type of quadratic forms conditioned on $\mc F_{n-1}$, by writing $\bA$ for the sandwiched matrices between $\bz_n^*$ and $\bz_n$ in any of the four quadratic forms in Eq.~\eqref{eq:mixed-terminal-restore}. Indeed, this allows us to invoke Hanson-Wright (cf.~Proposition~\ref{prop:hanson-wright}). Integrating the conditional Hanson-Wright tail as in Part I of the proof of Lemma~\ref{lem:one-step-covariance-bounds}, Eq.~\eqref{eq:one-step-conditional-quadratic-moments}, implies $\|\bz_n^*\bA\bz_n-\Tr\bA\|_{L^p}=O(p) \cdot \|\|\bA\|_F\|_{L^p}$ for all of the four terms above. Finally, by the moment bounds at $(s,j)=(n-1,n)$ (cf.~Eq.~\eqref{eq:trace-conditional-hanson-wright} and the paragraph after Eq.~\eqref{eq:trace-growing-resolvent-moments}), we have the following bounds, where the first line applies to the first two numerators in Eq.~\eqref{eq:mixed-terminal-restore}, and the second to each quadratic-form factor of its third numerator.
\begin{align}
    \|\Tr\bA\|_{L^p}&=O\prn{\frac{(\Lambda_n\vee\Lambda_n')^6}{n^2}},\quad
    \|\bz_n^*\bA\bz_n-\Tr\bA\|_{L^p}=O\prn{\frac{p(\Lambda_n\vee\Lambda_n')^6}{n^{5/2}}}
    ,\nonumber\\
    \|\Tr\bA\|_{L^{2p}}&=O\prn{\frac{(\Lambda_n\vee\Lambda_n')^4}{n}},\quad
    \|\bz_n^*\bA\bz_n-\Tr\bA\|_{L^{2p}}=O\prn{\frac{p(\Lambda_n\vee\Lambda_n')^4}{n^{3/2}}}
    .\label{eq:mixed-terminal-form-moments}
\end{align}
We point out that these bounds hold for $p=O(\log n)$. Combining the displays above and substituting them into Eq.~\eqref{eq:mixed-terminal-restore} with an additional application of H\"older, we conclude that
\begin{align}
    \left\|\Tr \prn{\bQ(\wt{\bR}_n\wt{\bR}_n'-\wt{\bR}_{n-1}\wt{\bR}_{n-1}')}\right\|_{L^p}
    \!=O\!\left(\frac{(\Lambda_n+\Lambda_n')(\Lambda_n\vee\Lambda_n')^6+\Lambda_n\Lambda_n'(\Lambda_n\vee\Lambda_n')^8}{n^2}\right)
    =O\prn{\frac{(\Lambda_n\vee\Lambda_n')^{10}}{n^2}}.\label{eq:mixed-trace-delete-moment}
\end{align}
Taking the expectation proves $|\xi_1|=O((\Lambda_n\vee\Lambda_n')^{10}/n^2)$.

\subparagraph{Part I-ii: Bounding $\chi_1$'s fluctuation $\xi_2$.}
We will invoke Efron-Stein in Lemma~\ref{lem:moment-efron-stein}, which boils down to bounding each $\|\xi_2-\xi_2^{(i)}\|_{L^p}$ for $1\leq i\leq n$, where $\xi_2^{(i)}$ denotes $\xi_2$'s replacement of $\bz_i$ by its independent copy, using the convention in Lemma~\ref{lem:moment-efron-stein}. By exchangeability, we take $i=n$ and let $\what{\bz}_n$ be a copy of $\bz_n$ independent of the entire sample. Canceling the common expectation in $\xi_2-\xi_2^{(n)}$ and subtracting the same mixed trace based on the first $n-1$ observations yields
\begin{align*}
    \xi_2-\xi_2^{(n)} =\Tr\prn{\bQ\big(\wt{\bR}_n\wt{\bR}_n'-\wt{\bR}_{n-1}\wt{\bR}_{n-1}'\big)}
    -\Tr\prn{\bQ\big(\wt{\bR}_n^{(n)}(\wt{\bR}_n')^{(n)}-\wt{\bR}_{n-1}\wt{\bR}_{n-1}'\big)},
\end{align*}
and this allows us to reuse the three terms analysis from Eq.~\eqref{eq:mixed-terminal-restore}. Let $U_n^{(n)}=\what{\bz}_n^*\wt{\bR}_{n-1}\what{\bz}_n$ and $(U_n')^{(n)}=\what{\bz}_n^*\wt{\bR}_{n-1}'\what{\bz}_n$. It boils down to controlling $\bz_n^* \bA \bz_n /(1+U_n) - \what{\bz}_n^* \bA \what{\bz}_n / (1 + U_n^{(n)})$ (and the ones measured at $z'$), for $\bA$ the numerator matrices following the same bounds in Eq.~\eqref{eq:mixed-terminal-form-moments}. Since $\bA$ and $G_n=G_n(z,0)$ are $\mc F_{n-1}$-measurable, both ratios have the same conditional law. Subtracting $\Tr\bA/(1+G_n)$ from each and applying the triangle inequality yields
\begin{align*}
    \left\|\frac{\bz_n^*\bA\bz_n}{1+U_n}-\frac{\what{\bz}_n^*\bA\what{\bz}_n}{1+U_n^{(n)}}\right\|_{L^p} \leq 2 \norm{\frac{\bz_n^*\bA\bz_n}{1+U_n}-\frac{\Tr\bA}{1+G_n}}_{L^p}
    =2\norm{\frac{\bz_n^*\bA\bz_n-\Tr\bA}{1+U_n}-\frac{\Tr\bA\,(U_n-G_n)}{(1+U_n)(1+G_n)}}_{L^p}.
\end{align*}
See Eq.~\eqref{eq:trace-ratio-centered-identity} in the proof of Lemma~\ref{lem:trace-centered-increments} for a similar analysis already done, and use its denominator bound at $(s,j)=(n-1,n)$, which implies $|1+G_n|^{-1}\leq\Lambda_n$. Combining with the same upper bound for $|1+U_n|^{-1} \leq \Lambda_n$ shown in Part I-i, and $\|U_n-G_n\|_{L^{2p}}=O(p\Lambda_n^2/\sqrt n)$ (cf.~Eq.~\eqref{eq:trace-conditional-hanson-wright} and the paragraph after Eq.~\eqref{eq:trace-growing-resolvent-moments}), it follows that
\begin{align*}
    &\left\|\frac{\bz_n^*\bA\bz_n}{1+U_n}-\frac{\what{\bz}_n^*\bA\what{\bz}_n}{1+U_n^{(n)}}\right\|_{L^p}\leq O\prn{(\Lambda_n\vee\Lambda_n')\|\bz_n^*\bA\bz_n-\Tr\bA\|_{L^p}
    +\frac{p(\Lambda_n\vee\Lambda_n')^4}{\sqrt n}\|\Tr\bA\|_{L^{2p}}}.
\end{align*}
For the first two differences in Eq.~\eqref{eq:mixed-terminal-restore}, use the first line of Eq.~\eqref{eq:mixed-terminal-form-moments} at $p$ for the centered numerator and at $2p$ for the trace. For the third difference, change one ratio at a time and apply H\"older to the two factors in $L^{2p}$. The second line of Eq.~\eqref{eq:mixed-terminal-form-moments} bounds the centered numerator in $L^{2p}$ as written and the trace in $L^{4p}$ upon replacing $p$ by $2p$. Combining all these bounds, we conclude that
\begin{align*}
    \|\xi_2-\xi_2^{(n)}\|_{L^p}
    \leq O\prn{2\frac{p(\Lambda_n\vee\Lambda_n')^{10}}{n^{5/2}}+2\frac{p(\Lambda_n\vee\Lambda_n')^8}{n^{3/2}}\frac{(\Lambda_n\vee\Lambda_n')^5}{n}}
    =O\prn{\frac{p(\Lambda_n\vee\Lambda_n')^{13}}{n^{5/2}}}.
\end{align*}
We are now ready to substitute this bound for every coordinate into Lemma~\ref{lem:moment-efron-stein}, Eq.~\eqref{eq:moment-efron-stein}, and we have
\begin{align*}
    \|\xi_2\|_{L^p}
    &\leq O(\sqrt p) \cdot \left(\sum_{i=1}^n\|\xi_2-\xi_2^{(i)}\|_{L^p}^2\right)^{1/2}
    =O\prn{\sqrt{np}\cdot \frac{p(\Lambda_n\vee\Lambda_n')^{13}}{n^{5/2}}}
    =O\prn{\frac{p^{3/2}(\Lambda_n\vee\Lambda_n')^{13}}{n^2}}.
\end{align*}
Combining the bound above with Part I-i and using $\chi_1=\xi_1+\xi_2$, we finish the proof of Part I.

\paragraph{Part II: Bounding the term $\chi_2$.} By Part III-i of the proof of Theorem~\ref{thm:self-consistent-path-stability}, Eq.~\eqref{eq:path-terminal-reciprocal-bound}, $|\delta_n\delta_n'|^{-1}=O((\Lambda_n\vee\Lambda_n')^2)$. Thus to prove the asserted bound for $|\chi_2|$, it suffices to show $|\delta_n\delta_n'\chi_2|=O((\Lambda_n\vee\Lambda_n')^{14}/n^2)$. Take $\bA=\wt{\bR}_{n-1}'\bQ\wt{\bR}_{n-1}$, so $\Tr\bA=\Tr(\bQ\wt{\bR}_{n-1}\wt{\bR}_{n-1}')$. Applying the terminal identity and then centering the numerator, we have
\begin{align}
    & \delta_n\delta_n'\chi_2
    \stackrel{\mathrm{(i)}}{=}\delta_n\delta_n'\E\Tr\bA-\E\brk{\bz_n^*\wt{\bR}_n' \bQ \wt{\bR}_n \bz_n} \stackrel{\mathrm{(ii)}}{=}\delta_n\delta_n'\E\Tr\bA-\E \brk{\frac{\bz_n^*\bA\bz_n}{(1+U_n)(1+U_n')}}\nonumber\\
    &=-\underbrace{\E\Tr\bA\cdot \left(\E\brk{\frac1{(1+U_n)(1+U_n')}}-\delta_n\delta_n'\right)}_{:=\xi_3}
    -\underbrace{\E\brk{(\bz_n^*\bA\bz_n-\E\Tr\bA)\left(\frac1{(1+U_n)(1+U_n')}-\delta_n\delta_n'\right)}}_{:=\xi_4}. \hspace{1.5em} \label{eq:mixed-terminal-freeze}
\end{align}

In (i), apply Lemma~\ref{lem:mixed-endpoint-identities}, Eq.~\eqref{eq:mixed-terminal-T-insertion}, and use exchangeability in $\bT_n=\sum_{i=1}^n\bx_i\bx_i^*$. At $z'=z$, use the diagonal identity of the same lemma and the local resolvent bound to differentiate under the expectation. In (ii), we apply Sherman-Morrison to $\wt{\bR}_n\bz_n$ and $\bz_n^*\wt{\bR}_n'$.

\subparagraph{Part II-i: Bounding $\chi_2$'s bias $\xi_3$.}
In $\xi_3$ from Eq.~\eqref{eq:mixed-terminal-freeze}, the factor $|\E\Tr\bA|$ is $O((\Lambda_n\vee\Lambda_n')^4/n)$ by the second line of Eq.~\eqref{eq:mixed-terminal-form-moments}. It suffices to bound the mean error of the reciprocal product. At $p=2$, combining Part I-ii with Eq.~\eqref{eq:trace-centered-old-samples} in Part I of the proof of Lemma~\ref{lem:trace-drift-components} at $(k,j)=(n,n)$ yields $\|U_n-\E U_n\|_{L^2}=O(\Lambda_n^2/\sqrt n+\Lambda_n^{11}/n)=O(\Lambda_n^3/\sqrt n)$, since $n=\Omega(\Lambda_n^{16})$. Applying the reciprocal expansion immediately before Eq.~\eqref{eq:one-step-deterministic-denominator-expansion} in Part III of the proof of Lemma~\ref{lem:one-step-G-bound}, centered at $\E U_n=\wb G_n(z,0)$, with the denominator bounds of Lemma~\ref{lem:resolvent-bounds-strict-positivity}.(i)--(ii), Eqs.~\eqref{eq:helper-path-resolvent-quadratic-bounds} and~\eqref{eq:helper-denominator-separation}, yields $\|(1+U_n)^{-1}-\delta_n\|_{L^2}=O(\Lambda_n^5/\sqrt n)$ and $|\E(1+U_n)^{-1}-\delta_n|=O(\Lambda_n^9/n)$, since the linear term has zero mean.

As the same estimates hold at $z'$, we can expand the product in Eq.~\eqref{eq:mixed-terminal-freeze} and apply the mean bounds to the two linear terms. By Cauchy-Schwarz to the third term below, we finally have
\begin{align*}
    \xi_3
    &=\E\Tr\bA\cdot\E\brk{\delta_n'\left(\frac1{1+U_n}-\delta_n\right)+\delta_n\left(\frac1{1+U_n'}-\delta_n'\right)+\left(\frac1{1+U_n}-\delta_n\right)\left(\frac1{1+U_n'}-\delta_n'\right)}\\
    &=O\prn{\frac{(\Lambda_n\vee\Lambda_n')^4}{n}\cdot\frac{\Lambda_n'\Lambda_n^9+\Lambda_n(\Lambda_n')^9+\Lambda_n^5(\Lambda_n')^5}{n}}
    =O\prn{\frac{(\Lambda_n\vee\Lambda_n')^{14}}{n^2}}.
\end{align*}

\subparagraph{Part II-ii: Bounding $\chi_2$'s covariance $\xi_4$.} We bound the two terms (the numerator and the denominator fluctuations) in $\xi_4$ separately. For the numerator, the second line of Eq.~\eqref{eq:mixed-terminal-form-moments} at $p=2$ bounds $\|\bz_n^*\bA\bz_n-\Tr\bA\|_{L^2}$ by $O((\Lambda_n\vee\Lambda_n')^4/n^{3/2})$. To change the center to $\E\Tr\bA$, recall that $\Tr\bA=\Tr(\bQ\wt{\bR}_{n-1}\wt{\bR}_{n-1}')$. Inserting the full trace yields
\begin{align*}
    \Tr\bA-\E\Tr\bA=\prn{\Tr\bA-\Tr(\bQ\wt{\bR}_n\wt{\bR}_n')}+\prn{\Tr(\bQ\wt{\bR}_n\wt{\bR}_n')-\E\Tr\bA} = \prn{\Tr\bA-\Tr(\bQ\wt{\bR}_n\wt{\bR}_n')} + \chi_1.
\end{align*}
At $p=2$, applying Eq.~\eqref{eq:mixed-trace-delete-moment} to the first difference yields an $L^2$ bound of $O((\Lambda_n\vee\Lambda_n')^{10}/n^2)$. Combining this with $\|\chi_1\|_{L^2}=O((\Lambda_n\vee\Lambda_n')^{13}/n^2)$ from Part I yields $\|\Tr\bA-\E\Tr\bA\|_{L^2}=O((\Lambda_n\vee\Lambda_n')^{13}/n^2)$. Adding this to the centered bound yields $\|\bz_n^*\bA\bz_n-\E\Tr\bA\|_{L^2}=O((\Lambda_n\vee\Lambda_n')^{5}/n^{3/2})$, since $n=\Omega((\Lambda_n\vee\Lambda_n')^{16})$. Indeed, Assumption~\ref{assmp:spectral-stability} at both parameters implies $(\Lambda_n\vee\Lambda_n')^8=O(\sqrt n)$, so the term with exponent $13$ contributes at most $O((\Lambda_n\vee\Lambda_n')^5/n^{3/2})$.
For the denominator term, using Part II-i and changing the two factors one at a time, we have
\begin{align*}
    \left\|\frac1{(1+U_n)(1+U_n')}-\delta_n\delta_n'\right\|_{L^2}
    &=\left\|\left(\frac1{1+U_n}-\delta_n\right)\frac1{1+U_n'}+\delta_n\left(\frac1{1+U_n'}-\delta_n'\right)\right\|_{L^2}\\
    &\leq\Lambda_n'\left\|\frac1{1+U_n}-\delta_n\right\|_{L^2}+|\delta_n|\left\|\frac1{1+U_n'}-\delta_n'\right\|_{L^2}
    =O\prn{\frac{(\Lambda_n\vee\Lambda_n')^6}{\sqrt n}}.
\end{align*}
Applying Cauchy-Schwarz in Eq.~\eqref{eq:mixed-terminal-freeze} therefore implies $|\xi_4|=O(((\Lambda_n\vee\Lambda_n')^{5}/n^{3/2})((\Lambda_n\vee\Lambda_n')^6/\sqrt n))\allowbreak=O((\Lambda_n\vee\Lambda_n')^{11}/n^2)$.
Combining the bounds for $\xi_3$ and $\xi_4$ in Eq.~\eqref{eq:mixed-terminal-freeze} gives $|\delta_n\delta_n'\chi_2|=O((\Lambda_n\vee\Lambda_n')^{14}/n^2)$. Applying Eq.~\eqref{eq:path-terminal-reciprocal-bound} from Part III-i of the proof of Theorem~\ref{thm:self-consistent-path-stability} yields $|\chi_2|=O((\Lambda_n\vee\Lambda_n')^{16}/n^2)$.

\subsubsection{Proof of Lemma~\ref{lem:mixed-single-trace-comparison}}\label{proof:mixed-single-trace-comparison}
\paragraph{Part I: Bounding the term $\chi_3$.}
Using the definition of $\chi_2$ to substitute for the terminal divided difference in $\chi_3$ yields
\begin{align*}
    \chi_3=\left(1-\frac{\delta_n\delta_n'}{\delta\lin(z)\delta\lin(z')}\right)
    \left(\E\Tr(\bQ\wt{\bR}_{n-1}\wt{\bR}_{n-1}')-\chi_2\right).
\end{align*}
The second factor is $O((\Lambda_n\vee\Lambda_n')^4/n)$. Indeed, taking $\bA=\wt{\bR}_{n-1}'\bQ\wt{\bR}_{n-1}$ in the second line of Eq.~\eqref{eq:mixed-terminal-form-moments} in Part I-i of the proof of Lemma~\ref{lem:terminal-mixed-replacement} gives this bound for the expected trace. The bound $|\chi_2|=O((\Lambda_n\vee\Lambda_n')^{16}/n^2)$ in that lemma is absorbed since $n=\Omega((\Lambda_n\vee\Lambda_n')^{16})$.

For the first factor, Part IV-iii of the proof of Theorem~\ref{thm:self-consistent-path-stability} bounds $|\delta_n-\delta\lin(z)|$ by $O(\Lambda_n^{17}/n)$. By Lemma~\ref{lem:population-scalar-stability}, $|\delta\lin(z)|=|zs_n(z)|\leq |z|/|z_-|\leq\Lambda_n$, and $|\delta_n|\leq\Lambda_n$ by Lemma~\ref{lem:resolvent-bounds-strict-positivity}.(ii), Eq.~\eqref{eq:helper-denominator-separation}. Applying Eqs.~\eqref{eq:path-terminal-reciprocal-bound} and~\eqref{eq:path-linear-reciprocal-transfer} in the proof of Theorem~\ref{thm:self-consistent-path-stability} at $z$ and $z'$ yields $|\delta\lin(z)\delta\lin(z')|^{-1}=O((\Lambda_n\vee\Lambda_n')^2)$. Changing one increment at a time then yields
\begin{align*}
    \left|1-\frac{\delta_n\delta_n'}{\delta\lin(z)\delta\lin(z')}\right|
    \leq\frac{|\delta\lin(z)-\delta_n|\,|\delta\lin(z')|+|\delta_n|\,|\delta\lin(z')-\delta_n'|}{|\delta\lin(z)\delta\lin(z')|}
    =O\prn{\frac{(\Lambda_n\vee\Lambda_n')^{20}}{n}}.
\end{align*}
Multiplying the two factors yields $|\chi_3|=O((\Lambda_n\vee\Lambda_n')^{24}/n^2)$.

\paragraph{Part II: Bounding the term $\chi_4$.}
For the final term, it suffices to bound the divided difference $n\delta\lin(z)\delta\lin(z')\chi_4$ by $O((\Lambda_n\vee\Lambda_n')^{24}/n)$, because applying Eqs.~\eqref{eq:path-terminal-reciprocal-bound} and~\eqref{eq:path-linear-reciprocal-transfer} in the proof of Theorem~\ref{thm:self-consistent-path-stability} at $z$ and $z'$ bounds $|n\delta\lin(z)\delta\lin(z')|^{-1}$ by $O((\Lambda_n\vee\Lambda_n')^2/n)$, and the desired bound on $\chi_4$ then follows.

Let $\Gamma$ be the positively oriented boundary of the union of the disks centered at $z$ and $z'$ with radii $|z_-|/4$ and $|z'_-|/4$, respectively. We can apply Cauchy's formula on this union to obtain
\begin{align*}
    n\delta\lin(z)\delta\lin(z')\chi_4
    =\frac1{2\pi\iu}\oint_\Gamma
    \frac{\zeta\brk{\E G_{n+1}^{\bQ}(\zeta,0)-G_1^{\bQ}(\zeta,w_0\lin(\zeta))}}{(\zeta-z)(\zeta-z')}\,\de\zeta.
\end{align*}
This allows us to invoke Theorem~\ref{thm:terminal-weighted-traces} for the single resolvents. We bound the other terms first. Since $\zeta\in\Gamma$ lies outside the open disk centered at $z$, we have $|\zeta-z|\geq |z_-|/4$. By the triangle inequality for the distance to $\R_{\geq0}$,
\begin{align*}
    |\zeta_-|\leq |\zeta-z|+|z_-|\leq 5|\zeta-z|.
\end{align*}
The same argument at $z'$ implies $|\zeta-z'|\geq |z'_-|/4$ and $|\zeta_-|\leq |\zeta-z'|+|z'_-|\leq 5|\zeta-z'|$. On the boundary arc centered at $z$, the latter bound yields $|\zeta_-|/(|\zeta-z|\,|\zeta-z'|)\leq 5/|\zeta-z|$. Integrating along this arc gives at most $10\pi$, since $|\zeta-z|$ is its radius. Interchanging $z$ and $z'$ gives the same integral bound on the other arc, and adding the two contributions and dividing by $2\pi$ yields
\begin{align*}
    |n\delta\lin(z)\delta\lin(z')\chi_4|
    &\leq10\sup_{\zeta\in\Gamma}\frac{|\zeta|}{|\zeta_-|}\left|\E G_{n+1}^{\bQ}(\zeta,0)-G_1^{\bQ}(\zeta,w_0\lin(\zeta))\right|
    =O\prn{\frac{(\Lambda_n\vee\Lambda_n')^{24}}{n}}.
\end{align*}
We apply Theorem~\ref{thm:terminal-weighted-traces} with $p=2$ and $r=0$, followed by Jensen's inequality in the last step. As in Part III of the proof of Theorem~\ref{thm:terminal-weighted-traces}, on each circle $\Lambda_n(\zeta)=O(\Lambda_n\vee\Lambda_n')$ with the same minimizing cutoff $m$ from Assumption~\ref{assmp:spectral-stability} and distance exponent $\nu+1$, using reflection with $\bQ^*$ for potential values of $\zeta$ in the lower half-plane. The remaining factor satisfies $|\zeta|/|\zeta_-|\leq\Lambda_n(\zeta)$ by Assumption~\ref{assmp:spectral-stability}. The proof is complete.

    \section{Supplementary details for the numerical experiments}\label{app:numerical-details}

We provide further details on the setups, formulas, and numerical methods in experiments. The corresponding code is available in the \href{https://github.com/Moriartycc/dimension-free-rmt}{public repository}.

\paragraph{Figure~\ref{fig:intro-gram-spectra}.}
For CIFAR-10, we use a fixed ImageNet-pretrained feature extractor, center the features using the training-set mean, and scale them so that the training covariance has largest eigenvalue one. Each trial samples $n$ distinct test images, and histogram heights use all $n$ eigenvalues. The colors indicate pooled eigenvalue quantiles, and cropped insets preserve the full-distribution normalization.

For the polynomial model, we generate $d=\max\{2^{14},256n\}$ Gaussian coordinates and approximate the remaining contribution to $\bX\bX^*/n$ by its expectation, $\tau_1(d)\bI_n/n$. This differs from the fixed covariance truncation in Section~\ref{sec:numerical}.

\paragraph{Figure~\ref{fig:deterministic-bulk-laws}.}
We numerically solve for the boundary values of the self-consistent equation and recover the densities by Stieltjes inversion. We sum the population terms explicitly up to a cutoff, then separate the first-order tail and approximate the quadratic remainder by quadrature. The limiting densities use the parametrizations in Appendix~\ref{proof:canonical-bulk-laws}.

\paragraph{Figure~\ref{fig:directional-measures}.}
Each $\mu_{\omega,j}$ is a probability measure, including its atom at zero. We plot $\mu_\omega$ without renormalization, so that it has mass $\Tr(\bSigma^2)/n$ and no atom at zero.

\paragraph{Figure~\ref{fig:interpolation-endpoint-errors}.}
We use $n=400$ Gaussian observations under the power law and critical law with $\alpha=\beta=0.2$, truncated at $d=2^{16}$. With $a_n=n\sigma_n$ and $b_n=\sum_{j=n+1}^{d}\sigma_j$, we plot $z=a_n\zeta$ and $z=b_n+a_n\zeta$, respectively. Each grid has 25 equally spaced real coordinates on the displayed interval and 17 geometrically spaced heights.

For backward, linear, and forward interpolation in Eqs.~\eqref{eq:self-consistent-expected-trace-update}--\eqref{eq:self-consistent-forward-update}, respectively, the rows estimate these errors with $p=1$ and $p=2$:
\begin{equation*}
    \frac{\prn{\norm{G_n(z,0)-G_1(z,w_0)}_{L^p},\quad \norm{G_n(z,0)-G_1(z,w_0\lin)}_{L^p},\quad \norm{G_n(z,\mu_n)-G_1(z,w_0\lin)}_{L^p}}}{1+|G_1(z,w_0\lin)|}.
\end{equation*}
We compute the mean absolute and root mean squared errors over 400 trials shared across methods and grid points, using 1,024 independent draws to approximate the backward expectations. Here $G_n$ uses $n-1$ observations, and $\mu_n$ is the actual forward endpoint. Forward updates use fixed-point iteration with a root-solver fallback.

\subsection{Additional experiments for ridge-type filters}\label{app:additional-ridge-experiments}
We begin with computing the effective risk formulas when $k=2$. By Fa\`a di Bruno's formula~\citep{apostol2000calculating}, also used by~\citet{cheng2024dimensionfree} for their calculations, from Eq.~\eqref{eq:ridge-effective-regularization}, we can calculate the derivatives below, where $\lambda_\star=-\omega(-\lambda)$ and
\begin{align*}
    \lambda=\lambda_\star\prn{1-\frac1n\Tr\prn{\bSigma(\bSigma+\lambda_\star\bI)^{-1}}}
    \qquad\text{and}\qquad
    \mathsf{J}_n(\lambda)=\frac1n\Tr\prn{\bSigma^2(\bSigma+\lambda_\star\bI)^{-2}}.
\end{align*}
It follows from that formula:
\begin{align*}
    \omega'(-\lambda)&=\frac{1}{1-\mathsf{J}_n(\lambda)}, \qquad \qquad \omega''(-\lambda)=\frac{2\Tr\prn{\bSigma^2(\bSigma+\lambda_\star\bI)^{-3}}}{n(1-\mathsf{J}_n(\lambda))^3},\\
    \omega'''(-\lambda)&=\frac{6\Tr\prn{\bSigma^2(\bSigma+\lambda_\star\bI)^{-4}}}{n(1-\mathsf{J}_n(\lambda))^4}
    +\frac{12\prn{\Tr\prn{\bSigma^2(\bSigma+\lambda_\star\bI)^{-3}}}^2}{n^2(1-\mathsf{J}_n(\lambda))^5}.
\end{align*}
Substituting into Proposition~\ref{prop:ridge-type-risk}.(i), we have the explicit variance and squared bias for $k=2$:
\begin{align*}
    \effectiveVariance[\upphi_{2,\lambda}]
    &=\frac{\tau^2}{1-\mathsf{J}_n(\lambda)}\Bigg\{\mathsf{J}_n(\lambda)
    +\frac{2\lambda\Tr\prn{\bSigma^2(\bSigma+\lambda_\star\bI)^{-3}}}{n(1-\mathsf{J}_n(\lambda))^2}\\
    &\qquad+\frac{\lambda^2\Tr\prn{\bSigma^2(\bSigma+\lambda_\star\bI)^{-4}}}{n(1-\mathsf{J}_n(\lambda))^3}
    +\frac{2\lambda^2\prn{\Tr\prn{\bSigma^2(\bSigma+\lambda_\star\bI)^{-3}}}^2}{n^2(1-\mathsf{J}_n(\lambda))^4}\Bigg\},\\
    \effectiveBias[\upphi_{2,\lambda}]
    &=\frac{\lambda_\star^2\bbeta^*\bSigma(\bSigma+\lambda_\star\bI)^{-2}\bbeta}{1-\mathsf{J}_n(\lambda)}\Bigg\{1
    +\frac{2\lambda\Tr\prn{\bSigma^2(\bSigma+\lambda_\star\bI)^{-3}}}{n(1-\mathsf{J}_n(\lambda))^2} +\frac{\lambda^2\Tr\prn{\bSigma^2(\bSigma+\lambda_\star\bI)^{-4}}}{n(1-\mathsf{J}_n(\lambda))^3} \\
    &\qquad
    +\frac{2\lambda^2\prn{\Tr\prn{\bSigma^2(\bSigma+\lambda_\star\bI)^{-3}}}^2}{n^2(1-\mathsf{J}_n(\lambda))^4}\Bigg\} -\frac{2\lambda\lambda_\star}{(1-\mathsf{J}_n(\lambda))^2}
    \brc{1+\frac{\lambda\Tr\prn{\bSigma^2(\bSigma+\lambda_\star\bI)^{-3}}}{n(1-\mathsf{J}_n(\lambda))^2}}
    \\
    & \qquad \times \bbeta^*\bSigma^2(\bSigma+\lambda_\star\bI)^{-3}\bbeta +\frac{\lambda^2}{(1-\mathsf{J}_n(\lambda))^3}\bbeta^*\bSigma^3(\bSigma+\lambda_\star\bI)^{-4}\bbeta.
\end{align*}
Figures~\ref{fig:ridge-q24-k1-absolute} to \ref{fig:ridge-q24-k2-relative} repeat Figures~\ref{fig:ridge-q16-k1-absolute} to \ref{fig:ridge-q16-k2-relative} with $q=24$. We keep all other settings and plotting conventions.

\begin{figure}[!htbp]
    \renewcommand{\thefigure}{\arabic{figure}a}
    \renewcommand{\theHfigure}{\arabic{figure}a}
    \centering
    \includegraphics[width=\linewidth]{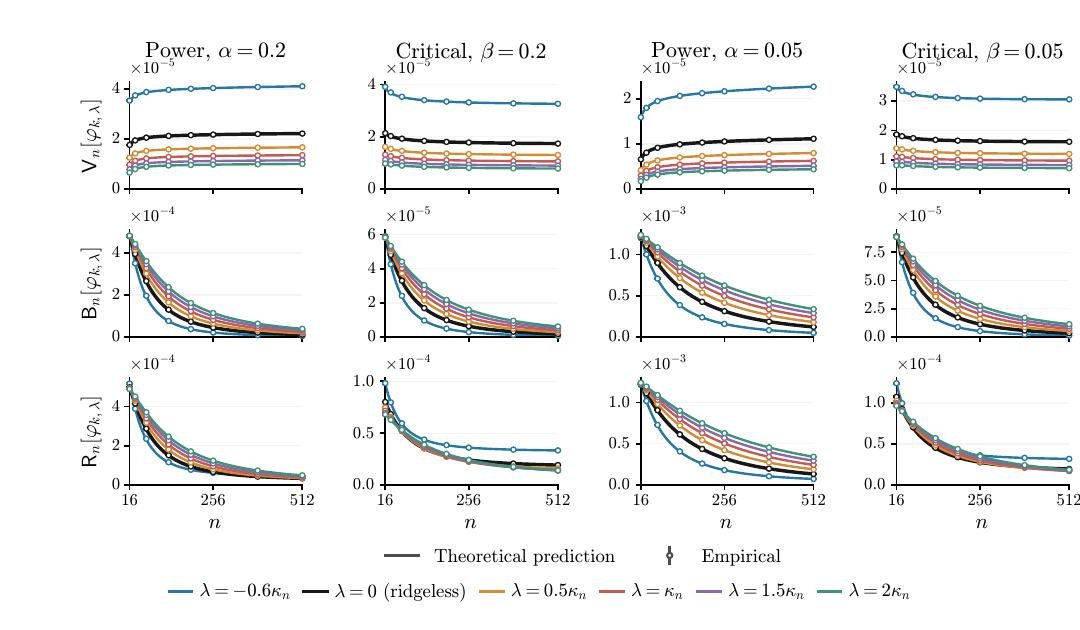}
    \caption{Ridge-type filters with $q=24$ and $k=1$. Rows show $\effectiveVariance{}$, $\effectiveBias{}$, and $\effectiveRisk{}$.}
    \label{fig:ridge-q24-k1-absolute}
\end{figure}

\begin{figure}[!htbp]
    \addtocounter{figure}{-1}
    \renewcommand{\thefigure}{\arabic{figure}b}
    \renewcommand{\theHfigure}{\arabic{figure}b}
    \centering
    \includegraphics[width=\linewidth]{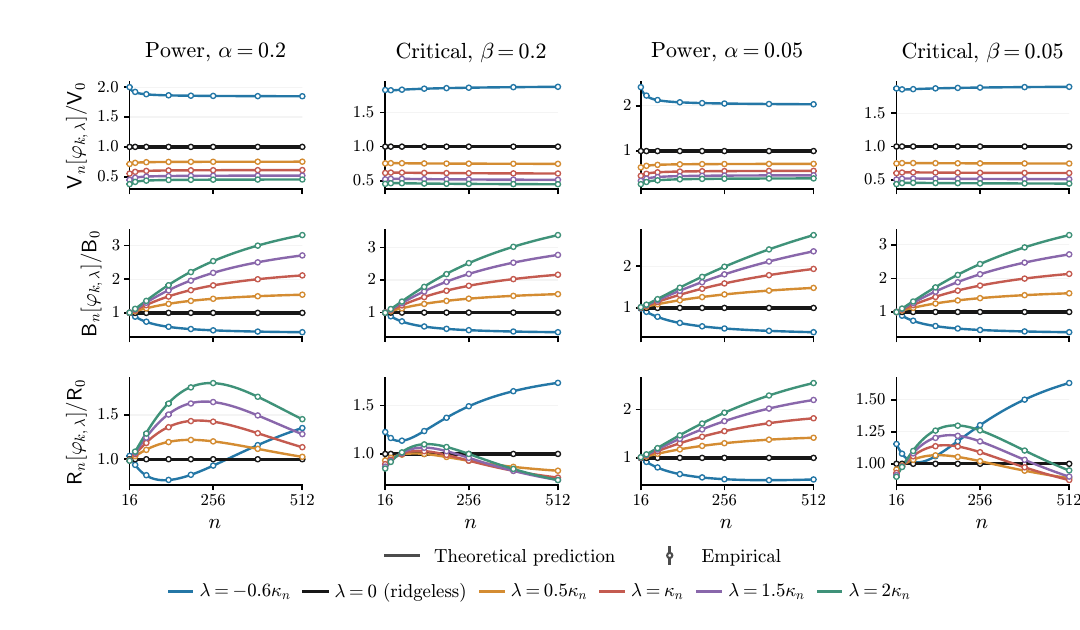}
    \caption{Ridge-type filters with $q=24$ and $k=1$, normalized by the respective ridgeless baselines.}
    \label{fig:ridge-q24-k1-relative}
\end{figure}

\begin{figure}[!htbp]
    \addtocounter{figure}{-1}
    \renewcommand{\thefigure}{\arabic{figure}c}
    \renewcommand{\theHfigure}{\arabic{figure}c}
    \centering
    \includegraphics[width=\linewidth]{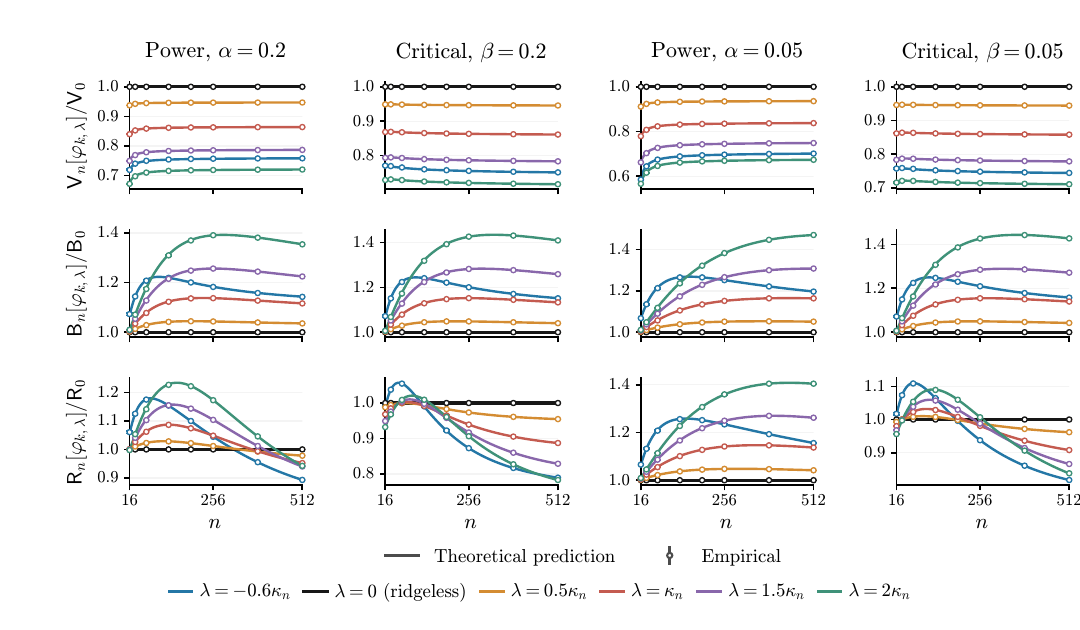}
    \caption{Iterated ridge with $q=24$ and $k=2$, normalized by the respective ridgeless baselines.}
    \label{fig:ridge-q24-k2-relative}
\end{figure}

\subsection{Additional experiments for gradient flow}\label{app:additional-gradient-flow-experiments}
We explain our numerical approach to find the $t_{\mathrm{opt}}$ curve first. We evaluate the risk in Proposition~\ref{prop:gradient-flow-risk} and its derivative using the formulas in Section~\ref{sec:gradient-flow}. We locate sign changes of the risk derivative on a logarithmic grid in $\wt{t}=t\tau_1(n)/n$, refine each detected root by Brent's method in logarithmic time, and compare the risks at these stationary points with those at $t=0$ and the ridgeless limit $t=\infty$. The smallest value among these candidates determines the plotted optimum. For Figure~\ref{fig:numerical-gradient-flow}, we do not assume a unique stationary point, and the search uses $151$ logarithmically spaced times in $[10^{-12},10^3]$, separately from the $91$ displayed times in $[0,4]$.

\begin{figure}[!htbp]
    \centering
    \includegraphics[width=\linewidth]{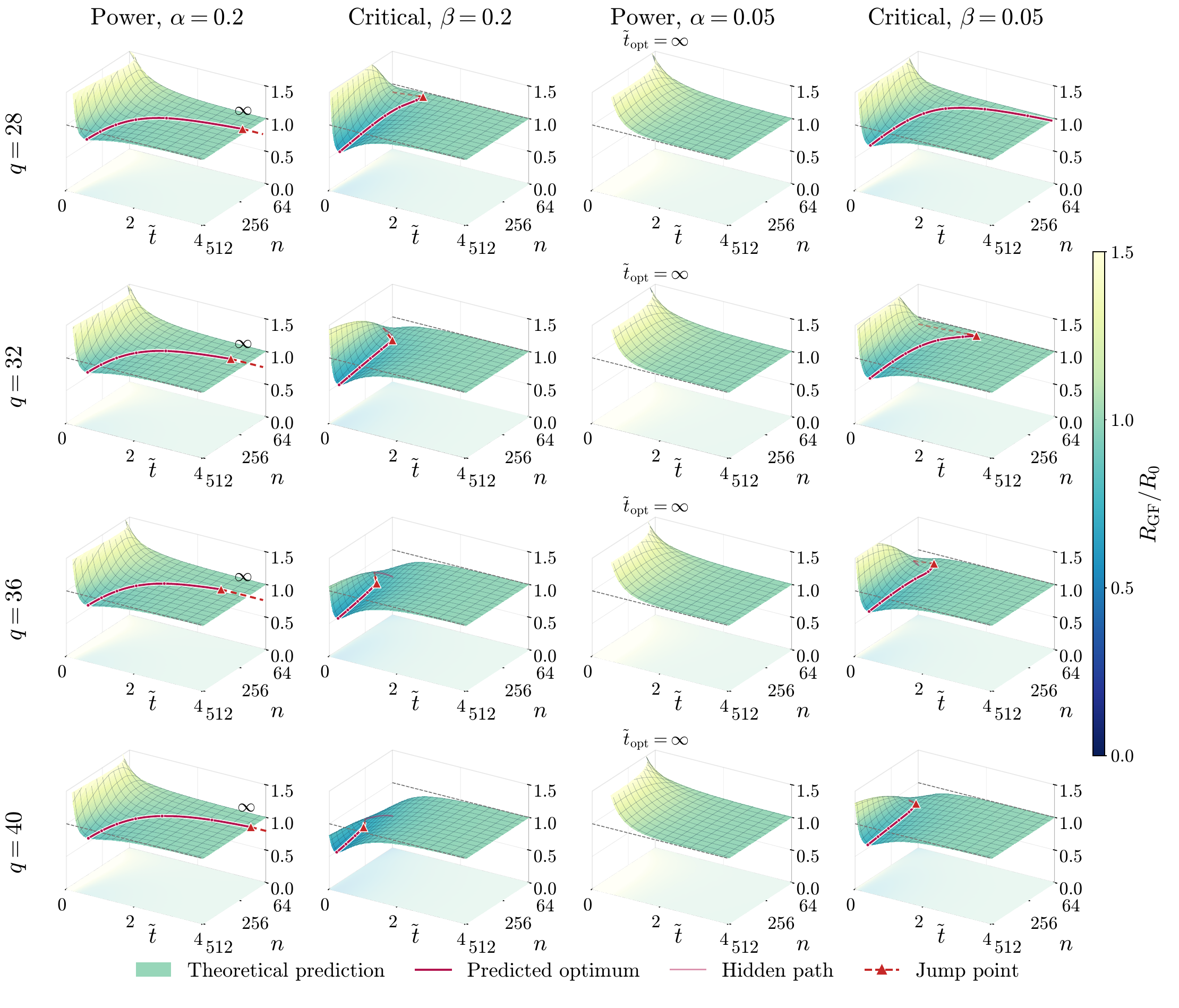}
    \caption{Gradient-flow theoretical risk surfaces relative to ridgeless for $q=28,32,36,40$ (top to bottom), with columns as in Figure~\ref{fig:numerical-gradient-flow}. \textcolor[HTML]{B3154F}{Curves and points} show predicted finite stopping optima. \textcolor[HTML]{C62828}{Triangles and dashed connectors} indicate jumps between sampled minimizing branches.}
    \label{fig:gradient-flow-source-transitions}
\end{figure}

Next, Figure~\ref{fig:gradient-flow-source-transitions} supplements these results with theoretical risk surfaces for $q\in\{28,32,36,40\}$, using $161$ search times in $[10^{-10},10^4]$ and otherwise the same settings. Relative improvements over ridgeless of at most $10^{-8}$ are treated as numerical ties. Among these risk surfaces, the critical law with $\beta=0.2$ and $q=32$ illustrates a jump between two finite stopping times, whereas the power law with $\alpha=0.2$ exhibits an apparent ridgeless-to-finite transition. It remains intriguing to investigate these complicated behaviors further.

\end{document}